\documentclass[11pt]{article}
\usepackage[T1]{fontenc}
\usepackage[utf8]{inputenc}
\usepackage{lmodern}
\usepackage[margin=1in]{geometry}
\usepackage{amsmath,amssymb,amsthm,mathtools,aliascnt}
\usepackage{microtype,graphicx,booktabs,tabularx,enumitem,mathrsfs}
\usepackage[colorlinks=true,linkcolor=blue,citecolor=blue,urlcolor=blue]{hyperref}
\usepackage[nameinlink,noabbrev]{cleveref}
\newtheorem{theorem}{Theorem}[section]
\newaliascnt{proposition}{theorem}
\newtheorem{proposition}[proposition]{Proposition}
\aliascntresetthe{proposition}
\crefname{proposition}{Proposition}{Propositions}
\newaliascnt{lemma}{theorem}
\newtheorem{lemma}[lemma]{Lemma}
\aliascntresetthe{lemma}
\crefname{lemma}{Lemma}{Lemmas}
\newaliascnt{corollary}{theorem}
\newtheorem{corollary}[corollary]{Corollary}
\aliascntresetthe{corollary}
\crefname{corollary}{Corollary}{Corollaries}
\theoremstyle{definition}
\newaliascnt{definition}{theorem}

\aliascntresetthe{definition}
\crefname{definition}{Definition}{Definitions}
\newaliascnt{problem}{theorem}
\newtheorem{problem}[problem]{Problem}
\aliascntresetthe{problem}
\crefname{problem}{Problem}{Problems}
\theoremstyle{remark}
\newaliascnt{remark}{theorem}
\newtheorem{remark}[remark]{Remark}
\aliascntresetthe{remark}
\crefname{remark}{Remark}{Remarks}
\newcommand{\R}{\mathbb R}
\newcommand{\C}{\mathbb C}
\newcommand{\Z}{\mathbb Z}
\newcommand{\N}{\mathbb N}
\newcommand{\A}{A}
\newcommand{\fc}[1]{\widehat{#1}_{\mathrm c}}
\newcommand{\PW}{\operatorname{PW}}
\newcommand{\supp}{\operatorname{supp}}
\newcommand{\conv}{\operatorname{conv}}
\newcommand{\dist}{\operatorname{dist}}
\newcommand{\vol}{\operatorname{vol}}
\newcommand{\sinc}{\operatorname{sinc}}
\newcommand{\norm}[1]{\lVert#1\rVert}
\newcommand{\dd}{\,\mathrm d}

\newcommand{\interior}{\operatorname{int}}

\setlist{itemsep=3pt,topsep=5pt}
\newcommand{\E}{E}
\newcommand{\Epol}{E^*}
\newcommand{\PiN}{\Pi_n^d}

\hypersetup{pdfauthor={Xinpeng Huang},pdftitle={Polynomial spaces are locally Paley\textendash Wiener: spectral geometry and sampling redundancy},pdfsubject={Local spectral geometry and sampling redundancy for polynomial spaces},bookmarksdepth=2}
\title{Polynomial spaces are locally Paley--Wiener:\\[4pt]
\large spectral geometry and sampling redundancy}
\author{Xinpeng Huang\\[3pt]
\small School of Geosciences and Info-Physics\\
\small Central South University, Changsha, China}
\date{29 November 2026}
\newif\ifinternaldiscussion
\internaldiscussionfalse 
\begin{document}
\maketitle
\begin{abstract}
We identify total-degree polynomial spaces on centrally symmetric convex
bodies locally with Paley--Wiener spaces. This correspondence makes
local spectral geometry a Fourier-analytic tool for studying polynomial
sampling and interpolation. For a Bernstein--Markov measure
with positive continuous density in the interior, the reproducing kernel
at scale $1/n$ about an interior point $x$, normalized by its diagonal
value at $x$, converges to the normalized Fourier transform of the
indicator of $A_x$. Here $A_x$ is the dual unit ball of the Baran
metric, which we call the local spectral body. Local recovery with
$L^2$ norm control gives bulk universality and Christoffel asymptotics
uniformly on compact interior subsets. The diagonal limit identifies
the Christoffel sampling density with the local density of degrees of
freedom, proportional to $|A_x|$. The off-diagonal limit retains the
shape of $A_x$.

Local spectral geometry complements Christoffel sampling theory by
relating the shape of the local spectrum to sampling redundancy.
Let $N_n$ denote the dimension of the polynomial space. Uniformly stable
sampling with positive weights and $N_n+o(N_n)$ nodes is possible
precisely when each $L^2(A_x)$ admits a normalized exponential Riesz
basis with bounds uniform in $x$. The regular hexagon admits such
sampling. On balls in dimensions $d\ge2$, the independent results of
Ortega-Cerd\`a and Wan on exponential Riesz bases imply a positive
relative redundancy at every fixed stability level.

For asymptotically isometric sampling, where the sampling bounds tend
to one, the criterion for vanishing relative redundancy is the existence of an exponential
orthonormal basis on every local spectrum. Euclidean packing bounds
quantify the optimal node ratio and determine it exactly on balls in
dimensions eight and twenty-four. Among centrally symmetric convex
bodies, precisely the intervals, parallelograms, and centrally symmetric
hexagons admit asymptotically isometric sampling with vanishing relative
redundancy. On every other body in this class, including all those in
dimensions $d\ge3$, a positive proportion of nodes beyond $N_n$ is
necessary for positive cubature exact through degree $2n$.
\end{abstract}



\clearpage
\setcounter{tocdepth}{1}
\begingroup
\fontsize{9}{10}\selectfont
\makeatletter
\let\paper@tocaddvspace\addvspace
\renewcommand{\addvspace}[1]{\paper@tocaddvspace{\glueexpr#1*2/5\relax}}
\makeatother
\tableofcontents
\endgroup
\clearpage
\section*{General Introduction}
\addcontentsline{toc}{section}{General Introduction}

Where should a function space be sampled, and how many observations
are needed for stable recovery? Christoffel sampling connects these
questions to the spatial distribution of degrees of freedom. For an
$N$-dimensional space $V\subset C(E)$ with its $L^2(\mu)$ norm and
orthonormal basis $\phi_1,\ldots,\phi_N$, the reproducing kernel satisfies
\[
 K_V(x,x)=\sum_{j=1}^N|\phi_j(x)|^2,
 \qquad \int_E K_V(x,x)\,d\mu(x)=N.
\]
The measure $K_V(x,x)\,d\mu(x)$ therefore distributes the dimension
over the domain. Its mass on a region is the trace of the corresponding
concentration operator, giving a precise sense in which it measures
local degrees of freedom. We use the classical convention that the
Christoffel function is $\lambda_V(x)=K_V(x,x)^{-1}$ wherever
$K_V(x,x)>0$.

Sampling independently from
\[
 d\rho_V(x)=\frac{K_V(x,x)}{N}\,d\mu(x)
\]
and assigning weights $N/(mK_V(x_j,x_j))$ to $m$ samples makes
the empirical Gram matrix unbiased. With $m$ of order $N\log N$,
this gives stable norm discretization with high probability.
This strategy appears as coherence-optimal sampling in
Hampton--Doostan \cite{HamptonDoostan2015} and as optimal weighted
least squares in Cohen--Migliorati \cite{CohenMigliorati2017};
Adcock \cite{Adcock2025} gives a systematic account under the name
Christoffel sampling. For polynomial spaces, Christoffel asymptotics
connect this sampling density to the equilibrium measure. This also
motivates the equilibrium-measure sampling and Christoffel weights
of Narayan, Jakeman, and Zhou \cite{NarayanJakemanZhou2017}.

These connections suggest a geometric approach to sampling:
understand the local degrees of freedom, distribute observations
accordingly, and determine the additional sampling needed for
stability. We develop this approach for total-degree polynomial spaces
by identifying their local Paley--Wiener models. The diagonal
Christoffel asymptotic and the full universality limit then describe
two aspects of the same geometry: spectral volume and spectral shape.
The latter leads to a criterion for stable sampling with asymptotically
one observation per coefficient.

\paragraph{Stable sampling and relative redundancy.}
Write $N_n=\dim\Pi_n^d=\binom{n+d}{d}$. We seek distinct nodes
$\Lambda_n\subset E$ and positive weights satisfying the
Marcinkiewicz--Zygmund inequalities
\[
 \alpha\|p\|_{L^2(\mu)}^2\le
 \sum_{\lambda\in\Lambda_n}\omega_{n,\lambda}|p(\lambda)|^2
 \le\beta\|p\|_{L^2(\mu)}^2,
 \qquad p\in\Pi_n^d,
\]
with fixed $0<\alpha\le\beta<\infty$. The relative redundancy is
$|\Lambda_n|/N_n-1$, and we call the cardinality
\emph{asymptotically critical} when this ratio tends to zero.
General finite-dimensional results provide upper bounds on sampling
budgets. Dolbeault--Cohen \cite{DolbeaultCohen2022} obtain near-best
approximation in expected $L^2$ error from a constant times the
dimension; the sparsification theorem of Batson, Spielman, and
Srivastava \cite{BSS2014} gives positive norm discretization with any
fixed relative excess and stability bounds depending on that excess.
Our question is when the relative excess can vanish while the
stability bounds remain uniform. We allow general positive weighted
arrays, so the question concerns the sampling capacity of the space
beyond any particular random sampling algorithm.

\paragraph{Christoffel asymptotics and universality through local spectra.}
Let $E=-E\subset\R^d$ be a convex body and let $\mu$ be a finite
positive Bernstein--Markov measure with support $E$. On an interior
open set $U$, assume $d\mu=w\,dx$ with $w$ positive and continuous.
If $E^*$ is the polar body, define
\[
 A_x=\conv\left\{
 \frac{y}{\sqrt{1-(x\cdot y)^2}}:y\in E^*\right\},
 \qquad q_E(x)=\frac{|A_x|}{(2\pi)^d}.
\]
The body $A_x$ is the dual unit ball of the Baran metric. Its role
in the equilibrium density is established in pluripotential theory
\cite{BedfordTaylor1986,BLMR2010}. Part~\ref{part:local} identifies
$\PW(A_x)$ as the microscopic polynomial space near $x$ and proves
\[
 n^{-d}K_n^\mu(x+u/n,x+v/n)
 \longrightarrow
 \frac{1}{w(x)(2\pi)^d}\int_{A_x}e^{i(u-v)\cdot\xi}\,d\xi
\]
uniformly for centers in compact subsets of $U$ and bounded real
$u,v$ (\cref{p1:thm:main}). Products of Chebyshev ridge polynomials
recover smooth interior spectral data with the full limiting norm
and spatial tail control (\cref{p1:thm:recovery}). This recovery
identifies the kernel and later permits the transfer of lower
sampling bounds.

On the diagonal, the limit gives
\[
 K_n^\mu(x,x)\sim\frac{n^d q_E(x)}{w(x)},
 \qquad
 \frac{K_n^\mu(x,x)w(x)}{N_n}
 \longrightarrow d!\,q_E(x).
\]
Thus the Christoffel sampling density converges locally in the
interior to the equilibrium density $d!q_E(x)$, expressed as a
local spectral volume. The normalized off-diagonal kernel instead
converges to $|A_x|^{-1}\int_{A_x}e^{i(u-v)\cdot\xi}\,d\xi$,
which retains the full spectral shape and describes correlations
between nearby point evaluations. These correlations enter the
Gram matrices that govern stability at small redundancy.

\paragraph{From the density law to a local basis condition.}
Uniformly stable, asymptotically critical arrays have the equilibrium
measure as their unweighted limiting node distribution
(\cref{p2:thm:critical-equidistribution}), by the convergence criterion
of Berman, Boucksom, and Witt Nystr\"om
\cite{BermanBoucksomWittNystrom2011}. This agrees with the limiting
Christoffel sampling density, but a macroscopic distribution does
not control the microscopic Gram matrices. Part~\ref{part:sampling}
extracts the additional condition from the full local model:
for every $x\in U$, $L^2(A_x)$ admits a normalized exponential Riesz
basis with bounds $\alpha/\beta,\beta/\alpha$ independent of $x$
(\cref{p2:thm:near-critical-transfer}). To obtain this pointwise
conclusion, we first use the fact that the evaluation relations have
total dimension $|\Lambda_n|-N_n=o(N_n)$. An averaged projection
trace estimate makes their local defects vanish at almost every
center, and exact-norm recovery supplies completeness in the limit.
The resulting bases have common bounds. Hausdorff continuity of
$A_x$ and compactness of exponential Riesz bases with these bounds
then extend their existence to every $x\in U$
(\cref{p2:lem:spectral-basis-compactness}).

Consequently, if $L^2(A_{x_0})$ admits no exponential Riesz basis at
even a single point $x_0\in U$, then every fixed ratio of MZ bounds
requires a positive relative redundancy gap
(\cref{p2:cor:fixed-condition-gap}). This obstruction applies to all
positive weighted arrays, including those with the equilibrium node
distribution, and supplies a structural lower bound complementary
to general sampling upper bounds.

\paragraph{A local-to-global criterion.}
Part~\ref{part:local-global} assumes positive continuous density on
the whole interior. If normalized exponential Riesz bases with common
bounds $a,b$ exist on $A_x$ for centers in a dense subset, there are
positive weighted rules with $N_n+o(N_n)$ nodes and bounds
$a-o(1),b+o(1)$ (\cref{p3:thm:local-global}). Uniform Fourier
localization assembles finitely many local sampling systems at each
stage. A positive correction restores the omitted energy, including
directions concentrated near the boundary, with a node count
controlled by a concentration trace. At $M_n=N_n+o(N_n)$ the two
sampling criteria can be read side by side
(\cref{p3:cor:existence-equivalence,p3:cor:isometric-equivalence}):
\[
\begin{aligned}
 \text{uniformly stable sampling}
   &\ \Longleftrightarrow\ \text{exponential Riesz bases on } A_x,\\
 \text{asymptotically isometric sampling}
   &\ \Longleftrightarrow\ \text{exponential orthonormal bases on } A_x .
\end{aligned}
\]
Both lines concern existence with positive sampling weights and require
the indicated basis at every interior point. In the first, the normalized
exponentials have Riesz bounds uniform in $x$. The second criterion has
a geometric interpretation through Fuglede's spectral set conjecture
\cite{Fuglede1974}, proved for convex bodies in every dimension by
Lev and Matolcsi \cite{LevMatolcsi2022}.
For a convex body $K$, the space $L^2(K)$ admits an orthonormal basis of
exponentials if and only if $K$ tiles $\mathbb R^d$ by translations.
Thus it is the local frequency body $A_x$ whose
translates must fit together at critical cardinality. Allowing extra nodes
leads to the packing estimates of Part~\ref{part:packing}.

\paragraph{Sampling redundancy on the disk, the regular hexagon, and higher-dimensional balls.}
The disk and the regular hexagon, with area measure, have the same
polynomial dimension $N_n=(n+1)(n+2)/2$. On the hexagon there are
$N_n+o(N_n)$ samples with frame operators converging to the
identity. On the disk, every fixed stability level requires a
positive proportion of additional nodes
(\cref{p3:cor:contrasting-domains}). Every local spectrum of the
regular hexagon is a hexagon or a parallelogram that tiles by translations.
For a ball $B^d$, $d\ge2$, the spectrum at its center is itself a ball.
The independent nonexistence results of Wan \cite{Wan2026} and
Ortega-Cerd\`a \cite{Ortega2026} exclude an exponential Riesz basis there,
so the fixed-stability gap persists in every such dimension.
For smooth positively curved bodies
in dimension $d\ge2$, the spectrum at the symmetry center is covered
by Ortega-Cerd\`a's extension to bounded convex domains with $C^2$
boundary. On the ball, Antezana, Marzo, and Ortega-Cerd\`a
\cite[Theorem~2]{AMO2021} used the rescaled kernel at the origin and the
obstruction to orthogonal exponentials on the ball to rule out orthogonal
bases of polynomial reproducing kernels. The fixed-stability gap answers
their question about simultaneous sampling and interpolation for the
weights considered there.

\paragraph{Necessary cubature redundancy in dimension \texorpdfstring{$d\ge3$}{d >= 3}.}
A positive cubature rule exact through degree $2n$ gives $S_n=I$ on
$\Pi_n^d$ and hence inherits the tiling requirement at vanishing relative
redundancy. In dimension $d\ge3$, every centrally symmetric convex body
has at least one interior point whose local spectrum fails to tile.
Consequently every such rule requires at least $(1+\delta)N_n$ nodes
for all sufficiently large $n$, with $\delta=\delta(E,\mu)>0$
(\cref{piv:classification,piv:cubature-gap}).
On the three-dimensional cube this excludes positive rules supported on
the cube of degree $2n-1$ with $n^3/6+O(n^2)$ nodes, as asked by
Li, Sun and Xu \cite[p.~121]{LiSunXu2009}.
Outside the class of centrally symmetric convex bodies, Gaussian
cubature attains the dimension bound on a nonconvex region in three
variables \cite{BerensSchmidXu1995}.

In the plane, the centrally symmetric bodies for which every local
spectrum tiles are precisely the parallelograms and centrally symmetric
hexagons (\cref{piv:plane}). These bodies admit asymptotically isometric
sampling at node ratio one. Exact moment constraints can still cost
additional nodes: on the square, a strictly positive analytic multiplier
of the product Chebyshev measure forces $\liminf M_n/N_n\ge4/3$ for
exact positive cubature (\cref{p2:thm:square-redundancy}).

Parts~\ref{part:local}--\ref{part:local-global} establish the local
spectral model, its necessary sampling conditions and the converse
construction. Part~\ref{part:packing} bounds the spatial distribution of
asymptotically isometric arrays by Euclidean packing densities
(\cref{piv:main}), evaluates the bounds on balls and cubes, and classifies
the centrally symmetric convex bodies that attain node ratio one
(\cref{piv:classification}). Appendix~\ref{p2:app:recovery} gives higher-order
recovery for quantitative aliasing. Appendices~\ref{p2:sec:hexagon-phase}
and~\ref{app:hexagon-numerics} treat the explicit phase grid and the
numerical kernel profiles on the hexagon. Appendix~\ref{p2:app:auxiliary}
collects $L^p$ sampling and least-squares estimates.
Appendix~\ref{sph:app} develops the spherical harmonic analogue.

\clearpage
\part{Bulk Universality and Local Spectral Geometry}\label{part:local}
\section{Polynomials at Scale \texorpdfstring{$1/n$}{1/n} and the Local Paley--Wiener Limit}\label{p1:sec:local-introduction}

Bulk universality describes the microscopic behavior of reproducing kernels of polynomial spaces. In one variable its basic profile is the sine kernel, $\sin\pi(a-b)/(\pi(a-b))$. Riemann--Hilbert methods provide precise asymptotics for varying analytic exponential weights \cite{DKMVZ1999}, while Lubinsky's comparison method establishes bulk universality for broad classes of fixed regular measures under local density assumptions \cite{Lubinsky2009}. Christoffel asymptotics relate the diagonal kernel to the equilibrium density; see Nevai \cite{Nevai1986} and Totik \cite{Totik2000}. Levin and Lubinsky \cite{LevinLubinsky2015} formulate one-dimensional Christoffel and universality problems in Paley--Wiener spaces. For a recent treatment through tangent measures and chains of de Branges spaces, see Eichinger, Luki\'c, and Woracek \cite{ELW2026}.

In several variables, Xu obtained Christoffel asymptotics for the classical weights on the ball and simplex \cite{XuBallAsymptotics1996,XuSimplexAsymptotics1999}. Kro\'o and Lubinsky \cite{KrooLubinsky2013} used Xu's kernel formulas \cite{XuBall1999,XuSimplex1998} and Jacobi Mehler--Heine asymptotics \cite{Szego1975} to identify the corresponding bulk limits as normalized Bessel functions. Their comparison theorems \cite[Theorems 1.1 and 1.5]{KrooLubinsky2013} then transfer reference asymptotics to other regular measures under a degree-stability hypothesis. For general background on multivariate Christoffel--Darboux kernels, see Dunkl and Xu \cite[Section 3.6]{DunklXu2014}.

\subsection{The Baran dual unit ball and the local kernel}

To describe the kernel on a general centrally symmetric convex body, we use the dual unit ball of the Baran metric. A translation preserves total degree, so we place the symmetry center at the origin. Let $E=-E\subset\R^d$ be the resulting body and let
\[
E^*=\{y\in\R^d:t\cdot y\leq1\text{ for every }t\in E\}
\]
be its polar. At $x\in\operatorname{int}E$, the Baran metric and its dual unit ball are
\[
F_E(x;v)=\sup_{y\in\partial E^*}
\frac{|y\cdot v|}{\sqrt{1-(x\cdot y)^2}},
\qquad
\A_x=\operatorname{conv}\left\{
\frac{y}{\sqrt{1-(x\cdot y)^2}}:y\in\partial E^*
\right\}.
\]
In particular, the support function of $\A_x$ satisfies
\begin{equation}\label{p1:eq:co-ball}
h_{\A_x}(v)=F_E(x;v).
\end{equation}
These formulas follow from the extremal-function representation of Lundin \cite{Lundin1985} and
Baran \cite{Baran1988}, in the form recorded by Bos, Calvi, and Levenberg \cite[Proposition 1.2]{BosCalviLevenberg2001}. Related geometric characterizations appear in Bos, Levenberg, and Waldron \cite{BosLevenbergWaldron2008} and in Kro\'o's survey \cite{KrooSurvey2003}. We call $\A_x$ the \emph{local spectral body} at $x$.

At the symmetry center, $\A_0=E^*$. The connection between multivariate polynomial approximation and entire functions of exponential type is classical; see Ganzburg \cite{Ganzburg1982}. His later sharp-constant theorem \cite[Theorem 1.2]{Ganzburg2020} gives the center Christoffel constant explicitly. Compactness and recovery methods for such extremal limits are developed further in \cite{Ganzburg2023}. Here the same perspective leads to exact-norm recovery at arbitrary interior centers and a position-dependent limiting kernel.

The volume of the dual unit ball also enters the equilibrium density formula of Bedford and Taylor \cite{BedfordTaylor1986}, extended to general convex bodies by Burns, Levenberg, Ma'u, and R\'ev\'esz \cite{BLMR2010}. On the ball and simplex, Piazzon \cite{Piazzon2019} identifies polynomial eigenspaces of the Baran Laplacian. Berman and Ortega-Cerd\`a \cite{BermanOrtega2018} relate real polynomial sampling to Paley--Wiener sampling and pluripotential equilibrium measures. In Section 1.7 they ask how the microscopic two-point kernel behaves in higher dimensions beyond the known examples.

Away from the symmetry center, the allowed frequency region $\A_x$ varies with $x$. For an arbitrary centrally symmetric convex body, assume global Bernstein--Markov control and write $d\mu=w\,dt$ on an interior open set $U$, with $w$ positive and continuous there. We prove
\[
n^{-d}K_n^\mu(x+u/n,x+v/n)
\longrightarrow \frac{1}{(2\pi)^d w(x)}
\int_{\A_x}e^{i(u-v)\cdot\xi}\,d\xi,
\]
uniformly for $x$ in compact subsets of $U$ and $u,v$ in bounded sets. At a general $x$, the translated body $E-x$ has its symmetry center at $-x$. We therefore construct polynomials centered at $x$ while keeping the original body and the total-degree constraint fixed.

The proof identifies both the functions and their norms in the limiting space. Products of Chebyshev ridge polynomials recover every function with smooth Fourier data compactly supported in $\operatorname{int}\A_x$. Their spatial decay gives convergence of the rescaled norms. In the other direction, the local expansion of the Siciak extremal function places the Fourier support of every bounded-energy limit in $\A_x$. These two facts determine the limiting norm of each finite linear combination of point evaluations and hence the two-point kernel. Its diagonal gives the Christoffel asymptotic. When the local spectral body is an ellipsoid, its Fourier transform yields the Bessel profile familiar from the ball and simplex.

The regular hexagon illustrates how the local spectral shape can vary within the bulk; see \cref{p1:sec:hexagon} and \cref{p1:fig:hexagon-coballs}. With the domain and the total-degree filtration fixed, the local spectral body is a hexagon inside the incircle and a parallelogram in each of the six outer corner regions. Along the interior arcs of the incircle, one pair of opposite vertices ceases to be extreme. The spectral body and its normalized Fourier transform remain continuous across this change of combinatorial type. The limiting kernel has explicit forms in both regimes, reducing to a product of two sinc functions in the corner regions. The Christoffel asymptotic records the area of the local spectral body, while the full universality limit retains its shape through the normalized Fourier transform.

\subsection{Standing assumptions and the local limit}\label{p1:sec:setting}

We use $\|\cdot\|$ for the Euclidean norm and $|A|$ for the Lebesgue volume of a measurable set $A\subset\R^d$.

Throughout, $\E\subset\R^d$ is a centrally symmetric convex body:
\[
\E=-\E,\qquad 0\in\E^\circ,\qquad \operatorname{int}\E\ne\varnothing.
\]
Its polar is $\Epol=\{y:x\cdot y\leq1\text{ for all }x\in\E\}$. We write
\[
\Pi_n^d=\{p\in\C[x_1,\ldots,x_d]:\deg p\leq n\},\qquad
N_n=\dim\Pi_n^d=\binom{n+d}{d}.
\]

Let $\mu$ be a finite positive measure with $\supp\mu=\E$. We assume the following $L^2$ Bernstein--Markov condition: there are constants $B_n^{\mathrm{BM}}\geq1$ such that
\begin{equation}\label{p1:eq:BM}
\|p\|_{L^\infty(\E)}\leq B_n^{\mathrm{BM}}\|p\|_{L^2(\mu)},\qquad p\in\Pi_n^d,
\qquad \log B_n^{\mathrm{BM}}=o(n).
\end{equation}
We also assume that on an open set $U\subset\operatorname{int}\E$,
\begin{equation}\label{p1:eq:density}
d\mu(t)=w(t)\,dt,
\qquad w\in C(U),\quad w(t)>0\quad(t\in U).
\end{equation}
Condition \eqref{p1:eq:density} specifies the restriction of $\mu$ to $U$ and allows $w$ to vary independently of the symmetry of $\E$. On the rest of $\E$ we use only the global measure assumptions above.

For $x\in\operatorname{int}\E$, define
\begin{equation}\label{p1:eq:gxy}
g_x(y)=\frac{y}{\sqrt{1-(x\cdot y)^2}},\qquad y\in\Epol,
\end{equation}
and let $\A_x$ be the convex hull defined in the introduction, with support function \eqref{p1:eq:co-ball}. Since $x$ is an interior point, $|x\cdot y|<1$ uniformly for $y\in\Epol$, so $\A_x$ is a compact, full-dimensional, centrally symmetric convex body.

Write $q_E(x)=|A_x|/(2\pi)^d$. In Parts~\ref{part:local-global} and~\ref{part:packing} the density hypothesis \eqref{p1:eq:density} holds throughout $\operatorname{int}E$.

We use the Fourier transform convention, initially for $f\in L^1\cap L^2$ and then extended to $L^2$ by Plancherel,
\[
\widehat f(\xi)=\int_{\R^d}f(u)e^{-iu\cdot\xi}\,du.
\]
For complex arguments, choose a real-valued orthonormal basis $\{p_{j,n}\}_{j=1}^{N_n}$ of $\Pi_n^d$ in $L^2(\mu)$ and use the holomorphic continuation
\[
K_n^\mu(z,w)=\sum_{j=1}^{N_n}p_{j,n}(z)p_{j,n}(w),\qquad z,w\in\C^d.
\]
On real arguments this is the usual reproducing kernel. On complex arguments we use the displayed holomorphic continuation, so the Hermitian kernel is $K_n^\mu(z,\overline w)$. The Christoffel function is
\[
\lambda_n^\mu(x)=\frac{1}{K_n^\mu(x,x)}
=\inf\left\{\int_E|p|^2\,d\mu:p\in\Pi_n^d,\ p(x)=1\right\}.
\]
For a compact convex body $A\subset\R^d$, the Paley--Wiener space is
\[
\PW(A)=\{f\in L^2(\R^d):\supp\widehat f\subseteq A\},
\]
where each member is identified with its canonical entire extension. We equip $\PW(A)$ with the $L^2(\R^d)$ norm. When we compare polynomial norms at $x\in U$ with their local limit, the change of variables $t=x+u/n$ gives the limiting squared norm $w(x)\|f\|_2^2$ on $\PW(\A_x)$. This factor accounts for the density in the kernel formula below.

\begin{theorem}[Local Paley--Wiener limit]\label{p1:thm:main}
Assume \eqref{p1:eq:BM}--\eqref{p1:eq:density}. Let $D\Subset U$ be compact and let $R>0$. Then, uniformly for $x\in D$ and $u,v\in\R^d$ with $\|u\|,\|v\|\leq R$,
\begin{equation}\label{p1:eq:main-kernel}
\frac{1}{n^d}K_n^\mu\left(x+\frac{u}{n},x+\frac{v}{n}\right)
\longrightarrow
\frac{1}{(2\pi)^d w(x)}
\int_{\A_x}e^{i(u-v)\cdot\xi}\,d\xi.
\end{equation}
The convergence also holds locally uniformly for the complexifications of $u$ and $v$.
\end{theorem}

\begin{corollary}[Bulk universality]\label{p1:cor:universality}
Under the hypotheses of \cref{p1:thm:main}, uniformly for $x$ in compact subsets of $U$ and for $u,v$ in bounded subsets of $\R^d$,
\begin{equation}\label{p1:eq:normalized-universality}
\frac{K_n^\mu(x+u/n,x+v/n)}{K_n^\mu(x,x)}
\longrightarrow
\frac{1}{|\A_x|}\int_{\A_x}e^{i(u-v)\cdot\xi}\,d\xi.
\end{equation}
\end{corollary}

\begin{proof}
Divide \eqref{p1:eq:main-kernel} by its value at $u=v=0$. The denominator is bounded away from zero on compact subsets of $U$. Indeed, $\A_x$ contains $\Epol$, since every $y\in\partial\Epol$ lies on the segment from $0$ to $g_x(y)$, and $w$ is bounded above on each compact set.
\end{proof}

\begin{corollary}[Christoffel asymptotics]\label{p1:cor:christoffel}
Under the hypotheses of \cref{p1:thm:main}, uniformly for $x$ in compact subsets of $U$,
\begin{equation}\label{p1:eq:christoffel-limit}
N_n\lambda_n^\mu(x)
\longrightarrow
\frac{(2\pi)^d}{d!}\frac{w(x)}{|\A_x|}.
\end{equation}
\end{corollary}

\begin{proof}
At $u=v=0$, \eqref{p1:eq:main-kernel} gives $n^{-d}K_n^\mu(x,x)\to |\A_x|/((2\pi)^dw(x))$. Since $N_n/n^d\to1/d!$, taking the reciprocal of the diagonal limit and multiplying by $N_n/n^d$ gives \eqref{p1:eq:christoffel-limit}.
\end{proof}

\begin{remark}[Equilibrium density]
For the probability-normalized pluripotential equilibrium measure, the formula of Bedford and Taylor \cite{BedfordTaylor1986} and Burns et al.\ \cite[Theorem 1.1]{BLMR2010} gives
\begin{equation}\label{p1:eq:equilibrium-density}
\rho_E(x)=\frac{d!}{(2\pi)^d}|\A_x|.
\end{equation}
Thus \cref{p1:cor:christoffel} reads $N_n\lambda_n^\mu(x)\to w(x)/\rho_E(x)$. The two-point limit in \cref{p1:thm:main} retains the shape of $\A_x$ through its Fourier transform, whose value at the origin is precisely the volume in \eqref{p1:eq:equilibrium-density}.
\end{remark}

\section{Extremal Growth and Chebyshev Ridge Wave Packets}

The proof of \cref{p1:thm:main} uses the same convex body $A_x$ twice.
Its support function bounds the growth of rescaled polynomials in
imaginary directions. Its interior supplies the frequencies for the
recovery construction. Compactness and recovery together determine the
limiting norms of point evaluations.

\subsection{Extremal growth and Fourier support}

To identify the possible Fourier support of a local polynomial limit, we first compute its exponential growth in imaginary directions. Write $V_\E$ for the Siciak extremal function. Central symmetry of $\E$ gives the following representation, valid also for nonsmooth bodies. The density $w$ enters only later, in the norm estimates.

\begin{lemma}[Lundin--Baran representation and local blow-up]\label{p1:lem:extremal}
For a centrally symmetric convex body $\E$,
\begin{equation}\label{p1:eq:BCL}
V_\E(z)=\sup_{y\in\partial\Epol}V_{[-1,1]}(y\cdot z),\qquad z\in\C^d.
\end{equation}
For every compact $D\Subset\operatorname{int}\E$ and every compact $C\subset\C^d$,
\begin{equation}\label{p1:eq:extremal-blowup}
nV_\E\left(x+\frac{\zeta}{n}\right)
=h_{\A_x}(\operatorname{Im}\zeta)+O_{D,C}(n^{-1}),
\qquad x\in D,\ \zeta\in C.
\end{equation}
\end{lemma}

\begin{proof}
The identity \eqref{p1:eq:BCL} is the consequence of Lundin's and Baran's results stated in Proposition 1.2 of Bos--Calvi--Levenberg \cite{BosCalviLevenberg2001}, after normalizing each nonzero real linear functional $y\cdot x$ by $\max_{t\in\E}|y\cdot t|=1$. This normalization is exactly $y\in\partial\Epol$.

For $s$ in a compact subinterval of $(-1,1)$ and $h\in\C$ sufficiently small, the interval extremal function satisfies
\begin{equation}\label{p1:eq:interval-expansion}
V_{[-1,1]}(s+h)
=\frac{|\operatorname{Im}h|}{\sqrt{1-s^2}}+O(|h|^2),
\end{equation}
uniformly in $s$ and $h$. Indeed, locally $V_{[-1,1]}(z)=|\operatorname{Im}\arccos z|$, and Taylor expansion of the analytic branch of $\arccos$ gives \eqref{p1:eq:interval-expansion}.

If $x\in D$, compactness of $D$ and $\partial\Epol$ gives $|x\cdot y|\leq1-\delta_D$ for some $\delta_D>0$. Apply \eqref{p1:eq:interval-expansion} with $s=x\cdot y$ and $h=(y\cdot\zeta)/n$. Uniformly in $x\in D$, $y\in\partial\Epol$, and $\zeta\in C$,
\[
nV_{[-1,1]}\left(x\cdot y+\frac{y\cdot\zeta}{n}\right)
=\frac{|y\cdot\operatorname{Im}\zeta|}{\sqrt{1-(x\cdot y)^2}}+O(n^{-1}).
\]
Taking the supremum over $y$ and using \eqref{p1:eq:gxy} yields \eqref{p1:eq:extremal-blowup}.
\end{proof}

\begin{lemma}[Interior ridge directions generate the co-ball]\label{p1:lem:dense-ridge}
For every $x\in\operatorname{int}\E$,
\[
\A_x=\overline{\operatorname{conv}\{g_x(y):y\in\operatorname{int}\Epol\}}.
\]
Consequently, if $C\Subset\operatorname{int}\A_x$, there exist finitely many points $y_1,\ldots,y_r\in\operatorname{int}\Epol$, including $d$ linearly independent directions and their negatives, such that
\[
C\subset\operatorname{int}P\subset P\subset\A_x,
\qquad P=\operatorname{conv}\{0,g_x(y_1),\ldots,g_x(y_r)\}.
\]
\end{lemma}

\begin{proof}
Fix $y_0\in\partial\Epol$ and $0<t<1$. Then $ty_0\in\operatorname{int}\Epol$, and
\[
g_x(ty_0)=\frac{t\sqrt{1-(x\cdot y_0)^2}}{\sqrt{1-t^2(x\cdot y_0)^2}}\,g_x(y_0).
\]
The scalar factor belongs to $(0,1)$ and tends to $1$ as $t\uparrow1$. Since $0\in\A_x$ and $\A_x$ is convex, $g_x(ty_0)\in\A_x$ and converges to $g_x(y_0)$. This proves the closure identity.

For the finite assertion, choose $0<\rho<\rho'<1$ such that $C\subset\rho\A_x$. For each unit vector $v$, the closure identity gives an interior ridge direction $y$ with $v\cdot g_x(y)>\rho' h_{\A_x}(v)$, since $h_{\A_x}(v)>0$. With $y$ fixed, the same strict inequality holds for all unit vectors in a neighborhood of $v$. A finite subcover therefore provides a polytope $P$ of the stated form with $h_P\geq\rho' h_{\A_x}$ and $P\subset\A_x$. It follows that $C\subset\rho\A_x\subset\operatorname{int}(\rho'\A_x)\subset\operatorname{int}P$. Enlarging this finite set by $d$ independent interior directions and their negatives preserves both inclusions. We order the directions so that the first $d$ are independent.
\end{proof}

\subsection{Exact-norm recovery by Chebyshev ridge wave packets}\label{p1:sec:recovery}

Fix $x\in\operatorname{int}\E$. We construct a polynomial whose profile on the scale $1/n$ tends to a prescribed Fourier integral. The basic observation is that $T_k(y\cdot t)=\cos(k\arccos(y\cdot t))$ has, near $x$, a phase with gradient proportional to $g_x(y)$. Products of these polynomials combine directions, while the sum of their degrees determines which frequencies can occur.

\subsubsection{A smooth lifting lemma}

Let $y_1,\ldots,y_r$ be chosen as in \cref{p1:lem:dense-ridge}, and abbreviate $g_j=g_x(y_j)$. Define
\[
B:\R^r\to\R^d,\qquad B\lambda=\sum_{j=1}^r\lambda_jg_j,
\]
and the open frequency-allocation simplex
\[
\Delta_r=\{\lambda\in\R^r:\lambda_j>0,\ \lambda_1+\cdots+\lambda_r<1\}.
\]
The coordinate $\lambda_j$ represents the fraction of the total degree assigned to direction $j$. The image of the closed simplex is $P$. Since $B$ is surjective, the image of its interior is $\operatorname{int}P$. The next lemma lets us prescribe a smooth frequency density on this image while keeping the fractions in a compact subset of $\Delta_r$.

\begin{lemma}[Smooth amplitude lifting]\label{p1:lem:lift}
For every $b\in C_c^\infty(\operatorname{int}P)$ there exists $a\in C_c^\infty(\Delta_r)$ such that
\begin{equation}\label{p1:eq:pushforward}
B_*(a(\lambda)\,d\lambda)=b(\xi)\,d\xi.
\end{equation}
Equivalently, for every continuous $\varphi$,
\[
\int_{\Delta_r}a(\lambda)\varphi(B\lambda)\,d\lambda
=\int_{\R^d}b(\xi)\varphi(\xi)\,d\xi.
\]
\end{lemma}

\begin{proof}
Choose a linear right inverse $S:\R^d\to\R^r$ of $B$ and a linear isomorphism $R:\R^{r-d}\to\ker B$. The map
\[
T(\xi,\eta)=S\xi+R\eta
\]
is an invertible linear map from $\R^d\times\R^{r-d}$ to $\R^r$, and $BT(\xi,\eta)=\xi$. Write $d\lambda=J\,d\xi\,d\eta$, where $J=|\det T|$.

For every $\xi_0\in\supp b\Subset\operatorname{int}P$, the fiber $B^{-1}(\xi_0)$ meets $\Delta_r$. Hence there are neighborhoods $U_{\xi_0}$ of $\xi_0$ and $V_{\xi_0}$ in $\R^{r-d}$ such that $T(U_{\xi_0}\times V_{\xi_0})\Subset\Delta_r$. By compactness, finitely many of the base neighborhoods $U_\ell$ cover $\supp b$. Retain the corresponding sets $V_\ell$ in the fiber coordinates. Take a smooth partition of unity $\{\psi_\ell\}$ near $\supp b$, with $\supp\psi_\ell\Subset U_\ell$, and choose $\chi_\ell\in C_c^\infty(V_\ell)$ with $\int\chi_\ell=1$.

Define, in the $(\xi,\eta)$ coordinates,
\[
a(T(\xi,\eta))=J^{-1}b(\xi)\sum_\ell\psi_\ell(\xi)\chi_\ell(\eta),
\]
and extend by zero. The support is compactly contained in $\Delta_r$, so $a\in C_c^\infty(\Delta_r)$. Integrating first in $\eta$ gives
\[
\int a(\lambda)\varphi(B\lambda)\,d\lambda
=\int b(\xi)\varphi(\xi)\sum_\ell\psi_\ell(\xi)\,d\xi
=\int b(\xi)\varphi(\xi)\,d\xi.
\]
\end{proof}

\subsubsection{The wave-packet construction}

For $y_j\in\operatorname{int}\Epol$, define on the real body $\E$
\[
\theta_j(t)=\arccos(y_j\cdot t).
\]
There is $\varepsilon_j>0$ such that
\begin{equation}\label{p1:eq:theta-range}
\varepsilon_j\leq\theta_j(t)\leq\pi-\varepsilon_j,
\qquad t\in\E.
\end{equation}
Moreover,
\begin{equation}\label{p1:eq:theta-gradient}
\nabla\theta_j(x)=-g_j,
\qquad
\theta_j\left(x+\frac{u}{n}\right)
=\theta_j(x)-\frac{g_j\cdot u}{n}+O(n^{-2}),
\end{equation}
uniformly for $u$ in bounded sets.

Given $b\in C_c^\infty(\operatorname{int}\A_x)$, choose the directions so that $\supp b\Subset\operatorname{int}P$ and lift $b$ to $a$ by \cref{p1:lem:lift}. Write $\Theta(t)=(\theta_1(t),\ldots,\theta_r(t))$. To center the packet at $x$, multiply each product by the compensating phase $e^{ik\cdot\Theta(x)}$ and set
\begin{equation}\label{p1:eq:wavepacket}
P_n(t)=2^rn^{-r}\sum_{\substack{k\in\Z^r\\a(k/n)\ne0}}a(k/n)e^{ik\cdot\Theta(x)}
\prod_{j=1}^rT_{k_j}(y_j\cdot t).
\end{equation}
The compact support of $a$ in $\Delta_r$ ensures that each term has positive integer indices $k_j$ and total degree at most $\sum_jk_j\leq(1-\eta)n$ for some fixed $\eta>0$. In particular,
\begin{equation}\label{p1:eq:degree}
P_n\in\Pi_n^d.
\end{equation}

\begin{lemma}[Poisson estimate]\label{p1:lem:poisson}
For $a\in C_c^\infty(\R^r)$, define
\[
S_n(t)=n^{-r}\sum_{k\in\Z^r}a(k/n)e^{ik\cdot t}.
\]
For every $m>0$ there is $C_m$ such that
\begin{equation}\label{p1:eq:poisson-est}
|S_n(t)|\leq C_m\bigl(1+n\dist(t,2\pi\Z^r)\bigr)^{-m},
\qquad t\in\R^r.
\end{equation}
\end{lemma}

\begin{proof}
Poisson summation gives
\[
S_n(t)=\sum_{\ell\in\Z^r}\widehat a\bigl(n(2\pi\ell-t)\bigr),
\]
with $\widehat a(\eta)=\int a(\lambda)e^{-i\lambda\cdot\eta}\,d\lambda$. By periodicity, take $t\in[-\pi,\pi]^r$. The term $\ell=0$ satisfies \eqref{p1:eq:poisson-est} by rapid decay of $\widehat a$. For $\ell\ne0$, $\|2\pi\ell-t\|\geq c_r\|\ell\|$, so summing a decay bound of order $L>\max(m,r)$ gives $O(n^{-L})$ uniformly on the cube. As the cube is bounded, this remainder also satisfies \eqref{p1:eq:poisson-est}.
\end{proof}

\begin{lemma}[Global localization of the wave packet]\label{p1:lem:localization}
For every $m>0$ there is $C_m$, independent of $n$ and $t\in\E$, such that
\begin{equation}\label{p1:eq:global-localization}
|P_n(t)|\leq C_m(1+n\|t-x\|)^{-m}+C_Mn^{-m}.
\end{equation}
\end{lemma}

\begin{proof}
Expanding each factor by $2T_k(s)=e^{ik\arccos s}+e^{-ik\arccos s}$ writes \eqref{p1:eq:wavepacket} as a sum of $2^r$ terms $S_n(t^\sigma(t))$, where
\[
t_j^\sigma(t)=\theta_j(x)+\sigma_j\theta_j(t),\qquad \sigma\in\{-1,1\}^r.
\]
For $\sigma=(-1,\ldots,-1)$, all phases vanish at $x$. Their coordinates are $\theta_j(x)-\theta_j(t)$. By \eqref{p1:eq:theta-range}, these differences lie strictly between $-\pi$ and $\pi$, so their nearest lattice point in $2\pi\Z$ is zero. For the first $d$ directions, which are linearly independent,
\[
|\theta_j(t)-\theta_j(x)|\geq |y_j\cdot(t-x)|.
\]
Indeed, $|\frac{d}{ds}\arccos s|\geq1$ on $(-1,1)$. If $Y$ is the matrix with rows $y_1^T,\ldots,y_d^T$, then $\|Y(t-x)\|\geq c\|t-x\|$. Hence
\[
\dist(t^{-}(t),2\pi\Z^r)\geq c\|t-x\|.
\]
Lemma~\ref{p1:lem:poisson} bounds the main branch by $C_m(1+n\|t-x\|)^{-m}$.

For every other branch, at least one coordinate has $\sigma_j=1$, so
\[
t_j^\sigma(t)=\theta_j(t)+\theta_j(x)\in[2\varepsilon_j,2\pi-2\varepsilon_j].
\]
That coordinate stays a fixed positive distance from $2\pi\Z$. Lemma~\ref{p1:lem:poisson} therefore bounds every non-main branch by $C_Mn^{-m}$. Summing the finitely many branches proves \eqref{p1:eq:global-localization}.
\end{proof}

\begin{theorem}[Ridge--Chebyshev recovery]\label{p1:thm:recovery}
Let $x\in\operatorname{int}\E$, let $b\in C_c^\infty(\operatorname{int}\A_x)$ and define
\[
f_b(u)=\int_{\A_x}b(\xi)e^{iu\cdot\xi}\,d\xi.
\]
There exist $P_n\in\Pi_n^d$ such that
\begin{equation}\label{p1:eq:recovery-local}
P_n\left(x+\frac{u}{n}\right)\longrightarrow f_b(u)
\end{equation}
locally uniformly for $u\in\R^d$. If in addition $x\in U$, the same polynomials satisfy
\begin{equation}\label{p1:eq:recovery-norm}
n^d\int_\E|P_n(t)|^2\,d\mu(t)
\longrightarrow
(2\pi)^dw(x)\int_{\A_x}|b(\xi)|^2\,d\xi.
\end{equation}
\end{theorem}

\begin{proof}
Take the polynomials in \eqref{p1:eq:wavepacket}. By \eqref{p1:eq:degree} they have degree at most $n$. The main term in the expansion above is
\[
S_n\bigl(\Theta(x)-\Theta(t)\bigr),
\]
where $S_n$ is as in \cref{p1:lem:poisson}. By \eqref{p1:eq:theta-gradient}, for $t=x+u/n$,
\[
\Theta(x)-\Theta\left(x+\frac{u}{n}\right)
=\frac1n(g_1\cdot u,\ldots,g_r\cdot u)+O(n^{-2}).
\]
Because $\|k\|=O(n)$ on the support of $a(k/n)$, the Taylor remainder changes the exponent by $O(n^{-1})$ uniformly for bounded $u$. The sum on the grid $n^{-1}\Z^r$ is therefore a Riemann sum converging locally uniformly to
\[
\int_{\Delta_r}a(\lambda)e^{i(B\lambda)\cdot u}\,d\lambda
=\int_{\A_x}b(\xi)e^{iu\cdot\xi}\,d\xi,
\]
where the last equality is \eqref{p1:eq:pushforward}. The other branches tend to zero locally uniformly by the proof of \cref{p1:lem:localization}. This proves \eqref{p1:eq:recovery-local}.

Now let $x\in U$ and choose a ball $U_x$ about $x$ with $\overline{U_x}\subset U$. Outside the ball, \eqref{p1:eq:global-localization} bounds $P_n$ by $C_Mn^{-M}$. Since $\mu$ is finite, this part of the rescaled energy is $O(n^{d-2M})$. Inside the ball, the change of variables $u=n(t-x)$ gives
\[
n^d\int_{U_x}|P_n(t)|^2w(t)dt
=\int_{n(U_x-x)}\left|P_n\left(x+\frac{u}{n}\right)\right|^2w\left(x+\frac{u}{n}\right)du.
\]
On this rescaled ball, $1+\|u\|\leq Cn$, so the residual $n^{-M}$ in \eqref{p1:eq:global-localization} is bounded by $C_M(1+\|u\|)^{-M}$. Extending the integrand by zero, we obtain the integrable majorant $C_M(1+\|u\|)^{-2M}$ for $2M>d$. By \eqref{p1:eq:recovery-local} and continuity of $w$, dominated convergence gives the limiting energy
\[
w(x)\int_{\R^d}|f_b(u)|^2du.
\]
Finally, Plancherel's theorem gives
\[
\int_{\R^d}|f_b(u)|^2du=(2\pi)^d\int_{\A_x}|b(\xi)|^2d\xi,
\]
which proves \eqref{p1:eq:recovery-norm}.
\end{proof}

\begin{remark}[Density of the recovery class]
The functions $f_b$, $b\in C_c^\infty(\operatorname{int}\A_x)$, form a dense subspace of $\PW(\A_x)$ in $L^2(\R^d)$. Indeed, $\partial\A_x$ has Lebesgue measure zero and $C_c^\infty(\operatorname{int}\A_x)$ is dense in $L^2(\A_x)$.
\end{remark}

\subsection{Compactness of rescaled polynomial spaces}

Recovery gives polynomials for each smooth spectral amplitude supported in $\A_x$. To identify the limiting space, we must also locate the frequencies of every possible limit. We do this by estimating the growth of the rescaled polynomials at complex arguments.

\begin{proposition}[Local compactness]\label{p1:prop:compactness}
Let $x\in U$, and let $p_n\in\Pi_n^d$ satisfy
\begin{equation}\label{p1:eq:bounded-energy}
\sup_n n^d\|p_n\|_{L^2(\mu)}^2<\infty.
\end{equation}
Then the rescaled functions $q_n(u)=p_n(x+u/n)$ have a subsequence converging locally uniformly on $\C^d$ to an entire function $f$. Every such limit satisfies
\begin{equation}\label{p1:eq:lsc}
w(x)\|f\|_{L^2(\R^d)}^2
\leq\liminf_{j\to\infty}n_j^d\|p_{n_j}\|_{L^2(\mu)}^2,
\end{equation}
where $q_{n_j}\to f$ is the chosen convergent subsequence, and
\begin{equation}\label{p1:eq:support}
\supp\widehat f\subseteq\A_x.
\end{equation}
The same assertion holds when $n$ is restricted at the outset to any increasing subsequence.
\end{proposition}

\begin{proof}
We first obtain compactness from a local evaluation estimate. We then identify the norm and the exponential type of each subsequential limit.

\smallskip
\noindent\emph{1. Local normality.} Choose nested cubes $Q_1\Subset Q_0\Subset U$ around $x$. Since $w$ is positive and continuous, $w\geq c>0$ on $Q_0$. By \eqref{p1:eq:bounded-energy},
\[
\int_{Q_0}|p_n(t)|^2dt\leq Cn^{-d}.
\]
To turn this integral bound into a pointwise bound, first work on the cube $I^d=[-1,1]^d$. For $p\in\Pi_n^d$, set $R(t)=p(t)\prod_j(1-t_j^2)$ and use the probability measure $d\nu=\pi^{-d}\prod_j(1-t_j^2)^{-1/2}dt$. Then
\[
\|R\|_{L^2(\nu)}^2\leq\pi^{-d}\|p\|_{L^2(I^d)}^2.
\]
The univariate orthonormal system is $1,\sqrt2T_1,\sqrt2T_2,\ldots$, whose diagonal kernel through degree $n+2$ is at most $2n+5$. Since $R$ has degree at most $n+2$ in each coordinate, expanding it in the tensor-product basis and applying Cauchy--Schwarz gives
\[
|p(s)|^2\leq\frac{(2n+5)^d}{\pi^d\prod_j(1-s_j^2)^2}
\|p\|_{L^2(I^d)}^2,\qquad s\in(-1,1)^d.
\]
Affine rescaling to the nested cubes, followed by the preceding integral bound, gives
\begin{equation}\label{p1:eq:nikolskii}
\|p_n\|_{L^\infty(Q_1)}\leq C.
\end{equation}
Applying the one-variable Bernstein--Walsh inequality successively in the coordinates of $Q_1$ bounds $|p_n(z)|$ by $C\exp(n\sum_j V_{I_j}(z_j))$, where $Q_1=\prod_j I_j$. The interval expansion \eqref{p1:eq:interval-expansion} therefore implies that
\[
\sup_{\zeta\in C}\left|p_n\left(x+\frac{\zeta}{n}\right)\right|\leq C_C
\]
for every compact $C\subset\C^d$. Thus $q_n$ is a normal family. Montel's theorem gives a locally uniformly convergent subsequence with entire limit $f$.

\smallskip
\noindent\emph{2. The limiting norm.} For every fixed $R>0$ and all sufficiently large $n$, the set $x+B_R/n$ lies in $U$. Changing variables gives
\[
n^d\int_{x+B_R/n}|p_n(t)|^2d\mu(t)
=\int_{B_R}|q_n(u)|^2w\left(x+\frac{u}{n}\right)du.
\]
Along the chosen convergent subsequence, local uniform convergence and continuity of $w$ imply
\[
w(x)\int_{B_R}|f(u)|^2du
\leq\liminf_{j\to\infty}n_j^d\|p_{n_j}\|_{L^2(\mu)}^2.
\]
Letting $R\to\infty$ proves \eqref{p1:eq:lsc}, and in particular $f|_{\R^d}\in L^2(\R^d)$.

\smallskip
\noindent\emph{3. Fourier support.} Bernstein--Walsh on $\E$ gives the directional growth estimate we need, provided the polynomials are uniformly bounded on $\E$. On $Q_1$ this follows from \eqref{p1:eq:nikolskii}, while on the entire body \eqref{p1:eq:BM} and \eqref{p1:eq:bounded-energy} give
\[
\|p_n\|_{L^\infty(\E)}\leq CB_n^{\mathrm{BM}} n^{-d/2}.
\]
We suppress the factor $B_n^{\mathrm{BM}}$ away from $x$ by using a multiplier of degree $o(n)$. Set
\[
m_n=\left\lceil\sqrt{n(1+\log B_n^{\mathrm{BM}})}\right\rceil.
\]
Then $m_n=o(n)$ and $\log B_n^{\mathrm{BM}}/m_n\to0$. Choose $\varepsilon>0$ so small that
\[
0<1-\varepsilon\|t-x\|^2\leq1,\qquad t\in\E,
\]
and define the polynomial localizer
\[
r_n(z)=\left(1-\varepsilon\sum_{j=1}^d(z_j-x_j)^2\right)^{m_n}.
\]
On the real body,
\[
0<r_n(t)\leq e^{-\varepsilon m_n\|t-x\|^2}.
\]
Choose $\rho>0$ with $B(x,\rho)\subset Q_1$. On this ball,
\eqref{p1:eq:nikolskii} bounds $p_n$ uniformly and $r_n\le1$.
On its complement in $\E$, the preceding estimates give
\[
 |p_n(t)r_n(t)|
 \le Cn^{-d/2}\exp\bigl(\log B_n^{\mathrm{BM}}
                         -\varepsilon m_n\rho^2\bigr)\le C
\]
for all sufficiently large $n$, since
$\log B_n^{\mathrm{BM}}/m_n\to0$. Combining the two regions gives
\begin{equation}\label{p1:eq:localized-sup}
\|p_nr_n\|_{L^\infty(\E)}\leq C.
\end{equation}
Also, $\deg(p_nr_n)\leq n+2m_n=n(1+o(1))$, and
\[
r_n\left(x+\frac{\zeta}{n}\right)
=\left(1-\frac{\varepsilon}{n^2}\sum_{j=1}^d\zeta_j^2\right)^{m_n}
\longrightarrow1
\]
locally uniformly on $\C^d$. Thus $p_nr_n$ has the same local limit $f$ as $p_n$.

The Bernstein--Walsh inequality now gives
\[
\left|p_nr_n\left(x+\frac{\zeta}{n}\right)\right|
\leq C\exp\left((n+2m_n)V_\E\left(x+\frac{\zeta}{n}\right)\right).
\]
By \cref{p1:lem:extremal},
\[
|f(\zeta)|\leq C\exp\bigl(h_{\A_x}(\operatorname{Im}\zeta)\bigr),
\qquad \zeta\in\C^d.
\]
Here the constant in \eqref{p1:eq:localized-sup} serves for every $\zeta$. Together with \eqref{p1:eq:lsc}, this global growth estimate is precisely the hypothesis of the $L^2$ Paley--Wiener theorem for compact convex spectral sets \cite{Hormander1983}, which yields \eqref{p1:eq:support}.
\end{proof}

\subsection{Convergence of the local reproducing kernels}

We compare the polynomial and Paley--Wiener kernels through the norms of point evaluations. Compactness gives an upper bound for these norms, and recovery gives the matching lower bound, as in the extremal argument of \cite{Ganzburg2023}. To obtain the kernel at two distinct points, we apply the argument to finite linear combinations of evaluations. For fixed $x\in U$, define
\[
\mathcal H_{n,x}=\{u\mapsto p(x+u/n):p\in\Pi_n^d\}
\]
with norm
\begin{equation}\label{p1:eq:scaled-norm}
\|q\|_{n,x}^2=n^d\int_\E|p(t)|^2d\mu(t),
\qquad q(u)=p(x+u/n).
\end{equation}
Its reproducing kernel is
\begin{equation}\label{p1:eq:scaled-kernel}
k_{n,x}(u,v)=n^{-d}K_n^\mu(x+u/n,x+v/n).
\end{equation}
The candidate limit is $\mathcal H_x=\PW(\A_x)$ with norm
\[
\|f\|_x^2=w(x)\int_{\R^d}|f(u)|^2du.
\]
Using the Fourier representation $f(u)=\int_{\A_x}b(\xi)e^{iu\cdot\xi}d\xi$, Plancherel gives
\[
\|f\|_x^2=(2\pi)^dw(x)\int_{\A_x}|b(\xi)|^2d\xi.
\]
Its reproducing kernel is
\begin{equation}\label{p1:eq:limit-kernel}
k_x(u,v)=\frac{1}{(2\pi)^dw(x)}\int_{\A_x}e^{i(u-v)\cdot\xi}d\xi.
\end{equation}

\begin{lemma}[Convergence of evaluation functionals]\label{p1:lem:evaluation}
Let $u_1,\ldots,u_m\in\R^d$ and $c_1,\ldots,c_m\in\C$. Define $L(f)=\sum_{j=1}^m c_jf(u_j)$. Then
\[
\|L\|_{\mathcal H_{n,x}^*}\longrightarrow\|L\|_{\mathcal H_x^*}.
\]
\end{lemma}

\begin{proof}
By the local evaluation bound in \cref{p1:prop:compactness}, the dual norms are uniformly bounded. Pass to a subsequence on which they approach their limsup, and for each degree choose $q_n$ in the unit ball of $\mathcal H_{n,x}$ with $|L(q_n)|=\|L\|_{\mathcal H_{n,x}^*}$. Such a maximizer exists in finite dimension. By \cref{p1:prop:compactness}, a further subsequence converges locally uniformly to $q\in\mathcal H_x$ with $\|q\|_x\leq1$. Since $L$ uses only finitely many values, $L(q_n)\to L(q)$ along this subsequence. Consequently,
\[
\limsup_{n\to\infty}\|L\|_{\mathcal H_{n,x}^*}\leq\|L\|_{\mathcal H_x^*}.
\]
For the lower bound, take any nonzero $f$ in the dense recovery class from \cref{p1:thm:recovery}. There are $q_n\in\mathcal H_{n,x}$ with $q_n\to f$ locally uniformly and $\|q_n\|_{n,x}\to\|f\|_x$. Thus
\[
\liminf_{n\to\infty}\|L\|_{\mathcal H_{n,x}^*}
\geq\frac{|L(f)|}{\|f\|_x}.
\]
Taking the supremum over the dense recovery class gives the reverse inequality.
\end{proof}

\begin{proof}[Proof of \cref{p1:thm:main}]
For any finite set of real points $u_1,\ldots,u_m$, the squared norm of a linear combination of evaluations is the corresponding quadratic form in the kernel matrix $(k_{n,x}(u_j,u_k))_{j,k}$. By \cref{p1:lem:evaluation}, these quadratic forms converge to those of $(k_x(u_j,u_k))_{j,k}$. Polarization therefore gives
\[
k_{n,x}(u,v)\longrightarrow k_x(u,v),\qquad u,v\in\R^d.
\]

To pass from real pointwise convergence to locally uniform convergence at complex arguments, we use the same local bounds once more. If $q\in\mathcal H_{n,x}$ has norm at most one, then its representing polynomial has $L^2(\mu)$ norm at most $n^{-d/2}$. The local Nikolskii estimate and the Bernstein--Walsh inequality on a fixed interior cube, as in Step 1 of \cref{p1:prop:compactness}, show that evaluation at any complex point in a fixed compact set has uniformly bounded norm on $\mathcal H_{n,x}$. Cauchy--Schwarz therefore gives local uniform bounds for the complexified kernels. Montel's theorem thus gives a normal family in $(u,v)\in\C^{2d}$. Every subsequential limit agrees with \eqref{p1:eq:limit-kernel} on $\R^d\times\R^d$. Applying the one-variable identity theorem in each coordinate identifies that limit on $\C^{2d}$ and proves local uniform convergence at the fixed center.

It remains to allow the center to vary. By compactness of $D$, it suffices to take arbitrary sequences $n_j\to\infty$ and $x_j\to x_\infty$ in $D$ and prove that $k_{n_j,x_j}$ and $k_{x_j}$ both converge to $k_{x_\infty}$ locally uniformly. First, $g_x(y)$ is jointly continuous on $D\times\partial\Epol$, so $\A_{x_j}\to\A_{x_\infty}$ in Hausdorff distance. Their support functions converge uniformly on bounded sets. Their indicators converge almost everywhere because their limiting boundary has measure zero, and all bodies lie in a common bounded set. Dominated convergence in \eqref{p1:eq:limit-kernel}, together with $w(x_j)\to w(x_\infty)$, therefore gives $k_{x_j}\to k_{x_\infty}$ locally uniformly, including at complex arguments.

For the polynomial kernels, we check compactness and recovery at the moving centers. Choose $\rho>0$ so that the cubes $Q_{0,j}=x_j+[-2\rho,2\rho]^d$ all lie in a common compact subset of $U$, and let $Q_{1,j}=x_j+[-\rho,\rho]^d$. The lower bound for $w$ on these cubes is independent of $j$. The evaluation estimate on $Q_{1,j}$, its complex extension and the bound for the polynomial multiplier centered at $x_j$ thus hold with uniform constants. On each fixed rescaled ball, $w(x_j+u/n_j)\to w(x_\infty)$ uniformly. Integrating there and then increasing its radius gives \eqref{p1:eq:lsc} with factor $w(x_\infty)$. The uniform expansion in \cref{p1:lem:extremal} and the convergence of the support functions give the growth bound $\lvert f(\zeta)\rvert\leq C\exp(h_{\A_{x_\infty}}(\operatorname{Im}\zeta))$ for every subsequential limit. Hence the compactness conclusion holds with limiting space $\mathcal H_{x_\infty}$.

For recovery, take $b\in C_c^\infty(\operatorname{int}\A_{x_\infty})$ and choose the directions $y_1,\ldots,y_r$ and amplitude $a$ at $x_\infty$. Keep these fixed and center \eqref{p1:eq:wavepacket} at $x_j$ by using the phases $\theta_\ell(x_j)$. In the local variables $x_j+u/n_j$, Taylor expansion of $\theta_\ell$ at $x_j$ gives $\theta_\ell(x_j)-\theta_\ell(x_j+u/n_j)=g_{x_j}(y_\ell)\cdot u/n_j+O(n_j^{-2})$ on bounded sets, with $g_{x_j}(y_\ell)\to g_{x_\infty}(y_\ell)$. Thus the principal Riemann sums converge locally uniformly to
\[
\int_{\Delta_r}a(\lambda)e^{i(B_{x_\infty}\lambda)\cdot u}\,d\lambda
=\int_{\A_{x_\infty}}b(\xi)e^{iu\cdot\xi}\,d\xi,
\]
where $B_{x_\infty}\lambda=\sum_\ell\lambda_\ell g_{x_\infty}(y_\ell)$. The localization estimate \eqref{p1:eq:global-localization} is uniform in $j$ because the directions are fixed in $\operatorname{int}\Epol$ and the first $d$ are independent. The change of variables in the proof of \cref{p1:thm:recovery} now gives convergence of the rescaled norms to the norm in $\mathcal H_{x_\infty}$.

We may therefore apply the evaluation-functional argument along $x_j$ and then polarize. This identifies the real pointwise limit of $k_{n_j,x_j}$ as $k_{x_\infty}$. The uniform complex evaluation bounds on the cubes give local uniform convergence by the same normal-family argument as at a fixed center. Since the conclusion holds along every convergent sequence of centers in $D$, the convergence in \eqref{p1:eq:main-kernel} is uniform on the stated compact sets.
\end{proof}

\section{Ellipsoids, Cross-Polytopes, and Classical Bulk Profiles}
\label{p1:sec:classical-profiles}

The main theorem gives the limiting normalized kernel in the form
\[
\mathcal U_x(z)=\frac{1}{|\A_x|}
\int_{\A_x}e^{iz\cdot\xi}\,d\xi,
\qquad z=u-v.
\]
At $z=0$ the unnormalized integral gives the volume of $\A_x$. At other displacements its normalization gives the characteristic function of the uniform probability measure on $\A_x$. We first evaluate this integral for ellipsoids, which include the local spectra on the interval and the ball. The simplex supplies a parallel example outside the centrally symmetric setting. On boxes the total-degree cutoff gives a cross-polytope, whose Fourier transform in dimension two factors along the two diagonals.

\subsection{Ellipsoids and normalized Bessel functions}

For $\nu>-1$, write
\[
J_\nu^*(r)=r^{-\nu}J_\nu(r),
\qquad
J_\nu^*(0)=\frac{1}{2^\nu\Gamma(\nu+1)},
\]
where the formula for positive real $r$ is extended by the even entire power series
\[
J_\nu^*(r)=2^{-\nu}\sum_{k=0}^{\infty}
\frac{(-r^2/4)^k}{k!\,\Gamma(k+\nu+1)}.
\]

\begin{proposition}[Fourier transform of an ellipsoid]\label{p1:prop:ellipsoid-bessel}
Let $G$ be a real symmetric positive-definite $d\times d$ matrix and let
\[
\mathcal E_G=\{\xi\in\R^d:\xi^{\mathsf T}G\xi\leq1\}.
\]
Then
\begin{equation}\label{p1:eq:ellipsoid-volume}
\vol(\mathcal E_G)=\frac{\omega_d}{\sqrt{\det G}},
\qquad
\omega_d=\frac{\pi^{d/2}}{\Gamma(d/2+1)},
\end{equation}
and, for $z\in\R^d$,
\begin{equation}\label{p1:eq:ellipsoid-fourier}
\frac{1}{\vol(\mathcal E_G)}
\int_{\mathcal E_G}e^{iz\cdot\xi}\,d\xi
=
\frac{J_{d/2}^*\!\left(\sqrt{z^{\mathsf T}G^{-1}z}\right)}
{J_{d/2}^*(0)}.
\end{equation}
The identity extends to complex $z$ by analytic continuation, with the right-hand side interpreted through the entire power series in $z^{\mathsf T}G^{-1}z$.
\end{proposition}

\begin{proof}
Substituting $\xi=G^{-1/2}\eta$ gives \eqref{p1:eq:ellipsoid-volume}. After division by the volume, the Jacobian cancels and the Fourier integral becomes
\[
\frac{1}{\omega_d}\int_{B_d}e^{i(G^{-1/2}z)\cdot\eta}\,d\eta.
\]
By rotational invariance this integral depends on $G^{-1/2}z$ through its length. In the present Fourier convention the unit-ball transform is
\[
\int_{B_d}e^{i\zeta\cdot\eta}\,d\eta
=(2\pi)^{d/2}\frac{J_{d/2}(\|\zeta\|)}{\|\zeta\|^{d/2}}.
\]
Dividing by $\omega_d$ and using
$J_{d/2}^*(0)=1/(2^{d/2}\Gamma(d/2+1))$
gives \eqref{p1:eq:ellipsoid-fourier}.
\end{proof}

The order $d/2$ of the Bessel function comes from the dimension of the ball. All dependence on the base point and the shape of the ellipsoid enters through the quadratic form $z^{\mathsf T}G^{-1}z$.

For the equilibrium weights on the ball and the simplex, one can also see this geometry in a differential operator. Piazzon \cite[Theorems 1 and 2]{Piazzon2019} identifies total-degree polynomial spaces with spectral subspaces of the corresponding Baran Laplacians. The degree-$n$ threshold is $n^2+O(n)$. H\"ormander's local spectral-function theorem \cite[Theorem 5.1]{Hormander1968} then gives the Fourier integral over the principal-symbol unit co-ball, with a lower-order remainder on compact interior sets. This recovers the Bessel profile in these Riemannian examples. The recovery construction of \cref{p1:sec:recovery} also applies to co-balls of other shapes.

\subsection{The interval and the sine kernel}

For $E=[-1,1]$, the Baran co-ball at $x\in(-1,1)$ is
\[
\A_x=\left[-\frac{1}{\sqrt{1-x^2}},
\frac{1}{\sqrt{1-x^2}}\right].
\]
Consequently,
\begin{equation}\label{p1:eq:interval-sine}
\mathcal U_x(z)
=\frac{1}{|\A_x|}\int_{\A_x}e^{iz\xi}\,d\xi
=\frac{\sin\!\left(z/\sqrt{1-x^2}\right)}
{z/\sqrt{1-x^2}}.
\end{equation}
For $d=1$, the Christoffel consequence of the main theorem gives
\[
w(x)K_n^\mu(x,x)\sim\frac{n}{\pi\sqrt{1-x^2}}.
\]
Thus the microscopic displacement $a/(w(x)K_n^\mu(x,x))$ is asymptotic to
$\pi a\sqrt{1-x^2}/n$. Substituting this scale into \eqref{p1:eq:interval-sine} yields
\[
\frac{\sin\pi(a-b)}{\pi(a-b)},
\]
which is the normalization in Lubinsky's one-dimensional bulk universality theorem \cite{Lubinsky2009}.

\subsection{The ball: the Bessel profile of Kro\'o--Lubinsky}\label{p1:subsec:ball}

For the unit ball $B^d$, the Baran co-ball is the ellipsoid
\[
\A_x^B=\{\xi:\xi^{\mathsf T}M_B(x)\xi\leq1\},
\qquad
M_B(x)=I-xx^{\mathsf T}.
\]
The matrix determinant lemma and the Sherman--Morrison formula give
\[
\det M_B(x)=1-\|x\|^2,
\qquad
M_B(x)^{-1}=I+\frac{xx^{\mathsf T}}{1-\|x\|^2}.
\]
Hence, for $z=u-v$,
\[
z^{\mathsf T}M_B(x)^{-1}z
=
\|u-v\|^2+\frac{(x\cdot(u-v))^2}{1-\|x\|^2}
=:G(x,u,v).
\]
Proposition~\ref{p1:prop:ellipsoid-bessel} therefore yields
\begin{equation}\label{p1:eq:ball-bessel}
\mathcal U_x^B(u-v)
=
\frac{J_{d/2}^*\!\left(\sqrt{G(x,u,v)}\right)}
{J_{d/2}^*(0)}.
\end{equation}
Kro\'o and Lubinsky obtain the same profile in Theorem~1.7 of \cite{KrooLubinsky2013} for regular measures on the ball with positive continuous density near the interior point under consideration. They first compute a reference limit from Xu's formula for Jacobi weights \cite{XuBall1999} and the Mehler--Heine formula, then transfer it to the other measures in their theorem. In \eqref{p1:eq:ball-bessel}, the quadratic form in their Bessel argument is the squared support function of the local ellipsoid.

\subsection{The simplex: a parallel ellipsoidal example}

The standard simplex provides a comparison outside the centrally symmetric setting. Write
\[
\Sigma^d=\{x_j>0,\ \sum_{j=1}^d x_j<1\},
\qquad
x_{d+1}=1-\sum_{j=1}^d x_j.
\]
Its Riemannian Baran metric and the corresponding polynomial eigenfunctions were studied by Piazzon \cite{Piazzon2019}. The inverse metric matrix, which defines the local co-ball, is
\[
M_T(x)=\operatorname{diag}(x_1,\ldots,x_d)-xx^{\mathsf T}.
\]
The Sherman--Morrison formula gives
\[
M_T(x)^{-1}
=
\operatorname{diag}(x_1^{-1},\ldots,x_d^{-1})
+\frac{\mathbf 1\mathbf 1^{\mathsf T}}{x_{d+1}}.
\]
Thus, with $z_j=u_j-v_j$ for $1\leq j\leq d$ and
$z_{d+1}=-\sum_{j=1}^d z_j$,
\[
z^{\mathsf T}M_T(x)^{-1}z
=
\sum_{j=1}^{d+1}\frac{z_j^2}{x_j}
=:H(x,u,v).
\]
The corresponding ellipsoid Fourier transform is therefore
\[
\frac{J_{d/2}^*\!\left(\sqrt{H(x,u,v)}\right)}
{J_{d/2}^*(0)},
\]
the profile in Theorem~1.8 of Kro\'o and Lubinsky \cite{KrooLubinsky2013}. Their reference-kernel representation uses Xu's Jacobi-weight Fourier-series formula on the simplex \cite{XuSimplex1998}.

\subsection{Local spectra of centrally symmetric polytopes}

The same integral gives different profiles when the local spectrum has flat faces. On a cube the co-ball is a weighted cross-polytope, and on a centrally symmetric polytope its shape can vary with position. The diagonal kernel records the volume of this body. Varying $z$ in the near-diagonal limit gives its full Fourier transform and thus determines the body itself.

For the polytopal examples here and in \cref{p1:sec:hexagon}, we use
\[
\operatorname{sinc}(t)=
\begin{cases}
\sin(t)/t,&t\ne0,\\
1,&t=0.
\end{cases}
\]

For the local spectral body $\A_x$, write
\[
\Phi_x(z)=\frac{1}{|\A_x|}
\int_{\A_x}e^{iz\cdot\xi}\,d\xi.
\]

For a centrally symmetric polytope with
$E^*=\operatorname{conv}\{\pm y_1,\ldots,\pm y_m\}$, the co-ball reduces to
\begin{equation}\label{p1:eq:finite-polar}
\A_x=\operatorname{conv}\{\pm g_x(y_1),\ldots,\pm g_x(y_m)\}.
\end{equation}
To see why only the vertices enter, fix $v$ and put $M=\max_j|v\cdot g_x(y_j)|$. The convex quadratic function
$(v\cdot y)^2+M^2(x\cdot y)^2$ is bounded by $M^2$ at each vertex and hence throughout $E^*$. Dividing by $1-(x\cdot y)^2>0$ gives $|v\cdot g_x(y)|\le M$. Since a vertex already attains $M$, the two sides of \eqref{p1:eq:finite-polar} have the same support function in every direction.

\subsection{Boxes, cubes, and parallelepipeds}\label{p1:subsec:cube}

Let
\[
B=\prod_{j=1}^d[\alpha_j,\beta_j],
\qquad
x\in\operatorname{int}B,
\]
and put
\[
\ell_j(x)=\sqrt{(x_j-\alpha_j)(\beta_j-x_j)}.
\]
The product property of the Siciak extremal function, together with the one-dimensional interval formula, gives
\[
F_B(x;v)=\max_{1\leq j\leq d}\frac{|v_j|}{\ell_j(x)}.
\]
Consequently the local spectral body is the weighted cross-polytope
\begin{equation}\label{p1:eq:box-spectral-body}
\A_x=
\left\{\xi\in\R^d:
\sum_{j=1}^d\ell_j(x)|\xi_j|\leq1
\right\},
\qquad
|\A_x|=\frac{2^d}{d!\prod_{j=1}^d\ell_j(x)}.
\end{equation}
To evaluate its Fourier transform, split the cross-polytope into orthants and set $t_j=\ell_j(x)|\xi_j|$. Summation over the signs of the coordinates gives
\begin{equation}\label{p1:eq:box-profile-general}
\Phi_x(z)
=d!\int_{\Delta_d}
\prod_{j=1}^d
\cos\left(\frac{z_jt_j}{\ell_j(x)}\right)\,dt,
\qquad
\Delta_d=\{t_j\geq0:\ t_1+\cdots+t_d\leq1\}.
\end{equation}
Thus the main theorem yields
\begin{equation}\label{p1:eq:box-christoffel}
N_n\lambda_n^\mu(x)
\longrightarrow
\pi^d w(x)\prod_{j=1}^d\ell_j(x).
\end{equation}

In dimension two, write $q_j=z_j/\ell_j(x)$. The diamond
$\{|\eta_1|+|\eta_2|\leq1\}$ becomes the square
$[-1,1]^2$ under
$s=\eta_1+\eta_2$, $t=\eta_1-\eta_2$, with
$d\eta_1d\eta_2=\frac12\,ds\,dt$. In these coordinates the Fourier integral separates, giving
\begin{equation}\label{p1:eq:rectangle-sinc}
\Phi_x(z)=
\operatorname{sinc}\left(\frac{q_1+q_2}{2}\right)
\operatorname{sinc}\left(\frac{q_1-q_2}{2}\right).
\end{equation}
In particular, for a rectangle
$[\alpha_1,\beta_1]\times[\alpha_2,\beta_2]$, the normalized bulk kernel has the explicit limit
\[
\frac{K_n^\mu(x+u/n,x+v/n)}{K_n^\mu(x,x)}
\longrightarrow
\operatorname{sinc}\left[
\frac12\left(
\frac{u_1-v_1}{\ell_1(x)}
+\frac{u_2-v_2}{\ell_2(x)}
\right)\right]
\operatorname{sinc}\left[
\frac12\left(
\frac{u_1-v_1}{\ell_1(x)}
-\frac{u_2-v_2}{\ell_2(x)}
\right)\right].
\]
The two diagonal directions in \eqref{p1:eq:rectangle-sinc} come from the diamond selected by the total-degree cutoff. A tensor-product cutoff selects a box and gives sinc factors along the coordinate axes. More generally, under an affine change of variables $x=c+Ls$, frequencies transform by $L^{-\mathsf T}$. The local spectral body on $c+L[-1,1]^d$ is therefore
\[
L^{-\mathsf T}
\left\{\eta:
\sum_{j=1}^d\sqrt{1-s_j^2}\,|\eta_j|\leq1
\right\},
\qquad s=L^{-1}(x-c),
\]
so the same calculation covers rotated rectangles and parallelepipeds.

For the unit cube, $\ell_j(x)=\sqrt{1-x_j^2}$. The product Chebyshev basis
gives another way to see the cross-polytope in \eqref{p1:eq:box-spectral-body}.
In angular coordinates each factor is a cosine. The total-degree constraint
restricts the nonnegative frequency indices to a simplex, and the signs
from expanding the cosines reflect it into all orthants. After local
coordinate scaling, the union of these simplices is $A_x$.

\section{Spectral Shape Transition on the Regular Hexagon}
\label{p1:sec:hexagon}

On a rectangle, changing the base point stretches the local spectrum
along its coordinate axes. On the regular hexagon, moving from the
center toward a corner changes the number of vertices of $A_x$.
The incircle marks the transition.

Let
\[
H=\operatorname{conv}
\{(\cos(k\pi/3),\sin(k\pi/3)):0\leq k<6\}.
\]
Its polar has vertices $\pm y_1,\pm y_2,\pm y_3$, where
\[
y_1=(1,1/\sqrt3),\qquad
y_2=(0,2/\sqrt3),\qquad
y_3=(-1,1/\sqrt3),
\qquad y_2=y_1+y_3.
\]
For an arbitrary interior base point $x=(x_1,x_2)$, set
\[
q_j=x\cdot y_j,\qquad
h_j=\sqrt{1-q_j^2},\qquad
b_j=\frac{y_j}{h_j},\qquad
S=h_1+h_2+h_3.
\]
The finite-polar formula \eqref{p1:eq:finite-polar} gives
\begin{equation}\label{p1:eq:hexagon-spectral-body}
\A_x=\operatorname{conv}\{\pm b_1,\pm b_2,\pm b_3\}.
\end{equation}
The vectors satisfy
\[
y_1=y_2-y_3,\qquad
y_2=y_1+y_3,\qquad
y_3=y_2-y_1.
\]
Each $b_j$ can thus be expressed in the basis $(b_k,b_\ell)$ of the other two directions. The hull of $\{\pm b_k,\pm b_\ell\}$ consists of vectors whose coefficients in that basis have absolute values summing to at most one. Consequently the pair $\pm b_j$ is redundant precisely when
\[
h_j\geq h_k+h_\ell,
\qquad \{j,k,\ell\}=\{1,2,3\}.
\]
For example, the identity
\[
b_2=\frac{h_1}{h_2}b_1+\frac{h_3}{h_2}b_3
\]
shows that the sum of the two coefficients is $(h_1+h_3)/h_2$.

To locate the transition in the physical hexagon, use the identity
\begin{equation}\label{p1:eq:hexagon-heron}
S(-h_1+h_2+h_3)(h_1-h_2+h_3)(h_1+h_2-h_3)
=3-4\|x\|^2,
\end{equation}
which follows by substituting
$q_1=x_1+x_2/\sqrt3$, $q_2=2x_2/\sqrt3$,
$q_3=-x_1+x_2/\sqrt3$ and using
$q_2=q_1+q_3$ and $q_1^2+q_2^2+q_3^2=2\|x\|^2$.
All $h_j$ are positive in the interior, so at most one triangle inequality can fail or become an equality. Reading off the sign in \eqref{p1:eq:hexagon-heron}, we obtain
\[
\begin{cases}
\|x\|<\sqrt3/2:&\text{all three strict triangle inequalities hold},\\
\|x\|=\sqrt3/2:&\text{exactly one triangle equality holds},\\
\|x\|>\sqrt3/2:&\text{exactly one of the three triangle inequalities fails}.
\end{cases}
\]
The circle $\|x\|=\sqrt3/2$ is the incircle of $H$, with six tangency points on $\partial H$. At interior points the central disk gives a hexagonal local spectrum, while the six corner regions and their open incircle arcs give parallelograms.

\subsection{Hexagonal spectra inside the incircle}
If $\|x\|<\sqrt3/2$, all six points in \eqref{p1:eq:hexagon-spectral-body} are vertices. Applying the shoelace formula in the cyclic order $b_1,b_2,b_3,-b_1,-b_2,-b_3$ gives
\begin{equation}\label{p1:eq:hexagon-area-central}
|\A_x|=
\frac{2S}{\sqrt3\,h_1h_2h_3}.
\end{equation}
For an explicit Fourier profile, define
\[
g_1=\frac{b_2-b_1}{2},\qquad
g_2=\frac{b_3-b_2}{2},\qquad
g_3=-\frac{b_1+b_3}{2}.
\]
The identities
\[
b_1=-g_1-g_2-g_3,\qquad
b_2=g_1-g_2-g_3,\qquad
b_3=g_1+g_2-g_3
\]
show that $\A_x\subseteq Z$, where
\[
Z=[-g_1,g_1]+[-g_2,g_2]+[-g_3,g_3].
\]
To prove equality and compute the Fourier transform, consider the three oriented areas
\[
\begin{aligned}
\det(g_1,g_2)&=
\frac{h_1+h_3-h_2}{2\sqrt3\,h_1h_2h_3},\\
\det(g_1,g_3)&=
\frac{h_2+h_3-h_1}{2\sqrt3\,h_1h_2h_3},\\
\det(g_2,g_3)&=
\frac{h_1+h_2-h_3}{2\sqrt3\,h_1h_2h_3}.
\end{aligned}
\]
They are positive by the strict triangle inequalities. The area of the zonotope is $4\sum_{i<j}\det(g_i,g_j)=2S/(\sqrt3\,h_1h_2h_3)$, which agrees with \eqref{p1:eq:hexagon-area-central}. The inclusion of full-dimensional compact convex sets with equal area gives $\A_x=Z$. Define the relative area weights
\[
\omega_{12}=\frac{h_1+h_3-h_2}{S},\qquad
\omega_{13}=\frac{h_2+h_3-h_1}{S},\qquad
\omega_{23}=\frac{h_1+h_2-h_3}{S}.
\]
These are the determinants above divided by their sum. To compute the transform, write the six translated parallelograms
\[
P_{ij}^{\pm}=\pm g_k+[-g_i,g_i]+[-g_j,g_j],
\qquad \{i,j,k\}=\{1,2,3\},
\]
as images of the six faces of the cube. Together they cover $Z$ exactly twice up to boundaries, since a generic fiber of the map
\[
[-1,1]^3\ni(t_1,t_2,t_3)\longmapsto t_1g_1+t_2g_2+t_3g_3
\]
is a segment with its two endpoints in the relative interiors of two faces. On each face the Fourier integral factors as
\[
\int_{P_{ij}^{\pm}}e^{iz\cdot\xi}\,d\xi
=
4\det(g_i,g_j)e^{\pm iz\cdot g_k}
\operatorname{sinc}(z\cdot g_i)\operatorname{sinc}(z\cdot g_j),
\]
Summing the two signs gives a cosine. Summing over the three pairs and dividing by twice the area of $Z$ yields
\begin{equation}\label{p1:eq:hexagon-fourier-central}
\begin{aligned}
\Phi_x(z)
={}&\omega_{12}\cos(z\cdot g_3)
\,\operatorname{sinc}(z\cdot g_1)
\,\operatorname{sinc}(z\cdot g_2)\\
&+\omega_{13}\cos(z\cdot g_2)
\,\operatorname{sinc}(z\cdot g_1)
\,\operatorname{sinc}(z\cdot g_3)\\
&+\omega_{23}\cos(z\cdot g_1)
\,\operatorname{sinc}(z\cdot g_2)
\,\operatorname{sinc}(z\cdot g_3).
\end{aligned}
\end{equation}
\subsection{Parallelogram spectra and the incircle interface}
Suppose $\|x\|>\sqrt3/2$, and let $j$ be the unique index for which
$h_j>h_k+h_\ell$, with $\{j,k,\ell\}=\{1,2,3\}$.
Then $\pm b_j$ lie inside the parallelogram generated by the other two pairs, and
\[
\A_x=\operatorname{conv}\{\pm b_k,\pm b_\ell\}.
\]
Put
\[
a=\frac{b_k+b_\ell}{2},\qquad
c=\frac{b_k-b_\ell}{2}.
\]
Since
\[
\A_x=[-a,a]+[-c,c],
\]
its area and normalized Fourier transform are
\begin{equation}\label{p1:eq:hexagon-area-parallelogram}
|\A_x|=\frac{4}{\sqrt3\,h_kh_\ell},
\qquad
\Phi_x(z)=
\operatorname{sinc}(z\cdot a)\operatorname{sinc}(z\cdot c).
\end{equation}
At an interior point of the incircle we have $h_j=h_k+h_\ell$. Here $b_j$ and $-b_j$ lie on opposite edges of the parallelogram. Both area formulas agree, and continuity of the vertices gives the same Fourier profile from either side.

Under the hypotheses of \cref{p1:thm:main}, the Christoffel limit is, in both regimes,
\begin{equation}\label{p1:eq:hexagon-christoffel}
N_n\lambda_n^\mu(x)\longrightarrow
\frac{2\pi^2w(x)}{|\A_x|},
\qquad N_n=\binom{n+2}{2}.
\end{equation}
The normalized local kernel converges to $\Phi_x(u-v)$, with $\Phi_x$ given by \eqref{p1:eq:hexagon-fourier-central} or \eqref{p1:eq:hexagon-area-parallelogram} according to the position of $x$.

For a point $x=(r,0)$ with $0\le r<1$, put $t=(1-r^2)^{-1/2}$. Then $h_1=h_3=1/t$ and $h_2=1$, so the triangle condition changes at $t=2$. For $t\ge2$ the spectral body is $[-t,t]\times[-t/\sqrt3,t/\sqrt3]$ and
\[
\Phi_{(r,0)}(z)=
\operatorname{sinc}(tz_1)\operatorname{sinc}(tz_2/\sqrt3),
\qquad t\ge2.
\]
Along the whole axis under consideration its area is
\[
|\A_{(r,0)}|=
\begin{cases}
2(t^2+2t)/\sqrt3,&1\leq t\leq2,\\
4t^2/\sqrt3,&t\geq2.
\end{cases}
\]

\begin{figure}[tbp]
\centering
\includegraphics[width=0.94\textwidth]{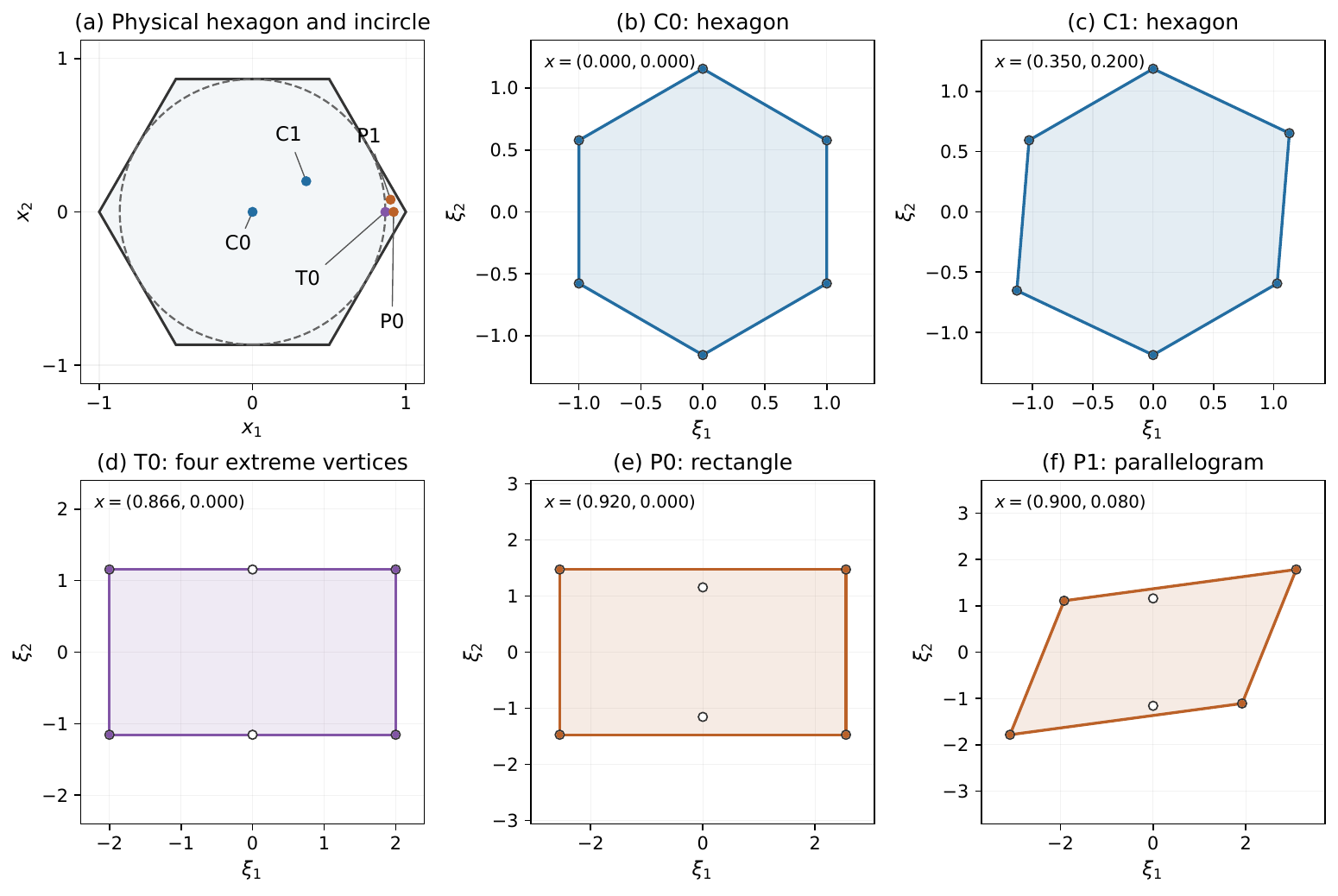}
\caption{The regular hexagon and local spectral bodies at five base points used in Appendix~\ref{app:hexagon-numerics}. C0 and C1 have hexagonal spectra. At T0, on the incircle, two of the six candidate vertices lie on edges and only four remain extreme. In the outer regions the spectrum is a rectangle at P0 and an oblique parallelogram at P1. Filled markers indicate extreme vertices; hollow markers indicate the remaining candidates. Spectral panels have independent axis ranges.}
\label{p1:fig:hexagon-coballs}
\end{figure}

Figure~\ref{p1:fig:hexagon-coballs} illustrates the central, transition, and corner regimes at both axial and off-axis points. Appendix~\ref{app:hexagon-numerics} compares the resulting Fourier profiles with the centered polynomial kernel at degree $64$.

Thus even on a fixed physical hexagon, the total-degree spaces have two different local Fourier profiles in the bulk. The three-term zonotope transform on the central disk becomes the product of two sinc functions on each corner region.

Central symmetry entered the proof of the local limit through
\eqref{p1:eq:BCL}, which identifies the extremal growth directions
with the directions used in the Chebyshev recovery. The Baran derivative
and equilibrium-volume formula also exist for general convex bodies
\cite{BLMR2010}. Extending the kernel limit requires a polynomial
realization of the full co-ball with the correct limiting energy, since
the ridge representation can fail \cite{BosCalviLevenberg2001}.

The recovered packets concentrate on the scale $1/n$ with controlled
tails. Testing a global sampling inequality on them gives the local
restrictions on an array that we study in Part~\ref{part:sampling}.

\clearpage
\part{Sampling, Interpolation, and Spectral Obstructions}\label{part:sampling}
\section{Macroscopic Equidistribution at Asymptotically Critical Cardinality}\label{p2:sec:sampling-frames}

Let $N_n=\dim\Pi_n^d$. Sampling asks whether point values determine a
polynomial with uniform norm control; interpolation asks whether arbitrary
nodal data can be realized with such control. At exactly $N_n$ nodes these
questions meet. We begin with their MZ formulation and the normalization
of point evaluations. We identify the equilibrium distribution forced by
uniformly stable, asymptotically critical sampling, then pass to tangent
sampling measures and derive geometric restrictions on node counts.
The Padua points give a critical-cardinality example on the square.

\subsection{Setting and the local kernel limit}\label{p2:sec:setting}

Under the standing assumptions of \cref{p1:sec:setting}, write
$K_n=K_n^\mu$ and $\|p\|_\mu=\|p\|_{L^2(\mu)}$. For a frequency body
$A$ the unweighted Paley--Wiener kernel is
\begin{equation}\label{p2:eq:local-body}
 k_A(u)=\frac1{(2\pi)^d}\int_A e^{iu\cdot\xi}\dd\xi .
\end{equation}
In particular, $k_{A_x}(0)=q_E(x)$.

By \cref{p1:thm:main},
\begin{equation}\label{p2:eq:kernel-limit}
 n^{-d}K_n(x+u/n,x+v/n)
 \longrightarrow \frac{k_{A_x}(u-v)}{w(x)},
\end{equation}
uniformly for $x$ in compact subsets of $U$ and for $u,v$ in bounded
real sets. In particular,
\begin{equation}\label{p2:eq:diagonal}
 K_n(x,x)\sim \frac{n^dq_E(x)}{w(x)}.
\end{equation}
Uniformity in the center will let us bound the mass of every rescaled
unit ball in a fixed interior neighborhood. By \eqref{p1:eq:equilibrium-density}, the
probability equilibrium density is $d!q_E$, whereas the phase density at
the $n^d$ normalization is $q_E$.

\subsection{Sampling, interpolation, and exact criticality}

At a node $\lambda\in E$, put
\[
 \kappa_{n,\lambda}(t)=\frac{K_n(t,\lambda)}{\sqrt{K_n(\lambda,\lambda)}}.
\]
For a finite set $\Lambda_n$, the Christoffel-normalized evaluation map is
\begin{equation}\label{p2:eq:sampling-operator}
 T_np=\left(\frac{p(\lambda)}{\sqrt{K_n(\lambda,\lambda)}}\right)_{\lambda\in\Lambda_n}.
\end{equation}
A uniform sampling array satisfies
\begin{equation}\label{p2:eq:normalized-frame}
 a\norm{p}_\mu^2\le \norm{T_np}_{\ell^2}^2\le b\norm{p}_\mu^2
 \qquad(p\in\Pi_n^d)
\end{equation}
with $0<a\le b<\infty$ independent of $n$. A uniform interpolation array is one for which the normalized kernels form a Riesz sequence, with bounds independent of $n$. This is equivalent to uniform solvability of normalized interpolation data together with a uniform upper Bessel bound for the normalized kernels.

We call an array \emph{uniformly completely interpolating} if it has both properties. Equivalently, it has exactly $N_n$ nodes and its Gram matrices satisfy
\begin{equation}\label{p2:eq:gram-bounds}
 aI\le G_n\le bI,\qquad
 G_n(\lambda,\lambda')=
 \frac{K_n(\lambda,\lambda')}{\sqrt{K_n(\lambda,\lambda)K_n(\lambda',\lambda')}}.
\end{equation}
Indeed, sampling requires at least $N_n$ vectors and interpolation allows at most $N_n$. At equality $T_n$ is square, so $T_n^*T_n$ and $G_n=T_nT_n^*$ have the same eigenvalues. When $|\Lambda_n|>N_n$, the extra dimensions instead appear in $\ker T_n^*$, a space that will enter the near-critical argument.

Finally, define the sup-norm norming constant
\begin{equation}\label{p2:eq:beta}
 \beta_n(\Lambda_n)=\sup_{0\ne p\in\Pi_n^d}
 \frac{\norm{p}_{L^\infty(E)}}{\max_{\lambda\in\Lambda_n}|p(\lambda)|},
\end{equation}
with value $+\infty$ when the nodes do not determine $\Pi_n^d$. This constant measures the worst loss in the supremum norm when one passes from the full domain to the nodes.

We work in $L^2(\mu)$ throughout the main sampling theory.
Appendix~\ref{p2:subsec:lp-mz} records the corresponding $L^p$
definitions and the relation between exact cardinality and the
Lebesgue constant.

\subsection{Christoffel weights and arbitrary weights at exact cardinality}\label{p2:sec:critical-weights}

Fix the continuous reference measure $\mu$. The variational formula for its Christoffel function is
\begin{equation}\label{p2:eq:christoffel-variational}
 \frac{1}{K_n(x,x)}
 =\min\left\{\norm{p}_\mu^2:p\in\Pi_n^d,\ p(x)=1\right\}.
\end{equation}
Indeed, the reproducing property and Cauchy--Schwarz give
$|p(x)|^2\le K_n(x,x)\norm{p}_\mu^2$, and equality in
\eqref{p2:eq:christoffel-variational} is attained by
$p=K_n(\cdot,x)/K_n(x,x)$. Thus $K_n(x,x)^{-1}$ is the minimum energy of a polynomial taking the value one at $x$. Normalizing evaluation by $\sqrt{K_n(x,x)}$ makes the functional have norm one. Its inner products with the other normalized evaluations form the Gram matrix in \eqref{p2:eq:gram-bounds}.

At exact cardinality, arbitrary positive sampling weights are comparable to the Christoffel weights. The following argument identifies where equality of the dimensions enters; Appendix~\ref{p2:subsec:lp-mz} gives the corresponding $L^p$ statement.

\begin{proposition}[Reduction of arbitrary weights at exact cardinality]\label{p2:prop:critical-weight-reduction}
Let $|\Lambda_n|=N_n=\dim\Pi_n^d$, and let
$\omega_{n,\lambda}>0$ be arbitrary weights. Suppose
\begin{equation}\label{p2:eq:critical-arbitrary-mz}
 a\norm{p}_\mu^2
 \le\sum_{\lambda\in\Lambda_n}\omega_{n,\lambda}|p(\lambda)|^2
 \le b\norm{p}_\mu^2,
 \qquad p\in\Pi_n^d,
\end{equation}
with $0<a\le b<\infty$ independent of $n$. Then, at every node,
\begin{equation}\label{p2:eq:critical-weight-comparison}
 \frac{a}{K_n(\lambda,\lambda)}
 \le\omega_{n,\lambda}
 \le\frac{b}{K_n(\lambda,\lambda)}.
\end{equation}
Consequently the Christoffel-normalized samples satisfy
\begin{equation}\label{p2:eq:critical-canonical-mz}
 \frac{a}{b}\norm{p}_\mu^2
 \le\sum_{\lambda\in\Lambda_n}
       \frac{|p(\lambda)|^2}{K_n(\lambda,\lambda)}
 \le\frac{b}{a}\norm{p}_\mu^2.
\end{equation}
In particular, a sequence of exact-cardinality node sets admits some positive weights satisfying a uniform two-sided $L^2$ MZ inequality if and only if it is uniformly completely interpolating in the normalization of \eqref{p2:eq:sampling-operator}.
\end{proposition}

\begin{proof}
Fix a node $\lambda\in\Lambda_n$. Substituting
$p=\kappa_{n,\lambda}$ in the upper inequality of
\eqref{p2:eq:critical-arbitrary-mz} and retaining the term at $\lambda$ gives
$\omega_{n,\lambda}K_n(\lambda,\lambda)\le b$.
The lower inequality makes the evaluation map injective. Since its domain and codomain have the same dimension, it is bijective, and there is a Lagrange polynomial $\ell_{n,\lambda}$ equal to one at $\lambda$ and zero at the other nodes. Therefore
\[
 \omega_{n,\lambda}
 \ge a\norm{\ell_{n,\lambda}}_\mu^2
 \ge \frac{a}{K_n(\lambda,\lambda)},
\]
where the last step follows from \eqref{p2:eq:christoffel-variational}.
It follows that the Christoffel-weighted sum lies between $1/b$ and $1/a$ times the original weighted sum. Combining these comparisons with \eqref{p2:eq:critical-arbitrary-mz} proves \eqref{p2:eq:critical-canonical-mz}. Since evaluation is square, Proposition~\ref{p2:prop:lp-critical} gives uniform complete interpolation in the normalized coordinates. Conversely, a uniformly completely interpolating array satisfies \eqref{p2:eq:critical-arbitrary-mz} with the Christoffel weights.
\end{proof}

\begin{corollary}[Exact norm identities at critical cardinality]\label{p2:cor:critical-isometry}
Suppose $|\Lambda_n|=N_n$ and positive weights satisfy
\[
 \norm{p}_\mu^2
 =\sum_{\lambda\in\Lambda_n}\omega_{n,\lambda}|p(\lambda)|^2
 \qquad(p\in\Pi_n^d).
\]
Then $\omega_{n,\lambda}=K_n(\lambda,\lambda)^{-1}$ at every node, and
$\{\kappa_{n,\lambda}:\lambda\in\Lambda_n\}$ is an orthonormal basis of $\Pi_n^d$. Conversely, such an orthonormal kernel basis gives this identity with the Christoffel weights.
\end{corollary}

\begin{proof}
Set $a=b=1$ in Proposition~\ref{p2:prop:critical-weight-reduction}. The normalized sampling operator is then an isometry between spaces of the same dimension. It is therefore unitary, and its Gram matrix is the identity. The converse is Parseval's identity.
\end{proof}

\begin{remark}[Scope of the reduction]\label{p2:rem:critical-weights-scope}
At critical cardinality the existence of uniform two-sided bounds on a sequence of node sets can therefore be tested with the Christoffel weights. The comparison in \eqref{p2:eq:critical-canonical-mz} allows the condition number of an individual rule to change under this normalization.

For a redundant array, testing the normalized kernel still gives $\omega_{n,\lambda}K_n(\lambda,\lambda)\le b$. The lower comparison used a Lagrange polynomial for each node. To see what changes in its absence, consider the space of constants with a probability measure. Two distinct nodes with weights $\varepsilon$ and $1-\varepsilon$ give an exact norm identity for every $0<\varepsilon<1$, even as the first weight tends to zero. This freedom to distribute sampling mass among redundant nodes will matter in the near-critical argument.
\end{remark}

A positive weighted MZ rule also gives stable least-squares
reconstruction. With weighted evaluation operator $T_n^\omega$,
the reconstruction operator is
$((T_n^\omega)^*T_n^\omega)^{-1}(T_n^\omega)^*$, and its norm is
at most $a^{-1/2}$. The precise noise and approximation estimates
are given in Appendix~\ref{p2:subsec:reconstruction}.

\subsection{Christoffel sampling and local phase volume}
\label{p2:subsec:leverage}

There is another way to select nodes from the continuous reference measure. Integrating the identity $K_n(t,t)=\sum_{k=1}^{N_n}|\phi_k(t)|^2$ over $E$ gives $\int_E K_n(t,t)\dd\mu(t)=N_n$. Hence
\[
 \dd\rho_n(t)=\frac{K_n(t,t)}{N_n}\dd\mu(t)
\]
is a probability measure. This is the Christoffel, or leverage-score, sampling measure used
in the coherence-optimal and weighted least-squares constructions of
Hampton--Doostan and Cohen--Migliorati
\cite{HamptonDoostan2015,CohenMigliorati2017}; see also
Adcock \cite{Adcock2025}. For $m$ independent draws
from $\rho_n$, the sample weights $N_n/(mK_n(t,t))$ make the empirical
Gram matrix unbiased. Their analysis gives high-probability stability
with a sampling budget of order $N_n\log N_n$, with constants depending on the prescribed success probability.

By \cref{p1:cor:christoffel} and $N_n\sim n^d/d!$, its density satisfies the compact-uniform limit
\begin{equation}\label{p2:eq:leverage-phase}
 \frac{\dd\rho_n}{\dd t}(x)
 =\frac{K_n(x,x)w(x)}{N_n}
 \longrightarrow d!\,q_E(x).
\end{equation}
Thus the limiting density of these random nodes is proportional to the
volume of the local spectral body. The same macroscopic distribution is
forced on every uniformly stable deterministic array whose cardinality
is asymptotically critical, as we now explain.

\subsection{Equilibrium distribution at asymptotically critical cardinality}
\label{p2:subsec:equidistribution}

An array is \emph{uniformly stable and asymptotically critical} if it
satisfies a weighted $L^2$ MZ inequality with bounds independent of $n$
and has $M_n=N_n+o(N_n)$ nodes. Its empirical distribution is determined
under the Bernstein--Markov condition alone. The following consequence
of \cite[Theorem~C]{BermanBoucksomWittNystrom2011} therefore applies
more generally than the local kernel limit: central symmetry and a
positive continuous density are not needed here. Write $\nu_E$ for the
probability-normalized pluripotential equilibrium measure of $E$.

\begin{theorem}[Equidistribution of uniformly stable, asymptotically critical arrays]
\label{p2:thm:critical-equidistribution}
Let $E\subset\R^d$ be a compact convex body with nonempty interior,
and let $\mu$ be a finite positive Bernstein--Markov measure on $E$.
Suppose that distinct nodes $x_{n,j}\in E$ and positive weights
$\omega_{n,j}$ satisfy
\begin{equation}\label{p2:eq:equidist-mz}
 a\norm{p}_\mu^2\le\sum_{j=1}^{M_n}\omega_{n,j}|p(x_{n,j})|^2
 \le b\norm{p}_\mu^2,\qquad p\in\Pi_n^d,
\end{equation}
with fixed $0<a\le b<\infty$. Set
$\sigma_n=\sum_j\omega_{n,j}\delta_{x_{n,j}}$, and let
$L_n=K_n^{\sigma_n}$ be the reproducing kernel for the discrete
inner product on $\Pi_n^d$. Then
\begin{equation}\label{p2:eq:discrete-bergman-limit}
 \beta_n:=\frac{1}{N_n}\sum_{j=1}^{M_n}
 \omega_{n,j}L_n(x_{n,j},x_{n,j})\delta_{x_{n,j}}
 \rightharpoonup\nu_E.
\end{equation}
If in addition $M_n=N_n+o(N_n)$, then
\begin{equation}\label{p2:eq:critical-counting-limit}
 \frac1{M_n}\sum_{j=1}^{M_n}\delta_{x_{n,j}}
 \rightharpoonup\nu_E,\qquad
 \frac1{N_n}\sum_{j=1}^{M_n}\delta_{x_{n,j}}
 \rightharpoonup\nu_E.
\end{equation}
All convergences are weak convergences of finite measures on $E$.
\end{theorem}

\begin{proof}
Scaling $\mu$ and all $\omega_{n,j}$ by the same constant leaves
$\beta_n$ unchanged, so assume $\mu(E)=1$. The Bernstein--Markov
condition supplies constants $C_n\ge1$ with $C_n^{1/n}\to1$ such that
$\norm{p}_E\le C_n\norm{p}_\mu$, where $\norm{\cdot}_E$ is the
supremum norm on $E$. Applying \eqref{p2:eq:equidist-mz} to $1$ gives
$a\le s_n:=\sigma_n(E)\le b$. For the probability measure
$\widehat\sigma_n=\sigma_n/s_n$, we consequently have
\[
 \norm{p}_{L^2(\widehat\sigma_n)}\le\norm{p}_E
 \le C_n\sqrt{b/a}\,\norm{p}_{L^2(\widehat\sigma_n)}.
\]
Let $B_\infty$ and $B_{2,n}$ denote the corresponding unit balls in the
complex polynomial space, of real dimension $2N_n$. Their inclusions
give
\[
 0\le\log\frac{\vol B_{2,n}}{\vol B_\infty}
 \le 2N_n\log(C_n\sqrt{b/a})=o(nN_n).
\]
This is the norm-volume criterion accompanying Theorem~C of
\cite{BermanBoucksomWittNystrom2011}. Its polynomial specialization
therefore yields
$N_n^{-1}K_n^{\widehat\sigma_n}(x,x)\dd\widehat\sigma_n(x)
\rightharpoonup\nu_E$.
Since $K_n^{\widehat\sigma_n}=s_nL_n$, this is precisely
\eqref{p2:eq:discrete-bergman-limit}.

To pass from kernel-weighted mass to node counts, choose an orthonormal
basis $(\psi_k)_{k=1}^{N_n}$ for $\Pi_n^d$ in $L^2(\sigma_n)$ and form
$V_{jk}=\sqrt{\omega_{n,j}}\psi_k(x_{n,j})$. The lower MZ bound makes
this inner product nondegenerate, and $V^*V=I$. Thus $VV^*$ is an
orthogonal projection of rank $N_n$. Its diagonal entries satisfy
\begin{equation}\label{p2:eq:discrete-leverage-defect}
 \ell_{n,j}:=\omega_{n,j}L_n(x_{n,j},x_{n,j})\in(0,1],\qquad
 \sum_j\ell_{n,j}=N_n,\qquad
 \sum_j(1-\ell_{n,j})=M_n-N_n.
\end{equation}
In particular,
\[
 \frac1{N_n}\sum_j\delta_{x_{n,j}}-\beta_n
 =\frac1{N_n}\sum_j(1-\ell_{n,j})\delta_{x_{n,j}}
\]
is a positive measure of mass $(M_n-N_n)/N_n$. This mass tends to zero
under asymptotic criticality, proving the second limit in
\eqref{p2:eq:critical-counting-limit}. The first follows from
$M_n/N_n\to1$.
\end{proof}

\begin{corollary}[Complete sampling and interpolation]
\label{p2:cor:complete-equidistribution}
Under the measure and domain assumptions of
\cref{p2:thm:critical-equidistribution}, every uniformly stable complete
sampling and interpolation array is equidistributed with respect to
$\nu_E$.
\end{corollary}

\begin{proof}
Such an array has $M_n=N_n$ and satisfies
\eqref{p2:eq:equidist-mz}, so the theorem applies. In this case $V$ is
square and unitary, and every $\ell_{n,j}$ equals one.
\end{proof}

For the centrally symmetric bodies of Part~\ref{part:local},
\eqref{p1:eq:equilibrium-density} identifies this common limit as
\begin{equation}\label{p2:eq:counting-phase-density}
 \dd\nu_E(x)=\rho_E(x)\dd x
 =\frac{d!}{(2\pi)^d}|A_x|\dd x.
\end{equation}
The reference measure affects the sampling norm and the weights, while
this necessary limiting node distribution depends only on $E$.
For general MZ arrays the kernel in
\eqref{p2:eq:discrete-leverage-defect} is the discrete kernel $L_n$.
For cubature exact on $\Pi_{2n}^d$, the discrete and continuous inner
products coincide on $\Pi_n^d$, so $L_n=K_n^\mu$; we return to this
case in \cref{p2:cor:cubature-equidistribution}.

The equilibrium law specifies how many nodes should lie in each fixed
macroscopic region at critical cardinality. Uniform stability also
requires that the samples resolve all local oscillations. The next
section retains the full spectral body $A_x$ in the rescaled MZ
inequality. Section~\ref{p2:sec:spectral-geometry} then shows that,
for some bodies, no uniformly stable asymptotically critical array
exists, including arrays with the equilibrium distribution.

\subsection{Relation to earlier work}

For real multivariate polynomials, Berman and Ortega-Cerd\`a
\cite{BermanOrtega2018} derived sampling bounds in terms of the
equilibrium measure. Antezana, Marzo, and Ortega-Cerd\`a \cite{AMO2021}
obtained interpolation density conditions and used a local kernel limit
on the ball to exclude orthogonal reproducing-kernel bases. The application of
\cref{p2:cor:smooth-physical} below treats the same family of ball measures.
Li \cite{Li2025} studies weighted $L^p$ sampling and interpolation
families on the ball, including microscopic Beurling--Landau densities
and degree-enlarged Fekete arrays.

The density estimates originate in Landau's theorem for Paley--Wiener
spaces \cite{Landau1967}. Related results for arrays include spherical
harmonics \cite{Marzo2007}, compact manifolds \cite{OrtegaPridhnani2012},
and localized reproducing-kernel spaces \cite{FuhrEtAl2017}. A closer
precedent for the transfer step is the work of Lev and Ortega-Cerd\`a
\cite[Theorem~6 and Lemma~16]{LevOrtega2016}, who use localized
recovery to pass from arrays of line-bundle sections to sampling and
interpolation in the Bargmann--Fock space. For polynomials the limiting
space varies with $x$, so this argument relates a single array to a
whole family of Fourier sampling problems.

On the constructive side, Gr\"ochenig \cite{Grochenig2020} relates
Marcinkiewicz--Zygmund inequalities to reconstruction and quadrature.
Cohen and Migliorati \cite{CohenMigliorati2017} use the kernel diagonal
to design weighted random sampling, and Dolbeault and Cohen
\cite{DolbeaultCohen2022} obtain near-best approximation in expected
$L^2$ error with a sampling budget proportional to dimension.
The finite-frame sparsification of Batson, Spielman, and Srivastava
\cite{BSS2014} gives deterministic positive norm discretization.
The lower bounds in Section~\ref{p2:sec:spectral-geometry} address the
complementary question of whether relative redundancy can vanish at
a fixed stability level. In the present polynomial setting, the bulk
kernel diagonal identifies the sampling density with local phase volume,
while the full local spectrum supplies a further structural constraint.

There is a related frame-theoretic picture near critical density.
Bownik and van Velthoven \cite{BownikVelthoven2025,
BownikVelthovenRedundancy2025,BownikVelthovenRiesz2026}
construct frames and Riesz sequences arbitrarily close to critical
density, in particular in general reproducing-kernel settings under
localization hypotheses. Enstad and van Velthoven
\cite[Theorem~1.1]{EnstadVelthoven2025} characterize separated exponential frames at critical density by the existence of Riesz-basis weak limits of
translates. Their use of averaged frame defects is closely related
to the finite-dimensional projection argument below. Here polynomial
localization connects global point counts to spectral bases at almost
every physical center; Part~\ref{part:local-global} supplies the
converse construction under common local bounds.

For prescribed nodes, the one-dimensional results of Coppersmith and
Rivlin \cite{CoppersmithRivlin1992}, Platte, Trefethen, and Kuijlaars
\cite{PlatteTrefethenKuijlaars2011}, and Adcock, Platte, and Shadrin
\cite{AdcockPlatteShadrin2019} describe how endpoint resolution affects
stability and approximation rates. Section~\ref{p2:sec:runge} examines a
multivariate source of instability in the interior and its consequences
for polynomial reconstruction of continuous functions.

\section{Tangent Sampling and Local Landau Density Bounds}\label{p2:sec:transfer}

The local kernel limit controls finite Gram matrices. Passing a lower
sampling bound to an entire Paley--Wiener space also requires control of
samples outside each fixed window. We first obtain this control by
localized recovery. The resulting transfer yields necessary densities
and, for prescribed lattices, tests involving the shape of the spectrum.

\subsection{Localized recovery and tangent sampling}

Consider the nodes in a window of radius $R/n$ about $x$. We wish to pass the polynomial sampling inequalities to the rescaled limit. After rescaling, compact convergence identifies the samples in each fixed window. To recover the whole norm in $\PW(A_x)$, we must also estimate the sampled energy outside the window as $R$ increases. The spatial decay of the wave packets in \cref{p1:sec:recovery} provides this estimate. Lev--Ortega-Cerd\`a \cite[Theorem~6 and Lemma~16]{LevOrtega2016} use a related argument to pass from sampling and interpolation for line-bundle sections to their Bargmann--Fock limits.

We record a stronger approximation estimate for the aliasing problem in \cref{p2:sec:runge}, where the rate of convergence of the samples will matter.

\begin{theorem}[Recovery to every prescribed algebraic order]\label{p2:thm:recovery}
Let $E=-E$ be a convex body and $x\in\interior E$. Suppose $\widehat f\in C_c^\infty(\interior A_x)$. For every integer $J\ge1$ there are $P_n^{[J]}\in\Pi_n^d$ such that, for every $L\ge0$,
\begin{equation}\label{p2:eq:all-orders}
 \sup_{u\in n(E-x)}(1+\norm{u})^L
 \left|P_n^{[J]}(x+u/n)-f(u)\right|
 \le C_{J,L}n^{-J}.
\end{equation}
The same sequence works for all $L$ at a fixed $J$. It also satisfies, for every $M$,
\begin{equation}\label{p2:eq:packet-decay}
 |P_n^{[J]}(t)|\le C_{J,M}
 \big((1+n\norm{t-x})^{-M}+n^{-M}\big),\qquad t\in E.
\end{equation}
If $\mu$ is a finite positive measure supported on $E$ and $\mu|_V=w(t)\dd t$ on a neighborhood $V$ of $x$, with $w$ positive and continuous, then
\begin{equation}\label{p2:eq:packet-norm}
 n^d\norm{P_n^{[J]}}_{L^2(\mu)}^2\longrightarrow w(x)\norm{f}_2^2.
\end{equation}
The construction and estimates are uniform on bounded subsets of $C_K^\infty$ for a fixed compact $K\Subset\interior A_x$.
\end{theorem}

The proof in \cref{p2:app:recovery} corrects the amplitude in the Chebyshev ridge construction order by order while keeping its spectral support inside $A_x$. For fixed orders, the approximation and decay estimates depend only on $E$, $x$, and the spectral data. The measure enters through the energy limit \eqref{p2:eq:packet-norm}.

\begin{theorem}[Weighted tangent-sampling principle]\label{p2:thm:tangent-sampling}
Assume the standing hypotheses of \cref{p1:sec:setting}, and fix $x\in U$. Let $\Lambda_n\subset E$ be finite, let $\omega_{n,\lambda}>0$, and suppose
\begin{equation}\label{p2:eq:weighted-frame}
 a\norm{p}_\mu^2\le
 \sum_{\lambda\in\Lambda_n}\omega_{n,\lambda}|p(\lambda)|^2
 \le b\norm{p}_\mu^2\qquad(p\in\Pi_n^d)
\end{equation}
with constants independent of $n$. The measures
\begin{equation}\label{p2:eq:tangent-measure}
 \sigma_{n,x}=\frac{n^d}{w(x)}
 \sum_{\lambda\in\Lambda_n}\omega_{n,\lambda}\delta_{n(\lambda-x)}
\end{equation}
are locally bounded uniformly in $n$. Every vague subsequential limit $\sigma$ satisfies
\begin{equation}\label{p2:eq:tangent-frame}
 a\norm{f}_2^2\le\int_{\R^d}|f|^2\dd\sigma\le b\norm{f}_2^2,
 \qquad f\in\PW(A_x).
\end{equation}
\end{theorem}

\begin{proof}
We first bound the rescaled sampling measures, then use localized
recovery to pass both inequalities to a dense class. The final step
identifies the continuous extension with actual point evaluation.

\smallskip\noindent\emph{1. Local and total sampling mass.}
Choose $\rho>0$ with $\overline{B(x,3\rho)}\subset U$. The limiting kernel in \eqref{p2:eq:kernel-limit} is continuous and strictly positive on the diagonal. Compact-uniform convergence therefore gives $r,c,C>0$ such that, for all sufficiently large $n$ and $z\in B(x,2\rho)$,
\[
 |K_n(t,z)|\ge cn^d\quad(\norm{t-z}\le r/n),\qquad
 K_n(z,z)\le Cn^d.
\]
Insert $p=K_n(\cdot,z)$ into the upper inequality of \eqref{p2:eq:weighted-frame}. Since $\norm{K_n(\cdot,z)}_\mu^2=K_n(z,z)$,
\[
 \sum_{\lambda\in B(z,r/n)}\omega_{n,\lambda}\le Cn^{-d}.
\]
Covering rescaled unit balls by a bounded number of $r$-balls gives
\begin{equation}\label{p2:eq:local-mass}
 \sigma_{n,x}(B(u,1))\le C
\end{equation}
whenever $\norm{u}\le n\rho$. In particular, the measures are uniformly bounded on each compact set and admit vague subsequential limits. Testing the same upper bound with $p=1$ gives a second estimate, valid on the whole rescaled domain:
\begin{equation}\label{p2:eq:total-mass}
 \sigma_{n,x}(\R^d)\le Cn^d.
\end{equation}

\smallskip\noindent\emph{2. Recovery and a uniform sampled tail.}
Fix a vaguely convergent subsequence with limit $\sigma$ and take $\widehat f\in C_c^\infty(\interior A_x)$. With the Fourier normalization of \cref{p2:sec:setting}, \cref{p1:thm:recovery} gives polynomials $P_n$ such that $F_n(u)=P_n(x+u/n)$ converges to $f(u)$ locally uniformly and
\[
 n^d\norm{P_n}_\mu^2/w(x)\longrightarrow\norm{f}_2^2.
\]
We show that their sampled energies also converge. By \eqref{p1:eq:global-localization},
$|F_n(u)|\le C_M((1+\norm{u})^{-M}+n^{-M})$ on $n(E-x)$. The squared contribution of the constant term is $O(n^{d-2M})$ by \eqref{p2:eq:total-mass}. For the decaying term, cover the portion $R<\norm{u}\le n\rho/2$ with unit balls of bounded overlap. Each has mass bounded by \eqref{p2:eq:local-mass}, and each shell $j\le\norm{u}<j+1$ meets $O((1+j)^{d-1})$ such balls. Summing over $j$ gives $O(R^{d-2M})$ for $R\ge1$ and $2M>d$. On $\norm{u}>n\rho/2$ the decaying term itself is $O(n^{-M})$, so the total mass bound again applies. Together these estimates give
\begin{equation}\label{p2:eq:sample-tail}
 \int_{\norm{u}>R}|F_n(u)|^2\dd\sigma_{n,x}(u)
 \le C_M\big(R^{d-2M}+n^{d-2M}\big).
\end{equation}
\smallskip\noindent\emph{3. Passing the entire sampling inequality.}
Choose $\eta_R\in C_c(\R^d)$ with $0\le\eta_R\le1$, equal to one
on $B_R$ and zero outside $B_{2R}$. For fixed $R$,
\[
\begin{aligned}
 &\int\eta_R|F_n|^2\,d\sigma_{n,x}
       -\int\eta_R|f|^2\,d\sigma\\
 &=\int\eta_R(|F_n|^2-|f|^2)\,d\sigma_{n,x}
   +\int\eta_R|f|^2\,d(\sigma_{n,x}-\sigma)
 \longrightarrow0.
\end{aligned}
\]
The first term tends to zero by local uniform convergence and bounded
mass on $B_{2R}$; the second tends to zero by vague convergence.
The limit measure inherits a uniform bound on the mass of unit balls
from \eqref{p2:eq:local-mass}. Since $f$ is Schwartz, its squared
integral against $\sigma$ outside $B_R$ also tends to zero.
We may therefore first let $n\to\infty$ and then $R\to\infty$ in
the cutoff identity, using \eqref{p2:eq:sample-tail}, to conclude that
\[
 \int |F_n|^2\dd\sigma_{n,x}\longrightarrow\int|f|^2\dd\sigma.
\]
Multiplying \eqref{p2:eq:weighted-frame} by $n^d/w(x)$ now gives \eqref{p2:eq:tangent-frame} for smooth interior spectral data.

\smallskip\noindent\emph{4. Extension to all bandlimited functions.}
Such data are dense in $\PW(A_x)$ because the boundary of $A_x$ has zero measure. The upper inequality extends the restriction map continuously into $L^2(\sigma)$. Moreover,
\[
 \norm{f-g}_\infty\le (2\pi)^{-d/2}|A_x|^{1/2}\norm{f-g}_2
\]
for $f,g\in\PW(A_x)$. If a sequence of smooth spectral data converges
to $f$ in $\PW(A_x)$, it consequently converges uniformly to $f$.
On each compact set, local finiteness of $\sigma$ identifies its
$L^2(\sigma)$ limit with $f$. Thus the continuous extension is
indeed restriction of the canonical continuous representative of $f$.
The lower inequality then passes by continuity, proving the theorem.
\end{proof}

\subsection{Tangent frames and tangent Riesz sequences}

For $x\in U$, consider the rescaled counting measures
\[
 \nu_{n,x}=\sum_{\lambda\in\Lambda_n}\delta_{n(\lambda-x)}.
\]
Different sampling nodes can approach the same point under rescaling. The limit of their counting measures records the number of nodes that merge there. Write such a limit as
\[
 \nu=\sum_{\gamma\in\Gamma}m_\gamma\delta_\gamma,
 \qquad m_\gamma\in\{1,2,\ldots\}.
\]
All sums over a multiset count these multiplicities.

\begin{proposition}[Separate sampling and interpolation transfer]\label{p2:prop:tangent-frame-riesz}
Assume the standing hypotheses of \cref{p1:sec:setting}, and fix $x\in U$.
\begin{enumerate}
\item If $\{\Lambda_n\}$ satisfies the normalized sampling inequality \eqref{p2:eq:normalized-frame} with bounds $a,b$, then its rescaled counting measures are locally bounded and every vague subsequential limit satisfies
\begin{equation}\label{p2:eq:multiset-tangent-frame}
 a\norm{f}_2^2\le
 \frac1{q_E(x)}\sum_{\gamma\in\Gamma}m_\gamma|f(\gamma)|^2
 \le b\norm{f}_2^2,
 \qquad f\in\PW(A_x).
\end{equation}
Equivalently, the exponentials $|A_x|^{-1/2}e^{-i\gamma\cdot\xi}$, repeated $m_\gamma$ times, form a frame for $L^2(A_x)$.
\item If $\{\Lambda_n\}$ is a uniform interpolation array with Riesz bounds $a,b$, its nodes are separated at scale $1/n$ on compact subsets of $U$. Every local limit is a simple separated set $\Gamma$ for which
\begin{equation}\label{p2:eq:interpolation-tangent-riesz}
 \big\{|A_x|^{-1/2}e^{-i\gamma\cdot\xi}:\gamma\in\Gamma\big\}
\end{equation}
is a Riesz sequence in $L^2(A_x)$ with bounds $a,b$.
\end{enumerate}
\end{proposition}

\begin{proof}
For sampling, the upper inequality tested on $\kappa_{n,z}$, together with the compact-uniform near-diagonal kernel limit, gives
\[
 \#\bigl(\Lambda_n\cap B(z,r/n)\bigr)\le C
\]
for some $r,C>0$, uniformly for $z$ in a smaller compact bulk neighborhood and large $n$. Indeed, every normalized kernel correlation in that ball is bounded away from zero. Covering bounded rescaled sets by finitely many such balls proves local boundedness and vague compactness. The limits are locally finite integer-valued atomic measures, with coalescing nodes represented by their multiplicities.

Apply \cref{p2:thm:tangent-sampling} with $\omega_{n,\lambda}=1/K_n(\lambda,\lambda)$. On every bounded rescaled set,
\[
 \frac{n^d}{w(x)K_n(\lambda,\lambda)}\longrightarrow\frac1{q_E(x)}
\]
uniformly. Thus the associated weighted tangent measure is $q_E(x)^{-1}\nu$, and \eqref{p2:eq:multiset-tangent-frame} follows. The stated frame interpretation follows from the Fourier convention in \cref{p2:sec:setting}.

For interpolation, the lower bound on every $2\times2$ principal Gram matrix gives
\[
 1-|G_n(\lambda,\lambda')|\ge a\qquad(\lambda\ne\lambda').
\]
If two nodes in a fixed compact subset of $U$ had rescaled distance tending to zero, \eqref{p2:eq:kernel-limit} would force $|G_n(\lambda,\lambda')|\to1$. The displayed inequality therefore gives separation and ensures that every limit point has multiplicity one. For a finite subset of $\Gamma$, take the nodes converging to its respective points. Their Gram matrices converge entrywise to the exponential Gram matrix, so the upper and lower inequalities hold for every finite coefficient vector. This gives the Riesz sequence in \eqref{p2:eq:interpolation-tangent-riesz}.
\end{proof}

\begin{theorem}[Complete interpolation forces local exponential Riesz bases]\label{p2:thm:riesz-transfer}
If $\{\Lambda_n\}$ is uniformly completely interpolating, then, for each $x\in U$, every local subsequential limit $\Gamma$ is separated and \eqref{p2:eq:interpolation-tangent-riesz} is a Riesz basis of $L^2(A_x)$ with the same bounds. Consequently, nonexistence of an exponential Riesz basis for one $A_x$, $x\in U$, excludes a uniformly completely interpolating polynomial array.
\end{theorem}

\begin{proof}
By \cref{p2:prop:tangent-frame-riesz}, interpolation gives the Riesz bounds for the exponentials on $\Gamma$. Sampling gives the lower frame bound for that same family. A function orthogonal to every exponential has zero frame coefficients, so this bound forces the function to vanish. The Riesz sequence thus spans $L^2(A_x)$.
\end{proof}

\subsection{A uniform Landau counting estimate}\label{p2:sec:landau-density}

Write $B_R=B(0,R)$ and let $A\subset\R^d$ be a compact set of positive volume. Put $q=|A|/(2\pi)^d$, and let $P_A$ be the orthogonal projection onto $\PW(A)$. Its kernel is $k_A(u-v)$, with
\[
 k_A(0)=q,\qquad \int_{\R^d}|k_A(u)|^2\dd u=q.
\]
The spatial concentration operator
\[
 C_R=P_A1_{B_R}P_A\big|_{\PW(A)}
\]
is positive, trace class, and has spectrum in $[0,1]$. Its trace and trace defect are
\begin{align}
 \operatorname{tr}C_R&=q|B_R|,\\
 V_A(R):=\operatorname{tr}(C_R-C_R^2)
 &=\int_{B_R}\int_{\R^d\setminus B_R}|k_A(u-v)|^2\dd v\dd u
 =o(|B_R|).\label{p2:eq:landau-trace-defect}
\end{align}
To obtain the last estimate, write the double integral as
\[
 \int_{\R^d}|k_A(t)|^2
 \bigl(|B_R|-|B_R\cap(B_R-t)|\bigr)\dd t.
\]
For each fixed $t$, the fraction of $B_R$ lost under translation by $t$ tends to zero. After division by $|B_R|$, dominated convergence applies since $|k_A|^2$ is integrable. For a fixed $\tau\in(0,1)$ and $0\le s\le1$, the difference $|\mathbf1_{[\tau,1]}(s)-s|$ is bounded by $\max(\tau^{-1},(1-\tau)^{-1})s(1-s)$. Summing over the eigenvalues shows that the dimension of the concentration space above threshold $\tau$ is $q|B_R|+o(|B_R|)$. Its leading term $(2\pi)^{-d}|B_R\times A|$ is the volume in position and frequency available within the window, in the Fourier normalization of \cref{p2:sec:setting}.

\begin{lemma}[Uniform finite-radius Landau estimates]\label{p2:lem:uniform-landau}
Fix $A$ and positive bounds $a,b$.
\begin{enumerate}
\item There are $r>0$ and $C>0$, depending only on $A,a,b$, such that every locally finite multiset $\Gamma$ satisfying
\[
 a\norm{f}_2^2\le q^{-1}\sum_{\gamma\in\Gamma}|f(\gamma)|^2
 \le b\norm{f}_2^2\qquad(f\in\PW(A))
\]
obeys, for every $y\in\R^d$ and $R>0$,
\begin{equation}\label{p2:eq:uniform-landau-lower}
 \#\bigl(\Gamma\cap B(y,R+r)\bigr)
 \ge q|B_R|-C V_A(R).
\end{equation}
The count includes multiplicities.
\item Define
\[
 \epsilon_A(r)=q^{-1}\int_{\norm{u}\ge r}|k_A(u)|^2\dd u
 \longrightarrow0.
\]
If the normalized kernels $k_A(\cdot-\gamma)/\sqrt q$, $\gamma\in\Gamma$, form a Riesz sequence with lower bound $a$, then, whenever $\epsilon_A(r)<a$,
\begin{equation}\label{p2:eq:uniform-landau-upper}
 \#\bigl(\Gamma\cap\overline{B(y,R)}\bigr)
 \le \frac{q|B_{R+r}|}{1-\epsilon_A(r)/a}
 \qquad(y\in\R^d,\ R>0).
\end{equation}
\end{enumerate}
\end{lemma}

\begin{proof}
For the first part, the upper sampling bound tested on $k_A(\cdot-z)/\sqrt q$ gives a uniform local counting bound. In fact, continuity at zero provides $s>0$ such that $|k_A(u)|\ge q/2$ for $\norm{u}\le s$, and then
\[
 \#\bigl(\Gamma\cap B(z,s)\bigr)\le4b.
\]
Covering a unit ball by balls of radius $s$ gives $\sup_z\#(\Gamma\cap B(z,1))\le M$, with $M$ depending only on $A,b$.

Choose a Schwartz function $\psi$ with $\widehat\psi=1$ on a neighborhood of $A$. Since $f=f*\psi$, Cauchy--Schwarz gives
\[
 |f(\gamma)|^2\le\norm{\psi}_1
 \int_{\R^d}|f(u)|^2|\psi(\gamma-u)|\dd u.
\]
The uniform local counting bound and Schwartz decay imply
\[
 \sup_u\sum_{\gamma\in\Gamma}|\psi(\gamma-u)|\le C_0,
 \qquad
 \sup_u\sum_{\norm{\gamma-u}\ge r}|\psi(\gamma-u)|
 \le\eta(r),\quad\eta(r)\to0,
\]
with constants uniform over the whole class of multisets. Splitting the integral into $B(y,R)$ and its complement therefore gives
\begin{equation}\label{p2:eq:uniform-discrete-tail}
 q^{-1}\sum_{\gamma\notin B(y,R+r)}|f(\gamma)|^2
 \le C_1\int_{\R^d\setminus B(y,R)}|f(u)|^2\dd u
       +\eta_1(r)\norm{f}_2^2,
 \qquad \eta_1(r)\to0.
\end{equation}
Choose $r$ with $\eta_1(r)<a/4$, and then $\tau\in(0,1)$ with $C_1(1-\tau)<a/4$. These choices depend only on $A,a,b$.

Let $H_{R,y}$ be the span of eigenvectors with eigenvalue at least $\tau$ for concentration on $B(y,R)$. Every $f\in H_{R,y}$ has at most $(1-\tau)\norm{f}_2^2$ of its energy outside this ball. Equation~\eqref{p2:eq:uniform-discrete-tail} therefore bounds its exterior sampled energy by $a\norm{f}_2^2/2$. The sample restriction to $\Gamma\cap B(y,R+r)$ is injective on $H_{R,y}$: otherwise the sampling lower bound would be contradicted. Hence the number of samples there is at least $\dim H_{R,y}$. If $\{\lambda_j\}$ are the concentration eigenvalues, then
\[
 \dim H_{R,y}\ge\sum_j\lambda_j
       -\frac1{1-\tau}\sum_j\lambda_j(1-\lambda_j)
 =q|B_R|-\frac{V_A(R)}{1-\tau}.
\]
Translation invariance removes $y$ from the trace terms, proving \eqref{p2:eq:uniform-landau-lower}.

For the second part, the Riesz lower bound first implies uniform separation:
its $2\times2$ Gram blocks give
$1-|k_A(\gamma-\gamma')|/q\ge a$ for distinct nodes, while $k_A$ is
continuous and $k_A(0)=q$. Thus bounded sets contain finitely many nodes.
Take the finite subfamily indexed by $\Gamma\cap\overline{B(y,R)}$, of size $m$, and let $V:\C^m\to\PW(A)$ synthesize its normalized kernels. Its orthogonal range projection is $Q=V(V^*V)^{-1}V^*$, so $Q\le a^{-1}VV^*$. Each kernel has at most $\epsilon_A(r)$ of its squared norm outside $B(y,R+r)$. Positivity and the trace inequality give
\[
 \operatorname{tr}\bigl[(I-C_{R+r,y})Q\bigr]
 \le a^{-1}\sum_{\gamma\in\Gamma\cap\overline{B(y,R)}}
 \int_{\R^d\setminus B(y,R+r)}\frac{|k_A(u-\gamma)|^2}{q}\dd u
 \le m\epsilon_A(r)/a,
\]
where $C_{R+r,y}$ is the translated concentration operator. On the other hand, $Q\le I$ gives
\[
 \operatorname{tr}(C_{R+r,y}Q)\le\operatorname{tr}C_{R+r,y}=q|B_{R+r}|.
\]
Adding these estimates and using $\operatorname{tr}Q=m$ proves \eqref{p2:eq:uniform-landau-upper}.
\end{proof}

\subsection{The polynomial array density theorem}

For a fixed $x\in U$, define the local array densities by the indicated order of limits:
\begin{align}
 \mathcal D_x^-(\Lambda)
 &=\liminf_{R\to\infty}\liminf_{n\to\infty}
   \frac{\#\bigl(\Lambda_n\cap B(x,R/n)\bigr)}{q_E(x)|B_R|},\\
 \mathcal D_x^+(\Lambda)
 &=\limsup_{R\to\infty}\limsup_{n\to\infty}
   \frac{\#\bigl(\Lambda_n\cap B(x,R/n)\bigr)}{q_E(x)|B_R|}.
 \label{p2:eq:local-array-densities}
\end{align}

\begin{theorem}[Local Landau necessary density]\label{p2:thm:local-landau}
Under \cref{p1:thm:main}, the following hold at every fixed $x\in U$.
\begin{enumerate}
\item Every normalized sampling array satisfies $\mathcal D_x^-(\Lambda)\ge1$.
\item Every uniform interpolation array satisfies $\mathcal D_x^+(\Lambda)\le1$.
\item A uniformly completely interpolating array has $\mathcal D_x^-(\Lambda)=\mathcal D_x^+(\Lambda)=1$.
\end{enumerate}
For every sampling tangent multiset, its lower uniform Beurling density is at least $q_E(x)$; for every interpolation tangent set, its upper uniform Beurling density is at most $q_E(x)$.
\end{theorem}

\begin{proof}
Fix $x$ and put $A=A_x$, $q=q_E(x)$. Proposition~\ref{p2:prop:tangent-frame-riesz} supplies locally bounded rescaled counting measures. For sampling and a fixed radius $R>r$, select a subsequence realizing the lower limit of the counts in $B_R$, and then a vaguely convergent subsequence with limit $\nu$. The open-set inequality for vague convergence and \eqref{p2:eq:uniform-landau-lower} give
\[
 \liminf_{n\to\infty}\#\bigl(\Lambda_n\cap B(x,R/n)\bigr)
 \ge\nu(B_R)
 \ge q|B_{R-r}|-C V_A(R-r).
\]
The constants $r,C$ depend only on $A,a,b$ and hence are common to the tangents selected at all radii $R$. After dividing by $q|B_R|$, use $|B_{R-r}|/|B_R|\to1$ and \eqref{p2:eq:landau-trace-defect} to obtain the lower density bound.

For interpolation, choose a subsequence realizing the upper limit of the counts and pass to a vague limit $\nu=\sum_{\gamma\in\Gamma}\delta_\gamma$. The compact-set inequality and \eqref{p2:eq:uniform-landau-upper} yield, for every fixed $r$ with $\epsilon_A(r)<a$,
\[
 \limsup_{n\to\infty}\#\bigl(\Lambda_n\cap B(x,R/n)\bigr)
 \le\nu(\overline{B_R})
 \le\frac{q|B_{R+r}|}{1-\epsilon_A(r)/a}.
\]
After division by $q|B_R|$, first let $R\to\infty$ and then $r\to\infty$. This proves the second assertion. The third follows from the first two and $\mathcal D_x^-\le\mathcal D_x^+$.

For a tangent counting measure, define uniform Beurling densities using the infimum or supremum over centers $y$ of $\nu(B(y,R))/|B_R|$, followed by the appropriate lower or upper limit as $R\to\infty$. The uniform-in-$y$ estimates of \cref{p2:lem:uniform-landau} give the final statement directly.
\end{proof}

\paragraph{Normalization and scope.}
The denominator has an intrinsic reproducing-kernel interpretation. For each fixed $R>0$,
\begin{equation}\label{p2:eq:local-trace-density}
 \int_{B(x,R/n)}K_n(t,t)\dd\mu(t)
 \longrightarrow q_E(x)|B_R|.
\end{equation}
Thus the density theorem compares the number of observations with the trace of the polynomial concentration operator in the same window. The factor $w$ cancels between the Christoffel asymptotic and $\dd\mu=w\dd t$. Equivalently, using the equilibrium density $\dd\nu_E=d!q_E(t)\dd t$ from \eqref{p1:eq:equilibrium-density} and $N_n\sim n^d/d!$, the denominator may be written as $N_n\nu_E(B(x,R/n))$ without changing the iterated densities.

Both densities use a fixed interior center and take $n\to\infty$ before $R\to\infty$. The lower count concerns Christoffel-normalized sampling. If a positive-weight rule satisfies $0<c\le\omega_{n,\lambda}K_n(\lambda,\lambda)\le C<\infty$, it gives the same normalized inequality with adjusted bounds. At exact cardinality this comparison follows from \cref{p2:prop:critical-weight-reduction}. For general redundant weighted rules, \cref{p2:thm:tangent-sampling} measures sampling through the mass attached to the nodes.

\subsection{Connections with polynomial, spherical, and RKHS density}

The trace in \eqref{p2:eq:local-trace-density} is the natural local dimension. In the spherical-harmonic setting of Marzo \cite{Marzo2007}, rotational symmetry makes this density constant, and the concentration operator supplies the corresponding Landau count. In orthonormal tangent coordinates its frequency body is the Euclidean unit ball. Dividing counts by Euclidean ball volume gives the threshold $|B_1|/(2\pi)^d$; Marzo's convention divides by $R^d$, so its constant contains an additional factor $|B_1|$. The compact-manifold theorem of \cite{OrtegaPridhnani2012} uses an analogous frequency-first, large-window-second limit. Here $A_x$ can be anisotropic and vary with $x$.

The equilibrium normalization agrees with the necessary polynomial sampling densities of Berman--Ortega-Cerd\`a \cite{BermanOrtega2018} and interpolation densities of Antezana--Marzo--Ortega-Cerd\`a \cite{AMO2021}. The latter authors count nodes in fixed physical balls. In \cref{p2:thm:local-landau}, the physical window $B(x,R/n)$ shrinks with the degree, so the proof first estimates counts at fixed rescaled radius $R$. On the ball, Li \cite{Li2025} obtains microscopic density conditions for weighted $L^p$ families in the intrinsic ball metric. The Euclidean windows used here express the threshold in terms of the volume of $A_x$ for each interior center.

The RKHS density theory in \cite{FuhrEtAl2017} also compares node counts with the integrated kernel diagonal, using weak localization and homogeneous approximation. In the polynomial array setting the spaces change with $n$. Passing to their common tangent space at $x$ makes the concentration estimates uniform over the subsequences used in the count.

\clearpage

\subsection{Packing, covering, and tiling for a lattice}

Petersen and Middleton \cite{PetersenMiddleton1962} characterized
alias-free lattice sampling by disjoint spectral translates along the
reciprocal lattice. The same periodization also gives the interpolation
criterion below.

Let $\Gamma=L\Z^d$ be a full-rank spatial lattice and
\[
 \Gamma^\#=2\pi L^{-T}\Z^d
\]
its reciprocal lattice. Let $Q$ be a fundamental cell of $\Gamma^\#$ and $A\subset\R^d$ be bounded, measurable, and of positive volume. For this generality, spectral containment in $A$ means that $\widehat f=0$ almost everywhere outside $A$. Define
\[
 m_A(\xi)=\sum_{\eta\in\Gamma^\#}\mathbf1_A(\xi+\eta),\qquad\xi\in Q.
\]

\begin{proposition}[Exact lattice classification]\label{p2:prop:lattice}
The lattice $\Gamma$ samples $\PW(A)$ if and only if $m_A\le1$ almost everywhere. It interpolates $\PW(A)$ if and only if $m_A\ge1$ almost everywhere. Both properties hold precisely when $m_A=1$ almost everywhere, that is, when $A$ tiles by $\Gamma^\#$. In the latter case the exponentials form an orthogonal basis of $L^2(A)$.
\end{proposition}

\begin{proof}
For $h\in L^2(A)$, periodize its zero extension over $\Gamma^\#$:
\[
 H(\xi)=\sum_{\eta\in\Gamma^\#}h(\xi+\eta),\qquad \xi\in Q.
\]
Fourier-series Parseval on $Q$ gives
\begin{equation}\label{p2:eq:periodization}
 \sum_{\gamma\in\Gamma}
 \left|\int_A h(\xi)e^{i\gamma\cdot\xi}\dd\xi\right|^2
 =|Q|\int_Q|H(\xi)|^2\dd\xi.
\end{equation}
Because $A$ is bounded, the multiplicities are bounded. When $m_A\le1$, at most one term of $H(\xi)$ is nonzero for almost every $\xi$, and the right side is $|Q|\norm{h}_2^2$. If $m_A\ge2$ on a set of positive measure, there are distinct $\eta_1,\eta_2\in\Gamma^\#$ and a positive-measure subset $S\subset Q$ for which both $S+\eta_j$ lie in $A$. The function $h=\mathbf1_{S+\eta_1}-\mathbf1_{S+\eta_2}$ has positive norm and vanishing periodization. Its Fourier coefficients in \eqref{p2:eq:periodization} all vanish, proving necessity of the packing condition.

For a finitely supported coefficient family $c_\gamma$, let $P(\xi)=\sum c_\gamma e^{i\gamma\cdot\xi}$. It is $\Gamma^\#$-periodic and
\begin{equation}\label{p2:eq:multiplicity-riesz}
 \int_A|P|^2=\int_Qm_A|P|^2,\qquad
 \int_Q|P|^2=|Q|\sum|c_\gamma|^2.
\end{equation}
The Riesz-sequence lower bound therefore holds when $m_A\ge1$. If $m_A=0$ on a set of positive measure, approximation of its indicator in $L^2(Q)$ by trigonometric polynomials shows that the lower bound fails. Bounded multiplicity supplies the upper bound in every case. Interpolation is equivalent to this Riesz-sequence property. The simultaneous case follows, and \eqref{p2:eq:multiplicity-riesz} then gives orthogonality.
\end{proof}

When $m_A=1$, integration over $Q$ gives $|A|=|Q|=(2\pi)^d/|\det L|$. Thus a simultaneous lattice has spatial density $1/|\det L|=|A|/(2\pi)^d$.

\section{Local Aliasing and Runge Instability}\label{p2:sec:runge}

\Cref{p2:prop:lattice} gives a geometric test for sampling on a
lattice. We now examine the effect of a reciprocal-lattice overlap
on polynomial reconstruction and compare it with the classical
endpoint instability on an interval.

\subsection{Classical sampling-rate benchmarks}\label{p2:sec:classical-rates}

First consider how node spacing affects uniform control of a polynomial even on an interval. For $M+1$ nodes $-1=x_0<\cdots<x_M=1$ and $n\le M$, set
\[
 B(M,n)=\sup\{\norm{p}_{[-1,1]}:p\in\Pi_n^1,\ \max_j|p(x_j)|\le1\}.
\]
This is the norming constant $\beta_n$ for the indicated nodes. For equispaced nodes, Coppersmith and Rivlin \cite{CoppersmithRivlin1992} identify its exponential scale. There are absolute positive constants such that
\begin{equation}\label{p2:eq:coppersmith-rivlin}
 c_0e^{c_1n^2/M}\le B(M,n)\le C_0e^{C_1n^2/M}.
\end{equation}
Consequently, uniform norming requires $M\gtrsim n^2$, and $M\ge c n^2$ suffices for each fixed $c>0$. Earlier proofs of quadratic-rate sufficiency are discussed in \cite[Remark~2.3]{AdcockPlatteShadrin2019}.

For quantile nodes generated by the fixed probability density
\[
 \rho(t)=g(t)(1-t)^\alpha(1+t)^\beta,\qquad
 \alpha,\beta>-1,\quad 0<c\le g\le C,
 \qquad \int_{-1}^{x_j}\rho(t)\dd t=j/M,
\]
Adcock, Platte and Shadrin \cite[Corollary~3.2 and Section~5]{AdcockPlatteShadrin2019} obtain a threshold exponent $2(\gamma_*+1)$, with $\gamma_*=\max\{\alpha,\beta,-1/2\}$. A sufficiently large $M\ge C_\rho n^{2(\gamma_*+1)}$ ensures bounded $B(M,n)$. Their lower estimate gives the necessary power when $\max\{\alpha,\beta\}>-1/2$. In the remaining case, the dimension of $\Pi_n^1$ gives the necessary linear power. These are uniform-norm estimates for the nodes prescribed by $\rho$. An $L^2(\mu)$ sampling inequality refers to the separate choice of orthogonality measure, and a uniform bound for a least-squares operator in the sup norm also requires control of that operator.

The exponent reflects the endpoint spacing in the angular variable. Near $t=1$, the last quantile gap has length comparable to $M^{-1/(\alpha+1)}$. Since $1-\cos\theta$ is comparable to $\theta^2$ near zero, the corresponding angular gap has length comparable to $M^{-1/[2(\alpha+1)]}$. Requiring it to resolve the angular scale $1/n$ gives the exponent $2(\alpha+1)$. At the other endpoint the exponent is $2(\beta+1)$, and the dimension supplies the third exponent $1$. Their maximum is precisely $2(\gamma_*+1)$.

For equispaced nodes with $n\ll M\ll n^2$, every fixed interior window of width $1/n$ contains a growing number of samples, yet \eqref{p2:eq:coppersmith-rivlin} forces $B(M,n)\to\infty$. In this example the loss occurs near the endpoints. We turn to a different source of instability in several dimensions, where aliasing on one fixed interior neighborhood can force the norming constant to diverge.

\subsection{A multivariate local aliasing mechanism}

On a rescaled lattice, frequencies that differ by a reciprocal lattice vector have the same nodal oscillation up to a constant phase. If both frequencies lie strictly inside the local spectrum, the recovery theorem transfers the resulting cancellation to polynomials.

\begin{theorem}[Superalgebraic norming obstruction]\label{p2:thm:aliasing}
Let $E=-E$ be a convex body, $x\in\interior E$, and $V$ a fixed neighborhood of $x$. Suppose a full-rank lattice $\Gamma$ has a nonzero $\eta\in\Gamma^\#$ such that
\[
 \interior A_x\cap(\interior A_x-\eta)\ne\varnothing.
\]
Suppose every node in $\Lambda_n\cap V$ belongs to
\[
 x+n^{-1}(\tau_n+\Gamma),
\]
with arbitrary offsets $\tau_n$. The exterior nodes may be placed arbitrarily. Then, for every integer $J\ge1$,
\begin{equation}\label{p2:eq:beta-rate}
 \beta_n(\Lambda_n)\ge c_Jn^J
\end{equation}
for all sufficiently large $n$. Equivalently, $\beta_n/n^s\to\infty$ for every fixed $s\ge0$.
\end{theorem}

\begin{proof}
By strict overlap we may choose a ball $B$ such that $\overline B$ and $\overline B+\eta$ are disjoint compact subsets of $\interior A_x$. Take $g$ with $0\ne\widehat g\in C_c^\infty(B)$ and, for an offset $\tau$, put
\[
 f_\tau(u)=(1-e^{i\eta\cdot(u-\tau)})g(u).
\]
For $u=\tau+\gamma$, the defining relation $\eta\cdot\gamma\in2\pi\Z$ makes the two terms cancel. Their Fourier supports, however, are disjoint. Hence $f_\tau$ is nonzero for every $\tau$. All these functions are Schwartz, and their Fourier transforms depend continuously on the single phase $e^{-i\eta\cdot\tau}$ in a fixed smooth spectral class. By compactness of the circle there are a common radius $R$ and a constant $c>0$ such that $\sup_{\norm{u}\le R}|f_\tau(u)|\ge c$ for every offset.

Apply \cref{p2:thm:recovery} at order $J$ to $f_{\tau_n}$ and denote the resulting polynomial by $P_n$. The constants are uniform over the compact family of phases. At every node $\lambda\in V$ we have $f_{\tau_n}(n(\lambda-x))=0$, so \eqref{p2:eq:all-orders} gives $|P_n(\lambda)|=O(n^{-J})$. For $\lambda\notin V$, the rescaled distance is at least a fixed multiple of $n$. Uniform Schwartz decay gives the same bound there. The fixed ball of radius $R$ eventually lies in $n(E-x)$, so $\norm{P_n}_E\ge c/2$ by uniform approximation on that ball. Dividing by the largest sampled value proves \eqref{p2:eq:beta-rate}, with the convention $\beta_n=\infty$ if the samples all vanish. For real polynomials, take a phase-rotated real part at a point attaining $\norm{P_n}_E$.
\end{proof}

For convex bodies, overlap of positive volume is equivalent to overlap of the interiors, since boundaries have measure zero and any nonempty open intersection contains a ball. Translates meeting only at the boundary therefore have disjoint spectra up to sets of measure zero in $L^2$.

\begin{remark}[Perturbations and the order of limits]\label{p2:rem:jitter}
If the rescaled local nodes are within $\varepsilon_n$ of $\tau_n+\Gamma$, uniform bounds on $\nabla f_\tau$ give
\[
 \beta_n\gtrsim_J(n^{-J}+\varepsilon_n)^{-1}.
\]
With agreement restricted to a rescaled window of radius $R_n\to\infty$, Schwartz decay contributes the further term $C_MR_n^{-M}$ to the denominator. Thus this lower bound diverges when the rescaled perturbation tends to zero and the window expands. The exponent $J$ is fixed throughout each estimate. An exponential lower bound would require control of the constants as $J$ increases.
\end{remark}

\begin{corollary}[Loss of the weighted sampling lower bound]\label{p2:cor:energy-aliasing}
Let $\mu$ have positive continuous density near $x$, and let positive sample weights satisfy $\sum_\lambda \omega_{n,\lambda}\le Cn^s$ for a fixed $s$. Under the aliasing hypotheses of \cref{p2:thm:aliasing}, the optimal lower weighted sampling constant obeys
\[
 \inf_{0\ne p\in\Pi_n^d}
 \frac{\sum_\lambda \omega_{n,\lambda}|p(\lambda)|^2}{\norm{p}_\mu^2}
 \le C_Jn^{s+d-2J}\qquad\text{for every fixed }J.
\]
\end{corollary}
\begin{proof}
For the polynomials $P_n$ constructed above, the bound on the largest sample gives
\[
 \sum_\lambda \omega_{n,\lambda}|P_n(\lambda)|^2
 \le C_J n^{-2J}\sum_\lambda \omega_{n,\lambda}
 \le C'_J n^{s-2J}.
\]
For the denominator, the family $f_\tau$ has a common positive lower bound on its $L^2$ norm over a sufficiently large fixed ball by compactness and the nonvanishing of each member. Uniform approximation on that ball and the positive continuous density at $x$ give $\norm{P_n}_\mu^2\ge c n^{-d}$ after the change of variables $u=n(t-x)$. Taking the quotient yields the asserted power.
\end{proof}
In particular, suppose the chosen positive weights satisfy the upper $L^2(\mu)$ MZ bound. Testing with $p=1$ gives $\sum_\lambda \omega_{n,\lambda}\le b\mu(E)$. The corollary then applies with $s=0$, so strict local lattice aliasing forces the lower MZ bound to vanish asymptotically, regardless of the total number of nodes.

\subsection{Two lattices of equal density on the square}

We now rotate a lattice inside a small square while retaining the same exterior nodes. This preserves both the macroscopic counting density and polynomial uniqueness, so the two orientations can be compared on those terms. Set
\[
 E=[-1,1]^2,\quad x_0=(\sqrt3/2,0),\quad
 h=2\pi/3,\quad\tau=(\pi/3,\pi/3),
\]
and let
\[
 U=x_0+(-1/128,1/128)^2,\quad
 W=x_0+(-1/256,1/256)^2,\quad
 V=x_0+(-1/512,1/512)^2.
\]
Use the tensor Chebyshev grid outside $V$ and insert one of two lattices inside $U$. More precisely, put
\[
 C_n=\{(\cos(j\pi/(2n)),\cos(k\pi/(2n))):0\le j,k\le2n\},
\]
\begin{equation}\label{p2:eq:hybrid-arrays}
 \Lambda_n^s=\big((x_0+n^{-1}(\tau+\Gamma_s))\cap U\big)
 \cup(C_n\setminus V),\qquad
 \Gamma_{\rm ax}=h\Z^2,\quad\Gamma_{\rm rot}=hR_{\pi/4}\Z^2.
\end{equation}

\begin{proposition}[Same density, different norming behavior]\label{p2:prop:hybrid}
Both arrays in \eqref{p2:eq:hybrid-arrays} determine $\Pi_n^2$ and have
$O(n^2)$ nodes. The counting measures
\[
 \nu_n^s=n^{-2}\sum_{\lambda\in\Lambda_n^s}\delta_\lambda,
 \qquad s\in\{\mathrm{ax},\mathrm{rot}\},
\]
have the same weak limit $\nu$, with
$\dd\nu\ge(256/121)q_E(t)\dd t$ in the bulk. Their uniform-norm sampling constants satisfy
\[
 \sup_n\beta_n(\Lambda_n^{\rm rot})<\infty,
 \qquad
 \frac{\beta_n(\Lambda_n^{\rm ax})}{n^s}\longrightarrow\infty
 \quad\text{for every fixed }s\ge0.
\]
\end{proposition}

\begin{proof}
The shared exterior set contains the $n+1$ nonpositive first-coordinate columns of $C_n$, each with $2n+1$ distinct second coordinates. If $p\in\Pi_n^2$ vanishes at these nodes, its restriction to each column is a univariate polynomial of degree at most $n$ with $2n+1$ zeros. It vanishes identically on each column. For every fixed second coordinate, $p$ now vanishes at $n+1$ distinct first coordinates, so $p=0$ on $E$. This proves uniqueness for both arrays.

For the cube, the Dubiner metric is
\[
 \delta_E(t,z)=\max_j|\arccos t_j-\arccos z_j|.
\]
A covering of $E$ by Dubiner balls of radius $\theta/n$, $\theta<\pi/2$, gives the norming bound $\sec\theta$ \cite{Vianello2018,PiazzonVianello2019}. This also applies to complex polynomials by taking a phase-rotated real part at a point of maximum modulus. To obtain such a covering for the rotated array, put $a=h/\sqrt2$. The columns of the rescaled lattice have horizontal spacing $a$ and vertical spacing $2a$ within each column. Rounding first to the nearest column and then along that column gives coordinate errors at most $a/2$ and $a$, respectively. On $U$, the infinitesimal metric has the bound
\[
 \max_j\frac{|v_j|}{\sqrt{1-t_j^2}}
 \le\frac{33}{32}\max(2|v_1|,|v_2|).
\]
Indeed, with $r=1/128$,
\[
 1-(\sqrt3/2+r)^2>3871/16384>(16/33)^2,
 \qquad 1-r^2>(32/33)^2.
\]
For $n\ge512$, rounding a point of $W$ to the physical rotated lattice stays inside $U$, since $a/n<1/256$. Integrating the metric bound along the segment to the rounded point gives Dubiner distance at most $\theta/n$, where
\[
 \theta=\frac{33a}{32}=\frac{11\sqrt2\pi}{32}<\frac\pi2.
\]
For a point outside $W$, round both angular coordinates to the nearest multiple of $\pi/(2n)$. The resulting $C_n$ node is at Dubiner distance at most $\pi/(4n)$. Since cosine is $1$-Lipschitz, its coordinate displacement is at most $\pi/(4n)<1/512$, so it stays outside $V$. The retained nodes therefore cover $E$ at the required radius. The finitely many smaller degrees have finite norming constants by uniqueness.

At $x_0$,
\[
 A_{x_0}=\{\xi:\tfrac12|\xi_1|+|\xi_2|\le1\}.
\]
The reciprocal lattice of $\Gamma_{\rm ax}$ contains $3e_1$. Since $A_{x_0}$ extends to $\pm2e_1$, the point $3e_1/2$ belongs to the interiors of both $A_{x_0}$ and $A_{x_0}+3e_1$. Only axis-lattice nodes occur inside $V$, so \cref{p2:thm:aliasing} gives the stated growth of $\beta_n(\Lambda_n^{\rm ax})$.

Finally, all coordinates of $C_n$ are algebraic. A local lattice coordinate is an algebraic number plus a nonzero algebraic multiple of $\pi$; in the rotated case the multiple is $(1+\sqrt2 k)/(3n)$, and in the axis case it is $(1+2k)/(3n)$. Hence the local and exterior arrays are disjoint. Lattice Riemann sums and the Chebyshev angular Riemann sum give the common limit
\begin{equation}\label{p2:eq:hybrid-density}
 \dd\nu(t)=\left[
 \frac9{4\pi^2}\mathbf1_U(t)+
 \frac{4\mathbf1_{E\setminus V}(t)}
 {\pi^2\sqrt{1-t_1^2}\sqrt{1-t_2^2}}
 \right]\dd t.
\end{equation}
Here $q_E(t)=1/(2\pi^2\sqrt{1-t_1^2}\sqrt{1-t_2^2})$. On $U$ the metric bound gives $q_E(t)\le1089/(1024\pi^2)$, so the lattice term alone is at least $(256/121)q_E(t)$. Outside $U$ the Chebyshev term equals $8q_E(t)$. This proves the claimed bound for the common density.
\end{proof}

In this oversampled example the two lattices have equal covolume, while only the axis lattice has the strict reciprocal overlap used in the proof. The difference in norming constants also has consequences for reconstruction of continuous functions.

\subsection{Divergence for a fixed continuous target}

Platte, Trefethen and Kuijlaars \cite{PlatteTrefethenKuijlaars2011} relate rapid approximation of analytic functions from equispaced samples to large condition numbers. Adcock, Platte and Shadrin \cite{AdcockPlatteShadrin2019} treat nonuniform sampling. For continuous functions, the uniform boundedness principle gives a consequence of norming loss for each prescribed sequence of polynomial-exact reconstruction operators.

\begin{proposition}[Runge consequence for polynomial-exact reconstruction]\label{p2:prop:runge}
Let each $\Lambda_n$ determine $\Pi_n^d$, and let
\[
 R_n:\ell^\infty(\Lambda_n)\longrightarrow\Pi_n^d
\]
be a fixed sequence of linear maps satisfying $R_n(p|_{\Lambda_n})=p$ for every $p\in\Pi_n^d$. Put $\mathcal I_nf=R_n(f|_{\Lambda_n})$. Then
\begin{equation}\label{p2:eq:reconstruction-norm}
 \norm{\mathcal I_n}_{C(E)\to C(E)}
 =\norm{R_n}_{\ell^\infty\to C(E)}\ge\beta_n(\Lambda_n).
\end{equation}
If $\sup_n\beta_n=\infty$, there is a single $f\in C(E)$ such that
\[
 \sup_n\norm{\mathcal I_nf}_\infty=\infty.
\]
In particular, its reconstruction errors are unbounded along a subsequence. If there are exactly $N_n$ unisolvent nodes, $\mathcal I_n$ is ordinary interpolation and its Lebesgue constant equals $\beta_n$.
\end{proposition}

\begin{proof}
By polynomial exactness, applying $R_n$ to the nodal values of any $p\in\Pi_n^d$ gives $\norm{R_n}\ge\beta_n$. Write $Q_n:C(E)\to\ell^\infty(\Lambda_n)$ for evaluation. To construct a norm-one right inverse $E_n$ of this contraction, choose continuous bumps supported in pairwise disjoint neighborhoods of the nodes, with each bump valued in $[0,1]$ and equal to one at its node. Multiplying by the prescribed values and summing extends each nodal vector with its norm preserved. From $Q_nE_n=I$ we obtain
\[
 \norm{R_nQ_n}\le\norm{R_n}\le\norm{R_nQ_n}\norm{E_n}
 =\norm{R_nQ_n}.
\]
This proves \eqref{p2:eq:reconstruction-norm}. Each $\mathcal I_n$ is bounded since $Q_n$ is a contraction and $R_n$ acts on a finite-dimensional space. If $\sup_n\beta_n=\infty$, uniform boundedness gives a single $f\in C(E)$ with $\sup_n\norm{\mathcal I_nf}_\infty=\infty$. The triangle inequality then gives unbounded errors along a subsequence. In the unisolvent square case every nodal vector corresponds to a unique polynomial, so the supremum defining $\norm{R_n}$ is exactly \eqref{p2:eq:beta}.
\end{proof}

Here the target $f$ can depend on the chosen reconstruction sequence. An analytic target would require estimates on a space carrying an analytic norm.

\section{Smooth Boundary Curvature and the Sampling Redundancy Gap}\label{p2:sec:spectral-geometry}

The equilibrium law in \cref{p2:thm:critical-equidistribution} fixes the
macroscopic distribution of any uniformly stable asymptotically critical
array. The tangent-sampling results of \cref{p2:sec:transfer} impose
an additional requirement on its microscopic sampling geometry. We now ask what happens when the
number of nodes approaches $N_n$. With only
$o(N_n)$ extra nodes, the relations among evaluations disappear locally
at almost every center, forcing exponential Riesz bases there with
common bounds. Continuity of the spectra then gives a basis at every
center, so one obstructed spectrum forces a positive redundancy gap.
Boundary curvature gives a geometric condition under which this obstruction
occurs. We compare the resulting gap with the sampling rules supplied by sparsification.

\subsection{Asymptotically critical sampling and the redundancy gap}

With $M_n>N_n$, relations among the weighted evaluation functionals form a space of dimension $M_n-N_n$. Such relations can survive at an individual center. When the excess is $o(N_n)$, however, their orthogonal projection has small trace on average over windows of radius $R/n$. This provides enough local independence to recover a basis at almost every center. The mechanism is related to the critical-density frame argument of Enstad--van Velthoven \cite{EnstadVelthoven2025} discussed above.

\begin{lemma}[Compactness of exponential Riesz bases under spectral limits]
\label{p2:lem:spectral-basis-compactness}
Let $A_j$ and $A$ be convex bodies in $\R^d$ with nonempty interior,
and suppose $A_j\to A$ in Hausdorff distance. For every $j$, assume
that
\[
 \bigl\{|A_j|^{-1/2}e^{-i\gamma\cdot\xi}:\gamma\in\Gamma_j\bigr\}
\]
is a Riesz basis of $L^2(A_j)$ with the same bounds
$0<\alpha\le\beta<\infty$. Then $L^2(A)$ admits a normalized
exponential Riesz basis with bounds $\alpha,\beta$.
\end{lemma}

\begin{proof}
We must preserve both finite linear independence and completeness.
Finite Gram matrices will give the former; uniform control of
coefficient tails will give the latter.

\smallskip\noindent\emph{1. Convergence of domains and frequency sets.}
Translate all spectral bodies by the same vector so that
$r\overline B\subset A$ for some $r>0$. This only changes the
normalized exponentials by constant phases. Set
$\varepsilon_j=d_H(A_j,A)$. Their support functions differ by at
most $\varepsilon_j$, and $h_A\ge r$ on the unit sphere. Thus,
for all large $j$,
\[
 (1-\varepsilon_j/r)A\subset A_j
       \subset(1+\varepsilon_j/r)A.
\]
These inclusions imply $\mathbf1_{A_j}\to\mathbf1_A$ in $L^1$,
$|A_j|\to|A|>0$, and eventual inclusion in $A_j$ of every compact
subset of $\operatorname{int}A$. They also place all the bodies,
after finitely many indices have been discarded, in a common ball
$B(0,R)$. Applying the Riesz lower bound to coefficients $1,-1$ at
distinct $\gamma,\gamma'\in\Gamma_j$ gives
\[
 2\alpha\le\frac1{|A_j|}\int_{A_j}
 |e^{-i\gamma\cdot\xi}-e^{-i\gamma'\cdot\xi}|^2\,d\xi
 \le R^2|\gamma-\gamma'|^2.
\]
Thus all $\Gamma_j$ have a common separation constant. In particular,
$\#(\Gamma_j\cap B(z,r))\le C(1+r)^d$, uniformly in $j,z,r$.
A diagonal extraction gives vague convergence of their counting
measures to the counting measure of a separated set $\Gamma$.

\smallskip\noindent\emph{2. The limiting Riesz inequalities.}
Every finite subset of $\Gamma$ can be matched with distinct points
of $\Gamma_j$ converging to it. The corresponding normalized Gram
matrices converge, by the $L^1$ convergence of the indicators and
uniform convergence of the finitely many exponentials on $B(0,R)$.
Their bounds $\alpha,\beta$ pass to the limit, so the limiting
exponentials form a Riesz sequence.

\smallskip\noindent\emph{3. Completeness of the limit.}
Take $f\in C_c^\infty(\operatorname{int}A)$,
extended by zero, and put
$F(\gamma)=\int f(\xi)e^{i\gamma\cdot\xi}\,d\xi$.
For all large $j$, the frame inequalities on $A_j$ give
\[
 \alpha\|f\|_2^2\le
 \frac1{|A_j|}\sum_{\gamma\in\Gamma_j}|F(\gamma)|^2
 \le\beta\|f\|_2^2.
\]
For every $m$, $|F(\gamma)|\le C_m(1+|\gamma|)^{-m}$.
The common counting bound applied to dyadic annuli yields, for
$T\ge1$ and $2m>d$,
\[
 \sum_{\substack{\gamma\in\Gamma_j\\|\gamma|>T}}|F(\gamma)|^2
 \le C_m\sum_{\ell\ge0}(2^{\ell+1}T)^d(2^\ell T)^{-2m}
 \le C_m'T^{d-2m}.
\]
The estimate is uniform in $j$ and tends to zero with $T$.
For a fixed continuous compact cutoff, vague convergence passes the
coefficient sum to $\Gamma$; the tail estimate then removes that
cutoff. Together with $|A_j|\to|A|$, this passes both frame bounds
to $\Gamma$ with normalization $|A|^{-1}$.
Since $C_c^\infty(\operatorname{int}A)$ is dense in $L^2(A)$, the
upper bound extends the coefficient map continuously. Each coordinate
of this extension remains the inner product with the corresponding
normalized exponential, because that inner product is a bounded
functional on $L^2(A)$. The lower bound therefore gives completeness
on all of $L^2(A)$. The limiting system is both a frame and a Riesz
sequence, hence a Riesz basis with the stated bounds. In particular,
the frequency sets cannot all escape to infinity: an empty limit
would contradict the lower frame bound for any nonzero test $f$.
\end{proof}

\begin{theorem}[Asymptotically critical sampling forces local bases at every center]
\label{p2:thm:near-critical-transfer}
Assume the standing hypotheses on $E$, $\mu$, and $U$. Suppose that finite
sets $\Lambda_n\subset E$ and arbitrary positive weights
$\omega_{n,\lambda}$ satisfy
\[
 a\norm{p}_\mu^2\le
 \sum_{\lambda\in\Lambda_n}\omega_{n,\lambda}|p(\lambda)|^2
 \le b\norm{p}_\mu^2,
 \qquad p\in\Pi_n^d,
\]
with fixed $0<a\le b<\infty$, and
\[
 M_n:=|\Lambda_n|=N_n+o(N_n).
\]
Then, for every $x\in U$, the space $L^2(A_x)$ admits
an exponential Riesz basis. The normalized exponentials may be chosen
with Riesz bounds $a/b$ and $b/a$.
The same conclusion holds when the hypotheses are imposed only along
an arbitrary sequence of degrees tending to infinity.
\end{theorem}

\begin{proof}
The proof has two stages. We first use the small global excess of
nodes to obtain local bases on a full-measure set of centers. We then
use the preceding compactness lemma to include the remaining centers.

\smallskip\noindent\emph{1. Averaging the relations among evaluations.}
Define $(T_n^\omega p)_\lambda=\sqrt{\omega_{n,\lambda}}\,p(\lambda)$ and let $H_n=T_n^\omega(T_n^\omega)^*$ be the Gram matrix of the weighted evaluation kernels. The sampling inequalities give $aI\le (T_n^\omega)^*T_n^\omega\le bI$, so $T_n^\omega$ has rank $N_n$. The orthogonal projection $P_n$ onto $\ker (T_n^\omega)^*$ has rank
\[
 r_n=M_n-N_n=o(n^d).
\]
The nonzero eigenvalues of $H_n$ lie in $[a,b]$, so
\begin{equation}\label{p2:eq:near-critical-projection}
 a(I-P_n)\le H_n\le b(I-P_n).
\end{equation}
Write $p_{n,\lambda}=(P_n)_{\lambda\lambda}\ge0$. Then
$\sum_\lambda p_{n,\lambda}=r_n$.

Pass to a degree subsequence on which
$\sum_n r_n/n^d<\infty$; we keep the same indices for simplicity.
For a positive integer $R$ define
\[
 e_{n,R}(x)=\sum_{\substack{\lambda\in\Lambda_n\\
                         |\lambda-x|<R/n}}p_{n,\lambda}.
\]
Integration over $U$ gives
\[
 \int_U e_{n,R}(x)\dd x
 \le |B(0,R)|\,n^{-d}r_n.
\]
Since the right-hand side is summable in $n$, Tonelli's theorem gives $\sum_n e_{n,R}(x)<\infty$ for almost every $x\in U$. Intersect the resulting full-measure sets over $R\in\mathbb N$ and fix $x$ in the intersection. Along the chosen degree subsequence, $e_{n,R}(x)\to0$ for every fixed window radius $R$.

\smallskip\noindent\emph{2. Local independence and a discrete limit.}
For nodes in $J_{n,R}=\Lambda_n\cap B(x,R/n)$, the principal block
$(P_n)_{J_{n,R}}$ is positive semidefinite and has trace $e_{n,R}(x)$.
Its operator norm is at most this trace. Compressing \eqref{p2:eq:near-critical-projection} to the same coordinate subspace gives
\begin{equation}\label{p2:eq:near-critical-local-block}
 a(1-e_{n,R}(x))I\le (H_n)_{J_{n,R}}\le bI.
\end{equation}
In particular, for
$c_{n,\lambda}=\omega_{n,\lambda}K_n(\lambda,\lambda)$,
\[
 a(1-e_{n,R}(x))\le c_{n,\lambda}\le b
 \qquad(\lambda\in J_{n,R}).
\]
Thus the local Gram estimate controls the weights of every node in a fixed rescaled window once $n$ is sufficiently large.

Consider the weighted measures from Theorem~\ref{p2:thm:tangent-sampling},
\[
 \sigma_{n,x}=\frac{n^d}{w(x)}
       \sum_{\lambda\in\Lambda_n}\omega_{n,\lambda}
                    \delta_{n(\lambda-x)}.
\]
They have uniformly bounded mass on every compact set. On each fixed
rescaled ball the diagonal kernel limit and the preceding estimates give
\[
 \frac{n^d\omega_{n,\lambda}}{w(x)}
   =\frac{c_{n,\lambda}}{q_E(x)}+o(1),
\]
uniformly over its nodes. The atom masses are therefore bounded above
and bounded away from zero there. In particular the local numbers of
nodes are uniformly bounded.

Dividing the blocks in \eqref{p2:eq:near-critical-local-block} by their
diagonal square roots gives the normalized Gram bounds
\[
 \frac{a(1-e_{n,R}(x))}{b}I
 \le (G_n)_{J_{n,R}}
 \le \frac{b}{a(1-e_{n,R}(x))}I
\]
for all sufficiently large $n$. The $2\times2$ lower bounds and the
compact-uniform kernel limit exclude collisions of rescaled nodes.
Indeed, suppose two sequences of distinct nodes converge after rescaling to $u$ and $v$ in the same window. The two-node Gram matrix has least eigenvalue $1-|G_n(\lambda,\lambda')|$, so the limiting correlation satisfies $|k_{A_x}(u-v)|/q_E(x)\le 1-a/b$. By continuity at zero this gives $\norm{u-v}\ge s$ for some $s>0$ independent of the window. The local counting bounds allow us to extract a vague limit of the counting measures. Each limit point is simple by separation, so the node approaching it can be identified in a small ball. A diagonal subsequence then gives convergence of the coefficients $c_{n,\lambda}$ at every such point. We obtain a separated set $\Gamma$ with coefficients $c_\gamma\in[a,b]$ and
\[
 \sigma_{n,x}\longrightarrow
 \sigma=\frac1{q_E(x)}\sum_{\gamma\in\Gamma}c_\gamma\delta_\gamma
 \quad\hbox{vaguely}.
\]
\smallskip\noindent\emph{3. Completeness and the Riesz inequalities.}
This last extraction is made at the fixed center $x$. To see the
normalization explicitly, let $\mathcal Uf=(2\pi)^{-d/2}\widehat f$
be the unitary Fourier identification with $L^2(A_x)$, and put
$e_\gamma=|A_x|^{-1/2}e^{-i\gamma\cdot\xi}$. Then
\[
 \langle\mathcal Uf,e_\gamma\rangle
       =q_E(x)^{-1/2}f(\gamma).
\]
The weighted tangent-sampling theorem therefore implies that
\[
 \left\{\sqrt{c_\gamma}\,|A_x|^{-1/2}
                 e^{-i\gamma\cdot\xi}:\gamma\in\Gamma\right\}
\]
is a frame for $L^2(A_x)$ with bounds $a,b$, and hence is complete.
To obtain the Riesz inequalities, take any finite subset of $\Gamma$
and a window containing it in its interior. The corresponding
principal block of $H_n$ converges entrywise to the Gram matrix of
the displayed weighted exponentials: its entries are
$\sqrt{c_{n,\lambda}c_{n,\lambda'}}G_n(\lambda,\lambda')$.
By \eqref{p2:eq:near-critical-local-block}, the limiting block has
spectrum in $[a,b]$. Thus the synthesis operator has a positive
lower bound and closed range. Completeness, obtained separately from
the tangent-sampling theorem, makes its range dense. It is therefore
onto, and the weighted exponentials form a Riesz basis.
Since $a\le c_\gamma\le b$, removing the factors $\sqrt{c_\gamma}$
gives normalized Riesz bounds $a/b,b/a$.

\smallskip\noindent\emph{4. Extension to every center.}
We have proved the basis property with common normalized bounds on a
full-measure subset $D\subset U$. Every nonempty open subset of $U$
has positive Lebesgue measure, so $D$ is dense in $U$.
To include an arbitrary $x_0\in U$,
choose $x_j\in D$ with $x_j\to x_0$. The generating map
$g_x(y)=y/\sqrt{1-(x\cdot y)^2}$ is jointly continuous on a compact
interior neighborhood of $x_0$ times $\partial E^*$, so
$A_{x_j}\to A_{x_0}$ in Hausdorff distance. Applying
\cref{p2:lem:spectral-basis-compactness} with bounds $a/b,b/a$ gives
the required basis at $x_0$. The initial summable subsequence can be
chosen within any prescribed sequence of degrees tending to infinity,
which proves the final assertion.
\end{proof}

The common bounds, rather than a uniform geometric bound on all
spectra throughout $U$, are essential in the last step. The centers
$x_j\to x_0$ lie eventually in one compact interior neighborhood,
where a common outer radius is available. Also, if the global
squared-norm bound ratio is $\kappa=b/a$, the displayed normalized
local bounds are $\kappa^{-1},\kappa$; no equality of the optimal
global and local condition numbers is asserted.

\begin{corollary}[A positive redundancy gap at fixed conditioning]
\label{p2:cor:fixed-condition-gap}
Assume that $L^2(A_{x_0})$ admits no exponential Riesz basis at some
point $x_0\in U$. Fix $1\le\kappa<\infty$, where $\kappa$ denotes the
ratio of upper to lower squared-norm MZ bounds. There are
$\delta_\kappa>0$ and $n_\kappa$ such that every positive-weight MZ rule
on $\Pi_n^d$ with bound ratio at most $\kappa$ satisfies
\[
 |\Lambda_n|\ge(1+\delta_\kappa)N_n
 \qquad(n\ge n_\kappa).
\]
\end{corollary}

\begin{proof}
The lower sampling bound implies $|\Lambda_n|\ge N_n$. If the asserted gap failed, we could choose $n_j\to\infty$ and rules with $1\le|\Lambda_{n_j}|/N_{n_j}<1+1/j$ and MZ bound ratio at most $\kappa$. Dividing the weights at each degree by the lower bound gives uniform bounds $1,\kappa$. By \cref{p2:thm:near-critical-transfer}, $L^2(A_{x_0})$ would then admit an exponential Riesz basis, a contradiction.
\end{proof}

\begin{remark}
One obstructed interior center suffices for both exact cardinality and
the relative gap. Equivalently, any sequence of positive weighted rules
with $|\Lambda_n|/N_n\to1$ must have MZ bound ratios tending to infinity.
The conclusion at every center concerns existence of a local basis;
it does not assert that every tangent of the original sampling array
is a basis. The size of $\delta_\kappa$ and its dependence on the domain
remain open. Here $\kappa$ controls the frame-operator condition number, and $\sqrt{\kappa}$ controls the condition number of the weighted evaluation map, with equality for optimal MZ bounds.
\end{remark}

\subsection{The geometric origin of the gap}

The volume $|A_x|$ gives the number of local degrees of freedom per
unit spatial volume, up to the Fourier factor $(2\pi)^{-d}$. Sampling
at that density also requires control of the correlations between
evaluations. By the off-diagonal kernel limit of Part~\ref{part:local},
the correlation between the normalized evaluations at $x+u/n$ and
$x+v/n$ tends to
\[
 \frac{1}{|A_x|}\int_{A_x}e^{i(u-v)\cdot\xi}\dd\xi .
\]
The density $w(x)$ cancels in this normalization. The shape of $A_x$
therefore controls the Gram matrices of finite configurations of
rescaled nodes. Localized recovery transfers the lower sampling bound
to the same Paley--Wiener space.

The role of the excess $M_n-N_n$ can be seen directly in the proof of
\cref{p2:thm:near-critical-transfer}. For the weighted evaluation
operator $T_n^\omega$, the projection onto $\ker (T_n^\omega)^*$ has trace $M_n-N_n$.
Its range consists of the linear relations among the weighted evaluations.
If $M_n-N_n=o(N_n)$, the trace contributed by the nodes in a ball of
radius $R/n$ tends to zero on average over its center, for each fixed
$R$. One can thus pass to a subsequence on which every bounded window
at almost every center inherits a lower Gram bound. Together with
sampling, this gives exponential Riesz bases with common bounds.
Spectral continuity and \cref{p2:lem:spectral-basis-compactness}
extend their existence to every center. Consequently one obstructed
interior spectrum forces a positive fraction of extra nodes globally (\cref{p2:cor:fixed-condition-gap}).
The argument applies to every choice of positive discrete weights for
the given reference measure.

\subsection{How boundary curvature enters the interior spectrum}

For balls, the required nonexistence result was obtained independently
by Wan \cite[Theorem~1.1]{Wan2026} and Ortega-Cerd\`a
\cite[Theorem~1.1]{Ortega2026}; invertible affine changes of variables
give the same obstruction for ellipsoids. Wan uses translated windows
and energy estimates, whereas Ortega-Cerd\`a uses stationary systems
and projections. Either theorem, combined with our transfer
argument, yields the ball and spherical applications here. For general smooth convex
spectra we use Ortega-Cerd\`a's extension
\cite[Corollary~8.1]{Ortega2026}, which excludes exponential Riesz
bases on every nonempty bounded convex open set with $C^2$ boundary
in dimension $d\ge2$. A single such local spectrum suffices.

\begin{corollary}[Smooth spectral obstruction for arbitrary positive weights]\label{p2:cor:smooth-obstruction}
Assume the standing hypotheses of \cref{p2:sec:setting}. If $d\ge2$
and $\partial A_{x_0}$ is $C^2$ at some center $x_0\in U$, then for
each fixed MZ bound ratio $\kappa$ there is a positive relative
redundancy gap as in \cref{p2:cor:fixed-condition-gap}.
In particular, there is no uniformly stable positive weighted rule
with $N_n+o(N_n)$ nodes, and no uniformly completely interpolating
polynomial array. If $0\in U$, the condition $\partial E^*\in C^2$
gives the same conclusions.
\end{corollary}

\begin{proof}
The boundary of a convex body has measure zero. Thus
\cite[Corollary~8.1]{Ortega2026} excludes an exponential Riesz basis
on $L^2(A_{x_0})$, and \cref{p2:cor:fixed-condition-gap} applies.
At the symmetry center $A_0=E^*$.
\end{proof}

The next lemma describes how positive boundary curvature controls
the full local spectrum near the symmetry center. The redundancy
obstruction already follows from the center itself; persistence in
a neighborhood is an additional geometric conclusion.

\begin{lemma}[Persistence of smooth positively curved local spectra]
\label{p2:lem:smooth-spectra-persist}
Let $d\ge2$ and let $E=-E\subset\R^d$ be a convex body with $C^2$ boundary and everywhere positive principal curvatures. There is $\delta>0$ with $B(0,\delta)\subset\interior E$ such that, for every $\norm{x}<\delta$, the entire boundary $\partial A_x$ is $C^2$ and has everywhere positive principal curvatures. More precisely, writing
\[
 H_E(\xi)=\sup_{z\in E}z\cdot\xi
\]
for the positively homogeneous support function, one has
\begin{equation}\label{p2:eq:smooth-local-spectral-sublevel}
 A_x=\{\xi\in\R^d:H_E(\xi)^2-(x\cdot\xi)^2\le1\}
 \qquad(\norm{x}<\delta).
\end{equation}
\end{lemma}

\begin{proof}
Let $r(u)$ be the inverse of the Gauss map of $\partial E$, for $u\in S^{d-1}$. Positive curvature and compactness make the Gauss map a $C^1$ diffeomorphism. The identity $\nabla H_E(u)=r(u)$ shows that $H_E$ is $C^2$ away from the origin. Moreover,
\[
 D^2H_E(u)u=0,
 \qquad
 D^2H_E(u)|_{u^\perp}>0.
\]
Indeed, on $u^\perp$, the derivative $Dr(u)$ is the inverse of the positive shape operator of $\partial E$ at $r(u)$. Set $Q_0=H_E^2/2$. Then
\[
 D^2Q_0
 =\nabla H_E\otimes\nabla H_E+H_E D^2H_E
\]
is positive definite at every nonzero point. The second term is strictly positive on every vector with a nonzero component perpendicular to $u$; on the remaining radial line, the first term is strictly positive because $\nabla H_E(u)\cdot u=H_E(u)>0$. Continuity on $S^{d-1}$ and degree-zero homogeneity of the Hessian yield a constant $\eta>0$ such that
\[
 D^2Q_0(\xi)\ge\eta I\qquad(\xi\ne0).
\]

Choose $0<\delta<\min\{\dist(0,\partial E),\sqrt\eta\}$. For $\norm{x}<\delta$, define
\[
 Q_x(\xi)=\tfrac12\bigl(H_E(\xi)^2-(x\cdot\xi)^2\bigr).
\]
Then
\[
 D^2Q_x(\xi)=D^2Q_0(\xi)-xx^T
 \ge(\eta-\delta^2)I\qquad(\xi\ne0).
\]
Since $Q_x$ is homogeneous of degree two, $Q_x(\xi)=O(\norm{\xi}^2)$ and $\nabla Q_x(\xi)=O(\norm{\xi})$ at zero. It therefore extends as a $C^1$ function there, with $\nabla Q_x(0)=0$. Integrating the Hessian bound along line segments, splitting at zero for any segment through the origin, shows that $Q_x$ is globally strongly convex. In particular $Q_x(\xi)\ge(\eta-\delta^2)\norm{\xi}^2/2$, so
$B_x=\{Q_x\le1/2\}$ is a compact convex body with the origin in its interior. On its boundary, Euler's identity gives
$\nabla Q_x(\xi)\cdot\xi=2Q_x(\xi)=1$, so the gradient does not vanish. The boundary is $C^2$, and its second fundamental form on a tangent vector $v$ is
\[
 \frac{v^TD^2Q_x(\xi)v}{\norm{\nabla Q_x(\xi)}}>0
 \qquad(v\ne0).
\]

It remains to identify this body with the local spectral body. If
$y\in\partial E^*$, then $H_E(y)=1$ and direct substitution gives
\[
 2Q_x\!\left(\frac{y}{\sqrt{1-(x\cdot y)^2}}\right)=1.
\]
Conversely, if $\xi\in\partial B_x$, put $y=\xi/H_E(\xi)$. Then $y\in\partial E^*$, and
$H_E(\xi)^2-(x\cdot\xi)^2=1$ implies $g_x(y)=\xi$. Hence
$g_x(\partial E^*)=\partial B_x$. Taking convex hulls gives $A_x=B_x$, proving \eqref{p2:eq:smooth-local-spectral-sublevel} and the claimed regularity.
\end{proof}

\begin{corollary}[Smooth physical domains require relative redundancy]
\label{p2:cor:smooth-physical}
Let $d\ge2$ and let $E=-E\subset\R^d$ have $C^2$ boundary and everywhere
positive principal curvatures. Assume the standing measure hypotheses,
with $0\in U$. For every fixed finite MZ bound
ratio $\kappa\ge1$, there exist $\delta_\kappa>0$ and $n_\kappa$ such
that every positive-weight sampling rule with ratio at most $\kappa$ has
\[
 |\Lambda_n|\ge(1+\delta_\kappa)N_n\qquad(n\ge n_\kappa).
\]
In particular, there is no uniformly stable rule with $N_n+o(N_n)$ nodes,
and no uniformly stable simultaneous sampling--interpolation array at
exact cardinality, for any choice of positive discrete weights.
\end{corollary}
\begin{proof}
At $x=0$, \cref{p2:lem:smooth-spectra-persist} gives
$\partial A_0=\partial E^*\in C^2$. Apply
\cref{p2:cor:smooth-obstruction} at this single center.
\end{proof}

\begin{remark}[Scope of the geometric obstruction]\label{p2:rem:smooth-scope}
The curvature hypothesis guarantees $C^2$ regularity of the local spectra near the symmetry center. The planar body $\{t_1^4+t_2^4\le1\}$ illustrates its role. This body has smooth strictly convex boundary, while its polar $A_0$ is the $\ell^{4/3}$ unit ball and fails to be $C^2$ at the coordinate axes. For a body of this kind, applying \cref{p2:cor:smooth-obstruction} requires locating a center whose full spectral boundary is $C^2$.

The size of the gap also depends on the stability one prescribes.
Proposition~\ref{p2:prop:sparse-sampling} provides arrays with any
fixed positive relative redundancy if the sampling bounds may depend
on that choice. In \cref{p2:cor:smooth-physical} the ratio $\kappa$
is fixed before the gap is obtained for $(E,\mu)$.
\end{remark}

\subsection{Sampling with any fixed relative redundancy}

Tchakaloff's theorem supplies a positive cubature rule exact through degree $2n$; its weighted evaluations form a finite Parseval frame. Spectral sparsification selects and reweights a subset of these nodes. The next proposition gives the resulting trade-off between the size of the subset and the sampling bounds.

\begin{proposition}[A consequence of spectral sparsification]
\label{p2:prop:sparse-sampling}
Let $\mu$ be a finite positive measure supported on the compact set $E$
for which $\norm{\cdot}_\mu$ is a norm on $\Pi_n^d$. For every $\rho>1$
there are distinct nodes $\lambda_j\in E$ and positive weights $\omega_j$,
with $N_n\le M_n\le\lceil\rho N_n\rceil$, such that
\begin{equation}\label{p2:eq:sparse-mz}
 (1-\rho^{-1/2})^2\norm{p}_\mu^2
 \le\sum_{j=1}^{M_n}\omega_j|p(\lambda_j)|^2
 \le(1+\rho^{-1/2})^2\norm{p}_\mu^2
 \qquad(p\in\Pi_n^d).
\end{equation}
\end{proposition}

\begin{proof}
Take a positive rule exact on $\Pi_{2n}^d$ from Tchakaloff's theorem \cite{BayerTeichmann2006},
with nodes $t_i$ and weights $a_i$, and a real orthonormal polynomial
basis $\phi_1,\ldots,\phi_{N_n}$. Its evaluation vectors
$v_i=\sqrt{a_i}(\phi_j(t_i))_{j=1}^{N_n}$ satisfy
$\sum_i v_iv_i^T=I$. The Batson--Spielman--Srivastava theorem
\cite[Theorem~3.1]{BSS2014}, with the ceiling convention explained in
\cite[Theorem~2.1 and Remark~2.2]{BownikVelthoven2025}, provides a
nonnegative reweighting supported on at most $\lceil\rho N_n\rceil$
vectors. The reweighted sum of the matrices $v_iv_i^T$ has eigenvalues between $(1-\rho^{-1/2})^2$ and $(1+\rho^{-1/2})^2$. Discard the terms with zero coefficient and multiply each remaining $a_i$ by its coefficient. Testing the matrix inequalities on the coefficient vector of a polynomial gives \eqref{p2:eq:sparse-mz}. For complex coefficients, apply the real inequalities separately to the real and imaginary parts. The lower bound implies injectivity, hence $M_n\ge N_n$.
\end{proof}

For fixed $\rho>1$, the displayed bounds are uniform in $n$. Their ratio grows as $\rho\downarrow1$. On the Paley--Wiener side, Bownik and van Velthoven \cite{BownikVelthoven2025} likewise obtain exponential frames arbitrarily close to critical density, with bounds depending on the excess density. The next subsection compares this upper bound with the geometric lower bound at a prescribed ratio of sampling bounds.

\subsection{The conditioning--redundancy law}\label{p2:subsec:redundancy-gap}

For a fixed domain and measure, let $m_n(\kappa)$ be the least number of
nodes in a positive-weight MZ rule on $\Pi_n^d$ whose upper-to-lower
squared-norm bound ratio is at most $\kappa$. Positive cubature makes
this quantity finite for every $\kappa\ge1$. Under the hypotheses of
Corollary~\ref{p2:cor:smooth-physical},
\[
 \Delta_{E,\mu}(\kappa):=
 \liminf_{n\to\infty}\frac{m_n(\kappa)}{N_n}-1>0
 \qquad(1\le\kappa<\infty).
\]
For $\kappa>1$, Proposition~\ref{p2:prop:sparse-sampling}, with
\[
 \rho_\kappa=
 \left(\frac{\sqrt\kappa+1}{\sqrt\kappa-1}\right)^2,
\]
gives $m_n(\kappa)\le\lceil\rho_\kappa N_n\rceil$. Positive cubature also gives $m_n(\kappa)\le\dim\Pi_{2n}^d$. Hence
\[
 0<\Delta_{E,\mu}(\kappa)\le\min\{2^d-1,\rho_\kappa-1\}.
\]
The condition number of the weighted evaluation map is at most $\sqrt\kappa$.

The positive lower bound comes from the local spectral obstruction,
whereas these upper bounds come from finite-dimensional discretization.
This places the polynomial problem within the general question of
sampling budgets studied, for example, in
\cite{DolbeaultCohen2022,BSS2014}: after fixing a stability level,
the structure of the function space can force a nonzero excess over
dimension. Corollary~\ref{p3:cor:contrasting-domains} will compare
two spaces of identical dimension for which this excess has different
asymptotic behavior. Determining the size of the positive gap requires
more quantitative information than the compactness argument supplies.
At $\kappa=1$, the sampling identity becomes exact cubature,
as developed in \cref{p2:sec:cubature}.

\subsection{Polygonal spectra and local admissibility}\label{p2:sec:polygonal-lattices}

For a centrally symmetric polygon, each spectral body is again a centrally symmetric polygon. These spectra have a different basis property.

\begin{proposition}[Every symmetric polygon is pointwise admissible]\label{p2:prop:polygon}
Suppose $d=2$ and $E$ is a symmetric convex polygon. For every $x\in\interior E$, $L^2(A_x)$ admits an exponential Riesz basis.
\end{proposition}

\begin{proof}
Write $E^*=\conv\{v_1,\ldots,v_m\}$ with the vertex set symmetric. Then
\begin{equation}\label{p2:eq:polar-vertices}
 A_x=\conv\left\{\frac{v_j}{\sqrt{1-(x\cdot v_j)^2}}:1\le j\le m\right\}.
\end{equation}
Indeed, $h(y)=\sqrt{1-(x\cdot y)^2}$ is positive and concave on $E^*$. If $y=\sum_j\lambda_jv_j$, then
\[
 \frac y{h(y)}=\sum_j\frac{\lambda_jh(v_j)}{h(y)}\frac{v_j}{h(v_j)},
 \qquad \sum_j\frac{\lambda_jh(v_j)}{h(y)}\le1.
\]
Adding the unused coefficient to the origin puts $y/h(y)$ in the convex hull on the right of \eqref{p2:eq:polar-vertices}, since that hull is symmetric. The reverse inclusion follows by taking $y=v_j$ in the definition of $A_x$. Thus $A_x$ is a centrally symmetric polygon. Every such polygon is a zonotope, and Debernardi--Lev \cite[Theorem~1.1]{DebernardiLev2022} construct an exponential Riesz basis on each zonotope.
\end{proof}

For each $x$ the proposition provides a basis on $A_x$. If these bases
have common normalized Riesz bounds, \cref{p3:thm:local-global} constructs
globally stable polynomial sampling with $N_n+o(N_n)$ nodes. Pointwise
existence alone does not establish that uniform hypothesis. For polygons
with at most six sides the bases are orthonormal, so
\cref{p3:cor:hexagon} applies. Appendix~\ref{p2:sec:hexagon-phase}
constructs a particular phase grid on the regular hexagon and identifies
its local bases. In higher dimensions the pointwise basis theorem also
applies when $A_x$ is a zonotope; centrally symmetric polytopes outside
this class require a different argument.

\subsection{Critical sampling on the square}\label{p2:subsec:rectangle-lattices}

For $E=\prod_{j=1}^2[c_j-r_j,c_j+r_j]$, let
\[
 \ell_j(x)=\sqrt{r_j^2-(x_j-c_j)^2}.
\]
Specializing \eqref{p1:eq:box-spectral-body} to $d=2$ gives
\[
 A_x=\{\xi:\ell_1(x)|\xi_1|+\ell_2(x)|\xi_2|\le1\}.
\]
Writing $a=1/\ell_1(x)$ and $b=1/\ell_2(x)$, this diamond tiles by the reciprocal lattice generated by $(a,b)$ and $(a,-b)$. The corresponding critical spatial lattice is generated by
\begin{equation}\label{p2:eq:rectangle-lattice}
 (\pi\ell_1,\pi\ell_2),\qquad(\pi\ell_1,-\pi\ell_2).
\end{equation}
The normalized kernel is explicitly
\begin{equation}\label{p2:eq:diamond-profile}
 \Phi_x(u)=\frac{k_{A_x}(u)}{k_{A_x}(0)}
 =\sinc\!\left(\frac{au_1+bu_2}{2}\right)
  \sinc\!\left(\frac{au_1-bu_2}{2}\right),\qquad
 \sinc t=\frac{\sin t}{t}.
\end{equation}
This follows by mapping the square $[-1/2,1/2]^2$ to the diamond with columns $(a,b),(a,-b)$.

For the lattice $\Gamma_x$ in \eqref{p2:eq:rectangle-lattice}, the local cardinal expansion is
\[
 f(u)=\sum_{\gamma\in\Gamma_x}f(\gamma)\Phi_x(u-\gamma),
 \qquad
 \norm{f}_2^2=\frac1{q_E(x)}\sum_{\gamma\in\Gamma_x}|f(\gamma)|^2.
\]
Orthogonal expansion gives convergence in $\PW(A_x)$, and the bounded evaluation norm gives locally uniform convergence. For the rectangle, a global polynomial construction realizing these local lattices is provided by the Padua points.

On $[-1,1]^2$ with
\[
 \dd\mu_C=\frac{\dd x_1\dd x_2}{\pi^2\sqrt{1-x_1^2}\sqrt{1-x_2^2}},
\]
the Padua nodes have exactly $N_n$ points, are unisolvent, and satisfy
\begin{equation}\label{p2:eq:padua-mz}
 \norm{p}_{L^2(\mu_C)}^2\asymp
 \frac1{N_n}\sum_{\lambda\in\mathrm{Pad}_n}|p(\lambda)|^2.
\end{equation}
The inequalities in \cite[Lemma~4.3 and Proposition~4.4]{Padua2007} hold for $1<p<\infty$ and give \eqref{p2:eq:padua-mz} at $p=2$. For an orthonormal polynomial basis $(p_k)$, the normalized square evaluation matrix $V_{\lambda k}=N_n^{-1/2}p_k(\lambda)$ has singular values bounded above and bounded away from zero uniformly in $n$. Its row norms imply $K_n(\lambda,\lambda)\asymp N_n$, giving uniform complete interpolation in our normalization. Affine transport gives rectangles.

The angular description, up to the usual choice of Padua family, is
\[
 \theta_1=\frac{k\pi}{n},\qquad \theta_2=\frac{l\pi}{n+1},
 \qquad k+l\text{ odd}.
\]
At a fixed interior point, the coordinate derivatives of $c_j+r_j\cos\theta_j$ are $-\ell_j(x)$. On the scale $n(\lambda-x)$, increments $(1,1)$ and $(1,-1)$ of the angular indices therefore generate the lattice in \eqref{p2:eq:rectangle-lattice}, up to a simultaneous change of sign. Along a subsequence, the angular offset of $x$ fixes the translation. Thus every tangent of this Padua family realizes the lattice predicted by the diamond $A_x$. The global estimate \eqref{p2:eq:padua-mz} concerns the Chebyshev $L^2$ norm. The associated sup-norm Lebesgue constants can grow with $n$.

At near-critical cardinality the equilibrium measure fixes the macroscopic
node distribution. At each interior center the microscopic sampling
inequality also forces an exponential Riesz basis on $A_x$.
The Padua array realizes both requirements on the square. For a general
body, can a family of admissible local bases be assembled into a single
global array? Part~\ref{part:local-global} gives a construction under
common Riesz bounds and controls the number of additional nodes needed
near the boundary.

\clearpage
\part{The Local-to-Global Problem}\label{part:local-global}

Part~\ref{part:sampling} gives necessary local conditions for stable
sampling with $N_n+o(N_n)$ nodes. We now prove a sufficient condition:
normalized exponential Riesz bases on the local spectral bodies at
a dense set of centers, with one common pair of positive, finite
bounds. Under positive continuous density on the whole interior,
this yields the existence criterion in
\cref{p3:cor:existence-equivalence}. The construction preserves
these bounds asymptotically and needs no continuous choice of the local
bases. For centrally symmetric polygons with at most six sides, all local
bases can be chosen orthonormal, and the resulting frame operators $S_n$
converge to the identity.

The proof combines finite interior patches with a positive correction
for the omitted region. Uniform Fourier localization controls every
polynomial in each patch, while the correction retains polynomial
directions concentrated near the boundary. Its node count is controlled
by the trace of a concentration operator. This establishes existence
with vanishing relative redundancy; exact cardinality and stability of
a prescribed array remain separate questions.

\section{Local-to-Global Sampling with Uniform Riesz Bounds}\label{p3:sec:theorem}

We use the setting of \cref{p1:sec:setting} with the additional
requirement that \eqref{p1:eq:density} hold on the entire interior:

\begin{equation}\label{p3:eq:2}
d\mu(x)=w(x)\,dx,\qquad
w\in C(\operatorname{int}E),\qquad w(x)>0.
\end{equation}
The restriction of $\mu$ to the boundary may be singular and may have positive mass.

Write $\nu_E$ for the probability-normalized equilibrium measure.

\begin{theorem}[Local-to-global sampling with uniform Riesz bounds]\label{p3:thm:local-global}

Under \eqref{p1:eq:BM} and \eqref{p3:eq:2}, let
$D\subset\operatorname{int}E$ be dense. Suppose that, for every
$x\in D$, there is a countable set $\Gamma_x\subset\mathbb R^d$ such that
\[
e_{x,\gamma}(\xi)=|A_x|^{-1/2}e^{-i\gamma\cdot\xi},
\qquad \gamma\in\Gamma_x,
\]
is a complete Riesz basis of $L^2(A_x,d\xi)$. Assume that the same constants $0<a\le b<\infty$ satisfy
\begin{equation}\label{p3:eq:4}
a\sum_\gamma|c_\gamma|^2
\le
\left\|\sum_\gamma c_\gamma e_{x,\gamma}\right\|_{L^2(A_x)}^2
\le
b\sum_\gamma|c_\gamma|^2
\end{equation}
for every $x\in D$ and every finitely supported coefficient family.
The optimal bounds at different centers may vary; only one common
positive lower bound and one common finite upper bound are required.

Then there exist finite sets of distinct nodes $\Lambda_n\subset E$, positive weights $\omega_{n,\lambda}$, and numbers $\epsilon_n\to0$ such that
\begin{equation}\label{p3:eq:5}
|\Lambda_n|=N_n+o(N_n)
\end{equation}
and
\begin{equation}\label{p3:eq:6}
(a-\epsilon_n)\|p\|_{L^2(\mu)}^2
\le
\sum_{\lambda\in\Lambda_n}\omega_{n,\lambda}|p(\lambda)|^2
\le
(b+\epsilon_n)\|p\|_{L^2(\mu)}^2,
\quad p\in\Pi_n^d.
\end{equation}
Equivalently, the associated frame operators $S_n=(T_n^\omega)^*T_n^\omega$ satisfy
\[
(a-\epsilon_n)I_{\Pi_n^d}\le S_n\le(b+\epsilon_n)I_{\Pi_n^d}.
\]
If $a=b=1$, then $\|S_n-I_{\Pi_n^d}\|\to0$.

\end{theorem}

Because each basis vector has norm one, the assumptions imply
$a\le1\le b$. The construction selects only finitely many reference
centers from $D$ at each fixed stage. Neither a continuous or measurable
choice of $\Gamma_x$ nor a uniform counting rate for the family is
required. In particular, a full-measure set of admissible centers
is sufficient.

\begin{corollary}[An existence criterion for asymptotically critical sampling]
\label{p3:cor:existence-equivalence}
Assume \eqref{p1:eq:BM} and \eqref{p3:eq:2}. The following assertions
are equivalent:
\begin{enumerate}[label=\textup{(\roman*)}]
\item There are positive weighted sampling rules on $\Pi_n^d$ with
$N_n+o(N_n)$ distinct nodes and one pair of positive finite MZ bounds
independent of $n$.
\item For every $x\in\operatorname{int}E$, $L^2(A_x)$ admits a
normalized exponential Riesz basis, with one pair of positive finite
bounds valid at all these centers.
\item The same local basis property, with common bounds, holds for
Lebesgue-almost every $x\in\operatorname{int}E$.
\item The same local basis property, with common bounds, holds on
some dense subset of $\operatorname{int}E$.
\end{enumerate}
\end{corollary}

\begin{proof}
Apply \cref{p2:thm:near-critical-transfer} with
$U=\operatorname{int}E$ to obtain \textup{(i)}$\Rightarrow$\textup{(ii)}.
The implication \textup{(ii)}$\Rightarrow$\textup{(iii)} is immediate,
and a full-measure subset of an open set is dense, giving
\textup{(iii)}$\Rightarrow$\textup{(iv)}.
Theorem~\ref{p3:thm:local-global} gives
\textup{(iv)}$\Rightarrow$\textup{(i)}; its bounds are uniformly
positive and finite for all sufficiently large degrees, and finitely
many initial degrees can be handled by Tchakaloff's theorem \cite{BayerTeichmann2006}.
\end{proof}

The three local conditions are also directly equivalent, without using
the global construction: for fixed bounds, the set of admissible centers
is relatively closed in $\operatorname{int}E$ by spectral continuity
and \cref{p2:lem:spectral-basis-compactness}, so a dense set suffices.
No continuous or measurable choice of the frequency sets is required.

The equivalence concerns existence with some common bounds. Global
bounds $\alpha,\beta$ give local normalized bounds
$\alpha/\beta,\beta/\alpha$, whereas local bounds $a,b$ give global
bounds approaching $a,b$. No equality of the optimal constants is
asserted.

\begin{corollary}[The asymptotically isometric case]
\label{p3:cor:isometric-equivalence}
Under the same assumptions, positive weighted rules satisfying
\[
|\Lambda_n|=N_n+o(N_n),\qquad \|S_n-I\|\longrightarrow0
\]
exist if and only if $L^2(A_x)$ admits a normalized exponential
orthonormal basis at every interior center. Equivalently, such
bases exist almost everywhere, or on a dense subset of the interior.
\end{corollary}

\begin{proof}
For sufficiency apply \cref{p3:thm:local-global} with $a=b=1$.
For necessity set $\delta_n=\|S_n-I\|\to0$ and use the proof of
\cref{p2:thm:near-critical-transfer}. With the notation of that proof,
\[
(1-\delta_n)(I-P_n)\le H_n\le(1+\delta_n)(I-P_n).
\]
The same summable-subsequence argument gives
$e_{n,R}(x)\to0$ for every integer $R$ at almost every $x$.
On each such window the compressed Gram matrix satisfies
\[
(1-\delta_n)(1-e_{n,R}(x))I
\le (H_n)_{J_{n,R}}\le(1+\delta_n)I.
\]
Use fixed bounds, for example $1/2,3/2$ for large $n$, in the
separation and compactness steps. The limiting coefficients
$c_\gamma$ are all one, and every finite Gram matrix of the limiting
normalized exponentials is the identity. The weighted tangent-sampling
theorem supplies completeness, so the limiting system is an
orthonormal basis. A full-measure set of such centers is dense.
Spectral continuity and \cref{p2:lem:spectral-basis-compactness},
with bounds $1,1$, extend the conclusion to every interior center
and prove the remaining equivalences.
\end{proof}

We retain the Fourier convention and the canonical representatives
of $\PW(A)$ from \cref{p1:sec:setting}; in particular,
$\|f\|_2^2=(2\pi)^{-d}\|\widehat f\|_2^2$. For the localization
argument define the scaled Sobolev norm
\begin{equation}\label{p3:eq:7}
\|f\|_{H_n^s}^2
=(2\pi)^{-d}\int_{\mathbb R^d}
\left(1+\frac{\norm{\zeta}^2}{n^2}\right)^s|\widehat f(\zeta)|^2\,d\zeta.
\end{equation}
The multiplier $\theta(D/n)$ has Fourier symbol $\theta(\zeta/n)$.

\section{Continuity, Separation, and Local Counting Bounds}\label{p3:sec:local}

\begin{lemma}[Continuity and total spectral volume]\label{p3:lem:1}
The map $x\mapsto A_x$ is Hausdorff continuous on $\operatorname{int}E$. On each compact interior set the bodies have a common positive inradius and a finite common outer radius. The function $q_E$ is positive and continuous, and the probability-normalized equilibrium measure satisfies
\begin{equation}\label{p3:eq:8}
d\nu_E(x)=d!\,q_E(x)\,dx\quad\text{on }\operatorname{int}E,\qquad
\nu_E(\partial E)=0,\qquad
\int_{\operatorname{int}E}q_E(x)\,dx=\frac1{d!}.
\end{equation}

\end{lemma}

\begin{proof} On $K\times E^*$, $K\Subset\operatorname{int}E$, the denominators in \eqref{p1:eq:gxy} are uniformly positive. The generating map is uniformly continuous, and convexification preserves its Hausdorff bound. Moreover, $E^*\subset A_x$: every $y\in E^*$ belongs to the segment from zero to $y/\sqrt{1-(x\cdot y)^2}$. Compactness gives the outer radius. Volume continuity follows, for example, by converting small Hausdorff errors into small dilations using the common inradius.

The density formula is \eqref{p1:eq:equilibrium-density}. By \cite[Theorem~1.1]{BLMR2010}, the unnormalized Monge--Amp\`ere density is $d!|A_x|$, its total mass is $(2\pi)^d$, and it puts no mass on the real boundary. Dividing by the total mass gives all of \eqref{p3:eq:8}. This argument does not require a smooth boundary. \end{proof}

\begin{lemma}[Sampling, dilation, and separation]\label{p3:lem:2}
Suppose the normalized exponentials $|A|^{-1/2}e^{-i\gamma\cdot\xi}$, $\gamma\in\Gamma$, form a complete Riesz basis of $L^2(A)$ with bounds $a,b$. Write $q_A=|A|/(2\pi)^d$. Then
\begin{equation}\label{p3:eq:9}
a\|f\|_2^2
\le q_A^{-1}\sum_{\gamma\in\Gamma}|f(\gamma)|^2
\le b\|f\|_2^2,\qquad f\in \PW(A).
\end{equation}
For $t>0$, the spectrum $tA$ and set $t^{-1}\Gamma$ have the same normalized bounds. If $A\subset R\overline B$, then
\begin{equation}\label{p3:eq:10}
\operatorname{sep}\Gamma\ge\frac{\sqrt{2a}}R.
\end{equation}

\end{lemma}

\begin{proof} The synthesis operator $T:\ell^2(\Gamma)\to L^2(A)$ is invertible. Its polar decomposition shows that $TT^*$, as well as $T^*T$, has spectrum in $[a,b]$. Hence the analysis inequality has the same bounds $a,b$. Since
\[
\langle\widehat f,|A|^{-1/2}e^{-i\gamma\cdot\xi}\rangle
=\frac{(2\pi)^d}{\sqrt{|A|}}f(\gamma),
\]
Plancherel gives \eqref{p3:eq:9}. The bounds are not replaced by their reciprocals; reciprocal bounds concern the dual basis.

Unitary dilation on $L^2$ proves the second statement. Finally, for $\gamma\ne\lambda$, use the two coefficients $1,-1$ and $|e^{it}-1|\le|t|$:
\[
2a\le\|e_\gamma-e_\lambda\|_{L^2(A)}^2
\le R^2\norm{\gamma-\lambda}^2.
\]
This proves separation. \end{proof}

\subsection{An upper counting bound from kernel tails}

The proof needs only an upper counting estimate. The following argument supplies it without an appeal to a density theorem or a lattice structure.

\begin{lemma}[Counting from a Riesz lower bound]\label{p3:lem:3}
Let $A$ be a compact spectral body and suppose its normalized exponentials indexed by $\Gamma$ have Riesz lower bound $a>0$. Completeness is not needed here. For each fixed cube $Q$ of positive volume,
\begin{equation}\label{p3:eq:11}
\#\bigl(\Gamma\cap(z+nQ)\bigr)
\le q_A n^d|Q|+o_{A,a,Q}(n^d),
\end{equation}
uniformly over $z\in\mathbb R^d$.

\end{lemma}

\begin{proof} The reproducing kernel of $\PW(A)$ is
\[
k_A(u)=(2\pi)^{-d}\int_Ae^{iu\cdot\xi}\,d\xi,\qquad
k_A(0)=\|k_A\|_2^2=q_A.
\]
The normalized kernels $\kappa_\gamma=k_A(\,\cdot-\gamma)/\sqrt{q_A}$ have the same Riesz bounds as the normalized exponentials under the unitary Fourier identification. Set
\[
\tau_A(t)=q_A^{-1}\int_{\norm{u}>t}|k_A(u)|^2\,du.
\]
Then $\tau_A(t)\to0$. Separation follows as in Lemma~\ref{p3:lem:2}, so bounded sets contain finitely many points.

Let $B=z+nQ$, $D=B+B(0,t)$, and $m=\#(\Gamma\cap B)$. If $m=0$, there is nothing to prove. Let $V:\mathbb C^m\to \PW(A)$ have columns $\kappa_\gamma$, $\gamma\in\Gamma\cap B$, and put $G=V^*V\ge aI$. The orthogonal projection onto $W=\operatorname{ran}V$, viewed in ambient $L^2(\mathbb R^d)$, is
\[
P_W=VG^{-1}V^*.
\]
Multiplication by $\mathbf1_D$ need not preserve $\PW(A)$; all traces below are taken in the ambient space. Finite-rank cyclicity and positivity give
\[
\begin{aligned}
\operatorname{Tr}(\mathbf1_{D^c}P_W)
&=\operatorname{tr}\bigl(G^{-1}V^*\mathbf1_{D^c}V\bigr)\\
&\le a^{-1}\sum_{\gamma\in\Gamma\cap B}
\|\mathbf1_{D^c}\kappa_\gamma\|_2^2
\le m\,\tau_A(t)/a.
\end{aligned}
\]
No commutation of the two positive matrices is required for this trace inequality.

For an orthonormal basis $h_1,\ldots,h_m$ of $W$, the evaluation functional on $\PW(A)$ has squared norm $q_A$. Bessel's inequality therefore gives
\[
\sum_{j=1}^m|h_j(u)|^2\le q_A.
\]
Consequently,
\begin{equation}\label{p3:eq:12}
\operatorname{Tr}(\mathbf1_D P_W)\le q_A|D|,
\qquad
m\bigl(1-\tau_A(t)/a\bigr)\le q_A|D|.
\end{equation}
If $Q$ has side length $\ell$, choose $t=\sqrt n$. Then
\[
|D|\le(n\ell+2\sqrt n)^d,\qquad
1-\tau_A(\sqrt n)/a\longrightarrow1.
\]
Both estimates are independent of $z$, proving \eqref{p3:eq:11}. \end{proof}

Landau's density theorems are consistent with this estimate and give the stronger equality of upper and lower uniform densities for a complete Riesz basis. That stronger result is not required below. In particular, no lattice-type $O(n^{d-1})$ counting remainder is asserted for a general $\Gamma$.

\section{Uniform Localization and Finite Patch Samples}\label{p3:sec:localization}
The recovery results of Part~\ref{part:local} approximate each fixed
function with smooth compactly supported spectral data by polynomials.
For the converse construction
we need a bound uniform over all polynomials in $\Pi_n^d$.
We obtain it in three steps: control complex growth with a polynomial
prefactor, suppress Fourier frequencies outside a fixed local spectral
envelope, and convert the resulting Sobolev error into an error in
point samples. All windows are fixed throughout these steps.

\subsection{Interior complex growth}

We use the Bernstein--Walsh inequality and the extremal representation \eqref{p1:eq:BCL}:
\begin{equation}\label{p3:eq:13}
|r(z)|\le\|r\|_E e^{(\deg r)V_E(z)},
\qquad
V_E(z)=\sup_{y\in\partial E^*}V_{[-1,1]}(y\cdot z).
\end{equation}
The interval expansion, uniform away from its endpoints, implies
\begin{equation}\label{p3:eq:14}
V_E(x+iy)\le h_{A_x}(y)+C_K\norm{y}^2,
\quad x\in K,\quad \norm{y}\le t_K,
\quad K\Subset\operatorname{int}E.
\end{equation}
The remainder is bounded before taking the polar supremum; no derivative of a maximizing polar direction is used.

\begin{lemma}[A polynomial complex-growth prefactor]\label{p3:lem:4}
For every compact $K\Subset\operatorname{int}E$, there are integers $k_n=o(n)$ and $C_K<\infty$ such that, for $n\ge1$,
\begin{equation}\label{p3:eq:15}
|p(x+iy)|
\le C_K n^{d/2}
e^{(n+2k_n)V_E(x+iy)}\|p\|_\mu,
\qquad p\in\Pi_n^d,\quad x\in K,\quad y\in\mathbb R^d.
\end{equation}

\end{lemma}

\begin{proof} First, for every compact interior $K'$,
\begin{equation}\label{p3:eq:16}
\sup_{K'}|p|\le C_{K'}n^{d/2}\|p\|_\mu.
\end{equation}
This is the unscaled interior estimate underlying
\eqref{p1:eq:nikolskii}. Apply the tensor Chebyshev evaluation bound on
nested interior cubes, use the positive lower bound for $w$, and take
a finite cover of $K'$. In particular, it is independent of any
sampling theorem.

Choose $r_0>0$ such that $K+2r_0\overline B\Subset\operatorname{int}E$, and choose $c_0>0$ with
\[
\tfrac12\le1-c_0\norm{t-x}^2\le1\quad(t\in E,\ x\in K).
\]
Set
\[
k_n=\left\lceil\frac{\log B_n^{\mathrm{BM}}}{c_0r_0^2}\right\rceil+1,\qquad
R_{n,x}(z)=\left(1-c_0\sum_{\ell=1}^d(z_\ell-x_\ell)^2\right)^{k_n}.
\]
Then $k_n=o(n)$, and on the real body
\[
0<R_{n,x}(t)\le e^{-c_0k_n\norm{t-x}^2}.
\]
For $\norm{t-x}\le r_0$, use \eqref{p3:eq:16}. Outside this ball, use \eqref{p1:eq:BM} and $B_n^{\mathrm{BM}}e^{-c_0k_nr_0^2}\le1$. Thus
\[
\|pR_{n,x}\|_E\le C_Kn^{d/2}\|p\|_\mu.
\]
At a pure imaginary displacement from its own real center,
\[
R_{n,x}(x+iy)=(1+c_0\norm{y}^2)^{k_n}\ge1.
\]
Apply Bernstein--Walsh to the polynomial $pR_{n,x}$, whose degree is at most $n+2k_n$, and divide by this last factor. \end{proof}

The multiplier uses $\sum(z_\ell-x_\ell)^2$, not a Hermitian squared distance. It is used only for a pointwise bound and is never inserted into the contour calculation. Suppressing $B_n^{\mathrm{BM}}$ is essential: $e^{o(n)}n^{-M}$ need not tend to zero for any fixed $M$.

\subsection{Fourier leakage outside the local spectra}

\begin{lemma}[Localization outside a convex spectral envelope]\label{p3:lem:5}
Let $\chi\in C_c^\infty(\operatorname{int}E)$. Let $C$ be a compact convex set containing every $A_x$ for $x\in\operatorname{supp}\chi$, and let $\theta\in C_c^\infty(\mathbb R^d)$ equal one on a neighborhood of $C$. For every $s\ge0$ and $L>0$,
\begin{equation}\label{p3:eq:17}
\|(I-\theta(D/n))(\chi p)\|_{H_n^s}
\le C_{\chi,\theta,s,L}n^{-L}\|p\|_\mu,
\qquad p\in\Pi_n^d,
\end{equation}
for all sufficiently large $n$. All windows and spectral margins are fixed before taking this limit.

\end{lemma}

\begin{proof}
The case $\chi=0$ is trivial. Otherwise we estimate each Fourier
coefficient outside the envelope by a complex displacement. A finite
Taylor polynomial of the cutoff will cancel the lower-order terms
without requiring an analytic extension of that cutoff.

\smallskip\noindent\emph{1. A uniform spectral gap.}
Choose $\delta>0$ such that $\theta=1$ on $C+\delta\overline B$.
For $\theta(\xi)\ne1$, let $c_\xi$ be the metric projection of
$\xi$ onto $C$, and set
\[
v=\frac{\xi-c_\xi}{\norm{\xi-c_\xi}}.
\]
The projection property gives
\[
\xi\cdot v-h_C(v)=\operatorname{dist}(\xi,C).
\]
This distance is at least $\delta$; boundedness of $C$ also gives $\operatorname{dist}(\xi,C)\ge\norm{\xi}-R$ for some $R$. Therefore
\begin{equation}\label{p3:eq:18}
\xi\cdot v-h_C(v)\ge c_1\langle\xi\rangle,\qquad
\langle\xi\rangle=(1+\norm{\xi}^2)^{1/2},
\end{equation}
on the entire unbounded set $\{\theta\ne1\}$.

\smallskip\noindent\emph{2. Decay after a complex displacement.}
Apply \eqref{p3:eq:14}--\eqref{p3:eq:15} with $y=-tv$, using the
compact set $K=\operatorname{supp}\chi$. Since $A_x=-A_x$,
$h_{A_x}(-v)=h_{A_x}(v)$, and $h_{A_x}(v)\le h_C(v)$.
Put $H=\max_{\|v\|=1}h_C(v)$ and $\kappa_n=2k_n/n\to0$.
After including the Fourier exponential, the exponent is at most
\[
 -nt\big[\xi\cdot v-h_C(v)-\kappa_nH
                 -(1+\kappa_n)C_Kt\big].
\]
First choose $0<t_0\le t_K$ with $2C_Kt_0\le c_1/4$.
Then take $n$ large enough that $\kappa_n\le1$ and
$\kappa_nH\le c_1/4$. These choices are independent of $p$ and
of the frequency $\xi$. Since
$\langle\xi\rangle\ge1$, \eqref{p3:eq:18} then gives
\begin{equation}\label{p3:eq:19}
\left|e^{-in\xi\cdot(x-itv)}p(x-itv)\right|
\le Cn^{d/2}e^{-c_2n\langle\xi\rangle t}\|p\|_\mu
\end{equation}
for $x\in\operatorname{supp}\chi$ and $0\le t\le t_0$.

\smallskip\noindent\emph{3. Cancellation for a smooth cutoff.}
For an integer $M\ge0$, write $D_v=v\cdot\nabla_x$, and define
\[
a_M(x,t)=\sum_{j=0}^M\frac{(-it)^j}{j!}D_v^j\chi(x),
\qquad
F_M(t)=\int_{\mathbb R^d}
e^{-in\xi\cdot(x-itv)}p(x-itv)a_M(x,t)\,dx.
\]
Every derivative of $\chi$ is supported in $\operatorname{supp}\chi$,
so all integrals have fixed compact real support. Denote the
holomorphic factor by
$E_{n,\xi}(x,t)=e^{-in\xi\cdot(x-itv)}p(x-itv)$.
Direct differentiation gives
\[
 \partial_tE_{n,\xi}=-iD_vE_{n,\xi},\qquad
 \partial_ta_M+iD_va_M
   =\frac{i(-it)^M}{M!}D_v^{M+1}\chi.
\]
Indeed, the terms of orders $0,\ldots,M-1$ in the second identity
cancel in pairs. Integration by parts in $x$, with no boundary term,
then gives
\[
 F_M'(t)=\int E_{n,\xi}(x,t)
                 \bigl(\partial_ta_M+iD_va_M\bigr)(x,t)\,dx.
\]
Thus we have the exact identity
\begin{equation}\label{p3:eq:20}
F_M'(t)=\frac{i(-it)^M}{M!}
\int_{\mathbb R^d}
e^{-in\xi\cdot(x-itv)}p(x-itv)D_v^{M+1}\chi(x)\,dx.
\end{equation}
This uses no holomorphic extension of the smooth cutoff. The direction $v$ is held fixed in $t$; it need not depend smoothly on $\xi$.

\smallskip\noindent\emph{4. Integrating the remainder and rescaling.}
For each fixed $j$, the directional derivative norms satisfy
\[
 \sup_{\|v\|=1}\|D_v^j\chi\|_1
 \le\sum_{|\alpha|=j}\frac{j!}{\alpha!}\|\partial^\alpha\chi\|_1.
\]
Consequently, \eqref{p3:eq:19} gives
\[
|F_M(t_0)|\le C_Mn^{d/2}e^{-c_2n\langle\xi\rangle t_0}\|p\|_\mu
\]
and
\[
\begin{aligned}
\int_0^{t_0}|F_M'(t)|\,dt
&\le C_M n^{d/2}\|p\|_\mu
       \int_0^{t_0}t^M e^{-c_2n\langle\xi\rangle t}\,dt\\
&\le C_M' n^{d/2-M-1}\langle\xi\rangle^{-M-1}\|p\|_\mu.
\end{aligned}
\]
The last step uses $\int_0^\infty t^Me^{-at}dt=M!a^{-M-1}$.
The endpoint satisfies the same polynomial bound because $t_0$ is
fixed and $e^{-c_2t_0X}\le C_M X^{-M-1}$ for $X\ge1$.
Since $F_M(0)=\widehat{\chi p}(n\xi)$, the identity
$F_M(0)=F_M(t_0)-\int_0^{t_0}F_M'(t)dt$ yields
\begin{equation}\label{p3:eq:21}
|\widehat{\chi p}(n\xi)|
\le C_Mn^{d/2-M-1}\langle\xi\rangle^{-M-1}\|p\|_\mu,
\qquad \theta(\xi)\ne1.
\end{equation}
Finally, the substitution $\zeta=n\xi$, with Jacobian $n^d$, yields
\[
\begin{aligned}
\|(I-\theta(D/n))(\chi p)\|_{H_n^s}^2
&\le C_Mn^{2d-2M-2}\|p\|_\mu^2
\int_{\mathbb R^d}\langle\xi\rangle^{2s-2M-2}\,d\xi.
\end{aligned}
\]
Choose $M+1>s+d/2$ and $M+1\ge d+L$. The integral is finite and the resulting norm is $O(n^{d-M-1})\|p\|_\mu$, proving \eqref{p3:eq:17}. \end{proof}

\subsection{Finite patch sampling by point evaluations}

\begin{lemma}[Discrete Sobolev evaluation]\label{p3:lem:6}
If $\Gamma\subset\mathbb R^d$ is separated and $s>d/2$ is an integer, then, for every $e\in H^s(\mathbb R^d)$,
\begin{equation}\label{p3:eq:22}
\sum_{\gamma\in\Gamma}|e(\gamma)|^2
\le C_{\Gamma,s}\|e\|_{H^s}^2.
\end{equation}
Equivalently, for every fixed center $x_0$,
\begin{equation}\label{p3:eq:23}
n^{-d}\sum_{\gamma\in\Gamma}|e(x_0+\gamma/n)|^2
\le C_{\Gamma,s}\|e\|_{H_n^s}^2.
\end{equation}

\end{lemma}

\begin{proof}
Choose $r>0$ below half the separation of $\Gamma$. Local Sobolev
embedding on $B(\gamma,r)$ gives, with a constant independent of
$\gamma$,
\[
 |e(\gamma)|^2\le C\sum_{|\alpha|\le s}
                 \int_{B(\gamma,r)}|\partial^\alpha e(u)|^2\,du.
\]
The balls are disjoint, so summing proves \eqref{p3:eq:22}.
Scaling by $u=n(x-x_0)$ proves the second assertion, using
\[
\|n^{-d/2}e(x_0+\cdot/n)\|_{H^s}=\|e\|_{H_n^s}.
\]
\end{proof}

\begin{lemma}[A frozen patch sampling inequality]\label{p3:lem:7}
Fix a real smooth window $\chi_j$ with compact interior support, a center $x_j$, and $w_j>0$. Let
\[
\Omega_j=(1+\delta)A_{x_j},\qquad q_j=\frac{|\Omega_j|}{(2\pi)^d}.
\]
Suppose a convex envelope of the local spectra on the window support lies compactly inside $\Omega_j$. Choose $\theta_j\in C_c^\infty(\operatorname{int}\Omega_j)$, equal one near that envelope. Use the dilated Riesz sampling set $\Gamma_j$ supplied by Lemma~\ref{p3:lem:2}, and define
\begin{equation}\label{p3:eq:24}
U_{j,n}p(u)=n^{-d/2}\sqrt{w_j}\,
(\chi_jp)(x_j+u/n).
\end{equation}
Then, for every $L>0$,
\begin{equation}\label{p3:eq:25}
\begin{aligned}
a\|U_{j,n}p\|_2^2-C_{j,L}n^{-L}\|p\|_\mu^2
&\le q_j^{-1}\sum_{\gamma\in\Gamma_j}|U_{j,n}p(\gamma)|^2\\
&\le b\|U_{j,n}p\|_2^2+C_{j,L}n^{-L}\|p\|_\mu^2.
\end{aligned}
\end{equation}

\end{lemma}

\begin{proof}
For this fixed patch, abbreviate $U=U_{j,n}p$ and set
$g=\theta_j(D)U\in\PW(\Omega_j)$, $e=U-g$.
The cutoff makes $g$ exactly bandlimited; the preceding leakage
estimate controls the difference from the actual windowed polynomial.
Lemma~\ref{p3:lem:5} and the unitary rescaling in
Lemma~\ref{p3:lem:6} give
\[
\|e\|_{H^s}\le C_{j,L}'n^{-L}\|p\|_\mu.
\]
Define $T_jh=q_j^{-1/2}(h(\gamma))_{\gamma\in\Gamma_j}$.
For $g$, Lemma~\ref{p3:lem:2} gives
$a\|g\|_2^2\le\|T_jg\|_{\ell^2}^2\le b\|g\|_2^2$.
For the error $e$, which need not be bandlimited,
Lemma~\ref{p3:lem:6} gives
$\|T_je\|_{\ell^2}\le C_{j,L}n^{-L}\|p\|_\mu$.
Moreover,
\[
\|U_{j,n}p\|_2^2=w_j\int\chi_j^2|p|^2\,dx\le C_j\|p\|_\mu^2.
\]
It follows that $\|g\|_2+\|T_jg\|_{\ell^2}\le C_j\|p\|_\mu$.
The difference-of-squares estimate now gives
\[
\begin{aligned}
 \bigl|\|T_jU\|_{\ell^2}^2-\|T_jg\|_{\ell^2}^2\bigr|
 &\le\|T_je\|_{\ell^2}
          \bigl(2\|T_jg\|_{\ell^2}+\|T_je\|_{\ell^2}\bigr)\\
 &\le C_{j,L}n^{-L}\|p\|_\mu^2.
\end{aligned}
\]
The same argument in $L^2$ bounds
$\bigl|\|U\|_2^2-\|g\|_2^2\bigr|$ by the same quantity.
Combining these two estimates with the sampling inequalities for $g$
proves \eqref{p3:eq:25}.

The sample sum for $U$ is finite because $U$ has compact support and
$\Gamma_j$ is separated. At reference nodes outside that support,
$U=0$ and $g=-e$, so their entire tail is already controlled by
the estimate for $T_je$. Thus restricting to the finite set of
nonzero samples introduces no additional error. Each retained term
is an actual point evaluation at
\begin{equation}\label{p3:eq:26}
\lambda_{j,\gamma}=x_j+\gamma/n,\qquad
\omega_{j,\gamma}=\frac{w_j}{n^dq_j}\chi_j(\lambda_{j,\gamma})^2.
\end{equation}
Zero weights are omitted. \end{proof}

The factor $w_j$ is constant. Thus high-order localization never differentiates the merely continuous function $w$.

\section{Low-Trace Boundary Correction and Global Assembly}\label{p3:sec:assembly}
The finite patch inequalities sample the part of the norm selected by
an interior cutoff. The remaining part may still have operator norm
close to one, even when its normalized trace is small. We therefore
retain it through a positive correction. Its small trace controls the
number of correction nodes, while an additive approximation controls
the error in the full polynomial norm.

\subsection{Christoffel mass and small-trace correction}

\begin{lemma}[Weak convergence and concentration trace]\label{p3:lem:8}
Let $K_n^\mu$ be the reproducing kernel of $\Pi_n^d$, and define
\[
d\beta_n(x)=\frac{K_n^\mu(x,x)}{N_n}\,d\mu(x).
\]
These probability measures converge weakly on $E$ to $\nu_E$. If $0\le\rho\le1$ is continuous, define
\[
\langle C_{\rho,n}p,p\rangle_\mu=\int_E\rho|p|^2\,d\mu.
\]
Then
\begin{equation}\label{p3:eq:27}
0\le C_{\rho,n}\le I,\qquad
\frac{\operatorname{tr}C_{\rho,n}}{N_n}
=\int_E\rho\,d\beta_n
\longrightarrow\int_E\rho\,d\nu_E.
\end{equation}

\end{lemma}

\begin{proof} The compact-interior diagonal asymptotic in \cref{p1:thm:main,p1:cor:christoffel} gives
\[
\frac{w(x)K_n^\mu(x,x)}{N_n}\longrightarrow d!q_E(x)
\]
uniformly on each compact interior subset, since $N_n/n^d\to1/d!$. Every weak subsequential limit of $\beta_n$ therefore restricts to $d!q_E\,dx$ on the interior. By Lemma~\ref{p3:lem:1}, this restriction already has total mass one; no boundary mass remains in the limit. Compactness proves weak convergence of the whole sequence. Summing the concentration quadratic form over an orthonormal polynomial basis gives the trace identity, and weak convergence gives its limit. \end{proof}

Small normalized trace does not mean small operator norm. The next lemma retains directions that may concentrate almost entirely in the omitted region.

\begin{lemma}[Positive point approximation of a low-trace operator]\label{p3:lem:9}
Let $H$ be a finite-dimensional real subspace of $C(E)$ with its $L^2(\mu)$ inner product. For $0\le\rho\le1$ continuous, define $C_\rho$ as above. For $0<\tau<1$, there is a positive point sampling operator $S_\rho$ such that
\begin{equation}\label{p3:eq:28}
\|S_\rho-C_\rho\|\le\tau,\qquad
\#\operatorname{nodes}(S_\rho)
\le C_{\mathrm{sp}}\tau^{-2}(1+\operatorname{tr}C_\rho).
\end{equation}
Its nodes lie in $\operatorname{supp}(\rho\mu)$. The same statement holds on the complexification of $H$.

\end{lemma}

\begin{proof} Choose a real orthonormal basis and let $v(x)$ be its evaluation column. Then
\[
C_\rho=\int_E\rho(x)v(x)v(x)^T\,d\mu(x).
\]
If this is zero, use no nodes. Otherwise, finite-dimensional convexity on the compact support of $\rho\mu$ gives an exact finite positive representation
\begin{equation}\label{p3:eq:29}
C_\rho=\sum_{i=1}^J\alpha_i v(z_i)v(z_i)^T,\qquad \alpha_i>0.
\end{equation}
Indeed, normalize $\rho\mu$ to a probability measure and apply the finite-dimensional convex-hull theorem to the continuous matrix-valued function $vv^T$. The initial number of nodes is immaterial.

Let $A$ have rows $\sqrt{\alpha_i}v(z_i)^T$. Then
\[
A^TA=C_\rho,\qquad
\|A\|^2\le1,\qquad
\|A\|_F^2=\operatorname{tr}C_\rho.
\]
Apply the additive row-sparsification theorem \cite[Theorem~5, full version]{CNW2016} with
\[
k=\max\{1,\lceil\operatorname{tr}C_\rho\rceil\}.
\]
In the form needed here, that theorem says: if $\|A\|\le1$ and $\|A\|_F^2\le k$, positive reweighting of at most $Ck/\tau^2$ rows gives a Gram matrix within $\tau$ of $A^TA$ in operator norm. Retained rows are still evaluation vectors, and their squared reweighting factors multiply $\alpha_i$. This proves \eqref{p3:eq:28}, including the harmless integer term.

For a complex function $u+iv$, add the real estimates for $u$ and $v$. \end{proof}

The approximation is additive relative to the ambient norm. This
distinction is essential for the node budget. A full-rank $C_\rho$
can have small trace, but a multiplicative lower bound
$(1-\tau)C_\rho\le S_\rho$ with $\tau<1$ would force $S_\rho$
to have full rank and hence to use at least $\dim H$ point evaluations.
The additive estimate \eqref{p3:eq:28} permits a much smaller correction.

\subsection{A finite cover with small counting overhead}

\begin{lemma}[Choice of a fixed stage]\label{p3:lem:10}
Let $D\subset\operatorname{int}E$ be dense.
For $0<\eta<1/100$, one can choose nonnegative windows
$\chi_0,\chi_1,\ldots,\chi_J\in C_c^\infty(\operatorname{int}E)$,
closed base cubes $Q_j$ with pairwise disjoint interiors, slightly expanded concentric cubes $Q_j^+\Subset\operatorname{int}E$, centers $x_j\in D\cap\operatorname{int}Q_j$, and $\delta>0$, so that, with $w_j=w(x_j)$, $\Omega_j=(1+\delta)A_{x_j}$, and $q_j=|\Omega_j|/(2\pi)^d$,
\begin{equation}\label{p3:eq:30}
0\le\chi_0\le1,\qquad
\rho=1-\chi_0^2,\qquad
\int_E\rho\,d\nu_E\le\eta^4,
\end{equation}
\begin{equation}\label{p3:eq:31}
\sum_{j=1}^J\chi_j^2=\chi_0^2,\qquad
\operatorname{supp}\chi_j\subset\operatorname{int}Q_j^+,
\end{equation}
\begin{equation}\label{p3:eq:32}
|w_j/w(x)-1|\le\eta\quad(x\in Q_j^+),
\end{equation}
\begin{equation}\label{p3:eq:33}
\bigcup_{x\in Q_j^+}A_x
\subset(1+\delta/4)A_{x_j}
\Subset(1+\delta/2)A_{x_j}
\Subset\Omega_j,
\end{equation}
and
\begin{equation}\label{p3:eq:34}
\sum_{j=1}^Jq_j|Q_j^+|
\le\frac{1+C_d\eta}{d!}.
\end{equation}
Every object is fixed within the stage.

\end{lemma}

\begin{proof} Choose a compact interior set of $\nu_E$-mass at least $1-\eta^4$, and a smooth $\chi_0$ equal one there. Set $\delta=\alpha=\eta/(20d)$. Keep the finitely many cubes of a sufficiently fine grid meeting $\operatorname{supp}\chi_0$, and expand their side lengths by $1+2\alpha$. They remain inside a fixed compact interior neighborhood.
In each base cube choose $x_j\in D\cap\operatorname{int}Q_j$.
The following mesh bounds hold uniformly for any such choices:
all points in a cube and its small expansion are uniformly close.

Positivity and continuity of $w,q_E$ make \eqref{p3:eq:32} and
\[
q_E(x_j)\le(1+\eta)q_E(x)\quad(x\in Q_j)
\]
hold when the mesh is small. For \eqref{p3:eq:33}, use
\[
A_x\subset A_{x_j}+\omega\overline B
\subset(1+\omega/r)A_{x_j},
\]
where $r>0$ is a common inradius, and reduce the mesh so that $\omega/r\le\delta/4$.

Choose nonnegative smooth $\psi_j$, equal one on a neighborhood of $Q_j$, supported inside $Q_j^+$. On a neighborhood of $\operatorname{supp}\chi_0$, $Z=\sum_j\psi_j^2\ge1$. Define
\[
\chi_j=\chi_0\psi_j/Z^{1/2}
\]
there and extend by zero. This is smooth because $\operatorname{supp}\chi_0\Subset\{Z>0\}$, and it proves \eqref{p3:eq:31}.

Finally, disjointness of the base interiors and Lemma~\ref{p3:lem:1} give
\[
\begin{aligned}
\sum_jq_j|Q_j^+|
&=(1+\delta)^d(1+2\alpha)^d\sum_jq_E(x_j)|Q_j|\\
&\le(1+\delta)^d(1+2\alpha)^d(1+\eta)\int_Eq_E(x)\,dx\\
&\le(1+C_d\eta)/d!.
\end{aligned}
\]
\end{proof}

Thin overlap is used only to limit the number of nodes. No estimate asserts that all polynomials have little energy in the overlap. The square partition is an exact pointwise identity.

\subsection{Assembly and diagonalization}

\begin{proof}[Proof of \cref{p3:thm:local-global}]

We first fix an accuracy parameter and construct a rule for all
sufficiently large degrees at that accuracy. We then compare its
operator bounds and node count before choosing a sequence of stages.
This order allows every cutoff and basis to remain fixed while the
localization estimates are applied.

\smallskip\noindent\emph{1. A finite collection of interior patches.}
Fix a stage $\eta$ and apply Lemma~\ref{p3:lem:10} to the dense
set $D$ in the hypothesis. At each of its finitely many centers
$x_j\in D$, choose a Riesz basis allowed by \eqref{p3:eq:4}, and
dilate it to $\Omega_j$. Choose $\theta_j$ equal one near $(1+\delta/2)A_{x_j}$ and compactly supported inside $\Omega_j$. This verifies every spectral-envelope hypothesis of Lemma~\ref{p3:lem:7}.

Define the bulk quadratic form
\begin{equation}\label{p3:eq:35}
\mathcal S_{\mathrm{bulk},n}(p)
=\sum_{j=1}^Jq_j^{-1}
\sum_{\gamma\in\Gamma_j}|U_{j,n}p(\gamma)|^2.
\end{equation}
Its weights are the positive point weights \eqref{p3:eq:26}. Because the stage contains finitely many fixed patches, Lemma~\ref{p3:lem:7} gives $r_{\eta,n}\to0$ and
\[
a\,\mathcal Q_n(p)-r_{\eta,n}\|p\|_\mu^2
\le\mathcal S_{\mathrm{bulk},n}(p)
\le b\,\mathcal Q_n(p)+r_{\eta,n}\|p\|_\mu^2,
\]
where
\[
\mathcal Q_n(p)=\sum_jw_j\int\chi_j^2|p|^2\,dx.
\]
Let $F_n=I-C_{\rho,n}$. Equations \eqref{p3:eq:31}--\eqref{p3:eq:32} imply
\[
(1-\eta)\langle F_np,p\rangle_\mu
\le\mathcal Q_n(p)
\le(1+\eta)\langle F_np,p\rangle_\mu.
\]
There is no missing singular boundary term: $\chi_0$ vanishes near $\partial E$. Since $0\le F_n\le I$ and $a\le b$, the operator $S_{\mathrm{bulk},n}$ satisfies
\begin{equation}\label{p3:eq:36}
aF_n-e_{\eta,n}I
\le S_{\mathrm{bulk},n}
\le bF_n+e_{\eta,n}I,
\qquad e_{\eta,n}=b\eta+r_{\eta,n}.
\end{equation}

\smallskip\noindent\emph{2. Restoring the omitted part of the norm.}
Apply Lemma~\ref{p3:lem:9} with $\tau=\eta$, obtaining $S_{\partial,n}$ with
\[
\|S_{\partial,n}-C_{\rho,n}\|\le\eta.
\]
Scale its positive weights by $a$, and set
\begin{equation}\label{p3:eq:37}
S_n=S_{\mathrm{bulk},n}+aS_{\partial,n}.
\end{equation}
The lower bound is $a(F_n+C_{\rho,n})-(e_{\eta,n}+a\eta)I$. For the upper bound, use
\[
bF_n+aC_{\rho,n}
=bI-(b-a)C_{\rho,n}\le bI.
\]
Thus
\begin{equation}\label{p3:eq:38}
(a-\varepsilon_{\eta,n})I
\le S_n\le
(b+\varepsilon_{\eta,n})I,
\qquad
\varepsilon_{\eta,n}=(a+b)\eta+r_{\eta,n}.
\end{equation}
In particular, the boundary correction does not add $a$ to the limiting upper frame bound.

\smallskip\noindent\emph{3. The node budget at a fixed stage.}
We count nodes before merging coincidences. Lemma~\ref{p3:lem:3} and \eqref{p3:eq:34} imply
\begin{equation}\label{p3:eq:39}
\begin{aligned}
M_{\mathrm{bulk},n}
&\le\sum_j\#\bigl(\Gamma_j\cap n(Q_j^+-x_j)\bigr)\\
&\le n^d\sum_jq_j|Q_j^+|+o_\eta(n^d)\\
&\le(1+C_d\eta)\frac{n^d}{d!}+o_\eta(n^d).
\end{aligned}
\end{equation}
Here and below, dependence on $\eta$ permits dependence on all fixed-stage choices. No rate uniform over all centers is being asserted.

By Lemma~\ref{p3:lem:8} and \eqref{p3:eq:30}, for all sufficiently large $n$,
\[
\operatorname{tr}C_{\rho,n}\le2\eta^4N_n.
\]
The number of correction nodes therefore satisfies
\begin{equation}\label{p3:eq:40}
M_{\partial,n}
\le C_{\mathrm{sp}}\eta^{-2}
+2C_{\mathrm{sp}}\eta^2N_n.
\end{equation}
This explains the choice of $\eta^4$ in
\eqref{p3:eq:30}: the trace mass is divided by the square of the
accuracy $\tau=\eta$, leaving a relative correction cost of order
$\eta^2$. The remaining integer term is absorbed by increasing the
degree threshold.
Choose a threshold for this stage so that, for every larger $n$,
\[
N_n\ge\eta^{-4},\qquad
r_{\eta,n}\le\eta,\qquad
o_\eta(n^d)\le\eta N_n.
\]
Since $n^d/d!\le N_n$, the whole stage satisfies
\begin{equation}\label{p3:eq:41}
\varepsilon_{\eta,n}\le(a+b+1)\eta,\qquad
M_n\le
\bigl(1+(C_d+1)\eta+3C_{\mathrm{sp}}\eta^2\bigr)N_n.
\end{equation}

\smallskip\noindent\emph{4. Choosing rules for every degree.}
Choose $\eta_k\downarrow0$ with
\[
\eta_k<\min\left\{\frac1{100},\frac{a}{2(a+b+1)}\right\}.
\]
For each $k$, first fix every object of that finite stage. Then choose strictly increasing integer thresholds $n_k$ such that \eqref{p3:eq:41} holds for every $n\ge n_k$ using stage $k$. For $n_k\le n<n_{k+1}$, use that stage and put
\[
\epsilon_n=(a+b+1)\eta_k<a/2.
\]
This is a construction for all sufficiently large integer degrees, not just a subsequence.

Remove zero weights and merge coincident physical nodes by adding their weights. The form is unchanged and the node count cannot increase. Its positive lower bound makes the evaluation map on $\Pi_n^d$ injective. Hence the number of distinct nodes is at least $N_n$; \eqref{p3:eq:41} supplies the matching upper asymptotic.

For the finitely many degrees below $n_1$, use Tchakaloff's theorem \cite{BayerTeichmann2006} on $\Pi_{2n}^d$ to represent the full identity by positive evaluations. Since $a\le1\le b$, these rules also satisfy \eqref{p3:eq:6}, with $\epsilon_n=0$ if desired. This completes the proof. \end{proof}

\section{Sampling on the Hexagon}\label{p3:sec:consequences}

On the regular hexagon, the spectrum changes shape across the incircle
(\cref{p1:sec:hexagon}). Both shapes admit orthonormal exponential
bases. We first record how Riesz bounds behave under affine changes of
spectrum, then use lattice tilings to obtain asymptotically isometric
sampling on every centrally symmetric hexagon.

\begin{corollary}[A common affine spectral model]\label{p3:cor:affine}

Assume \eqref{p1:eq:BM} and \eqref{p3:eq:2}. Suppose a fixed compact spectral body $A_0$ has a normalized complete exponential Riesz basis with bounds $a,b$, and for every interior $x$,
\[
A_x=L_xA_0+t_x,\qquad L_x\in GL(d,\mathbb R).
\]
Then Theorem~\ref{p3:thm:local-global} applies with the same $a,b$.

\end{corollary}

\begin{proof}The unitary map
\[
f\longmapsto |\det L_x|^{-1/2}
f\bigl(L_x^{-1}(\xi-t_x)\bigr)
\]
from $L^2(A_0)$ to $L^2(A_x)$ takes the normalized exponentials to normalized exponentials with frequencies $L_x^{-T}\Gamma_0$, multiplied by scalar phases of modulus one. It preserves completeness and both Riesz bounds. No differentiable choice of the affine maps, or integrability condition such as $L_x=DF(x)^T$, is needed by this patch construction. \end{proof}

\begin{corollary}[Centrally symmetric polygons with at most six sides]\label{p3:cor:hexagon}

Let $E\subset\mathbb R^2$ be a centrally symmetric convex polygon with nonempty interior and at most six sides. Under \eqref{p1:eq:BM} and \eqref{p3:eq:2}, there are positive weighted sampling rules such that
\begin{equation}\label{p3:eq:42}
|\Lambda_n|=\frac{(n+1)(n+2)}2+o(n^2),\qquad
\|S_n-I_{\Pi_n^2}\|\longrightarrow0.
\end{equation}
In particular, this holds for a regular hexagon, for area measure and for the reference measure
\begin{equation}\label{p3:eq:43}
d\mu_*(x)=\frac{|A_x|}{2\pi^2}\,dx=2q_E(x)\,dx
\end{equation}
used in Appendix~\ref{p2:sec:hexagon-phase}.

\end{corollary}

\begin{proof}\emph{Local orthogonal bases.} Let the vertices of $E^*$ be $\{\pm v_1,\ldots,\pm v_m\}$, $m\le3$. Then
\begin{equation}\label{p3:eq:44}
A_x=\operatorname{conv}\left\{
\frac{\pm v_j}{\sqrt{1-(x\cdot v_j)^2}}:1\le j\le m
\right\}.
\end{equation}
To verify this reduction, put $h(y)=\sqrt{1-(x\cdot y)^2}$, a positive concave function on $E^*$. If $y=\sum_i t_i z_i$ is a convex combination of polar vertices, then
\[
\frac{y}{h(y)}
=\sum_i\frac{t_i h(z_i)}{h(y)}\frac{z_i}{h(z_i)},
\qquad
\sum_i\frac{t_i h(z_i)}{h(y)}\le1.
\]
Place the remaining coefficient at zero, which belongs to the symmetric convex hull. This proves \eqref{p3:eq:44}, as in \cref{p2:prop:polygon}.

Because $E^*\subset A_x$, the body is full dimensional. After deleting redundant vertices, it is a parallelogram or a centrally symmetric hexagon. We give a complete lattice-tiling argument for the latter case.

Write its six vertices in counterclockwise order as
\[
u,v,w,-u,-v,-w,\qquad P=u+v,\quad Q=v+w,\quad L=[P\ Q].
\]
The shoelace identity gives
\begin{equation}\label{p3:eq:45}
|A_x|=\det(u,v)+\det(v,w)+\det(u,w)=\det(P,Q)>0.
\end{equation}
The three opposite edge pairs are matched by translations $P,Q,Q-P$. Set $\mathcal L=L\mathbb Z^2$, and let $\zeta=2\pi L^{-T}k$, $k\in\mathbb Z^2\setminus\{0\}$. Apply the divergence theorem to $\zeta e^{-i\zeta\cdot\xi}$. Opposite edge integrals cancel because the outward normals are opposite and the phases agree under the matching translations. Therefore
\[
\int_{A_x}e^{-i\zeta\cdot\xi}\,d\xi=0.
\]
The $\mathcal L$-periodization of $\mathbf1_{A_x}$ consequently has all nonconstant Fourier coefficients zero. Its constant coefficient is $|A_x|/|\det L|=1$. Uniqueness of Fourier series for integrable periodic functions gives
\begin{equation}\label{p3:eq:46}
\sum_{\ell\in\mathcal L}\mathbf1_{A_x}(\xi-\ell)=1
\quad\text{for almost every }\xi.
\end{equation}
Thus $A_x$ is a lattice fundamental domain, up to null sets.

For a parallelogram with cyclic vertices $u,v,-u,-v$, use $P=u+v$, $Q=v-u$; it is directly a fundamental domain for $[P\ Q]\mathbb Z^2$. In either case, folding onto a standard lattice cell preserves Lebesgue measure and the dual characters. Hence
\[
\left\{|A_x|^{-1/2}e^{-i\gamma\cdot\xi}:
\gamma\in2\pi L^{-T}\mathbb Z^2\right\}
\]
is a complete orthonormal basis. All local bounds are exactly $a=b=1$, including locations where a pair of hexagonal vertices becomes redundant. Theorem~\ref{p3:thm:local-global} now proves \eqref{p3:eq:42}.

\smallskip\noindent\emph{Verification for the specified measures.} Area measure $\mu_0=dx$, or its probability normalization, has positive continuous interior density. Lemma~\ref{p3:lem:1} shows that \eqref{p3:eq:43} is a probability measure. Since $A_x\supset E^*$, its density is bounded below by $|E^*|/(2\pi^2)>0$. Both measures dominate a positive constant multiple of area measure.

The required Bernstein--Markov constants have polynomial growth. Let $M=\|p\|_E$, choose a maximizing $x_0$, and rotate the phase so $p(x_0)=M$. For every $y\in E$, the real polynomial
\[
r_y(t)=\operatorname{Re}p((1-t)x_0+ty),\qquad 0\le t\le1,
\]
has degree at most $n$ and supremum norm at most $M$. Markov's inequality on $[0,1]$ gives $\|r_y'\|_\infty\le2n^2M$. Thus at $t_n=1/(4n^2)$, $|p|\ge M/2$ throughout the homothetic copy $(1-t_n)x_0+t_nE$. Convexity keeps it in $E$, and its area is $t_n^2|E|$. Therefore, for $n\ge1$,
\begin{equation}\label{p3:eq:47}
\|p\|_E\le\frac{8n^2}{\sqrt{|E|}}\|p\|_{L^2(dx)}.
\end{equation}
Domination of area measure gives $B_n^{\mathrm{BM}}\le C_E(1+n)^2$ for both $\mu_0$ and $\mu_*$, after adjusting the constant. \end{proof}

Corollary~\ref{p3:cor:hexagon} does not identify the original phase grid as the constructed array. It constructs independent patch samples and positive correction nodes; it neither proves exact cardinality $N_n$ nor preserves a prescribed weight rule.

\begin{corollary}[Identical dimensions and contrasting redundancy]
\label{p3:cor:contrasting-domains}
Let $B\subset\R^2$ be the unit disk and let $H\subset\R^2$ be a
regular hexagon, each equipped with area measure. Their total-degree
polynomial spaces have the same dimension
$N_n=(n+1)(n+2)/2$, but the following hold:
\begin{enumerate}[label=\textup{(\roman*)}]
\item On $H$ there are positive weighted rules with
$M_n=N_n+o(N_n)$ and $\|S_n-I\|\to0$.
\item On $B$, for each $1\le\kappa<\infty$, there exist
$\delta_\kappa>0$ and $n_\kappa$ such that every positive weighted
MZ rule with upper-to-lower squared-norm bound ratio at most $\kappa$
has at least $(1+\delta_\kappa)N_n$ nodes when $n\ge n_\kappa$.
\end{enumerate}
\end{corollary}

\begin{proof}
Part \textup{(i)} is \cref{p3:cor:hexagon}.
The area-measure Bernstein--Markov estimate
\eqref{p3:eq:47} applies to every planar convex body and hence to $B$.
At the center of $B$, the local spectrum is again the disk.
The nonexistence theorem of either Wan \cite[Theorem~1.1]{Wan2026}
or Ortega-Cerd\`a \cite[Theorem~1.1]{Ortega2026}, together with
\cref{p2:cor:fixed-condition-gap}, proves \textup{(ii)}.
\end{proof}

On the regular hexagon one can also prescribe the nodes using a single
continuous phase map. Appendix~\ref{p2:sec:hexagon-phase} constructs that
map and identifies its local lattices. At the incircle interface the
change of coordinates is a shear that preserves local completeness
(\cref{p2:prop:hex-interface}). Stability of this prescribed finite grid
is the question in \cref{p2:prob:hex-phase-sampling}.

\clearpage
\part{Cubature, Asymptotically Isometric Sampling, and Geometric Packing}\label{part:packing}

How many nodes does exact integration require beyond the dimension of the
polynomial space? Yuan Xu closes \emph{Minimal Cubature Rules} with two open
problems, under the headings \emph{Better Lower Bound for the Number of
Nodes} and \emph{A Family of Cubature Rules on the Cube}
\cite[Chapter~10, Sections~10.1--10.2, pp.~222--227]{Xu2025Minimal}.
The three-variable question of Li, Sun and Xu
\cite[p.~121]{LiSunXu2009} gives a concrete target. A rule of degree
$2n-1$ must have at least $\dim\Pi_{n-1}^3=\binom{n+2}{3}$ nodes.
Can a rule of this degree on the cube attain the leading term $n^3/6$ with an error
of order $n^2$?

Why should the shape of a frequency domain constrain a cubature rule?
Exactness through degree $2n$ recovers the squared norm of every
$p\in\Pi_n^d$ by sampling $|p|^2$. Zooming in about $x$ turns point
evaluation at a rescaled node $u$ into a Fourier coefficient against
the exponential $\xi\mapsto e^{-iu\cdot\xi}$ on $A_x$.
At $M_n=N_n+o(N_n)$, stable sampling requires exponential Riesz bases
on these bodies. Requiring the discrete norms to approach the continuous
norm, $\|S_n-I\|\to0$, strengthens the local requirement to orthonormal
bases (\cref{p3:cor:isometric-equivalence}). By Fuglede's conjecture for
convex bodies, proved in all dimensions by Lev and Matolcsi
\cite{LevMatolcsi2022}, this requirement is equivalent to the local
spectral body $A_x$ tiling frequency space by translations.
Positive cubature satisfies the exact identity $S_n=I$, so it inherits
the same obstruction. With excess nodes, the tangent measures still
satisfy a Parseval identity, and packing bounds quantify the node count.
After the estimates on balls and the classification of the critical
case, Section~\ref{piv:square-gap} returns to exactness. An analytic
change of weight on the square separates the two node counts.

\section{Positive Cubature and Isometric Sampling}\label{p2:sec:cubature}
We use the standing assumptions of \cref{p1:sec:setting} with
$U=\interior E$ as in \eqref{p3:eq:2}. Boundary mass is allowed.
For a positive rule the node set is a subset of $E$. We count distinct
nodes after adding weights at coincident nodes.

\subsection{Exactness and the frame identity}

Applying an exact rule to $|p|^2$ gives the norm identity on $\Pi_n^d$.
Conversely, the identity
\[
 \sum_j\omega_j|p(z_j)|^2=\int_E|p|^2\dd\mu,
 \qquad p\in\Pi_n^d,
\]
implies exactness through degree $2n$. Polarization gives exactness on
$p\overline q$ for $p,q\in\Pi_n^d$, and every monomial of degree at most
$2n$ factors into two monomials of degree at most $n$. After a common
rescaling of weights, the case $\kappa=1$ therefore describes precisely
positive cubature. In particular, \cref{p2:cor:smooth-physical} forces a
positive proportion of extra cubature nodes on the smooth bodies covered
by that result.

Tchakaloff's theorem \cite{BayerTeichmann2006} gives a general upper bound. The associated kernel identities also quantify the excess node count.

\begin{proposition}[Weighted sampling with a bounded redundancy factor]\label{p2:prop:cubature-frame}
For each $n$, there are distinct nodes $\lambda_j\in E$ and positive weights $\omega_j$ such that
\begin{equation}\label{p2:eq:parseval-cubature}
 \norm{p}_\mu^2=\sum_{j=1}^{M_n}\omega_j|p(\lambda_j)|^2
 \qquad(p\in\Pi_n^d),\qquad
 N_n\le M_n\le N_{2n}.
\end{equation}
Thus the weighted kernels form a Parseval frame. Their Christoffel-normalized coefficients $c_j=\omega_jK_n(\lambda_j,\lambda_j)$ satisfy
\begin{equation}\label{p2:eq:cubature-christoffel}
 0<c_j\le1,\qquad \sum_{j=1}^{M_n}c_j=N_n,\qquad
 \sum_{j=1}^{M_n}(1-c_j)=M_n-N_n.
\end{equation}
In particular $N_{2n}/N_n\to2^d$ for fixed $d$.
\end{proposition}

\begin{proof}
By Tchakaloff's theorem there exists a positive cubature formula on $E$ exact on $\Pi_{2n}^d$ with at most $N_{2n}$ nodes. Applying it to $|p|^2$ gives \eqref{p2:eq:parseval-cubature}, and the resulting evaluation map is injective, so $M_n\ge N_n$. To bound $c_j$, apply the identity to the unit vector $\kappa_{n,\lambda_j}$ and retain the contribution of its own node. To compute the sum of the $c_j$, apply the same identity to each member of an orthonormal basis $(\phi_k)_{k=1}^{N_n}$ and sum over $k$. The kernel identity $\sum_k|\phi_k(\lambda_j)|^2=K_n(\lambda_j,\lambda_j)$ gives \eqref{p2:eq:cubature-christoffel}. Finally, the dimension formula $N_n=\binom{n+d}{d}$ gives the stated limit.
\end{proof}

In terms of normalized kernels the Parseval identity reads
\[
 \norm{p}_\mu^2=\sum_j c_j|\langle p,\kappa_{n,\lambda_j}\rangle_\mu|^2.
\]
The coefficients $c_j$ describe the contribution of each normalized
evaluation to the exact norm identity. The trace identity fixes their
total mass, while the node count is bounded above by $N_{2n}$.

\subsection{Relaxing the exact norm identity}

Let $a_n,b_n$ be the optimal squared-norm bounds in

\[
 a_n\norm{p}_\mu^2\le\sum_{\lambda\in\Lambda_n}
 \omega_{n,\lambda}|p(\lambda)|^2\le b_n\norm{p}_\mu^2,\qquad p\in\PiN .
\]
The weighted evaluation map $T_n^\omega$ takes $p$ to
$(\sqrt{\omega_{n,\lambda}}\,p(\lambda))_{\lambda\in\Lambda_n}$.
Its frame operator $S_n=(T_n^\omega)^*T_n^\omega$ acts on $\PiN$ and has
extreme eigenvalues $a_n,b_n$. In particular,
\begin{equation}\label{piv:isometry}
 \norm{S_n-I_{\PiN}}=\max\{|a_n-1|,|b_n-1|\}.
\end{equation}
We call the array \emph{asymptotically isometric} when this quantity tends
to zero. If only $\kappa_n=b_n/a_n\to1$ is prescribed, multiplying all
weights at degree $n$ by $2/(a_n+b_n)$ gives operator error
$(\kappa_n-1)/(\kappa_n+1)$. This normalization preserves the nodes.

To record the distribution of the nodes themselves, let
\[
 \mathcal N_n(\varphi)=\frac1{N_n}\sum_{\lambda\in\Lambda_n}\varphi(\lambda),
 \qquad \varphi\in C(E).
\]
Each distinct node counts once, regardless of its sampling weight. In
particular $\mathcal N_n(1)=M_n/N_n$, the node ratio. Define
\begin{equation}\label{piv:rho}
 \rho_{\mathrm{iso}}(E,\mu)
 =\inf_{\norm{S_n-I}\to0}\limsup_{n\to\infty}\frac{M_n}{N_n}
 =\lim_{\kappa\downarrow1}\limsup_{n\to\infty}\frac{m_n(\kappa)}{N_n},
\end{equation}
where $m_n(\kappa)$ is the minimum from \cref{p2:subsec:redundancy-gap}.
For each $\kappa>1$ an asymptotically isometric array eventually obeys
that bound ratio. Conversely, choose $\kappa_j\downarrow1$ and minimizers
for $m_n(\kappa_j)$. Use the $j$th family beyond a threshold where its
node ratio is within $1/j$ of its limiting upper bound, switching at the
next threshold. The common rescaling above gives vanishing operator error.
The relative redundancy is $\rho_{\mathrm{iso}}-1$.

\subsection{Near-minimal cubature and equilibrium distribution}
\label{p2:subsec:cubature-equidistribution}

The trace defect also explains the distribution of positive cubature
nodes when their count approaches the dimension lower bound.

\begin{corollary}[Equidistribution of asymptotically minimal positive cubature]
\label{p2:cor:cubature-equidistribution}
Let $E$ and $\mu$ satisfy the assumptions of
\cref{p2:thm:critical-equidistribution}, and let
$\sigma_n=\sum_{j=1}^{M_n}\omega_{n,j}\delta_{x_{n,j}}$ be positive
cubature measures supported in $E$ and exact on $\Pi_{2n}^d$.
Then
\begin{equation}\label{p2:eq:cubature-two-measures}
 \sigma_n\rightharpoonup\mu,\qquad
 \frac1{N_n}\sum_{j=1}^{M_n}
 \omega_{n,j}K_n^\mu(x_{n,j},x_{n,j})\delta_{x_{n,j}}
 \rightharpoonup\nu_E.
\end{equation}
If $M_n=N_n+o(N_n)$, the unweighted empirical measures also satisfy
\[
 \frac1{M_n}\sum_{j=1}^{M_n}\delta_{x_{n,j}}
 \rightharpoonup\nu_E.
\]
\end{corollary}

\begin{proof}
Exactness on $\Pi_{2n}^d$ identifies the discrete and continuous inner
products on $\Pi_n^d$. Thus \eqref{p2:eq:equidist-mz} holds with
$a=b=1$, and $L_n=K_n^\mu$. The kernel-weighted and unweighted
conclusions follow from \cref{p2:thm:critical-equidistribution}.
For the first conclusion in \eqref{p2:eq:cubature-two-measures},
positivity and exactness on constants give $\sigma_n(E)=\mu(E)$.
Given $f\in C(E)$, approximate it uniformly by a fixed polynomial $q$.
For all sufficiently large $n$, exactness gives
\[
 \left|\int_E f\dd\sigma_n-\int_E f\dd\mu\right|
 \le 2\mu(E)\norm{f-q}_E.
\]
Polynomial density in $C(E)$ proves the claim.
\end{proof}

These limits distinguish the distribution of the nodes from the
measure integrated by the rule. The original cubature weights reproduce
$\mu$; multiplying each weight by $K_n^\mu(x_{n,j},x_{n,j})/N_n$
gives a probability measure converging to $\nu_E$. With only $o(N_n)$
extra nodes, the nonnegative defect in
\eqref{p2:eq:cubature-christoffel} has vanishing relative mass, so
unweighted node counts have the same equilibrium limit. Thus
near-minimal rules for different Bernstein--Markov measures on the
same body have a common limiting node distribution, even when their
cubature weights and exactness properties differ.

\section{Asymptotically Isometric Sampling and Euclidean Packing}\label{piv:general}

The upper bound comes from the lattice sampling picture of Petersen and
Middleton \cite{PetersenMiddleton1962}. Fix an interior point $x$ and
rescale the nodes by $u=n(\lambda-x)$. Sampling the local space
$\PW(A_x)$ on a lattice periodizes its Fourier data over a reciprocal
lattice cell $Q$. On an overlap of two spectral copies, opposite Fourier
amplitudes can cancel. When the interiors are disjoint, Fourier-series
Parseval gives exact isometry using $|Q|/|A_x|$ nodes per local degree
of freedom (\cref{p2:prop:lattice}). This is the reciprocal of the
packing fraction. Choosing a nearly optimal lattice on each patch and
controlling the overlap and boundary errors gives the global upper bound.

For the lower bound, the same function that certifies an upper bound on
packing density gives a lower bound on sampling density. For Euclidean
balls the packing bound was obtained independently by Gorbachev
\cite{Gorbachev2000} and by Cohn and Elkies \cite[Theorem~3.1]{CohnElkies2003}.
Xu \cite{Xu2003} used spherical design bounds to count the nodes of
positive cubature rules on the ball. Kuperberg \cite[Section~6]{Kuperberg2006}
obtained local node-density bounds for positive cubature on the simplex
and pointed out how a bound on the sum of squared weights gives a lower
bound on the node count. Gorbachev \cite[pp.~15--17]{Gorbachev2014}
applied a planar packing certificate to designs on $\mathbb S^2$,
including formulas with arbitrary positive weights.

Here this argument applies in a window about each interior point $x$.
The interior of $A_x-A_x$ consists of differences of interior
frequencies. Polarizing the tangent Parseval identity fixes the
Fourier transform of the sampling measure there (\cref{piv:gap}).
A packing certificate for $A_x$ then bounds the sum of squared weights
in the window. At fixed total weight, concentrating the mass on fewer
nodes increases this sum. This gives a lower count of distinct nodes
in windows of side $R/n$ about $x$, first as $n\to\infty$ and then
as $R\to\infty$ (\cref{piv:localcount}).

To express the node counts relative to the local degrees of freedom,
recall that $q_E(x)=|A_x|/(2\pi)^d$ and that the probability measure
$\dd\nu_E(x)=d!\,q_E(x)\dd x$ gives zero mass to $\partial E$
(\cref{p3:lem:1}). The bounds involve two affine invariants of $A_x$.
For a convex body $K$ put $D=K-K$ and write
$\det L$ for the covolume of a full-rank lattice $L$. Define
\begin{align}
 \Delta_{\mathrm{lat}}(K)
  &=\sup_{\substack{L\,\mathrm{lattice}\\
             (L\setminus\{0\})\cap\interior D=\varnothing}}
                  \frac{|K|}{\det L},\label{piv:latdef}\\
 \Delta_{\mathrm{LP}}(K)
  &=|K|\inf_{F\in\mathcal F(K)}\frac{F(0)}{\fc{F}(0)}.
                                                     \label{piv:lpdef}
\end{align}
The class $\mathcal F(K)$ consists of real even Schwartz functions with
$\fc{F}\ge0$, $F(0),\fc{F}(0)>0$, and $F\le0$ outside
$\interior D$. Here $(\mathcal F_{\mathrm c} f)(\xi)=\fc{f}(\xi)=\widehat f(2\pi\xi)$ uses cycles per
unit length, the convention in the translative linear program of Cohn
and Elkies \cite[Theorem~B.1]{CohnElkies2003}. An unmarked Fourier transform
continues to use the angular frequency of Part~\ref{part:local}.

\begin{theorem}[Distribution of asymptotically isometric samples]\label{piv:main}
For every asymptotically isometric positive weighted array on $E$ and
every nonnegative $\varphi\in C(E)$,
\begin{equation}\label{piv:lower}
 \liminf_{n\to\infty}\mathcal N_n(\varphi)\ \ge\
 \int_{\interior E}\frac{\varphi(x)}{\Delta_{\mathrm{LP}}(A_x)}\dd\nu_E(x).
\end{equation}
The same conclusion holds along any subsequence of degrees tending to
infinity on which $\norm{S_n-I}\to0$. There exists a single
asymptotically isometric positive weighted array for which every such
$\varphi$ satisfies
\begin{equation}\label{piv:upper}
 \limsup_{n\to\infty}\mathcal N_n(\varphi)\ \le\
 \int_{\interior E}\frac{\varphi(x)}{\Delta_{\mathrm{lat}}(A_x)}\dd\nu_E(x).
\end{equation}
In particular, taking $\varphi=1$ gives
\[
 \int_{\interior E}\frac{\dd\nu_E(x)}{\Delta_{\mathrm{LP}}(A_x)}
 \ \le\ \rho_{\mathrm{iso}}(E,\mu)\ \le\
 \int_{\interior E}\frac{\dd\nu_E(x)}{\Delta_{\mathrm{lat}}(A_x)} .
\]
\end{theorem}

The distinction between absolute and relative density is visible in the
two factors of the integrals. The equilibrium measure $\nu_E$ counts
local degrees of freedom. The reciprocal packing densities bound the
number of nodes per local degree of freedom. If all $A_x$ are affine
images of one body, the latter bounds are constant even where $q_E$ varies.
If an array attains the lower total count, its spatial distribution is
also determined. Put $h(x)=1/\Delta_{\mathrm{LP}}(A_x)$. Whenever
$M_n/N_n\to\int h\dd\nu_E$, applying \eqref{piv:lower} to $\varphi$
and $\norm{\varphi}_\infty-\varphi$ for $\varphi\ge0$ gives
\begin{equation}\label{piv:empirical}
 \frac1{M_n}\sum_{\lambda\in\Lambda_n}\delta_\lambda
 \ \rightharpoonup\ \frac{h\,\nu_E}{\int h\dd\nu_E}.
\end{equation}

\subsection{Comparison of packing densities}\label{piv:proofs}

We prove \cref{piv:main}. We first compare the two packing invariants,
then construct the array for the upper estimate and count the nodes of an
arbitrary asymptotically isometric array for the lower estimate.

\begin{lemma}[Comparison and continuity]\label{piv:continuity}
For every centrally symmetric spectral body $K$,
\[
 1/d!\ \le\ \Delta_{\mathrm{lat}}(K)
   \ \le\ \Delta_{\mathrm{LP}}(K)\ \le\ 1 .
\]
Both densities are affine invariant and continuous on full-dimensional
origin-symmetric convex bodies in the Hausdorff metric.
\end{lemma}
\begin{proof}
For an admissible lattice $L$, Poisson summation and the signs give
$\fc{F}(0)/\det L\le\sum_{\ell\in L}F(\ell)\le F(0)$.
For the upper bound $1$, take $F=\phi*\widetilde\phi$ with real
$\phi\in C_c^\infty(\interior K)$, $0\le\phi\le1$,
$\widetilde\phi(u)=\phi(-u)$, approximating $\mathbf1_K$ in $L^1$ and $L^2$.
Then $F(0)=\int\phi^2$ and $\fc{F}(0)=(\int\phi)^2$.
For the dimensional lower bound translate the center to zero and maximize
$|\det V|$ over matrices with columns in $K$. Replacing one column at a
time shows $V^{-1}K\subset[-1,1]^d$, while convexity gives
$V^{-1}K\supset\{u:\sum|u_j|\le1\}$. The box $V[-1,1]^d$ tiles by a
lattice and has volume at most $d!\,|K|$.

Under $K\mapsto TK+b$, admissible lattices transform by $L\mapsto TL$,
and test functions by $F\mapsto F(T^{-1}\cdot)$. For continuity put
\[
 u(K)=|K|/\Delta_{\mathrm{lat}}(K),\qquad
 \ell(K)=|K|/\Delta_{\mathrm{LP}}(K).
\]
The first is the infimum of admissible covolumes. The second is the
supremum of $\fc{F}(0)/F(0)$. Both increase under inclusion and are
homogeneous of degree $d$. When $K'$ is Hausdorff-close to $K$,
$(1-\epsilon)K\subset K'\subset(1+\epsilon)K$, with $\epsilon\to0$.
The squeeze proves continuity. The displayed dimensional bounds make both
integrals in Theorem~\ref{piv:main} finite.
\end{proof}

\subsection{Reciprocal lattices and the upper estimate}\label{piv:construction}
Let $L$ be a frequency lattice and let $\Gamma$ be the sampling lattice
with reciprocal $\Gamma^\#=L$ in the notation of \cref{p2:prop:lattice}.
Its point density is $D_\Gamma=\det L/(2\pi)^d$.
If translates of a spectral body $\Omega$ by $L$ have disjoint
interiors, then every $g\in\PW(\Omega)$ satisfies
\begin{equation}\label{piv:latticeparseval}
 D_\Gamma^{-1}\sum_{\gamma\in\Gamma}|g(\gamma)|^2=\|g\|_2^2 .
\end{equation}
Indeed, on any fundamental parallelepiped of $L$ the periodization of
$\widehat g$ has at most one nonzero summand almost everywhere. Its
$L^2$ norm squared equals $\|\widehat g\|_2^2$, and its Fourier
coefficients are $(2\pi)^d g(\pm\gamma)/\det L$ with the corresponding
choice of sign. Parseval's identity for Fourier series gives
\eqref{piv:latticeparseval}. An overlap of positive measure allows
two opposite Fourier amplitudes with zero periodization and therefore
precludes lattice sampling. This is the lattice criterion of
\cref{p2:prop:lattice} and Petersen--Middleton \cite{PetersenMiddleton1962}.
The least possible relative density of one exact sampling lattice for
$\Omega$ is consequently $1/\Delta_{\mathrm{lat}}(\Omega)$.
\begin{proof}[Proof of the upper bound in Theorem~\ref{piv:main}]
Fix a stage $0<\eta<1/100$. Use the cover of Lemma~\ref{p3:lem:10}. Its base cubes
$Q_j$ have disjoint interiors and their expansions $Q_j^+$ have side
lengths $(1+O_d(\eta))$ times those of $Q_j$. The windows satisfy
\[
 \sum_j\chi_j^2=\chi_0^2,\qquad
 \supp \chi_j\subset\interior Q_j^+,\qquad
 \int(1-\chi_0^2)\dd\nu_E\le\eta^4 .
\]
On $Q_j^+$ we have $|w(x_j)/w(x)-1|\le\eta$. A convex envelope of
$\{A_x:x\in Q_j^+\}$ lies strictly inside
$\Omega_j=(1+\delta)A_{x_j}$ with $\delta=O_d(\eta)$.
By Lemma~\ref{piv:continuity}, refining the mesh in the proof of
Lemma~\ref{p3:lem:10} also ensures $u(A_{x_j})\le(1+\eta)u(A_x)$ for $x\in Q_j$.
All refinements take place on a fixed compact interior neighborhood
at this stage. We choose the expanded cubes to have diameter at most $\eta$.

For each $j$ choose a packing lattice $L_j$ for $\Omega_j$ with
$\det L_j\le(1+\eta)(1+\delta)^d u(A_{x_j})$.
Choose $\Gamma_j$ with $\Gamma_j^\#=L_j$ and put $D_j=\det L_j/(2\pi)^d$. As in
Lemma~\ref{p3:lem:7}, set
\[
 U_{j,n}(u)=n^{-d/2}\sqrt{w(x_j)}\,(\chi_jp)(x_j+u/n).
\]
Choose a smooth multiplier $\theta_j$ compactly supported in
$\interior \Omega_j$, equal to one near the envelope.
For $g=\theta_j(D)U_{j,n}$, \cref{p3:lem:5,p3:lem:6} give
\[
 \|U_{j,n}-g\|_{H^k}
  +\bigg(\sum_{\gamma\in\Gamma_j}|(U_{j,n}-g)(\gamma)|^2\bigg)^{1/2}
 \le C_{j,k,m}n^{-m}\|p\|_\mu
\]
for integers $k>d/2$ and every fixed $m$. Here $D$ in the multiplier
denotes the angular Fourier variable. Apply \eqref{piv:latticeparseval}
to $g$ and use $\|U_{j,n}\|_2\le C_j\|p\|_\mu$. This gives
\[
 \left|D_j^{-1}\sum_{\gamma\in\Gamma_j}|U_{j,n}(\gamma)|^2
                  -w(x_j)\int\chi_j^2|p|^2\dd x\right|
            \le C_{j,m}n^{-m}\|p\|_\mu^2 .
\]
The nonzero summands are positive point evaluations at
$\lambda=x_j+\gamma/n$, with weights
$w(x_j)\chi_j(\lambda)^2/(n^dD_j)$. Summing over the finitely many
patches approximates $\int\chi_0^2|p|^2\dd \mu$ to additive error
$(\eta+o_\eta(1))\|p\|_\mu^2$, where $o_\eta(1)\to0$ with $n$ at a
fixed stage.

A fixed lattice of density $D_j$ contributes at most
$n^d D_j|Q_j^+|+O_j((n+1)^{d-1})$ points, by counting its fundamental
cells meeting the boundary of the dilated cube. Disjointness of the
base interiors and the choice of mesh give
\begin{align*}
 \sum_j D_j |Q_j^+|
 &\le \frac{1+C_d\eta}{(2\pi)^d}
                         \sum_j \int_{Q_j}u(A_x)\dd x\\
 &\le \frac{1+C_d\eta}{d!}
             \int_{\interior E}\frac{\dd\nu_E}{\Delta_{\mathrm{lat}}(A_x)} .
\end{align*}
For $0\le\varphi\le1$ of Lipschitz constant at most one, each patch
contributes at most its lattice count times $\sup_{Q_j^+}\varphi$.
This supremum is at most $\varphi(x)+\eta$ for every $x\in Q_j$.
Repeating the preceding sum with this factor bounds the normalized
interior count by $\int\varphi\,\Delta_{\mathrm{lat}}(A_x)^{-1}\dd\nu_E
+O_d(\eta)+o_\eta(1)$, uniformly over these tests. Here we used
$\Delta_{\mathrm{lat}}^{-1}\le d!$ and $N_n/n^d\to 1/d!$.
The complementary positive form $\int(1-\chi_0^2)|p|^2\dd \mu$ can be
large on individual unit polynomials. Lemma~\ref{p3:lem:8} bounds its trace by
$2\eta^4N_n$ for all sufficiently large $n$. Lemma~\ref{p3:lem:9} approximates
it to operator error $\eta$ by positive point evaluations at
$O(\eta^{-2}+\eta^2N_n)$ nodes. Add these evaluations and merge coincident
nodes by adding weights. The full operator has error at most
$2\eta+o_\eta(1)$.

Choose $\eta_j\downarrow0$ and increasing degree thresholds so that the
fixed-stage errors in the norm and the normalized count are at most
$\eta_j$ beyond the corresponding threshold, including the term
$O(\eta_j^{-2}/N_n)$. Use stage $j$ from its threshold to the next.
The correction contributes at most its normalized node count to each
of the tests. Thus the thresholds work simultaneously for all of them.
Scaling gives the bound for arbitrary nonnegative Lipschitz functions,
and uniform approximation gives \eqref{piv:upper} for nonnegative
continuous functions, since the total counts are bounded by the estimate
for $\varphi=1$. Finitely many initial rules may be chosen by positive
exact cubature.
\end{proof}

\subsection{Tangent measures and the lower estimate}\label{piv:measures}
For a body $K$ in cycles per unit length, the corresponding
bandlimited space is $\PW(2\pi K)$ in the notation of this paper.
Call a positive Radon measure $\sigma$ Parseval on that space if
\begin{equation}\label{piv:parseval}
 \int |f|^2\dd \sigma=\int |f(u)|^2\dd u\qquad(f\in\PW(2\pi K)).
\end{equation}
A window containing total weight $V$ on $m$ atoms has sum of squared
weights at least $V^2/m$. Thus to force many nodes it suffices to
bound that sum from above while keeping track of the mass in the window.
A nonnegative Fourier kernel supplies the upper bound. Its diagonal
sees the squared weights, and the Parseval identity controls its total
energy. This argument survives collisions of nodes and even a continuous
component in the limiting measure.
\begin{lemma}[A window count]\label{piv:window}
Let $K$ be a convex body, $P=(-1/2,1/2)^d$, $F\in\mathcal F(K)$ and
$0<s<1$. For $R\ge1$ and $W=z+RP$, every positive $\sigma$ satisfying
\eqref{piv:parseval} obeys
\begin{equation}\label{piv:windowcount}
 \#(\supp \sigma\cap W)
 \ \ge\ s^d\frac{\fc{F}(0)}{F(0)}R^d
             -C_{K,F,s}(R+1)^{d-1}.
\end{equation}
The constant is uniform in $z$ and $\sigma$.
\end{lemma}
\begin{proof}
Choose a Schwartz function in $\PW(2\pi K)$ nonzero at the origin. Testing its
translates in \eqref{piv:parseval}, followed by a finite covering, bounds
the $\sigma$-mass of each unit ball by a constant depending only on $K$.
In particular $\sigma$ is tempered. Polarization gives the distributional
identity
\begin{equation}\label{piv:gap}
 \fc{\sigma}=\delta_0\qquad\text{on }\interior (K-K).
\end{equation}
Indeed, $\sigma-du$ annihilates $f\overline g$ whenever $\fc{f}$ and
$\fc{g}$ are smooth and compactly supported in $\interior K$.
Every $t\in\interior (K-K)$ has the form $a-b$ with
$a,b\in\interior K$. Given a test function supported sufficiently
near $t$, translate it by $b$ and convolve it with a smooth approximate
point mass at $-b$. These convolutions are transforms of such products
and converge to the test function in $C_c^\infty$. A partition of unity
and the symmetry of $K-K$ give \eqref{piv:gap}.

Set $D=K-K$, $F_s(\xi)=F(\xi/s)$ and $H=\mathcal F_{\mathrm c}^{-1}F_s$. Then $H\ge0$,
$H(0)=s^d\fc{F}(0)$, and $|F_s|\le F(0)$. Choose $\epsilon>0$
such that $sD+\epsilon B^d$ is compactly contained in $\interior D$.
Take a nonnegative Schwartz function $g$ of integral one with
$\supp \fc{g}\subset\epsilon B^d$. The normalized
squared modulus of a Schwartz function with sufficiently small Fourier
support provides one. Put $h=\mathbf1_W*g$ and $\tau=h\sigma$. We have
$0\le h\le1$, $\int h=|W|$, $\int h^2\le |W|$. Testing \eqref{piv:gap}
against the corresponding modulates of $h$ gives
$\fc{\tau}=\fc{h}$ on $sD$ and $\tau(\R^d)=|W|$.
By Fourier inversion, the sign of $F_s$ off $sD$, and Plancherel,
\begin{equation}\label{piv:energy}
 \iint H(u-v)\dd \tau(u)\dd \tau(v)
   =\int F_s(\xi)|\fc{\tau}(\xi)|^2\dd \xi
   \le F(0)\int h^2\dd u\le F(0)|W|.
\end{equation}

The mass of $\tau$ outside $W$ is $O((R+1)^{d-1})$. To see this,
cover unit-width annuli about $y\in W$ by unit balls. The uniform mass
bound and Schwartz decay give
\[
 \int_{W^c}g(u-y)\dd \sigma(u)
 \le C_m(1+\operatorname{dist}(y,W^c))^{-m}
\]
for every fixed $m>1$. Integrate in $y$ over $W$. Each inner shell of
unit thickness has volume $O_d((R+1)^{d-1})$.

If $W$ meets the support in finitely many points, it contains only their
atoms. Nonnegativity of $H$ and $\tau$ allows us to keep all diagonal
terms at these points in \eqref{piv:energy}. Cauchy--Schwarz yields
\[
  H(0)\,\tau(W)^2
  \le \#(\supp \sigma\cap W)\,
         H(0)\sum_{u\in W}\tau(\{u\})^2
  \le \#(\supp \sigma\cap W)\,F(0)|W|.
\]
Substitution of the shell estimate proves \eqref{piv:windowcount}, with a
larger constant for bounded $R$. An infinite support count already
satisfies the assertion.
\end{proof}
\begin{proof}[Proof of the lower bound in Theorem~\ref{piv:main}]
Fix an unbounded sequence of degrees with $\|S_n-I\|\to0$. At each
interior $x$, \cref{p2:thm:tangent-sampling} gives vague precompactness of
\[
 \sigma_{n,x}=\frac{n^d}{w(x)}
       \sum_{\lambda\in\Lambda_n}\omega_{n,\lambda}\,
                         \delta_{n(\lambda-x)} .
\]
Every subsequential limit $\sigma$ is Parseval on $\PW(A_x)$, with cycles-frequency body $\widetilde A_x=A_x/(2\pi)$. Indeed,
on every sufficiently late tail the two sampling bounds are
$1-\epsilon$ and $1+\epsilon$. Apply \cref{p2:thm:tangent-sampling} on that tail and let
$\epsilon\downarrow0$.

For integer $R\ge1$ let
$b_{n,R}(x)=\#\{\lambda:n(\lambda-x)\in RP\}$.
For fixed $R,x$, pass to a subsequence realizing its lower limit and
then to a vague limit. Any finite list of distinct support points in
$RP$ admits disjoint relatively compact open neighborhoods there.
Positivity and vague convergence give a node in each neighborhood for
all sufficiently large indices. Thus the limiting support count is at
most $\liminf_n b_{n,R}(x)$. The constant in Lemma~\ref{piv:window}
is common to all possible limits at this $x$. First send integer $R$ to
infinity, next $s$ to one, and finally take the supremum over $F$. This gives
\begin{equation}\label{piv:localcount}
  \liminf_{R\to\infty}\liminf_n R^{-d}b_{n,R}(x)
  \ge \ell(\widetilde A_x)=q_E(x)/\Delta_{\mathrm{LP}}(A_x).
\end{equation}
Now take a relatively open $V\subset E$ and a compact
$J\subset V\cap\interior E$. For fixed $R$ and sufficiently large $n$,
every node counted over $J$ lies in $V$. Fubini gives
$\int_J b_{n,R}(x)\dd x\le(R/n)^d\#(\Lambda_n\cap V)$.
Fatou's lemma, first in $n$ and then in $R$, and
\eqref{piv:localcount} yield
\[
 \liminf_n\frac{\#(\Lambda_n\cap V)}{N_n}
 \ge d!\int_J\ell(\widetilde A_x)\dd x .
\]
Exhaust $V\cap\interior E$ and use $\nu_E(\partial E)=0$. Applying the
result to $V=\{x:\varphi(x)>t\}$, integrating in $t$, and using Fatou's
lemma once more proves \eqref{piv:lower}.
\end{proof}

\begin{remark}[Spherical moments and the Euclidean Fourier gap]\label{piv:rem:design-gap}
Gorbachev \cite[slide~8]{Gorbachev2026Seminar} recalls the zonal harmonic
criterion of Delsarte, Goethals, and Seidel
\cite[Theorem~5.5]{DelsarteGoethalsSeidel1977} for a spherical $s$-design $X$. Write $R_\ell$ for the zonal
polynomial of degree $\ell$ with $R_\ell(1)=1$. The criterion is
\[
 \sum_{u,v\in X} R_\ell(\langle u,v\rangle)=0,\qquad 1\le\ell\le s .
\]
By the spherical harmonic addition formula, each sum is a positive
multiple of the sum of squares of the degree-$\ell$ harmonic moments
of $X$. Thus the empirical probability measure and normalized area
agree on all harmonics through degree $s$. In the spherical design
linear program these identities remove the low positive degrees from
the pair sum, leaving the constant term and a tail with controlled sign
\cite[Theorem~5.10]{DelsarteGoethalsSeidel1977}
\cite[slide~14]{Gorbachev2026Seminar}. The identity
$\fc{\sigma}=\delta_0$ on $\interior(K-K)$ plays the same role in
\cref{piv:window}, with $K=A_x/(2\pi)$ in cycles per unit length.
It fixes the response on all interior difference frequencies before the
sign of the certificate controls the remaining frequencies.
Appendix~\ref{sph:sec:design-count} applies the window count on $\mathbb S^d$.
\end{remark}

\section{Optimal Ratios on Balls and Spherical Packing}\label{piv:balls}
The optimal cubature redundancy on the ball and sphere is of particular
interest, especially because of its connection with the spherical design
problem. An equal-weight rule on the sphere that is exact through
degree $s$ is a spherical $s$-design. Xu \cite[Corollary~4.2]{Xu2003} used the relation between ball and sphere
cubature to bound the node count on the disk with projected surface measure.
Gorbachev \cite{Gorbachev2014} sharpened the bound with a planar
Cohn--Elkies certificate. Building on the spherical design bounds of
Delsarte, Goethals, and Seidel \cite{DelsarteGoethalsSeidel1977} and
Yudin \cite{Yudin1997}, Gorbachev \cite{Gorbachev2026Seminar} used
Fourier--Gegenbauer expansions and the Euclidean packing LP to obtain
a stronger high-dimensional lower bound for spherical designs.

To see the role of the measure, let $D_s$ denote the node count of a
positive cubature rule of precision $s$ for
$\dd\mu(t)=(1-\norm{t}^2)^{-1/2}\dd t$ on the disk. Xu's lower bound is
$(2/j_{1,1}^2+o(1))s^2$, where $j_{1,1}$ is the first positive zero
of $J_1$. With $M_n=D_{2n}$ and $N_n=(n+1)(n+2)/2$, this gives
$\liminf M_n/N_n\ge16/j_{1,1}^2=1.08977\ldots$.
Normalize the disk and sphere measures to mass one and lift each node
$t$ to $(t,\pm\sqrt{1-\norm{t}^2})$, splitting its weight equally.
Averaging over the two lifts cancels odd powers of the last coordinate
and replaces its square by $1-\norm{t}^2$. Thus the lifted rule is
exact through degree $s$ on $\mathbb S^2$. If $S_s$ is the least node
count for positive spherical cubature of that precision, then $S_s\le2D_s$.
For $c_*=\liminf S_s/s^2$, Gorbachev's planar certificate gives
$c_*\ge1/(\pi r^2)$ with $2\pi r^2=7.25519\ldots<7.25520$
\cite[p.~17]{Gorbachev2014}. Hence
\begin{equation}\label{p2:eq:disk-cubature-gap}
 \liminf_{n\to\infty}\frac{M_n}{N_n}\ \ge\ 4c_*
 \ \ge\ \frac{8}{7.25520}=1.1026574\ldots .
\end{equation}
In the degree-$s$ normalization this gives $\liminf D_s/s^2\ge
c_*/2\ge 1/7.25520=0.137832\ldots$. The coefficient $0.144$ printed
at the end of \cite[Section~4, p.~24]{Gorbachev2014} would require a
stronger estimate than this certificate and lift provide.

The local spectra of the two polynomial spaces explain why Euclidean
ball packing also enters our sampling bounds. For spherical polynomials
of degree at most $n$ on $\mathbb S^d$, the local spectral body in tangent
coordinates is $B^d$ (\cref{sph:prop:kernel}). For total-degree polynomials
on the solid ball $E=B^d$, \cref{p1:subsec:ball} gives the ellipsoid
\[
 A_x=(I-xx^{\mathsf T})^{-1/2}B^d,\qquad
 q_E(x)=\frac{v_d}{(2\pi)^d\sqrt{1-\norm{x}^2}},\qquad
 v_d=\frac{\pi^{d/2}}{\Gamma(1+d/2)} .
\]
The radial semiaxis is $(1-\norm{x}^2)^{-1/2}$ and the transverse semiaxes
are one. The affine map stretches the ball and a packing cell by the
same volume factor. Both packing bounds in \cref{piv:main} are therefore
constant in $x$, giving for every admissible measure $\mu$
\[
 \frac1{\Delta_{\mathrm{LP}}(B^d)}\le
 \rho_{\mathrm{iso}}(B^d,\mu)\le
 \frac1{\Delta_{\mathrm{lat}}(B^d)}.
\]
For a radial Schwartz function $p$ normalized by $p(0)=\fc{p}(0)=1$,
with $\fc{p}\ge0$ and $p(y)\le0$ for $\norm{y}\ge r$, the rescaling
$F(\xi)=p(r\xi/2)$ gives
\begin{equation}\label{piv:cert}
 \rho_{\mathrm{iso}}(B^d,\mu)
       \ge \frac{(2/r)^d}{v_d}
       =\Gamma(1+d/2)\left(\frac8t\right)^{d/2},
       \qquad t=2\pi r^2 .
\end{equation}

Cohn and Elkies \cite[Section~7 and Table~4]{CohnElkies2003} give
$t=7.25520$ in dimension two and $t=8.19385$ in dimension three.
On the disk the same planar certificate as in \eqref{p2:eq:disk-cubature-gap}
thus gives the lower endpoint $8/7.25520$ for all positive, asymptotically
isometric arrays and every admissible measure, including area measure.
The triangular and face-centered cubic lattices give the upper endpoints
in dimensions two and three through their packing densities
$\pi/(2\sqrt3)$ and $\pi/(3\sqrt2)$, respectively.
In dimensions eight and twenty-four the sharp linear programming bounds
match the densities of the $E_8$ and Leech lattices
\cite{Viazovska2017,CohnEtAl2017}. Their reciprocals are the exact
sampling ratios in the table. In these two dimensions the arrays
constructed in \cref{piv:main} attain the stated ratios, and
\eqref{piv:empirical} gives the limiting node distribution $\nu_{B^d}$
for every array that attains them.

\begin{center}
\begin{tabular}{@{}rlll@{}}
\toprule
Dimension & Lower endpoint & Upper endpoint & Upper construction\\
\midrule
2  & $1.102657404344\ldots$ & $1.102657790843\ldots$ & triangular lattice\\
3  & $1.2824462330\ldots$   & $1.3504744742\ldots$   & face-centered cubic lattice\\
8  & \multicolumn{2}{c}{$384/\pi^4=3.9421371857\ldots$} & $E_8$ lattice\\
24 & \multicolumn{2}{c}{$12!/\pi^{12}=518.24902266\ldots$} & Leech lattice\\
\bottomrule
\end{tabular}
\end{center}
On the disk, choosing arbitrary positive weights and allowing all nodes to move can improve the leading ratio over the triangular-lattice construction by less than
$3.87\cdot10^{-7}$. The small interval comes from the Fourier certificate
in \eqref{piv:cert}. Optimality of triangular disk packing alone would
leave open how much an irregular sampling measure could save.
Whether positive exact cubature on the disk can attain the
triangular-lattice sampling ratio remains open.

On the unit ball, Corollary~\ref{p2:cor:smooth-physical} answers the simultaneous sampling and interpolation question posed in \cite{AMO2021} for nodes in the ball and all the weights $\dd\mu_a(t)=(1-\norm{t}^2)^{a-1/2}\dd t$, $a\ge0$, considered there. Indeed, \cite[Proposition~4]{AMO2021} gives $\sup_{t\in E} K_n(t,t)\lesssim n^{d+2a}$, which implies the Bernstein--Markov condition by the reproducing inequality. The density of each $\mu_a$ is positive and continuous near the origin. The earlier nonexistence result in \cite[Theorem~2]{AMO2021} excludes orthogonal reproducing-kernel bases at every sufficiently large degree. The corollary excludes uniform Riesz bounds across degrees and allows arbitrary positive sampling weights.

\section{Necessary Cubature Redundancy in Higher Dimensions and the Fuglede Conjecture}\label{piv:cubes}

On a centrally symmetric convex body, a positive cubature rule of degree
$2n$ with $N_n+o(N_n)$ nodes would force each local spectral body $A_x$
to admit an exponential orthonormal basis (\cref{p3:cor:isometric-equivalence}).
The solution of Fuglede's conjecture for convex bodies
\cite{LevMatolcsi2022} then requires each $A_x$ to tile $\R^d$ by
translations. We examine this condition first on cubes, then identify the
two planar classes that admit asymptotically isometric sampling at node
ratio one (\cref{piv:plane}). In dimensions $d\ge3$ the tiling condition
fails somewhere on every centrally symmetric convex body, giving the
necessary cubature redundancy in \cref{piv:cubature-gap}.

Write $Q_d=[-1,1]^d$ and $Q_d^*=\{\xi:\sum_j|\xi_j|\le1\}$.
This same cross-polytope appears in angle coordinates because the indices share the one
degree constraint $\norm{\alpha}_1\le n$. Indeed, with
$x_j=\cos\theta_j$, a product of Chebyshev polynomials of multidegree
$\alpha$ expands into exponentials with frequency vectors $(\pm\alpha_1,\ldots,\pm\alpha_d)$.
As shown in \cref{p1:subsec:cube}, the local spectral body is
\begin{equation}\label{piv:cubetangent}
 A_x=L_x Q_d^*,\qquad
 L_x=\operatorname{diag}\big((1-x_1^2)^{-1/2},\ldots,(1-x_d^2)^{-1/2}\big).
\end{equation}
The diagonal map is precisely the local change from Chebyshev angles to
physical coordinates. In particular,
$\dd\nu_{Q_d}=\prod_j\dd x_j/(\pi\sqrt{1-x_j^2})$. Every interior point sees
the same cross-polytope up to this change of coordinates.
\begin{corollary}[A change between dimensions two and three]\label{piv:cube}
For every admissible measure,
\[
 \rho_{\mathrm{iso}}(Q_2,\mu)=1,\qquad
 1<\rho_{\mathrm{iso}}(Q_3,\mu)\le\frac{19}{18}.
\]
More generally, $\rho_{\mathrm{iso}}(Q_d,\mu)>1$ for $d\ge3$, and
$1/\Delta_{\mathrm{lat}}(Q_d^*)$ is an upper bound.
\end{corollary}
\begin{proof}
The translates of $Q_2^*$ by the lattice generated by $(1,1),(1,-1)$
tile the plane, so \eqref{piv:upper} gives ratio one.
For $d\ge3$ every facet of $Q_d^*$ is a $(d-1)$-simplex and has no
center of symmetry. Greenfeld--Lev \cite{GreenfeldLev2017} prove that every facet of
a spectral convex polytope must be centrally symmetric. Thus $L^2(Q_d^*)$
has no exponential orthonormal basis. If $\rho_{\mathrm{iso}}=1$, the diagonal choice in \eqref{piv:rho}
would give an asymptotically isometric array with $M_n/N_n\to1$.
This contradicts \cref{p3:cor:isometric-equivalence} at $x=0$. The upper
estimate follows from \eqref{piv:cubetangent} and Theorem~\ref{piv:main}.
Minkowski's critical lattice for the regular octahedron has packing density
$18/19$ \cite{BetkeHenk2000}, giving the stated three-dimensional ratio.
\end{proof}

\subsection{The square and its Chebyshev weights}\label{p2:sec:square-redundancy}

The planar case of \cref{piv:cube} also appears in classical cubature
formulas. On the square $[-1,1]^2$ take the product Chebyshev measure
$\dd\mu_C=\pi^{-2}\prod_{j=1}^2(1-t_j^2)^{-1/2}\dd t_1\dd t_2$.
The positive formulas in \cite[Section~2.2, formulas (2.4)--(2.5)]{XuSquare1996},
including the minimal formula of Morrow and Patterson \cite{MorrowPatterson1978},
integrate through degree $2n+1$ with
$M_n=N_n+\lceil(n+1)/2\rceil$ nodes. In particular, $M_n/N_n\to1$.

The local diamond and its critical lattice were identified in \cref{p2:subsec:rectangle-lattices}.

For the Padua family in \cref{p2:subsec:rectangle-lattices}, positive
weights give the sharp squared-norm bound ratio $2$. The single mode
where exactness fails identifies the source of this factor.

For the displayed Padua family, the standard positive cubature weights are
\begin{equation}\label{p2:eq:padua-weights}
 \omega_\xi=\frac1{n(n+1)}
 \begin{cases}
  1/2,&\xi\text{ is a vertex},\\
  1,&\xi\text{ lies on an edge but is not a vertex},\\
  2,&\xi\text{ is an interior node}.
 \end{cases}
\end{equation}
The corresponding rule $Q_n$ is exact through degree $2n-1$ and on the
subspace of $\Pi_{2n}^2$ orthogonal to $T_{2n}(x_1)$ in
$L^2(\mu_C)$ \cite{Padua2007}. It satisfies the sharper norm identity
\begin{equation}\label{p2:eq:padua-sharp-mz}
 Q_n(|p|^2)=\norm{p}_{L^2(\mu_C)}^2+
 \left|\left\langle p,\sqrt2T_n(x_1)\right\rangle_{L^2(\mu_C)}\right|^2,
 \qquad p\in\Pi_n^2.
\end{equation}
Indeed, $Q_n(T_{2n}(x_1))=1$, while its integral is zero. In the
Chebyshev expansion of $|p|^2$, the coefficient of this sole exceptional
mode is the squared coefficient of $\sqrt2T_n(x_1)$ in $p$.
Consequently the optimal squared-norm bounds for these weights are
$a=1$ and $b=2$, attained by $1$ and $\sqrt2T_n(x_1)$, respectively.
Padua points therefore give complete stable interpolation on $\Pi_n^2$
without giving exact cubature through degree $2n$ for $\mu_C$.

\begin{remark}[Dependence on the weight]\label{p2:rem:weights}
If $c\,\dd\mu\le\dd\nu\le C\,\dd\mu$ on the entire domain, with fixed positive constants, the polynomial norms and Christoffel functions are comparable. Consequently uniform complete interpolation transfers between $\mu$ and $\nu$. For two measures with positive continuous densities in the interior, this criterion also requires control of their ratio near the boundary.
\end{remark}

\subsubsection{Zero-redundancy cubature and a common node distribution}
\label{p2:subsec:square-equidistribution}

Changing the target measure can make exact norm recovery possible at
exactly $N_n$ nodes. For the product Chebyshev probability measure of
the second kind,
\[
 \dd\mu_U=\frac4{\pi^2}\sqrt{1-x_1^2}\sqrt{1-x_2^2}\dd x_1\dd x_2,
\]
the Morrow--Patterson construction gives positive cubature exact on
$\Pi_{2n}^2$ with $N_n$ nodes \cite{MorrowPatterson1978,XuSquareSurvey2018}.
Its relation to Padua points is particularly direct
\cite{BeberokBialasCiezDeMarchi2026}:
\begin{equation}\label{p2:eq:morrow-patterson-padua}
 \mathrm{MP}_n=\mathrm{Pad}_{n+2}\cap(-1,1)^2,\qquad
 |\mathrm{MP}_n|=N_n.
\end{equation}
With the Padua family above, the weights are
\begin{equation}\label{p2:eq:morrow-patterson-weights}
 v_\xi=\frac8{(n+2)(n+3)}(1-\xi_1^2)(1-\xi_2^2),
 \qquad \xi\in\mathrm{MP}_n.
\end{equation}

For completeness, these statements follow directly from the Padua rule.
Put $m=n+2$ and $h(x)=(1-x_1^2)(1-x_2^2)$, so that
$\dd\mu_U=4h\dd\mu_C$. For $q\in\Pi_{2n}^2$, the polynomial
$hq$ has total degree at most $2m$ and degree in $x_1$ at most
$2n+2<2m$. Its coefficient in the exceptional mode $T_{2m}(x_1)$
is therefore zero, and the rule $Q_m$ integrates $hq$ exactly.
The boundary nodes contribute zero. Multiplying its interior weights
$2/[m(m+1)]$ by $4h$ yields
\eqref{p2:eq:morrow-patterson-weights}. There are $2m+1$ boundary
nodes in $\mathrm{Pad}_m$, so
$|\mathrm{MP}_n|=N_{n+2}-(2n+5)=N_n$.
In particular, the norm identity is exact on $\Pi_n^2$, and the
weighted kernel vectors form an orthonormal basis.

\begin{table}[htbp]
\centering
\small
\begin{tabular}{@{}lcccc@{}}
\toprule
Nodes & Measure & Node count & Exact through & Bound ratio $\kappa$\\
\midrule
Padua & $\mu_C$ & $N_n$ & $2n-1$ & $2$\\
Morrow--Patterson & $\mu_U$ & $N_n$ & $2n$ & $1$\\
First-kind square rules & $\mu_C$ & $N_n+O(n)$ & $2n+1$ & $1$\\
\bottomrule
\end{tabular}
\caption{Three square constructions, with the squared-norm bound ratio
on $\Pi_n^2$ for their standard positive weights. The first-kind rules
are those recalled at the start of Section~\ref{p2:sec:square-redundancy}.}
\label{p2:tab:square-rules}
\end{table}

Both Chebyshev measures are Bernstein--Markov, as follows for example
from polynomial bounds on the sup norms of their orthonormal product
bases. Their positive interior densities give the same local spectral
bodies. Moreover, the empirical node measures of all three constructions
in \cref{p2:tab:square-rules} converge to the equilibrium measure of
the square,
\[
 \dd\nu_{[-1,1]^2}
 =\frac{\dd x_1\dd x_2}{\pi^2\sqrt{1-x_1^2}\sqrt{1-x_2^2}}
 =\dd\mu_C.
\]
For Padua points this follows from
\cref{p2:cor:complete-equidistribution}; for the two exact cubature
constructions it follows from
\cref{p2:cor:cubature-equidistribution}. The Morrow--Patterson
cubature weights themselves converge to $\mu_U$. Thus a common
macroscopic node distribution is compatible with different target
measures and different exactness at critical cardinality.

\subsection{What an explicit cubature grid sees}\label{piv:grids}
The following family of grids exhibits the change from two to three
dimensions under the same choice of weight. Consider the total-degree
space $\Pi_n^d$ on $Q_d=[-1,1]^d$, with
product Chebyshev probability measure. Take $n=2r-1$. In each coordinate
mark the abscissae $t_j=\cos(j\pi/(2r))$, $0\le j\le2r$. Keep a grid point
if its indices are all even or all odd. There are $r+1$ even choices and
$r$ odd choices per coordinate. With suitable positive weights this is
one of the exact cubature rules of De Marchi, Vianello and Xu
\cite[Theorem~2.2]{DeMarchiVianelloXu2009}. Its node count and the dimension are
\[
 M_n=(r+1)^d+r^d,\qquad N_n=\binom{2r-1+d}{d}.
\]
On the square the two leading terms agree. On the three-dimensional cube the same rule
uses three nodes for every two coefficients in the limit.
\begin{figure}[htbp]
\centering
\includegraphics[width=.98\linewidth]{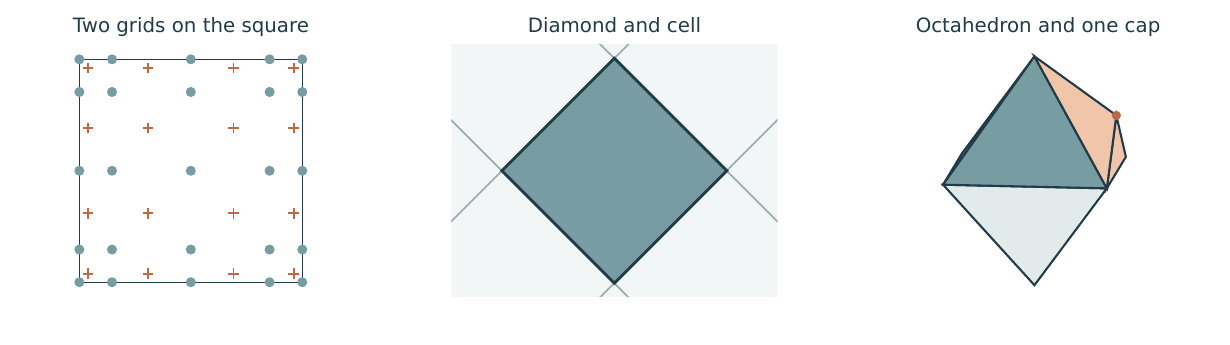}
\caption{Left, the all-even and all-odd nodes on the square for $r=4$.
Center, the diamond fills the reciprocal cell. Right, the octahedron
together with one of the eight tetrahedra needed to fill its reciprocal
cell. Each tetrahedron has volume $1/12$.}
\label{piv:figure}
\end{figure}

We can also see which frequencies the sampling grid identifies. In angle
coordinates on the torus, take the uniform average on the grid with
spacing $\pi/r$ and on its translate by $\pi(1,\ldots,1)/(2r)$, with
mass $1/2$ on each grid. Folding by the coordinate cosines gives exactly
the retained nodes and their positive weights. The first grid averages
$e^{ik\cdot\theta}$ to zero unless $k=2r\ell$. For such a $k$, the
translation multiplies the average by $(-1)^{\sum_j\ell_j}$. The two-grid
average is one on $2rD_d$ and zero elsewhere,
where
\[
 D_d=\{\ell\in\Z^d:\ \ell_1+\cdots+\ell_d\text{ is even}\}.
\]
The shortest nonzero $\ell^1$ length in $2rD_d$ is $4r$. This proves
exactness through degree $4r-1$, in particular for $|p|^2$ with
$p\in\Pi_{2r-1}^d$.

The lattice $D_d$ has covolume two. In two dimensions its Voronoi cell
is precisely the diamond $Q_2^*$, of area two. In three dimensions $Q_3^*$
is the octahedron with vertices $\pm e_1,\pm e_2,\pm e_3$. It has
volume $4/3$. The three-dimensional lattice still has a cell of volume
two. Its Voronoi cell extends past each triangular face of the octahedron.
For example, over the face joining $e_1,e_2,e_3$ it reaches
$(1/2,1/2,1/2)$. The tetrahedron between that face and the extra vertex
has volume $1/12$. Eight such tetrahedra account for all the unused
volume. The node ratio $3/2$ is visible in the cell as $2/(4/3)$.

The value $3/2$ belongs to this particular exact rule. For exact
trigonometric lattice cubature on the standard torus, Kuperberg
\cite[Section~3]{Kuperberg2006} expresses the leading term of the optimal
lattice node count in terms of the reciprocal lattice packing density of
the cross-polytope and explicitly uses the octahedral density $18/19$.
See also Cools and Lyness \cite{CoolsLyness2001} on trigonometric lattice
rules. On the algebraic cube $Q_3$, the octahedra are the local spectral
bodies of the total-degree polynomial space. The critical lattice gives
the optimal ratio $19/18$ among single sampling lattices in the tangent
model. By \cref{piv:cube} this is an upper bound for asymptotically
isometric sampling for every admissible $\mu$ on $Q_3$. The optimal
ratio for unrestricted arrays is strictly greater than one and its exact
value remains open. In dimension four, the explicit
packing of density $0.8248579\ldots$ recorded in \cite{Miller2025} gives
$\rho_{\mathrm{iso}}(Q_4,\mu)<1.213$. 

\subsection{Planar domains at critical density}\label{piv:planar}
The square belongs to a class that can be characterized completely.
\begin{corollary}[Planar classification]\label{piv:plane}
For a centrally symmetric planar convex body satisfying the standing
measure hypotheses, $\rho_{\mathrm{iso}}(E,\mu)=1$ if and only if $E$
is a parallelogram or a centrally symmetric hexagon.
\end{corollary}
\begin{proof}
The diagonal construction in \eqref{piv:rho} attains node ratio one whenever $\rho_{\mathrm{iso}}=1$. \Cref{p3:cor:isometric-equivalence} implies that
$A_0=E^*$ admits an orthonormal basis of exponentials. Iosevich, Katz
and Tao \cite{IosevichKatzTao2003} proved Fuglede's conjecture for
convex planar domains and classified the spectral convex bodies in the
plane. Applying their theorem to $E^*$ shows that it is a parallelogram
or a centrally symmetric hexagon. Polarity exchanges edges and vertices
and preserves central symmetry in the plane, so $E$ has the same
classification. Conversely, \cref{p3:cor:hexagon} constructs the required rules
for each of these two classes.
\end{proof}

\subsection{Necessary cubature redundancy in higher dimensions}\label{piv:all-dimensions}

Can a positive cubature rule exact through degree $2n$ use only
$N_n+o(N_n)$ nodes, where $N_n=\dim\Pi_n^d$? Gaussian quadrature on
an interval attains $N_n$ exactly. In dimensions $d\ge3$, every
centrally symmetric convex body requires a positive proportion of
additional nodes for all sufficiently large $n$
(\cref{piv:cubature-gap}). The obstruction already arises when the
sampling bounds merely tend to one.

Call a convex body $K$ \emph{spectral} if $L^2(K)$ admits an
orthonormal basis of normalized exponentials. By
\cref{p3:cor:isometric-equivalence}, attaining critical density requires
every interior $A_x$ to be spectral. The center alone can miss an
obstruction. For the regular octahedron $E=Q_3^*$, for example,
$A_0=Q_3$ is a cube and is spectral. We examine what happens to its
facets as the base point moves.

\begin{lemma}[Facets under displacement of the base point]\label{piv:displacement}
Suppose $d\ge3$ and the polar $E^*$ is a polytope with vertex pairs
$\{\pm v_j\}_{j=1}^m$. If the numbers $(u\cdot v_j)^2$ are pairwise
distinct, then for every sufficiently small $t>0$, no facet of $A_{tu}$
is centrally symmetric. The excluded directions $u$ lie in the finite
union of the hyperplanes $(v_i\pm v_j)^\perp$, $i<j$.
\end{lemma}
\begin{proof}
Fix such a direction $u$ and put $a_j=(u\cdot v_j)^2$. For $t$ near
zero, $tu\in\interior E$ and the finite-polar formula
\eqref{p1:eq:finite-polar} reads
\[
 A_{tu}=\conv\{\pm f_j(t)v_j:1\le j\le m\},
 \qquad f_j(t)=(1-a_jt^2)^{-1/2}.
\]
For four distinct indices $j_1,\ldots,j_4$, the corresponding functions
$f_{j_r}$ are linearly independent near zero. Indeed, if
$\sum_{r=1}^4 c_r f_{j_r}(t)=0$, comparing the coefficients of
$1,t^2,t^4,t^6$ gives $\sum_{r=1}^4 c_r a_{j_r}^k=0$ for $0\le k\le3$.
The Vandermonde matrix is invertible since the $a_{j_r}$ are distinct.
Consequently, for any signs $\epsilon_r\in\{-1,1\}$, the analytic vector
\[
 \sum_{r=1}^4 \epsilon_r f_{j_r}(t)v_{j_r}
\]
has a coordinate which is not identically zero, since each $v_{j_r}$ is nonzero.
The zeros of that coordinate are isolated. There are finitely many
index and sign choices, so on some common interval $0<t<t_0$ no four
vertices from distinct antipodal pairs have equal diagonal midpoints.

Every facet lies in a supporting hyperplane disjoint from the origin
and therefore contains at most one vertex from each ambient antipodal
pair. A centrally symmetric facet has dimension $d-1\ge2$, so its
vertices include two distinct pairs with the same midpoint, the center
of the facet. These four vertices would belong to four distinct
ambient pairs, which the preceding argument excludes.
\end{proof}

\begin{theorem}[Critical asymptotically isometric sampling]\label{piv:classification}
Under the standing hypotheses of Part~\ref{part:packing},
$\rho_{\mathrm{iso}}(E,\mu)=1$ exactly for intervals in dimension one
and for parallelograms and centrally symmetric hexagons in dimension
two. Every full-dimensional centrally symmetric convex body in
dimension $d\ge3$ satisfies $\rho_{\mathrm{iso}}(E,\mu)>1$.
\end{theorem}
\begin{proof}
On an interval the $(n+1)$-node Gaussian quadrature rule for $\mu$ has
positive weights and is exact through degree $2n+1$. It therefore gives
$S_n=I$ with $n+1=N_n$ nodes. The two-dimensional assertion is
\cref{piv:plane}.

Suppose $d\ge3$ and $\rho_{\mathrm{iso}}(E,\mu)=1$. The diagonal
choice following \eqref{piv:rho} gives an asymptotically isometric
array with $M_n/N_n\to1$, since injectivity gives $M_n\ge N_n$
eventually. By \cref{p3:cor:isometric-equivalence}, $A_x$ is spectral
at every interior point. At the origin $A_0=E^*$. In their proof of
Fuglede's conjecture for convex bodies, Lev and Matolcsi show that every
spectral convex body is a polytope \cite[Theorem~1.2]{LevMatolcsi2022}.
Hence $E^*$ is a polytope. Choose $u$ outside the hyperplanes in
\cref{piv:displacement}. For all sufficiently small $t>0$, the
polytope $A_{tu}$ has no centrally symmetric facet. Yet every facet of
a spectral convex polytope is centrally symmetric
\cite{GreenfeldLev2017}. This contradiction rules out ratio one.
\end{proof}

\begin{corollary}[A positive relative gap for cubature]\label{piv:cubature-gap}
Let $E\subset\R^d$ be a full-dimensional centrally symmetric convex
body, $d\ge3$, and let $\mu$ satisfy the standing measure hypotheses.
There are $\delta=\delta(E,\mu)>0$ and $n_0=n_0(E,\mu)$ such that,
for every $n\ge n_0$, the number $M_n$ of distinct nodes in any
positive cubature rule for $\mu$ on $E$ exact on $\Pi_{2n}^d$ satisfies
\[
 M_n\ge(1+\delta)N_n=(1+\delta)\binom{n+d}{d}.
\]
\end{corollary}
\begin{proof}
Exactness applied to $|p|^2$, $p\in\Pi_n^d$, gives sampling bounds
$a=b=1$. In particular, evaluation is injective, so $M_n\ge N_n$.
If the asserted gap failed, there would be degrees $n_j\to\infty$
and positive rules with $M_{n_j}/N_{n_j}\to1$. The subsequence
assertion of \cref{p2:thm:near-critical-transfer}, with $a=b=1$,
would give a normalized exponential orthonormal basis on every
interior $A_x$. By \cref{p3:cor:isometric-equivalence} this would
give $\rho_{\mathrm{iso}}(E,\mu)=1$, contrary to
\cref{piv:classification}.
\end{proof}

For fixed $d\ge3$ and precision $2n-1$, M\"oller's classical lower
bound \cite{Moller1979} adds a term of order $n^{d-1}$ to
$\dim\Pi_{n-1}^d$ (see also \cite[Chapter~5]{Xu2025Minimal}).
At the same precision, \cref{piv:cubature-gap} applied at $n-1$
requires at least $c\,n^d$ additional nodes for all sufficiently large $n$,
with $c>0$ independent of $n$.

\subsection{The three-variable cube question}\label{piv:cube-question}
For the product Chebyshev measure on $[-1,1]^3$, the cubature rules of
De Marchi, Vianello and Xu \cite{DeMarchiVianelloXu2009} use
$n^3/4+O(n^2)$ nodes at degree $2n-1$. Li, Sun and Xu
\cite[p.~121]{LiSunXu2009} recovered this count by discrete Fourier
analysis and asked whether $n^3/6+O(n^2)$ nodes could suffice.
Let $M_n$ count the nodes of a rule with their degree convention.
Exactness through degree $2n-1$ includes exactness on $\Pi_{2(n-1)}^3$.
Applying \cref{piv:cubature-gap} at $n-1$ therefore gives, for any
positive such rule supported on the cube and all sufficiently large $n$,
\[
 M_n\ \ge\ (1+\delta)\dim\Pi_{n-1}^3
          \ =\ (1+\delta)\binom{n+2}{3}.
\]
Thus the answer to their existence question is negative for positive
rules with nodes in the cube. The value of the optimal leading
coefficient remains undetermined. In particular, proving the
coefficient $1/4$ suggested in \cite[p.~121]{LiSunXu2009} would
require a quantitative strengthening for exact rules.

\subsection{Gaussian cubature outside the class of centrally symmetric convex bodies}\label{piv:gaussian-outside}
The next example, due to Berens, Schmid and Xu \cite{BerensSchmidXu1995},
attains the dimension bound at every degree $2n-1$ in three variables.
It is the family cited by Li, Sun and Xu immediately after their cube
question \cite[p.~121]{LiSunXu2009}.

Let $\Omega=F([-1,1]^3)$, where
\[
 F(x,y,z)=(x+y+z,\;xy+yz+zx,\;xyz),
\]
and let $\nu$ be the image under $F$ of product Lebesgue measure on the
cube. If $p\in\Pi_{2n-1}^3$, then $p\circ F$ has degree at most $2n-1$ in
each coordinate, since each component of $F$ is affine in each variable
separately. The tensor product of three $n$-point Gauss--Legendre rules
therefore integrates $p\circ F$ exactly. Two product nodes have the same
image precisely when their coordinates differ by a permutation: the
unordered coordinates are recovered as the roots of
\[
 t^3-u_1t^2+u_2t-u_3,\qquad u=F(x,y,z).
\]
The images are indexed by triples $1\le i\le j\le k\le n$. Summing
the positive product weights for each triple gives a positive rule on
$\Omega$ with exactly $\binom{n+2}{3}=\dim\Pi_{n-1}^3$ distinct nodes.

The support $\Omega$ is full dimensional, since the Jacobian of $F$ is
nonzero at triples of distinct interior coordinates. Its geometry differs
from that of the cube. For $0<a<1$, the two points $F(a,a,a)$ and
$F(-a,-a,-a)$ correspond to $(t-a)^3$ and $(t+a)^3$. Their midpoint
corresponds to $t^3+3a^2t$, whose nonreal roots exclude that point from
$\Omega$. The region also lacks central symmetry. The unique extrema of
$u_1$ force a putative center to be the midpoint of $F(1,1,1)$ and
$F(-1,-1,-1)$, namely $(0,3,0)$. The coordinate $u_2$ ranges from $-1$
to $3$, as one sees by evaluating its multiaffine expression at the
vertices of the cube, so its midpoint would have to be $1$.
This three-dimensional support admits exact Gaussian cubature despite
the positive relative gap on every centrally symmetric convex body in
the same dimension.

\section[A Redundancy Gap between Asymptotically Isometric Sampling and Exact Cubature]{A Redundancy Gap between Asymptotically Isometric \mbox{Sampling} and Exact Cubature}\label{piv:square-gap}

On the square $Q=[-1,1]^2$ every local spectral body is a diamond that
tiles the plane. The planar classification in \cref{piv:plane} therefore
allows asymptotically isometric sampling at node ratio one. To see what
exactness can cost on the same square, we perturb the product Chebyshev
probability measure $\mu_C$ of Section~\ref{p2:sec:square-redundancy}
by a positive analytic multiplier $h$.

Write $\dd\mu=h\dd\mu_C$ and $N_n=\dim\Pi_n^2=\binom{n+2}{2}$.
Let $m_n(\kappa)$ denote the least number of nodes in a positive-weight
sampling rule with squared-norm bound ratio at most $\kappa$.
At $\kappa=1$ a common rescaling of the weights gives exact cubature
through degree $2n$.

\begin{theorem}[Asymptotically isometric sampling and exact cubature]
\label{p2:thm:square-redundancy}
There is a strictly positive real analytic multiplier $h$ on a
neighborhood of $Q$, arbitrarily close to $1$ uniformly on $Q$, such
that $\mu=h\mu_C$ admits positive sampling rules with $M_n/N_n\to1$
and squared-norm bounds $1-\varepsilon_n,1+\varepsilon_n$, where
$\varepsilon_n\to0$. For this same measure the minimum node counts
satisfy
\begin{equation}\label{p2:eq:square-exact-gap}
 m_n(\kappa)=N_n+O_{h,\kappa}(n)\quad(\kappa>1),\qquad
 \liminf_{n\to\infty}\frac{m_n(1)}{N_n}\ge\frac43 .
\end{equation}
\end{theorem}

\begin{proof}
\emph{Sampling with a vanishing norm error.}
This part holds for every strictly positive continuous $h$ on $Q$.
Choose a positive polynomial $q$ of total degree at most $k$ with
$\|q/h-1\|_{L^\infty(Q)}\le\varepsilon<1$. Such an approximation
exists because $\min_Q h>0$. Put $s=n+\lceil k/2\rceil$. The positive
square rules of Morrow--Patterson and Xu in
Section~\ref{p2:sec:square-redundancy} are exact for $\mu_C$ through
degree $2s+1$ with $N_s+O(s)$ nodes. Denote their nodes and weights
by $z_j$ and $v_j$. Since $q|p|^2$ has degree at most $2n+k$, we
obtain, simultaneously for all $p\in\Pi_n^2$,
\begin{equation}\label{p2:eq:square-multiplier-rule}
 (1-\varepsilon)\norm{p}_\mu^2
 \le \sum_j v_jq(z_j)|p(z_j)|^2
 =\int_Q q|p|^2\dd\mu_C
 \le (1+\varepsilon)\norm{p}_\mu^2 .
\end{equation}
For fixed $\kappa>1$ choose $\varepsilon$ so that
$(1+\varepsilon)/(1-\varepsilon)\le\kappa$. The degree $k$ is then
fixed as $n$ grows, and $N_s+O(s)=N_n+O_{h,\kappa}(n)$. Injectivity
of evaluation gives the lower bound $N_n$. To make the norm error
tend to zero as well, take $k_n=\lfloor\sqrt n\rfloor$ and polynomial
approximants of degree at most $k_n$ converging uniformly to $h$.
They are eventually positive, their relative errors tend to zero,
and $s/n\to1$. This gives the first assertion.

\smallskip
\noindent\emph{Exactness and the number of independent moments.}
Normalize $\mu$ to mass one. For $a,b\ge0$ let
$\phi_{ab}(x,y)=T_a(x)T_b(y)$, with $T_a(\cos\theta)=\cos(a\theta)$.
The products with $a+b\le 2n$ form an orthogonal basis of $\Pi_{2n}^2$
in $L^2(\mu_C)$. Set $D_n=\dim\Pi_{2n}^2$ and omit the constant
mode, whose integral is already fixed. The $D_n-1$ remaining moments
of an exact $M$-node rule must be of the form
\begin{equation}\label{p2:eq:square-moment-map}
 \left(\sum_{j=1}^{M}\omega_j\phi_{ab}(z_j)\right)_{1\le a+b\le 2n},
 \qquad z_j\in Q,\quad \omega_j\ge0,\quad \sum_{j=1}^{M}\omega_j=1 .
\end{equation}
The nodes supply $2M$ real parameters and the weights supply $M-1$.
After eliminating the last weight, the map in
\eqref{p2:eq:square-moment-map} is polynomial on a compact subset of
$\R^{3M-1}$, hence Lipschitz there. If $3M<D_n$, its image has
Lebesgue measure zero in $\R^{D_n-1}$. This is the dimension
argument for finite atomic moment representations of di Dio and
Schm\"udgen \cite{diDioSchmuedgen2022}.

We now realize moments outside these images for every $n$ with one
analytic multiplier. Choose independent variables $U_{ab}$ uniformly
distributed on $[-1,1]$, one for each $a+b\ge1$. For a constant
$c>0$, put
\[
 c_{ab}=c\,e^{-(a+b)^2}U_{ab},\qquad
 h=1+\sum_{a+b\ge1} c_{ab}\phi_{ab}.
\]
Since $|T_a|\le1$ on $[-1,1]$, we have
$\|h-1\|_{L^\infty(Q)}
\le c\sum_{\ell\ge1}(\ell+1)e^{-\ell^2}$.
Choosing $c$ small makes this bound arbitrarily small and ensures
$h>0$ on $Q$. On each compact subset of $\C^2$, the Chebyshev
recurrence gives $|\phi_{ab}|\le C^{a+b}$ for some $C$ depending on
the compact set. The series therefore extends holomorphically to
$\C^2$. Orthogonality gives $\int_Q h\dd\mu_C=1$ and
\[
 \int_Q\phi_{ab}\,h\dd\mu_C
   =c_{ab}\norm{\phi_{ab}}_{L^2(\mu_C)}^2,\qquad a+b\ge1 .
\]
Thus the $D_n-1$ nonconstant moments have a joint Lebesgue density
for each $n$. With probability one their vector avoids the image of
\eqref{p2:eq:square-moment-map} for every $M$ with $3M<D_n$.
There are only countably many pairs $(n,M)$, so we can fix one
realization of $h$ for which $3m_n(1)\ge D_n$ for every $n$.
Finally $D_n/N_n=(2n+1)(2n+2)/[(n+1)(n+2)]\to4$, which gives
\eqref{p2:eq:square-exact-gap}.
\end{proof}

For $\mu_C$ itself the positive square formulas give
$m_n(1)/N_n\to1$. Multiplying $\mu_C$ by the $h$ of
\cref{p2:thm:square-redundancy} preserves the hypotheses of the bulk
limit. Indeed, each normalized Chebyshev product has supremum norm
at most $2$ on $Q$, so $\|p\|_{L^\infty(Q)}\le
2\sqrt{N_n}\,\|p\|_{L^2(\mu_C)}$. The positive lower bound on $h$
then gives Bernstein--Markov control for $\mu$, and its area density
is positive and continuous in the interior. The two measures have the
same local spectral bodies and normalized bulk kernel limits. The
multiplier changes each nonconstant Chebyshev moment by an
independently chosen amount. A uniform
approximation to $h$ controls the sampling error for every unit
polynomial by \eqref{p2:eq:square-multiplier-rule}. Exact cubature
must reproduce each of these moments at its prescribed value,
however small the change. Thus a positive fraction of extra nodes
can be required globally even when every tangent space admits
exact sampling at critical density.

\clearpage
\section*{Problems for Further Investigation}
\addcontentsline{toc}{section}{Problems for Further Investigation}

\paragraph{1. Local spectra beyond central symmetry.}
For which convex bodies $E$ without central symmetry do the
total-degree spaces $\Pi_n^d$, rescaled about an interior point $x$,
have a Paley--Wiener limit $\PW(A_x)$ as in \cref{p1:thm:main}?
The Baran dual unit ball is a natural candidate for $A_x$.
The bulk kernel on the simplex already gives an ellipsoidal example
\cite[Theorem~1.8]{KrooLubinsky2013}. If such a local model exists, does
the criterion of \cref{p3:cor:existence-equivalence} extend as well:
is uniformly stable positive sampling with $N_n+o(N_n)$ nodes equivalent
to the existence of exponential Riesz bases on all $A_x$ with common
normalized bounds? The Gaussian rule of Berens, Schmid and Xu in
\cref{piv:gaussian-outside} attains the dimension bound in three variables
on a nonconvex, non-centrally symmetric support. Keeping convexity
while relaxing central symmetry would isolate the role of the latter
hypothesis in the redundancy problem.

\paragraph{2. Stable sampling on the three-dimensional cube.}
Can positive sampling on $Q_3=[-1,1]^3$ attain $M_n=N_n+o(N_n)$
with a fixed finite ratio of sampling bounds? Asymptotically isometric
sampling on $Q_3$ requires a positive relative excess of nodes
(\cref{piv:cube}). Every interior spectrum is a linear image of the
regular octahedron $Q_3^*$. Thus \cref{p3:cor:existence-equivalence}
reduces the fixed-conditioning question to a single Fourier problem:
does $L^2(Q_3^*)$ admit a Riesz basis of exponentials?

\paragraph{3. Exact positive cubature on the regular hexagon.}
Let $H$ be the regular hexagon with area measure. Can positive cubature
rules supported on $H$ be exact on $\Pi_{2n}^2$ with $M_n=N_n+o(N_n)$?
Every local spectrum on $H$ tiles the plane, and
\cref{p3:cor:contrasting-domains} gives arrays at this cardinality with
$\|S_n-I\|\to0$. Exact cubature asks for $S_n=I$ at each degree, with
the same area measure and the same asymptotic node count. The analytic
weight in \cref{p2:thm:square-redundancy} shows on the square that
this change can increase the leading node count.

\paragraph{4. Can irregular samples on the disk improve on the triangular lattice?}
Take the unit disk $B^2$ with area measure. Its local spectra are
ellipses, so affine images of the triangular lattice give an
asymptotically isometric array with limiting node ratio $2\sqrt{3}/\pi$.
Is $\rho_{\mathrm{iso}}(B^2,dx)=2\sqrt{3}/\pi$, or can irregularly
placed nodes with freely chosen positive weights attain a smaller ratio?
The planar lower bound in \eqref{piv:cert} lies less than $3.87\times10^{-7}$
below the lattice value and applies to all such arrays. The remaining
question is whether the optimal lattice density also fixes the sampling
cost when no lattice is prescribed.

\clearpage
\appendix
\section{Higher-Order Chebyshev-Ridge Recovery}\label{p2:app:recovery}

Fix $x\in\interior E$ and write the function to be recovered as
\[
 f(u)=\int b(\xi)e^{iu\cdot\xi}\dd\xi,
 \qquad b=(2\pi)^{-d}\widehat f\in C_c^\infty(\interior A_x).
\]
We start with the ridge polynomials of \cref{p1:sec:recovery}. Their leading
term recovers $f$. To improve the accuracy to order $n^{-J}$, we expand the
error and show that each coefficient has spectral support in $\supp b$.
We can then use the same ridge directions to cancel successive coefficients.
After constructing the polynomials, we compute their rescaled norm at centers
where $\mu$ has a positive continuous density.

\paragraph{Finite directions and a fixed lift.}
Radial approximation of points of $\partial E^*$ gives the closure identity
\begin{equation}\label{p2:eq:app-closed-hull}
 A_x=\overline{\conv}\{g_x(y):y\in\interior E^*\}.
\end{equation}
Choose a compact $K\Subset\interior A_x$ containing $\supp b$.
By \cref{p1:lem:dense-ridge}, we can
choose $y_1,\ldots,y_r\in\interior E^*$ such that
\[
 K\Subset\interior P,\qquad
 P=\conv\{0,g_x(y_1),\ldots,g_x(y_r)\}.
\]
We may add directions so that the first $d$ are linearly independent. Set
\[
 B\lambda=\sum_{j=1}^r\lambda_jg_x(y_j),\qquad
 \Delta_r=\{\lambda\in\R^r:\lambda_j>0,\ \textstyle\sum_j\lambda_j<1\}.
\]
The surjective map $B$ sends the open simplex $\Delta_r$ onto
$\interior P$. We now choose a lift that will apply to every correction
amplitude supported in $K$. Take a right inverse $S$ of $B$ and an
isomorphism $R:\R^{r-d}\to\ker B$. In the coordinates
$T(\xi,\eta)=S\xi+R\eta$ we have $BT(\xi,\eta)=\xi$. Write
$J_T=|\det T|$. For each $\xi\in K$, some $T(\xi,\eta)$ lies in
$\Delta_r$. By openness, a product neighborhood of $(\xi,\eta)$ has
closure in $T^{-1}\Delta_r$. Compactness gives finitely many such
products $U_\ell\times V_\ell$ with the $U_\ell$ covering $K$.
Choose smooth functions $\psi_\ell$ whose sum is one near $K$,
with $\supp\psi_\ell\Subset U_\ell$, and functions
$\chi_\ell\in C_c^\infty(V_\ell)$ with integral one. For $c\in C_K^\infty$, define
\[
 (\mathcal Lc)(T(\xi,\eta))
 =J_T^{-1}c(\xi)\sum_\ell\psi_\ell(\xi)\chi_\ell(\eta).
\]
The products of cutoffs have support in a fixed compact set
$K_\lambda\Subset\Delta_r$ in the $\lambda$ coordinates. Integration
in $\eta$, using $\int\chi_\ell=1$, gives
\begin{equation}\label{p2:eq:app-lift}
 B_*((\mathcal Lc)(\lambda)\dd\lambda)=c(\xi)\dd\xi,
 \qquad B(\supp\mathcal Lc)\subseteq\supp c.
\end{equation}
The formula also bounds every smooth seminorm of $\mathcal Lc$ by finitely
many seminorms of $c$, so the lift is continuous. When $r=d$, integration
in the fiber amounts to evaluation at its single point.

\paragraph{A polynomial operator and its principal phase.}
Put $\theta_j(t)=\arccos(y_j\cdot t)$, $\Theta=(\theta_1,\ldots,\theta_r)$, and $a=\mathcal Lc$. Define
\begin{equation}\label{p2:eq:app-polynomial-operator}
 (\mathcal S_nc)(t)=2^rn^{-r}
 \sum_{\substack{k\in\Z^r\\a(k/n)\ne0}}
 a(k/n)e^{ik\cdot\Theta(x)}\prod_{j=1}^rT_{k_j}(y_j\cdot t).
\end{equation}
Since $K_\lambda$ is compactly contained in the open simplex, each
contributing $k$ satisfies $k_j\ge1$ and
$\sum_jk_j\le(1-\eta)n$ for a fixed $\eta>0$. The product in
\eqref{p2:eq:app-polynomial-operator} has total degree at most this
sum, so $\mathcal S_nc\in\Pi_n^d$. Let
$F_a(z)=\int a(\lambda)e^{i\lambda\cdot z}\dd\lambda$.
Write $T_{k_j}(y_j\cdot t)=\cos(k_j\theta_j(t))$. Expanding the product
of cosines gives one lattice sum for each choice of signs. Poisson
summation on each of these lattices yields
\begin{equation}\label{p2:eq:app-poisson}
 (\mathcal S_nc)(t)=
 \sum_{\sigma\in\{-1,1\}^r}\sum_{\ell\in\Z^r}
 F_a\bigl(n[\Theta(x)+\sigma\Theta(t)-2\pi\ell]\bigr),
\end{equation}
where $\sigma\Theta(t)$ denotes coordinatewise multiplication.
Only the all-minus sign with $\ell=0$ can give a vanishing phase.
Indeed, $y_j\in\interior E^*$ puts $\theta_j(E)$ in
$[\epsilon_j,\pi-\epsilon_j]$ for some $\epsilon_j>0$. A plus sign
therefore gives a coordinate $\theta_j(x)+\theta_j(t)$ uniformly
separated from $2\pi\mathbb Z$. On the all-minus branch each difference
lies in $[-\pi,\pi]$, so every nonzero lattice shift stays uniformly
away from the origin. For large $\norm{\ell}$ the shifted phase has
norm bounded below by a constant times $\norm{\ell}$, uniformly in $t$.
The Schwartz decay of $F_a$ thus makes the sum over these terms
$O_{M,c}(n^{-M})$ for every $M$. We obtain
\begin{equation}\label{p2:eq:app-principal}
 (\mathcal S_nc)(t)=F_a\bigl(n[\Theta(x)-\Theta(t)]\bigr)
       +O_{M,c}(n^{-M}),\qquad t\in E.
\end{equation}
By continuity of the lift and the common compact support of its images,
these constants are uniform for $c$ in bounded subsets of $C_K^\infty$.

\paragraph{A global Taylor expansion in a real parameter.}
For $u\in n(E-x)$ and $0\le h\le n^{-1}$ set
\[
 \Phi_h(u)=-\int_0^1D\Theta(x+shu)u\dd s.
\]
Convexity keeps $x+shu$ in $E$ throughout the integral.
For $h>0$ the fundamental theorem of calculus gives
$\Phi_h(u)=[\Theta(x)-\Theta(x+hu)]/h$, and at $h=0$ we have
$\Phi_0(u)=B^{\mathsf T}u$. To control the Taylor remainder on the
whole set $n(E-x)$, observe that
\[
 \Phi_{h,j}(u)=(y_j\cdot u)
 \int_0^1\frac{\dd s}{\sqrt{1-(y_j\cdot x+shy_j\cdot u)^2}}.
\]
The integral is positive and at least one. Since the directions span $\R^d$, there is $c_0>0$ such that
\begin{equation}\label{p2:eq:app-phase-estimates}
 \norm{\Phi_h(u)}\ge c_0\norm{u},\qquad
 \norm{\partial_h^m\Phi_h(u)}\le C_m\norm{u}^{m+1}
 \quad(m\ge0).
\end{equation}
For the derivative estimate, each differentiation under the integral
introduces a further factor of $u$. The derivatives of $\Theta$ are
bounded on a neighborhood of $E$ because all the $y_j$ lie in
$\interior E^*$.

Set $H_a(h,u)=F_a(\Phi_h(u))$. For $j\ge1$, each term in its $j$th
$h$ derivative contains a derivative of $F_a$ and $q$ factors of the
form $\partial_h^{m_i}\Phi_h$, with $m_i\ge1$ and $\sum_i m_i=j$.
These factors contribute at most $C(1+\norm{u})^{j+q}$, with $q\le j$.
The lower bound on $\norm{\Phi_h(u)}$ makes each derivative of $F_a$
at this point decay faster than any power of $\norm{u}$. Thus for every
$j,L\ge0$,
\[
 \sup_{\substack{u\in n(E-x)\\0\le h\le n^{-1}}}
 (1+\norm{u})^L|\partial_h^jH_a(h,u)|\le C_{j,L,a}.
\]
Taylor's theorem at $h=0$, evaluated at $h=n^{-1}$, gives
\begin{equation}\label{p2:eq:app-taylor}
 F_a(\Phi_{1/n}(u))=
 \sum_{j=0}^{J-1}n^{-j}H_{j,a}(u)+O_{J,L,a}\bigl(n^{-J}(1+\norm{u})^{-L}\bigr),
 \qquad H_{j,a}=\frac{\partial_h^jH_a(0,\cdot)}{j!}.
\end{equation}
Because $\norm{u}\le n\,\operatorname{diam}E$, choosing $M>J+L$ in
\eqref{p2:eq:app-principal} gives the same weighted remainder for $\mathcal S_nc$.

\paragraph{Spectral support and recursive correction.}
We next identify the spectral amplitudes of the Taylor coefficients.
At $h=0$, each $\partial_h^m\Phi_h(u)$ is a homogeneous polynomial of
degree $m+1$ in $u$, as follows by differentiating its integral formula.
Derivatives of $F_a$ insert powers of $\lambda$ in its Fourier integral.
Since $\Phi_0(u)=B^{\mathsf T}u$, each coefficient $H_{j,a}$ is a finite
linear combination of terms
\[
 u^\alpha\int\lambda^\beta a(\lambda)e^{iu\cdot B\lambda}\dd\lambda.
\]
In the coordinates $T(\xi,\eta)$, the pushforward
$B_*(\lambda^\beta\mathcal Lc\dd\lambda)$ has density
\[
 c_\beta(\xi)
 =c(\xi)\sum_\ell\psi_\ell(\xi)
     \int T(\xi,\eta)^\beta\chi_\ell(\eta)\dd\eta .
\]
Thus $c_\beta$ is the product of $c$ with a fixed smooth function.
Integration by parts in $\xi$ gives
\[
 u^\alpha\int c_\beta(\xi)e^{iu\cdot\xi}\dd\xi
   =i^{|\alpha|}\int\partial_\xi^\alpha c_\beta(\xi)e^{iu\cdot\xi}\dd\xi.
\]
It follows that $H_{j,a}$ has spectral amplitude $L_jc$, where $L_j$
is a finite sum of derivatives of smooth multiples of $c$.
In particular, $L_j:C_K^\infty\to C_K^\infty$ is continuous and
$\supp L_jc\subseteq\supp c$. At order zero, \eqref{p2:eq:app-lift}
gives $L_0c=c$. Substitution into \eqref{p2:eq:app-taylor} yields
\begin{equation}\label{p2:eq:app-expansion}
 (\mathcal S_nc)(x+u/n)=
 \int\left(\sum_{j=0}^{J-1}n^{-j}L_jc\right)(\xi)e^{iu\cdot\xi}\dd\xi
 +O_{J,L,c}\bigl(n^{-J}(1+\norm{u})^{-L}\bigr).
\end{equation}
Because the same $L_j$ apply to every amplitude in $C_K^\infty$, we can
cancel the error at order $n^{-m}$ by feeding its negative back into
the polynomial operator. More precisely, set
\[
 b_0=b,\qquad b_m=-\sum_{j=1}^mL_jb_{m-j},\qquad
 P_n^{[J]}=\sum_{m=0}^{J-1}n^{-m}\mathcal S_nb_m.
\]
All $b_m$ remain in $C_K^\infty$, and $P_n^{[J]}\in\Pi_n^d$.
Apply \eqref{p2:eq:app-expansion} to $b_m$ through order $J-m$.
The constant coefficient is $b$, while the coefficient of $n^{-r}$, $1\le r<J$, is
$b_r+\sum_{j=1}^rL_jb_{r-j}=0$. The remainders, after their prefactors, are all of order $n^{-J}$. Hence
\begin{equation}\label{p2:eq:app-recovered-estimate}
 \sup_{u\in n(E-x)}(1+\norm{u})^L
 |P_n^{[J]}(x+u/n)-f(u)|\le C_{J,L}n^{-J}.
\end{equation}
For fixed $J$ this is one sequence valid for every $L$.
At this order all constants depend on finitely many seminorms of $b$
and are therefore uniform on bounded subsets of $C_K^\infty$.

\paragraph{Localization and the exact norm.}
To estimate the polynomials on the physical domain, use the combined
amplitude $\mathcal L(\sum_{m<J}n^{-m}b_m)$ in
\eqref{p2:eq:app-principal}. It has fixed compact support and uniformly
bounded smooth seminorms. Taking $u=n(t-x)$ and $h=1/n$ in
\eqref{p2:eq:app-phase-estimates} gives
$\norm{\Theta(x)-\Theta(t)}\ge c_0\norm{t-x}$. The Schwartz bound for
the principal term and the estimate for the remaining phases then give,
for every $M$,
\begin{equation}\label{p2:eq:app-packet-decay}
 |P_n^{[J]}(t)|\le C_{J,M}
 \bigl[(1+n\norm{t-x})^{-M}+n^{-M}\bigr],\qquad t\in E.
\end{equation}
Suppose now that $\mu$ is finite and equals $w(t)\dd t$ near $x$,
with $w$ continuous and $w(x)>0$. Choose an open ball $V$ centered at
$x$ whose closure lies in this neighborhood and in $\interior E$.
On $E\setminus V$, \eqref{p2:eq:app-packet-decay} gives
$n^d\int_{E\setminus V}|P_n^{[J]}|^2\dd\mu=O(n^{d-2M})$.
Inside $V$, change variables $t=x+u/n$. On each fixed ball
$\norm{u}\le R$, the rescaled integrand converges uniformly to
$w(x)|f(u)|^2$ by \eqref{p2:eq:app-recovered-estimate} and continuity
of $w$. On the rest of $n(V-x)$ its integral is bounded by
\[
 C\int_{\norm{u}>R}(1+\norm{u})^{-2M}\dd u+C n^{d-2M}.
\]
Choose $2M>d$. Taking $n\to\infty$ and then $R\to\infty$ gives
\[
 n^d\norm{P_n^{[J]}}_{L^2(\mu)}^2
 \longrightarrow w(x)\int_{\R^d}|f(u)|^2\dd u.
\]
This completes the proof of Theorem~\ref{p2:thm:recovery}.\hfill$\square$

\clearpage
\section{A Phase Grid on the Regular Hexagon}
\label{p2:sec:hexagon-phase}

On the regular hexagon, the spectral body changes from a hexagon in the
central disk to a parallelogram in each corner. We construct a continuous
phase map whose inverse takes an ordinary lattice to a node array with
orthonormal exponential bases in both kinds of bulk limit and along their
common interface. For $n=6m$, the resulting node count lies between the dimensions of two
consecutive polynomial spaces and is asymptotically equivalent to both.

Li--Sun--Xu \cite{LiSunXu2008} develop discrete Fourier analysis for
trigonometric spaces on hexagons. Here the functions sampled on the
physical hexagon are ordinary algebraic polynomials of bounded total degree.
Phase coordinates will specify the nodes and their local limits.

We use the normalization
\[
 v_1=(1,0),\qquad v_2=(1/2,\sqrt3/2),\qquad
 v_3=(-1/2,\sqrt3/2),\qquad
 E=\{x:|v_j\cdot x|\leq1,\ 1\leq j\leq3\}.
\]
Thus $v_2=v_1+v_3$ and the incircle has radius one.
This is a rotated and dilated version of the hexagon in
\cref{p1:sec:hexagon}. If $H$ denotes that circumradius-one hexagon and
$\mathcal R_\theta$ is rotation through $\theta$, then
\[
 E=LH,\qquad L=\frac{2}{\sqrt3}\mathcal R_{-\pi/6},\qquad
 A^{E}_{Lx}=L^{-T}A^{H}_x.
\]
The last identity follows directly from the polar-body formula for the
local spectrum and identifies the coordinates used in the two hexagon figures. Put
\[
 t_j=v_j\cdot x,\qquad \theta_j=\arcsin t_j,\qquad
 d_j=(1-t_j^2)^{1/2},\qquad g_j=\nabla\theta_j=v_j/d_j.
\]
All differential formulas below concern the interior of $E$. Formula \eqref{p1:eq:hexagon-spectral-body}, in the present normalization, gives
\begin{equation}\label{p2:eq:hex-phase-body}
 A_x=\operatorname{conv}\{\pm g_1,\pm g_2,\pm g_3\}.
\end{equation}

\subsection{A phase map with the correct local Jacobian}

Let $C_0=\{x:\norm{x}\leq1\}$, and let $C_1$ be the closed corner between
the unit-circle arc from $v_1$ to $v_2$ and the two sides meeting at
$(1,1/\sqrt3)$.  Define
\begin{align}
 F_0(x)&=(\theta_1+\theta_2,\theta_2+\theta_3),
     &&x\in C_0, \label{p2:eq:hex-central-phase}\\
 F_1(x)&=(\theta_1+\theta_2,\theta_2-\theta_1+\pi/3),
     &&x\in C_1. \label{p2:eq:hex-corner-phase}
\end{align}
Let $R$ denote physical rotation through $\pi/3$ and put
\[
 S=\begin{pmatrix}1&-1\\1&0\end{pmatrix}.
\]
The remaining five corner formulas are
$F(R^kx)=S^kF_1(x)$ for $x\in C_1$ and $0\leq k<6$. On $C_0$ set $F=F_0$.

\begin{proposition}\label{p2:prop:hex-phase-map}
The preceding formulas define a homeomorphism $F:E\longrightarrow D$,
where $D$ is the simple, nonconvex dodecagon whose cyclically ordered
vertices are $\pi/3$ times
\[
\begin{gathered}
 (2,0),(3,1),(2,2),(2,3),(0,2),(-1,2),\\
 (-2,0),(-3,-1),(-2,-2),(-2,-3),(0,-2),(1,-2).
\end{gathered}
\]
The map is real analytic in each open bulk region, and
\begin{equation}\label{p2:eq:hex-phase-jacobian}
 \det DF(x)=|A_x|>0
\end{equation}
there.  In particular,
\begin{equation}\label{p2:eq:hex-phase-total-area}
 |D|=2\pi^2,\qquad \int_E |A_x|\,dx=2\pi^2.
\end{equation}
\end{proposition}

\begin{proof}
The identity
\[
 (d_1+d_2+d_3)(-d_1+d_2+d_3)
 (d_1-d_2+d_3)(d_1+d_2-d_3)=3(1-\norm{x}^2)
\]
follows from $t_2=t_1+t_3$ and
$t_1^2+t_2^2+t_3^2=3\norm{x}^2/2$.
At the origin the three triangle inequalities for $d_1,d_2,d_3$ are
strict. The displayed product stays positive in the open disk, so none of
the inequalities can become an equality there. To see their geometric
meaning, use $d_3g_3=d_2g_2-d_1g_1$. In the basis $(g_1,g_2)$ the sum
of the absolute coordinates of $g_3$ is $(d_1+d_2)/d_3$. Thus $g_3$ is
outside the convex hull of $\{\pm g_1,\pm g_2\}$ exactly when
$d_3<d_1+d_2$. The relations obtained by solving for $g_1$ and $g_2$
give the other two tests. All six points in \eqref{p2:eq:hex-phase-body}
are therefore extreme throughout the open disk.  In the open corner $C_1$,
$t_1,t_2\geq|t_3|$, so $d_3$ is the largest of the three $d_j$.
The identity therefore gives $d_3>d_1+d_2$, and $\pm g_3$ are redundant.
At its incircle interface,
$d_3=d_1+d_2$.

In the disk, the shoelace formula gives
\[
 |A_x|=\det(g_1,g_2)+\det(g_1,g_3)+\det(g_2,g_3)
       =\det(g_1+g_2,g_2+g_3)=\det DF_0(x).
\]
In the corner, $A_x=\operatorname{conv}\{\pm g_1,\pm g_2\}$, so
\[
 |A_x|=2\det(g_1,g_2)
       =\det(g_1+g_2,g_2-g_1)=\det DF_1(x).
\]
Since $\det R=\det S=1$, the same Jacobian identity holds in the five
rotated corners.

To prove injectivity of $F_0$, use $(t_1,t_3)$ as coordinates and choose
$\psi'=\arcsin$.  The map $F_0$ is the gradient, in these coordinates, of
\[
 \Psi(t_1,t_3)=\psi(t_1)+\psi(t_3)+\psi(t_1+t_3).
\]
Since $\psi''(t)=(1-t^2)^{-1/2}>0$ for $|t|<1$, this Hessian is
positive definite in the open disk in physical coordinates (an ellipse in
$(t_1,t_3)$ coordinates). The open segment between any two distinct points
of the closed disk lies in its interior. Along that segment the component
of $\nabla\Psi$ in its direction is strictly increasing, and the gradient
is continuous at the endpoints. It follows that $F_0$ is injective on
the closed disk.

On the arc $x=(\cos\varphi,\sin\varphi)$, $0\leq\varphi\leq\pi/3$,
\[
 \theta_1=\pi/2-\varphi,\qquad
 \theta_2=\pi/6+\varphi,\qquad
 \theta_3=\varphi-\pi/6.
\]
Both phase formulas therefore equal $(2\pi/3,2\varphi)$.
Moreover $F_0(Rx)=SF_0(x)$, because rotation transforms
$(\theta_1,\theta_2,\theta_3)$ into $(-\theta_3,\theta_1,\theta_2)$.
The six arcs map to the six sides of the hexagon $H_0$ with vertices
\[
 (2\pi/3,0),\ (2\pi/3,2\pi/3),\ (0,2\pi/3)
 \quad\hbox{and their negatives}.
\]
Since $F_0$ is continuous and injective on the closed disk, the Jordan
curve theorem identifies its image with $H_0$.

The map $F_1$ is injective as a composition of three injective maps,
namely $x\mapsto(t_1,t_2)$, the componentwise arcsine and the affine map
in \eqref{p2:eq:hex-corner-phase}.
Its three boundary pieces map to the sides of the triangle $T_0$ with
vertices
\[
 (2\pi/3,0),\quad (2\pi/3,2\pi/3),\quad(\pi,\pi/3).
\]
It therefore maps $C_1$ onto $T_0$.  Each triangle $S^kT_0$ lies outside
$H_0$ in the cone spanned by its two adjacent central vertices. Thus the
six triangles have disjoint interiors and attach to the respective sides of $H_0$.
The formulas agree on all shared boundary pieces, giving the asserted
homeomorphism onto their union.  Its outer boundary is the listed
dodecagon.  Finally $|H_0|=4\pi^2/3$ and $|T_0|=\pi^2/9$.
Change of variables on the seven open regions proves
\eqref{p2:eq:hex-phase-total-area}, since the intervening arcs have area zero.
\end{proof}

\subsection{The critical lattices generated by the phase map}

\begin{lemma}\label{p2:lem:hex-lattice-cell}
At every point in a smooth bulk region, $A_x$ tiles spectral space, up
to null sets, by the lattice $DF(x)^T\mathbb Z^2$.  Consequently
\begin{equation}\label{p2:eq:hex-critical-tangent-lattice}
 \Gamma_x=2\pi DF(x)^{-1}\mathbb Z^2
\end{equation}
is a critical sampling and interpolation lattice for $\PW(A_x)$, and
the normalized exponentials indexed by $\Gamma_x$ form an orthonormal
basis of $L^2(A_x)$.
\end{lemma}

\begin{proof}
In the disk put $P=g_1+g_2$ and $Q=g_2+g_3$.  The three pairs of opposite
edges of $A_x$ differ by translations $P$, $Q$, and $Q-P$, respectively,
and $|A_x|=\det(P,Q)$.

These are precisely the opposite-edge translations in the proof of
\cref{p3:cor:hexagon}. The determinant identity there shows that
$P\mathbb Z+Q\mathbb Z$ tiles spectral space by translates of $A_x$.

In the corner, set $P=g_1+g_2$ and $Q=g_2-g_1$.  Then
$A_x=[P\ Q][-1/2,1/2]^2$, which is a fundamental parallelogram for
$P\mathbb Z+Q\mathbb Z$.  The rotated corners are identical.
In each case $[P\ Q]=DF(x)^T$.
Folding a measurable lattice fundamental domain onto its standard
parallelogram is a unitary identification of the corresponding $L^2$
spaces and preserves dual-lattice exponentials.  The Fourier basis on the parallelogram therefore gives the claimed
orthonormal basis on $A_x$. Fourier transformation identifies its coefficients
with the samples of functions in $\PW(A_x)$, up to the fixed normalization.
\end{proof}

Consider the finite phase arrays
\begin{equation}\label{p2:eq:hex-phase-array}
 \Lambda_n^{\mathrm{ph}}
   =F^{-1}\bigl((2\pi/n)\mathbb Z^2\cap D\bigr).
\end{equation}
At a fixed point $x$ in a smooth region, both $F$ and its local inverse
are differentiable. Taylor expansion in the two coordinate systems shows
that the possible limits of $n(\Lambda_n^{\mathrm{ph}}-x)$ are translates of
\eqref{p2:eq:hex-critical-tangent-lattice}. The translation is determined
by the limit of $-nF(x)/(2\pi)$ modulo $\mathbb Z^2$.

\subsection{The incircle interface preserves local completeness}

Across the incircle the derivative of $F$ changes in its second row.
In tangent coordinates this shifts each lattice row by an amount depending
only on the row. Such shifts preserve the Fourier basis property on a square.

\begin{proposition}\label{p2:prop:hex-interface}
Let $x_*$ be on the open incircle arc separating $C_0$ and $C_1$.
Every tangent $\Gamma$ of \eqref{p2:eq:hex-phase-array} at $x_*$ makes
$\{|A_{x_*}|^{-1/2}e^{i\gamma\cdot\xi}:\gamma\in\Gamma\}$
an orthonormal basis of $L^2(A_{x_*})$.
\end{proposition}

\begin{proof}
At $x_*$, $d_3=d_1+d_2$.  Put
\[
 P=g_1+g_2,\qquad Q=g_2-g_1,\qquad
 s=\frac{d_2}{d_1+d_2}.
\]
Then $A_{x_*}=[P\ Q][-1/2,1/2]^2$ and the two one-sided derivatives are
\[
 DF_1=\begin{pmatrix}P^T\\Q^T\end{pmatrix},\qquad
 DF_0=\begin{pmatrix}P^T\\(Q+sP)^T\end{pmatrix}.
\]
The identity for the second row follows from
$d_3g_3=d_2g_2-d_1g_1$.
Use the coordinates
\[
 z_1=P\cdot u/(2\pi),\qquad z_2=Q\cdot u/(2\pi).
\]
The first component of $F$ is constant on the arc and increases into
$C_1$. Thus the tangent to the arc is $z_1=0$, with the disk on the
side $z_1<0$. Write $u=u(z)$ for the inverse linear coordinate map.
On every bounded set of $z$, the rescalings
\[
 \Phi_n(z)=\frac{n}{2\pi}\bigl(F(x_*+u(z)/n)-F(x_*)\bigr)
 \quad\hbox{converge uniformly to}\quad
 L(z)=\bigl(z_1,z_2+s\min\{z_1,0\}\bigr).
\]
Indeed, each phase formula is smooth up to this open arc, so its Taylor
remainder is uniform on bounded rescaled sets. The rescaled circular arc
converges to its tangent line. On the strip between the arc and that line
the two linear formulas differ by $s z_1=o(1)$.
Subtracting $s\min\{z_1,0\}$ from the second coordinate gives a continuous
inverse for $L$.

Pass to a subsequence on which the offsets of $-nF(x_*)/(2\pi)$ modulo
$\mathbb Z^2$ converge to $\tau=(\tau_1,\tau_2)$, choosing representatives
$\tau_n\to\tau$. Every limit of bounded node coordinates then lies in
$L^{-1}(\mathbb Z^2+\tau)$. To obtain the converse, fix
$v\in\mathbb Z^2+\tau$, put $z=L^{-1}(v)$ and take any sufficiently small
disk centered at $z$. The point $v$ has positive distance from the image
of its boundary under $L$, and that image has winding number one about
$v$. Uniform convergence on the boundary preserves this winding number
for $\Phi_n$ about $v_n=v+\tau_n-\tau$. Hence $v_n$ lies in the image
of the disk under $\Phi_n$ for all sufficiently large $n$. Each such
preimage is a rescaled node, and it is unique because $F$ is injective.
As the disk can be taken arbitrarily small, the preimages converge to $z$.
This establishes the full tangent set, including points on the interface.

Applying $L^{-1}$ to the translated lattice gives
\begin{equation}\label{p2:eq:hex-interface-row-shift}
 \{(k+\tau_1,\ell+a_k):k,\ell\in\mathbb Z\},
 \qquad a_k=\tau_2-s\min\{k+\tau_1,0\}.
\end{equation}
For any real sequence $(a_k)$ the exponentials
\[
 \exp\!\big(2\pi i[(k+\tau_1)t_1+(\ell+a_k)t_2]\big),
 \qquad k,\ell\in\mathbb Z,
\]
form an orthonormal basis of the unit square $[-1/2,1/2]^2$. Integration
in $t_1$ gives orthogonality for distinct $k$, and for fixed $k$ the
$t_2$ factors form a modulated Fourier basis. If a function is orthogonal
to the entire family, its $k$th modulated Fourier coefficient in $t_1$ is
orthogonal to this $t_2$ basis for every $k$. Thus all of these coefficient
functions vanish, and completeness in $t_1$ gives the result.
Under $\xi=[P\ Q]t$ the square becomes $A_{x_*}$ and
$\gamma\cdot\xi=2\pi(z_1t_1+z_2t_2)$ for a tangent node $\gamma=u(z)$.
The factor $|A_{x_*}|^{-1/2}$ normalizes the transported basis.
\end{proof}

For moving centers approaching a point of the open arc at finite rescaled
normal distance, the interface line acquires an offset. The two linear
branches still have the same first row and differ by a shear in the second,
so the resulting set has the same row form with modified shifts $a_k$.
At the endpoints of the arc the physical boundary also enters the scaling.

\subsection{Cardinality and a reference measure}

\begin{proposition}\label{p2:prop:hex-phase-count}
For $n=6m$, $m\ge1$, the array
\eqref{p2:eq:hex-phase-array} has cardinality
\begin{equation}\label{p2:eq:hex-phase-count}
 M_n=\frac{n^2}{2}+n+1,\qquad
 \dim\Pi_{n-1}^2<M_n<\dim\Pi_n^2.
\end{equation}
\end{proposition}

\begin{proof}
The polygon $nD/(2\pi)$ has integer vertices when $n=6m$.
Its area is $n^2/2$.  At $m=1$ each of its twelve edges has a primitive
integer displacement, so the boundary lattice-point count at general
$m$ is $12m=2n$.  Pick's theorem gives the total number of lattice
points, including boundary points, as area plus half the boundary
count plus one.  Since $F$ is injective, this is exactly $M_n$.
Finally $\dim\Pi_{n-1}^2=n(n+1)/2$ and
$\dim\Pi_n^2=(n+1)(n+2)/2$, which give the strict inequalities.
\end{proof}

In particular, along $n=6m$,
\[
 \frac{M_n}{\dim\Pi_{n-1}^2}\longrightarrow1,\qquad
 \frac{M_n}{\dim\Pi_n^2}\longrightarrow1.
\]
To formulate sampling and interpolation for this phase grid, fix the
reference measure
\begin{equation}\label{p2:eq:hex-phase-reference-measure}
 d\mu_*(x)=\frac{|A_x|}{2\pi^2}\,dx.
\end{equation}
The Jacobian and area formulas in Proposition~\ref{p2:prop:hex-phase-map}
show that $\mu_*$ is a probability measure and that $F_\#\mu_*$ is
normalized Lebesgue measure on $D$. Its density is positive and
continuous in the interior of $E$.

We verify the Bernstein--Markov condition for this choice of measure.
Each $y\in\partial E^*$ lies on the segment from $0$ to $g_x(y)$, so
$E^*\subseteq A_x$ and $\mu_*$ dominates a constant multiple of area
measure. Choose a fixed disk $B\subset\interior E$. For $p\in\Pi_n^2$,
$n\geq1$, let $x_0$ maximize $|p|$ and rotate the phase of $p$ so that
$p(x_0)=\norm{p}_E$. For each $y\in B$ the real polynomial
$q_y(s)=\operatorname{Re}p((1-s)x_0+sy)$ has supremum norm at most
$\norm{p}_E$ on $[0,1]$. Markov's inequality on this interval gives
$\norm{q_y'}_\infty\leq 2n^2\norm{p}_E$. With $s=1/(4n^2)$ it follows
that $|p|\geq\norm{p}_E/2$ on the entire disk $(1-s)x_0+sB$, whose area
is $s^2|B|$. Integration gives $\norm{p}_E\leq C(1+n)^2\norm{p}_{L^2(\mu_*)}$,
also valid for constants since $\mu_*$ is a probability measure.

\clearpage
\begin{figure}[!htbp]
 \centering
 \includegraphics[width=\textwidth]{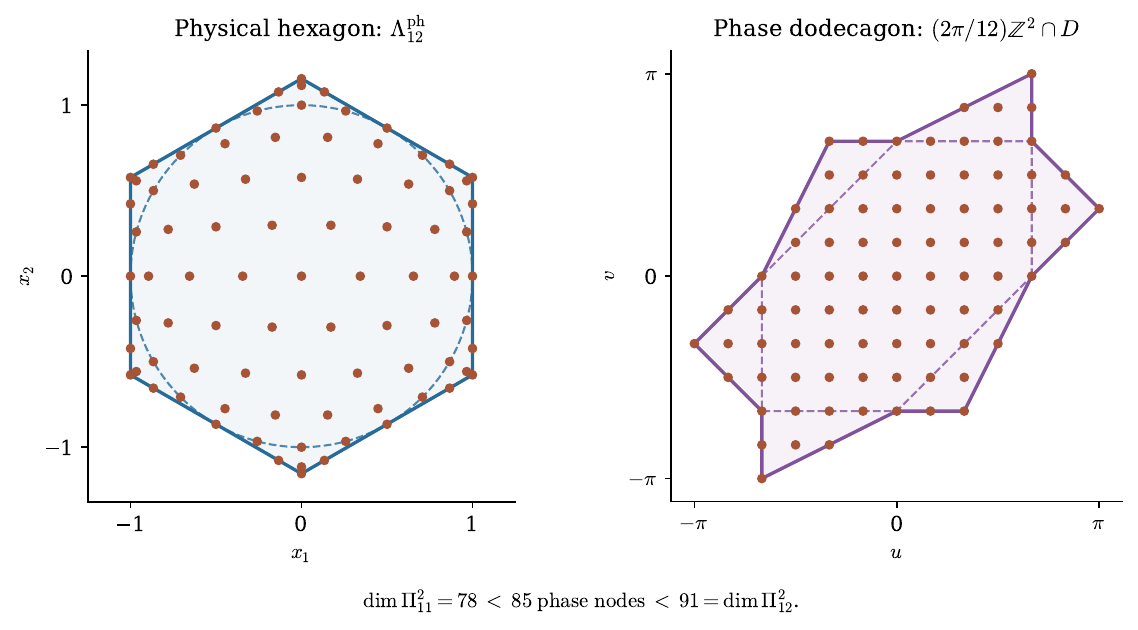}
 \caption{The unshifted phase grid for $n=12$ and its physical preimage.
 The dashed incircle maps to the dashed central hexagon. There are
 $85$ nodes, between $\dim\Pi_{11}^2=78$ and $\dim\Pi_{12}^2=91$.}
 \label{p2:fig:hexagon-phase-grid}
\end{figure}

\begin{problem}[Sampling and interpolation on the hexagonal phase grid]
\label{p2:prob:hex-phase-sampling}
For $n=6m$, does the phase grid $\Lambda_n^{\rm ph}$ provide uniformly
stable sampling for $\Pi_{n-1}^2$ and uniformly stable interpolation
for $\Pi_n^2$ in $L^2(\mu_*)$, with the Christoffel normalization of
the respective polynomial spaces?
\end{problem}

\clearpage
\section{Numerical Experiments on the Hexagon}
\label{app:hexagon-numerics}

We illustrate the position-dependent universality limit on the
circumradius-one hexagon $H$ of Section~\ref{p1:sec:hexagon}, with Lebesgue
measure $d\mu_0=dx$. The finite-dimensional space is the full total-degree
space $\Pi_{64}^2$, of dimension $2145$. For each fixed interior point $x$,
the centered kernel profile and its predicted limit are
\begin{equation}\label{eq:hexagon-numerical-profile}
 R_{n,x}(z)=
 \frac{K_n^{\mu_0}(x+z/(2n),x-z/(2n))}{K_n^{\mu_0}(x,x)},
 \qquad
 \Phi_x(z)=\frac{1}{|A_x|}\int_{A_x}e^{iz\cdot\xi}\,d\xi.
\end{equation}
Both equal one at $z=0$. The denominator is evaluated at the fixed base
point $x$, as in the universality theorem.

\begin{figure}[!htbp]
 \centering
 \includegraphics[width=0.64\textwidth]{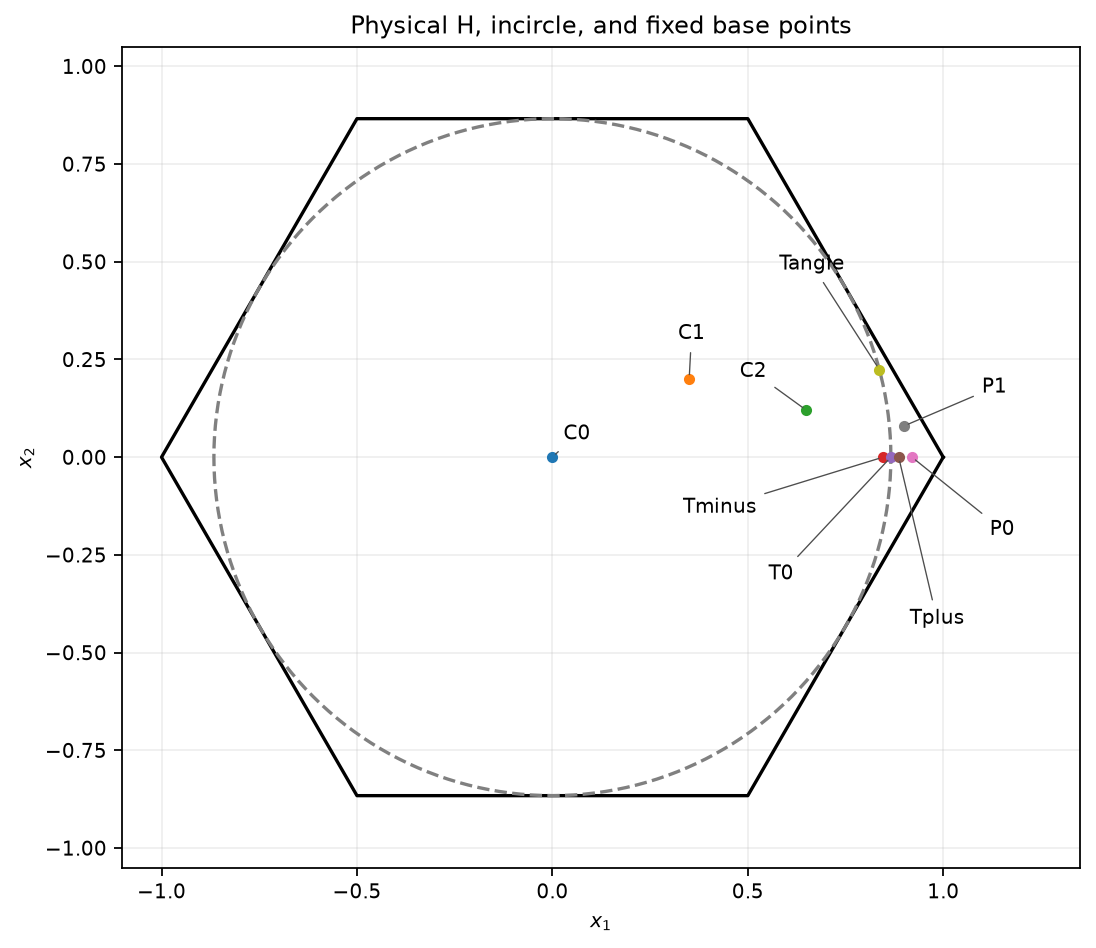}
 \caption{The circumradius-one physical hexagon, its incircle, and the
 nine fixed interior base points used in the numerical comparisons.
 The coordinates are also given in the corresponding profile panels.}
 \label{fig:hexagon-numerical-centers}
\end{figure}

Figure~\ref{fig:hexagon-numerical-centers} locates the nine base points.
The incircle has radius $r_* = \sqrt3/2$. In particular,
$\mathrm{T0}=(r_*,0)$ and
$\mathrm{Tangle}=r_*(\cos(\pi/12),\sin(\pi/12))$ lie on its interior arcs,
while $\mathrm{Tminus}=(r_*-0.02,0)$ and
$\mathrm{Tplus}=(r_*+0.02,0)$ straddle the transition. All points are
strictly inside $H$; this is a change of spectral shape within the bulk.

The maps use the uniform $41\times41$ displacement grid
$\mathcal G=\{-2+0.1j:0\le j\le40\}^2$. At $n=64$, both physical
arguments in \eqref{eq:hexagon-numerical-profile} remain inside $H$ for
every displayed point and displacement. The signed error is
$R_{64,x}-\Phi_x$, and the reported value
\[
 E_\infty(x)=\max_{z\in\mathcal G}|R_{64,x}(z)-\Phi_x(z)|
\]
is the maximum absolute error on this grid.

\clearpage
The polynomial kernel was computed from an orthonormal basis constructed
by degree-block multivariate Arnoldi orthogonalization. Integration on
$H$ used six triangles joining the origin to its edges and positive
tensor Gauss--Legendre quadrature of order $66$ on each triangle map.
This integrates the required polynomial products exactly in exact
arithmetic. The spectral formula was used to compute the comparison
target independently of the polynomial basis. Rebuilding the quadrature
with orders $70$ and $74$ at five representative centers changed the
normalized profiles by less than $10^{-7}$ on a separate $21\times21$
check grid, well below the discrepancies displayed here.

\begin{figure}[!htbp]
 \centering
 \includegraphics[width=\textwidth]{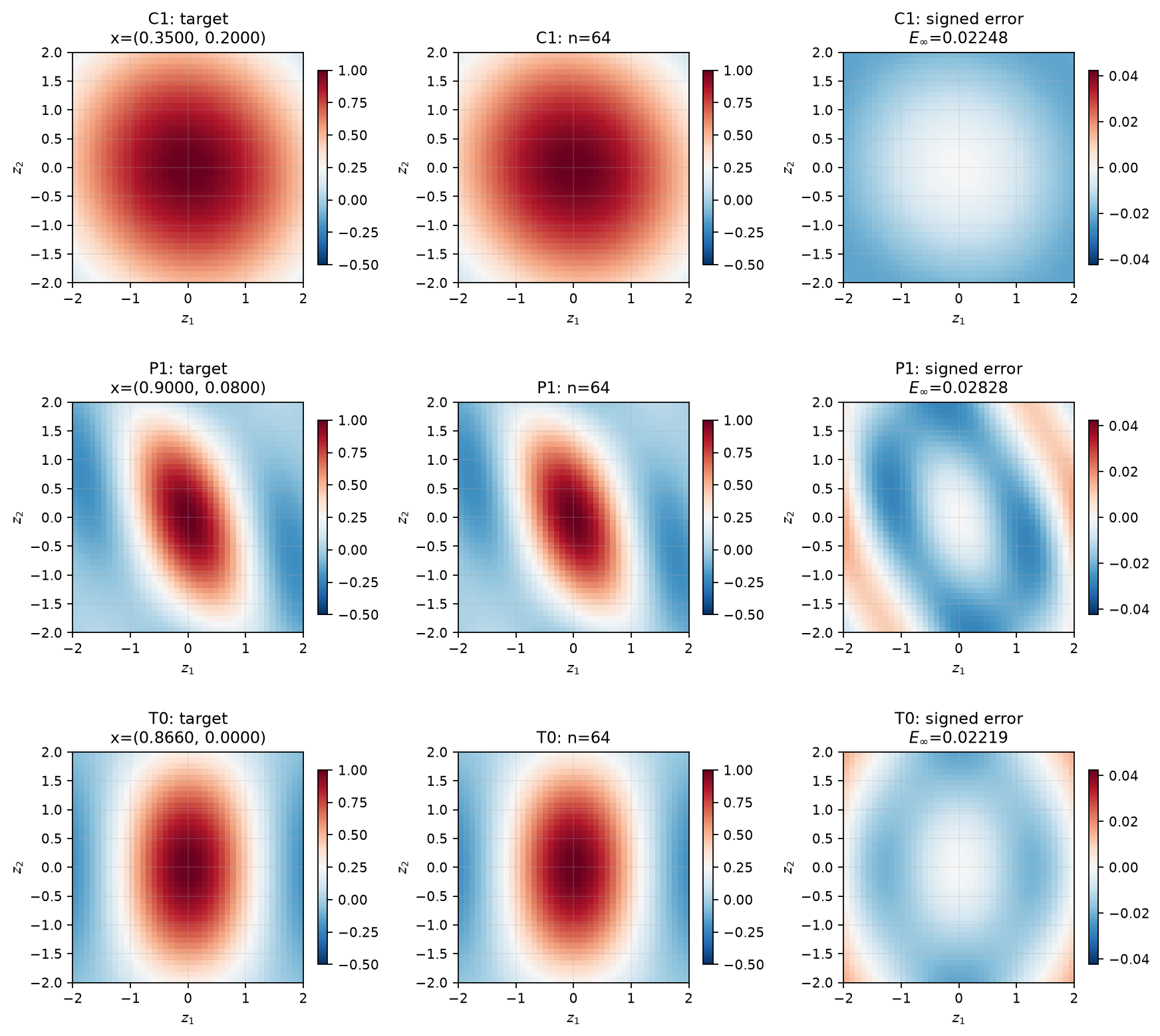}
 \caption{Theoretical profile $\Phi_x$ (left), centered polynomial kernel
 $R_{64,x}$ (middle), and signed difference (right), at C1, P1, and T0.
 These points represent a hexagonal spectrum, an off-axis parallelogram,
 and the exact transition, respectively. The first two columns share a
 common color scale; the error panels use a separate symmetric scale.}
 \label{fig:hexagon-numerical-reference}
\end{figure}

\clearpage
The comparisons show how the shape of $A_x$ appears in the local kernel.
C0 retains the sixfold symmetry of the central hexagon, while C1 and C2
show the directional distortion caused by moving the base point. At the axial corner point P0 the
spectral body is a rectangle and the target is a product of two sinc
functions. At P1 the corresponding parallelogram produces an oblique
profile. Across the nine centers, the grid errors range from $0.01175$
at C0 to $0.04241$ at P0. These are absolute errors for profiles normalized
to one at the origin.

\begin{figure}[!htbp]
 \centering
 \includegraphics[width=\textwidth]{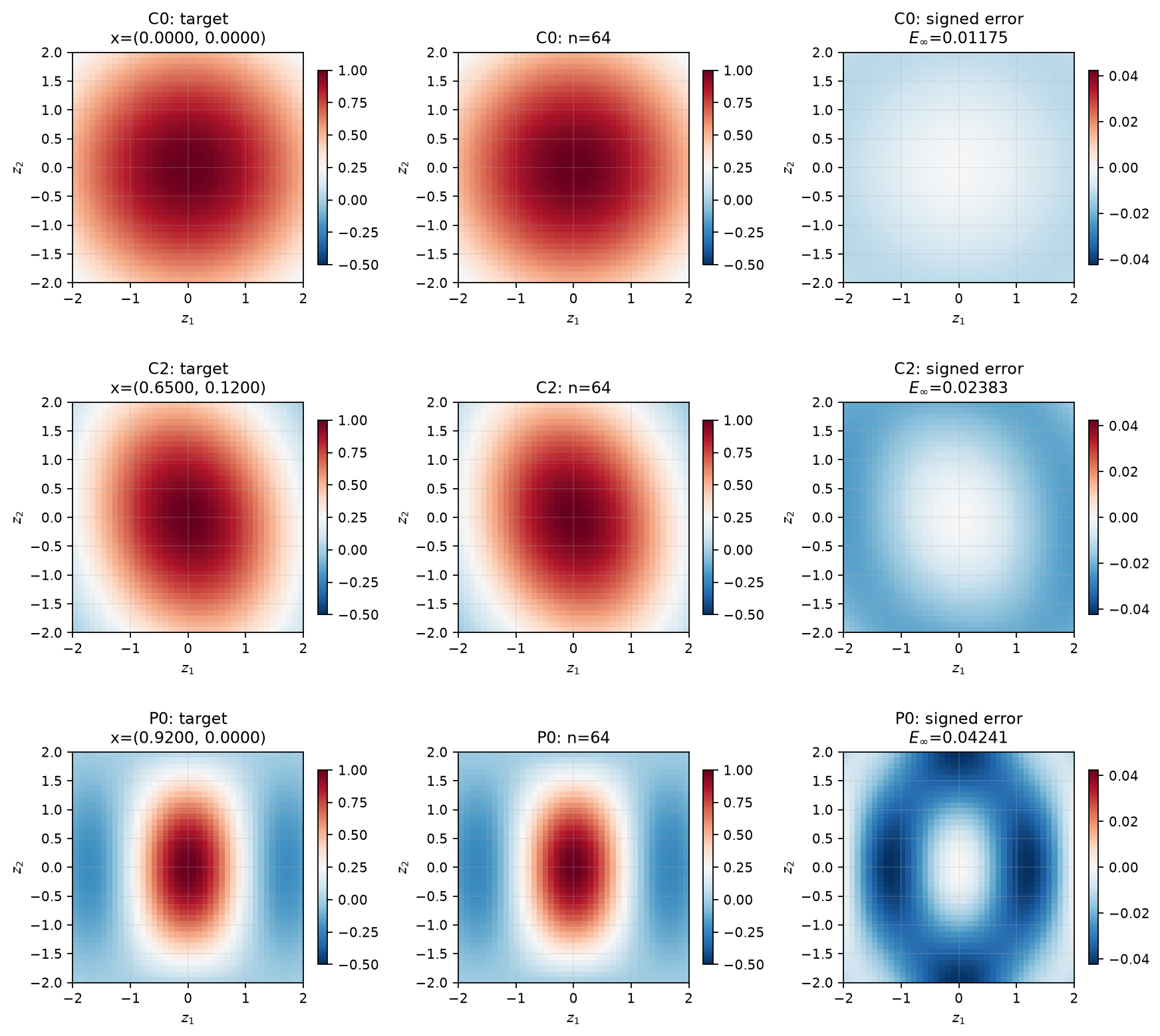}
 \caption{Comparison at the center C0, the noncentral inner point C2,
 and the axial outer point P0, with the same column conventions as in
 Figure~\ref{fig:hexagon-numerical-reference}. The finite-degree profiles
 reproduce the change from a hexagonal spectral transform to the
 rectangular product-of-sinc profile.}
 \label{fig:hexagon-numerical-inner-outer}
\end{figure}

\clearpage
The transition points also illustrate the continuity of the Fourier
profile when a pair of candidate spectral vertices ceases to be extreme.
At Tminus the local spectral body still has six extreme vertices; at T0
it has four, with the remaining two candidates lying on its edges; at
Tplus those candidates lie in its interior. Tangle gives the same
transition away from a symmetry axis. The degree-$64$ comparisons capture
the profiles on both sides and at the transition. They provide a
finite-degree illustration of the centered section of the two-point
universality limit; a convergence rate would require estimates beyond
this single-degree comparison.

\begin{figure}[!htbp]
 \centering
 \includegraphics[width=\textwidth]{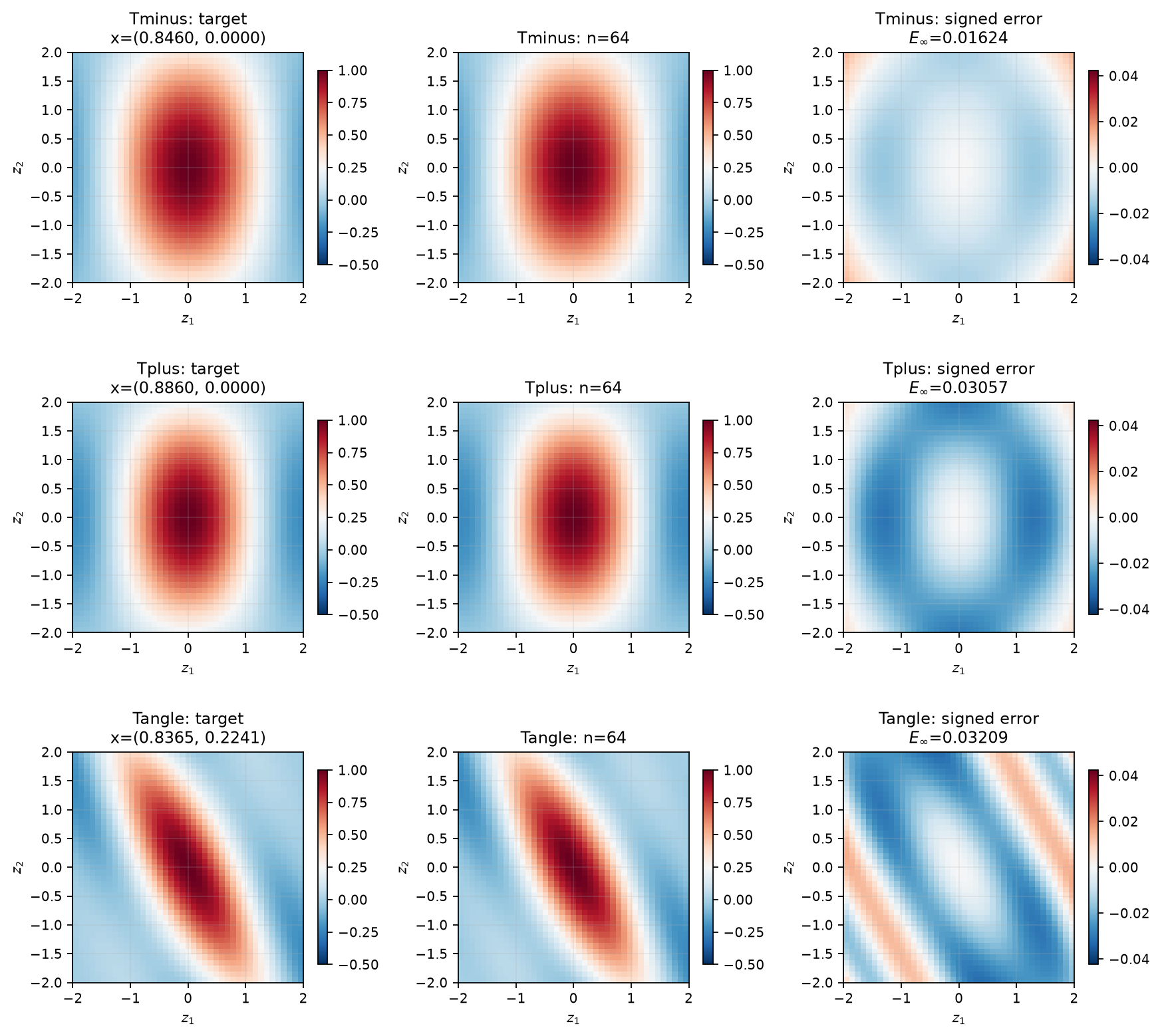}
 \caption{Profiles immediately inside and outside the incircle
 (Tminus and Tplus) and at the off-axis transition point Tangle.
 Together with the T0 row of
 Figure~\ref{fig:hexagon-numerical-reference}, these comparisons illustrate
 the spectral transition within the interior of the physical hexagon.}
 \label{fig:hexagon-numerical-transition}
\end{figure}

\clearpage

\clearpage
\section{\texorpdfstring{$L^p$}{Lp} Sampling and Least-Squares Estimates}\label{p2:app:auxiliary}

This appendix collects the general $L^p$ formulation and the
least-squares consequences of positive weighted sampling. The main
arguments use their $L^2$ setting and the Christoffel normalization
established in Section~\ref{p2:sec:sampling-frames}.

\subsection{\texorpdfstring{$L^p$}{Lp} Marcinkiewicz--Zygmund inequalities}\label{p2:subsec:lp-mz}

Fix $1\le p<\infty$ and a finite positive
measure $\mu$ for which $\norm{\cdot}_{L^p(\mu)}$ is a norm on
$\Pi_n^d$. For distinct nodes $\Lambda_n\subset E$ and positive weights
$\omega_{n,\lambda}$, independent of the polynomial being sampled, set
\[
 \nu_n=\sum_{\lambda\in\Lambda_n}\omega_{n,\lambda}\delta_\lambda,
 \qquad
 \norm{c}_{\ell^p(\omega_n)}
 =\left(\sum_{\lambda\in\Lambda_n}
              \omega_{n,\lambda}|c_\lambda|^p\right)^{1/p}.
\]
A weighted $L^p$ Marcinkiewicz--Zygmund (MZ) array satisfies
\begin{equation}\label{p2:eq:lp-mz}
 a_p\norm{q}_{L^p(\mu)}^p
 \le \sum_{\lambda\in\Lambda_n}
              \omega_{n,\lambda}|q(\lambda)|^p
 \le b_p\norm{q}_{L^p(\mu)}^p,
 \qquad q\in\Pi_n^d,
\end{equation}
with $0<a_p\le b_p<\infty$ independent of $n$. Thus evaluation preserves the norm up to fixed factors when the data are measured in $\ell^p(\omega_n)$. The weights may depend on $n$ and the whole node configuration.

For interpolation we require the upper inequality in \eqref{p2:eq:lp-mz}
and the following control on the interpolant. For every data vector
$c\in\C^{\Lambda_n}$ there is $q\in\Pi_n^d$ such that
\begin{equation}\label{p2:eq:lp-interpolation}
 q(\lambda)=c_\lambda\quad(\lambda\in\Lambda_n),\qquad
 \norm{q}_{L^p(\mu)}\le C_p\norm{c}_{\ell^p(\omega_n)},
\end{equation}
where $C_p$ is independent of $n$. In the Hilbert case the solution of minimum norm depends linearly on $c$. For general $p$, bijectivity of evaluation also gives a linear solution operator.

The endpoint $p=\infty$ uses the uniform norm on $E$ and
$\norm{c}_{\ell^\infty}=\max_{\lambda}|c_\lambda|$. Since every node has
positive mass, this maximum is also the $L^\infty(\nu_n)$ norm. Thus the MZ condition is
\begin{equation}\label{p2:eq:linfty-mz}
 a_\infty\norm{q}_{L^\infty(E)}
 \le \max_{\lambda\in\Lambda_n}|q(\lambda)|
 \le \norm{q}_{L^\infty(E)},
\end{equation}
with $a_\infty>0$ independent of $n$. This is exactly a uniform bound on
the norming constants in \eqref{p2:eq:beta}. Endpoint interpolation is defined
by \eqref{p2:eq:lp-interpolation} with the corresponding uniform norms.

\begin{proposition}[Exact cardinality and $L^p$ complete interpolation]
\label{p2:prop:lp-critical}
For each fixed $1\le p\le\infty$, sampling and interpolation together
force $|\Lambda_n|=N_n$. At this cardinality the MZ inequality is equivalent
to a uniformly bounded evaluation isomorphism with uniformly bounded
inverse. In particular, for $1\le p<\infty$,
\eqref{p2:eq:lp-mz} gives interpolation constant $C_p=a_p^{-1/p}$.
\end{proposition}

\begin{proof}
Let $V_n q=(q(\lambda))_{\lambda\in\Lambda_n}$. The lower MZ bound makes $V_n$ injective, so $|\Lambda_n|\ge N_n$. Interpolation makes $V_n$ surjective, so $|\Lambda_n|\le N_n$. When the dimensions agree, the lower MZ bound alone gives bijectivity. Substituting $q=V_n^{-1}c$ in \eqref{p2:eq:lp-mz} then gives $\norm{V_n^{-1}c}_{L^p(\mu)}\le a_p^{-1/p}\norm{c}_{\ell^p(\omega_n)}$. Conversely, the operator norms of $V_n$ and $V_n^{-1}$ give the upper and lower MZ bounds, respectively. At $p=\infty$ the same argument uses \eqref{p2:eq:linfty-mz}.
\end{proof}

\subsection{Lebesgue constants at exact cardinality}

At exact cardinality, the endpoint condition also identifies the
Lebesgue constant precisely. For a bijective evaluation map, let $\ell_{n,\lambda}$ denote the Lagrange
polynomials and put $I_nc=\sum_\lambda c_\lambda\ell_{n,\lambda}$. Then
\begin{equation}\label{p2:eq:critical-lebesgue}
 \beta_n(\Lambda_n)
 =\norm{I_n}_{\ell^\infty\to L^\infty(E)}
 =\sup_{t\in E}\sum_{\lambda\in\Lambda_n}
                         |\ell_{n,\lambda}(t)|.
\end{equation}
Indeed, each polynomial is the interpolant of its nodal values, which gives the first equality. For fixed $t$, choose the phases of the data so that all nonzero terms $c_\lambda\ell_{n,\lambda}(t)$ have the same argument. Taking the supremum over $t$ gives the second equality.

\begin{remark}[Intrinsic point-value scales]\label{p2:rem:lp-scales}
For $1\le p<\infty$ put
\begin{equation}\label{p2:eq:p-christoffel}
 \lambda_{n,p}^{\mu}(x)
 =\min\left\{\norm{q}_{L^p(\mu)}^p:
                     q\in\Pi_n^d,\ q(x)=1\right\}.
\end{equation}
By homogeneity this minimum equals $\|\delta_x\|^{-p}$, where $\delta_x:q\mapsto q(x)$ is point evaluation. At a node $\lambda$, every polynomial with $q(\lambda)=1$ contributes at least $\omega_{n,\lambda}$ to the discrete sum. Minimizing the upper MZ inequality over these polynomials yields the upper estimate below. At exact cardinality evaluation is bijective. The Lagrange polynomial at $\lambda$ has discrete norm to the $p$th power exactly $\omega_{n,\lambda}$. Applying the lower MZ inequality to this polynomial therefore gives the matching lower estimate, yielding
\begin{equation}\label{p2:eq:lp-weight-comparison}
 a_p\lambda_{n,p}^{\mu}(\lambda)
 \le\omega_{n,\lambda}
 \le b_p\lambda_{n,p}^{\mu}(\lambda).
\end{equation}
Comparing the corresponding discrete sums shows that the intrinsic weights $\lambda_{n,p}^{\mu}(\lambda)$ give MZ bounds $a_p/b_p$ and $b_p/a_p$.
\end{remark}

\subsection{Frame bounds, dual reconstruction, and noise}\label{p2:subsec:reconstruction}

For positive weights $\omega_{n,\lambda}$ define
\[
 T_n^\omega p=\big(\sqrt{\omega_{n,\lambda}}p(\lambda)\big)_{\lambda\in\Lambda_n},
 \qquad v_{n,\lambda}=\sqrt{\omega_{n,\lambda}}K_n(\cdot,\lambda).
\]
By the reproducing property, the entries of $T_n^\omega p$ are the frame coefficients $\langle p,v_{n,\lambda}\rangle_\mu$. Hence the weighted MZ inequality is equivalent to $\{v_{n,\lambda}\}$ being a frame for $\Pi_n^d$ with bounds independent of $n$. The normalization in \eqref{p2:eq:sampling-operator} is the special choice $\omega_{n,\lambda}=K_n(\lambda,\lambda)^{-1}$.

\begin{proposition}[Stable least-squares reconstruction]\label{p2:prop:frame-reconstruction}
Suppose
\begin{equation}\label{p2:eq:frame-omega}
 a\norm{p}_\mu^2\le\norm{T_n^\omega p}_2^2\le b\norm{p}_\mu^2
 \qquad(p\in\Pi_n^d)
\end{equation}
with $a,b>0$ independent of $n$. Then
\[
 S_n=(T_n^\omega)^*T_n^\omega,\qquad
 \mathcal R_n=S_n^{-1}(T_n^\omega)^*
\]
satisfy $aI\le S_n\le bI$, $\mathcal R_nT_n^\omega=I$, and
\begin{equation}\label{p2:eq:frame-noise}
 \norm{\mathcal R_n(T_n^\omega p+\eta)-p}_\mu
 \le a^{-1/2}\norm{\eta}_2.
\end{equation}
For any weighted data vector $y$, $\mathcal R_n y$ is the unique minimizer of $\norm{T_n^\omega q-y}_2$ over $q\in\Pi_n^d$. For exact data this gives the dual-frame formula
\[
 p=\sum_{\lambda\in\Lambda_n}\sqrt{\omega_{n,\lambda}}p(\lambda)
          S_n^{-1}v_{n,\lambda}.
\]
For $f\in C(E)$ and noisy weighted samples $T_n^\omega f+\eta$,
\begin{align}\label{p2:eq:frame-approximation}
 \norm{f-\mathcal R_n(T_n^\omega f+\eta)}_\mu
 &\le\inf_{p\in\Pi_n^d}\left\{\norm{f-p}_\mu
      +a^{-1/2}\norm{T_n^\omega(f-p)}_2\right\}
      +a^{-1/2}\norm{\eta}_2\notag\\
 &\le(1+\sqrt{b/a})\sqrt{\mu(E)}
       \inf_{p\in\Pi_n^d}\norm{f-p}_{L^\infty(E)}
       +a^{-1/2}\norm{\eta}_2.
\end{align}
\end{proposition}

\begin{proof}
Taking the squared norm of $T_n^\omega p$ gives the quadratic form of $S_n$, so \eqref{p2:eq:frame-omega} yields $aI\le S_n\le bI$. In particular $S_n$ is invertible. The normal equations for data $y$ are $S_nq=(T_n^\omega)^*y$, hence their unique solution is $\mathcal R_ny$. For $y=T_n^\omega p$ the solution equals $p$. Expanding $(T_n^\omega)^*y$ in the kernel vectors gives the dual-frame formula. Moreover $\mathcal R_n\mathcal R_n^*=S_n^{-1}$, so $\norm{\mathcal R_n}\le a^{-1/2}$. Applying this bound to a perturbation $\eta$ gives \eqref{p2:eq:frame-noise}.

For continuous $f$, subtract any $p\in\Pi_n^d$ and use $\mathcal R_nT_n^\omega p=p$. The triangle inequality gives the first line of \eqref{p2:eq:frame-approximation}. To control the discrete residual, apply the upper MZ bound to the constant polynomial $1$, which yields $\sum_\lambda\omega_{n,\lambda}\le b\mu(E)$. Consequently $\norm{T_n^\omega(f-p)}_2\le\sqrt{b\mu(E)}\norm{f-p}_\infty$. Together with $\norm{f-p}_\mu\le\sqrt{\mu(E)}\norm{f-p}_\infty$, this proves the second line.
\end{proof}

In \eqref{p2:eq:frame-noise} the data coordinates already include the factors $\sqrt{\omega_{n,\lambda}}$. Thus an error $e_\lambda$ in the raw value $p(\lambda)$ contributes $\omega_{n,\lambda}|e_\lambda|^2$ to the squared noise norm. For a continuous target, \eqref{p2:eq:frame-approximation} also includes the residual at the nodes, since the MZ bounds apply to polynomials.

\clearpage
\section{Spherical Harmonic Analogue}
\label{sph:app}

The spherical case separates spatial density from spectral shape particularly
clearly. The diagonal reproducing kernel is constant, while the local spectral
body at every point is a Euclidean ball. We verify the recovery and tail
estimates needed to transfer the projection-defect argument of
\cref{p2:thm:near-critical-transfer} to this setting. The nonexistence
theorems of Ortega-Cerd\`a \cite[Theorem~1.1]{Ortega2026} and Wan
\cite[Theorem~1.1]{Wan2026} then give a positive relative redundancy
gap at every fixed stability level. At asymptotic isometry the window
count gives an explicit lower bound in terms of packing in the tangent
dimension, with a consequence for spherical designs
(\cref{sph:sec:design-count}).

\subsection{The cumulative space and its local kernel}

Let $d\ge2$, let $dS$ be ordinary surface measure on the unit sphere
$\mathbb S^d\subset\R^{d+1}$, and write $s_d=|\mathbb S^d|$. Set
\[
 V_L=\bigoplus_{\ell=0}^L\mathcal H_\ell(\mathbb S^d),\qquad
 N_L=\dim V_L=\binom{L+d}{d}+\binom{L+d-1}{d}
       \sim\frac{2}{d!}L^d.
\]
Here $\mathcal H_\ell$ consists of spherical harmonics of degree $\ell$;
thus $V_L$ is the cumulative space, not a single eigenspace. Let $K_L$
be its reproducing kernel in $L^2(dS)$. Rotational invariance gives
\begin{equation}\label{sph:eq:diagonal}
 K_L(x,x)=k_L:=N_L/s_d,\qquad
 q_d:=\lim_{L\to\infty}L^{-d}k_L
       =\frac{|B^d|}{(2\pi)^d},
\end{equation}
where $B^d$ is the unit ball in $\R^d$.

\begin{proposition}[The local spectral ball]\label{sph:prop:kernel}
Identify $T_x\mathbb S^d$ with $\R^d$ by an orthonormal frame. Then
\begin{equation}\label{sph:eq:kernel}
 L^{-d}K_L\bigl(\exp_x(u/L),\exp_x(v/L)\bigr)
 \longrightarrow k_B(u-v):=
 \frac{1}{(2\pi)^d}\int_{B^d}e^{i(u-v)\cdot\xi}\dd\xi.
\end{equation}
The convergence is uniform for $u,v$ in compact sets and uniform over
centers and orthonormal frames.
\end{proposition}

\begin{proof}
The addition formula, in cumulative form, reads
\[
 K_L(x,y)=\frac{N_L}{s_d}
 \frac{P_L^{(d/2,d/2-1)}(x\cdot y)}
      {P_L^{(d/2,d/2-1)}(1)};
\]
see, for example, \cite{EhlerGrochenig2026}, after converting
from probability surface measure. The Mehler--Heine formula
\cite[Theorem~8.1.1]{Szego1975} gives the normalized limit
\[
 \Gamma(d/2+1)\left(\frac{2}{r}\right)^{d/2}J_{d/2}(r)
   =\frac{1}{|B^d|}\int_{B^d}e^{i(u-v)\cdot\xi}\dd\xi,
 \qquad r=\norm{u-v},
\]
with continuous value $1$ at $r=0$. Since
$L\,\dist(\exp_x(u/L),\exp_x(v/L))\to\norm{u-v}$ locally uniformly,
\eqref{sph:eq:diagonal} proves the assertion. Rotations make it uniform
in the center and frame.
\end{proof}

In particular, on $\mathbb S^2$ the normalized limit is $2J_1(r)/r$:
the local spectrum is a two-dimensional disk. The Christoffel function
is exactly $s_d/N_L$, and the probability measure
$N_L^{-1}K_L(x,x)dS(x)$ is exactly $dS/s_d$ at every degree.

\subsection{Localized recovery and weighted sampling limits}

Kernel convergence alone does not transfer the lower sampling bound:
the samples outside each fixed rescaled window must also be controlled.
Smooth spectral cutoffs provide that control.

\begin{lemma}[Spherical recovery with rapid spatial decay]
\label{sph:lem:recovery}
Suppose $\widehat f\in C_c^\infty(B^d)$, with Fourier convention as in
Part~\ref{part:local}. For every $x\in\mathbb S^d$ there are
$p_{L,x}\in V_L$ such that
\begin{align}
 p_{L,x}(\exp_x(u/L))&\longrightarrow f(u)
          &&\text{locally uniformly},\label{sph:eq:recovery}\\
 L^d\norm{p_{L,x}}_{L^2(dS)}^2&\longrightarrow\norm f_2^2,
          \label{sph:eq:norm}\\
 |p_{L,x}(y)|&\le C_m(1+L\dist(x,y))^{-m}
          &&(m\in\N).\label{sph:eq:decay}
\end{align}
The constants and convergence are uniform in $x$ and in orthonormal
tangent frames; in particular $x$ may depend on $L$.
\end{lemma}

\begin{proof}
We construct a smooth spectral cutoff whose local limit reproduces
the desired frequencies. Its off-diagonal decay will control both the
spatial tail and the rescaled norm of the recovered polynomial.

\smallskip\noindent\emph{1. A localized spectral cutoff.}
Choose $0<r<s<1$ with $\supp\widehat f\subset rB^d$, and an even
$h\in C_c^\infty((-s,s))$ with $h=1$ on $[-r,r]$. Put
\[
 c=(d-1)/2,\qquad \mathcal D=\sqrt{-\Delta_{\mathbb S^d}+c^2},
 \qquad Q_L=h(\mathcal D/L).
\]
On $\mathcal H_\ell$, the eigenvalue of $\mathcal D$ is $\ell+c$.
Consequently the range of $Q_L$ lies in $V_L$ for large $L$. Its kernel
$Q_L(y,z)$ satisfies, for every $m$,
\begin{equation}\label{sph:eq:smooth-kernel}
 |Q_L(y,z)|\le C_m L^d(1+L\dist(y,z))^{-m},
\end{equation}
by the smooth multiplier localization theorem
\cite[Theorem~3.5]{NPW2006}.
For completeness, the fixed shift $c$ in the eigenvalues does not
affect uniformity of this estimate. Written in the unshifted variable
$\ell/L$, the cutoff is $h_L(t)=h(|t|+c/L)$. For all large $L$ it
is constant on a fixed neighborhood of zero, so it is smooth there;
its supports and every fixed collection of derivative norms are
uniformly bounded. The localization theorem therefore applies with
constants independent of $L$.

\smallskip\noindent\emph{2. The local multiplier limit.}
With $K_j=0$ for $j<0$, summation by
parts gives the exact identity
\[
 Q_L(y,z)=-\int_0^\infty h'(t)
                    K_{\lfloor Lt-c\rfloor}(y,z)\dd t.
\]
Indeed, the coefficient of the degree-$\ell$ projection on the
right is $-\int_{(\ell+c)/L}^\infty h'(t)dt=h((\ell+c)/L)$.
The derivative $h'$ is supported away from zero. Applying
\cref{sph:prop:kernel} with bandwidth $\lfloor Lt-c\rfloor$ therefore
gives, locally uniformly,
\begin{equation}\label{sph:eq:smooth-limit}
 L^{-d}Q_L(\exp_x(u/L),\exp_x(v/L))
 \longrightarrow q_h(u-v):=
 (2\pi)^{-d}\int_{\R^d}h(|\xi|)e^{i(u-v)\cdot\xi}\dd\xi.
\end{equation}
Dominated convergence in the $t$ integral follows also from
$|K_j(y,z)|\le k_j$ and $j=O(L)$ on its support.
More explicitly, the local limit of $L^{-d}K_{\lfloor Lt-c\rfloor}$
is the inverse Fourier transform of $\mathbf1_{tB^d}$. Integrating
these indicators against $-h'(t)$ gives $h(|\xi|)$, which identifies
the multiplier in \eqref{sph:eq:smooth-limit}.

\smallskip\noindent\emph{3. Spatial decay and norm recovery.}
Fix $0<\rho<\pi/4$ and a smooth radial cutoff $\chi_x$, equal to one
near $x$ and supported in the cap $\dist(x,y)<\rho$. Define, by zero
outside that cap,
\[
 g_{L,x}(y)=\chi_x(y)f(L\exp_x^{-1}y),\qquad
 p_{L,x}=Q_Lg_{L,x}.
\]
The function $f$ is Schwartz, so
$|g_{L,x}(z)|\le C_m(1+L\dist(x,z))^{-m}$ for every $m$.
The triangle inequality gives
\[
 1+L\dist(x,y)
 \le(1+L\dist(y,z))(1+L\dist(x,z)).
\]
Using order $m+d+1$ in \eqref{sph:eq:smooth-kernel} and order $m$
for the decay of $g_{L,x}$ therefore gives
\[
\begin{aligned}
 &(1+L\dist(x,y))^m|p_{L,x}(y)|\\
 &\qquad\le C_m\int_{\mathbb S^d}
       L^d(1+L\dist(y,z))^{-d-1}\,dS(z)\le C_m'.
\end{aligned}
\]
The last integral is uniformly bounded by geodesic polar coordinates
and $\sin t\le t$. This proves \eqref{sph:eq:decay} on the whole
sphere, not only in the coordinate cap.

In normal coordinates the volume element is
\[
 dS(\exp_x(v/L))=L^{-d}j(v/L)\dd v,\qquad
 j(v/L)=\left(\frac{\sin(|v|/L)}{|v|/L}\right)^{d-1}.
\]
Using \eqref{sph:eq:smooth-limit}, the bound $L^{-d}|Q_L|\le C$, and
Schwartz decay of $f$, dominated convergence gives
$p_{L,x}(\exp_x(u/L))\to(q_h*f)(u)=f(u)$, locally uniformly.
The last equality follows from $h(|\xi|)\widehat f(\xi)=\widehat f(\xi)$.
Finally, normal coordinates cover the sphere except its antipode,
which has measure zero. Hence
\[
 L^d\|p_{L,x}\|_{L^2(dS)}^2
 =\int_{|v|<\pi L}|p_{L,x}(\exp_x(v/L))|^2j(v/L)\,dv.
\]
Since $0\le j(v/L)\le1$, \eqref{sph:eq:decay} bounds the integrand,
extended by zero, by $C_m(1+|v|)^{-2m}$. Choose $2m>d$ and apply
dominated convergence to obtain \eqref{sph:eq:norm}.
All constructions commute with rotations, so the same estimates and
limits hold uniformly for moving centers and tangent frames.
\end{proof}

Let $\Lambda_L\subset\mathbb S^d$ be finite and let $\omega_{L,\lambda}>0$.
We use the MZ convention
\begin{equation}\label{sph:eq:mz}
 a\norm p_{L^2(dS)}^2\le
 \sum_{\lambda\in\Lambda_L}\omega_{L,\lambda}|p(\lambda)|^2
 \le b\norm p_{L^2(dS)}^2,\qquad p\in V_L.
\end{equation}
The ratio $b/a$ concerns squared norms.

\begin{lemma}[Local mass and tangent sampling]\label{sph:lem:sampling}
Assume \eqref{sph:eq:mz} with fixed $a,b$. There is $C$, depending
only on $d,b$, such that for $R\ge1$ and all large $L$,
\begin{equation}\label{sph:eq:mass}
 \sum_{\dist(z,\lambda)<R/L}\omega_{L,\lambda}
       \le C L^{-d}(1+R)^d,\qquad z\in\mathbb S^d.
\end{equation}
Fix any centers $x_L$ and tangent frames. On the cap of fixed radius
$\rho<\pi/4$ consider the rescaled measures
\[
 \nu_L=L^d\sum_{\dist(x_L,\lambda)<\rho}
               \omega_{L,\lambda}\delta_{L\exp_{x_L}^{-1}\lambda}.
\]
They have vaguely convergent subsequences, and every vague limit $\nu$
satisfies
\begin{equation}\label{sph:eq:tangent-mz}
 a\norm f_2^2\le\int_{\R^d}|f|^2\dd\nu\le b\norm f_2^2,
                  \qquad f\in\PW(B^d).
\end{equation}
\end{lemma}

\begin{proof}
For some fixed $\rho_0>0$, \cref{sph:prop:kernel} gives
$|K_L(z,\lambda)|\ge k_L/2$ when
$\dist(z,\lambda)<\rho_0/L$. Apply the upper MZ bound to
$K_L(\cdot,z)/\sqrt{k_L}$ to obtain mass at most $4b/k_L$ in this cap.
Covering a cap of radius $R/L$ by at most $C_d(1+R)^d$ such caps proves
\eqref{sph:eq:mass}. For $R$ comparable to or larger than $L$, the same
bound follows by covering the whole sphere. This gives local boundedness
of $\nu_L$ and hence vague compactness.

For the packets in \cref{sph:lem:recovery}, dyadic annuli and
\eqref{sph:eq:mass} give, whenever $2m>d$ and $R\ge1$,
\begin{equation}\label{sph:eq:sample-tail}
 L^d\sum_{\dist(x_L,\lambda)>R/L}
       \omega_{L,\lambda}|p_{L,x_L}(\lambda)|^2
       \le C_{f,m}R^{d-2m}.
\end{equation}
Pass to the limit in \eqref{sph:eq:mz}, multiplied by $L^d$, first
with a compactly supported continuous window, then let its radius tend
to infinity. Local uniform recovery, vague convergence,
\eqref{sph:eq:norm}, and \eqref{sph:eq:sample-tail} prove
\eqref{sph:eq:tangent-mz} for $\widehat f\in C_c^\infty(B^d)$.
These functions are dense in $\PW(B^d)$. The upper bound extends the
restriction map continuously to that space. This extension agrees
with pointwise evaluation, since convergence in $\PW(B^d)$ implies
uniform convergence on $\R^d$, and hence convergence in $L^2(\nu)$
on compact sets. The lower bound extends by continuity as well.
\end{proof}

\subsection{Near-critical cardinality forces a basis on the ball}

\begin{theorem}[Spherical near-critical transfer]\label{sph:thm:transfer}
Suppose \eqref{sph:eq:mz} holds with fixed $0<a\le b<\infty$ along
any sequence $L\to\infty$, and
$M_L:=|\Lambda_L|=N_L+o(N_L)$. Then $L^2(B^d)$ admits an exponential
Riesz basis. Its normalized exponentials can be chosen with Riesz
bounds $a/b$ and $b/a$.
\end{theorem}

\begin{proof}
Let $T_Lp=(\sqrt{\omega_{L,\lambda}}p(\lambda))_\lambda$,
$H_L=T_LT_L^*$, and let $P_L$ be the orthogonal projection onto
$\ker T_L^*$. Its rank is $r_L=M_L-N_L=o(L^d)$, and
\begin{equation}\label{sph:eq:projection}
 a(I-P_L)\le H_L\le b(I-P_L),\qquad
 \sum_\lambda(P_L)_{\lambda\lambda}=r_L.
\end{equation}
The recovery lemma is uniform in the center, so we may select a
different favorable center at each bandwidth. To make all fixed
rescaled windows favorable simultaneously, choose
\[
 R_L=\left(\frac{L^d}{\max\{1,r_L\}}\right)^{1/(2d)}.
\]
Then $R_L\to\infty$, $R_L/L\le L^{-1/2}\to0$, and
$R_L^d r_L/L^d\to0$: for $r_L\ge1$ the latter quantity is
$\sqrt{r_L/L^d}$, and for $r_L=0$ it is zero. For
\[
 e_L(x)=\sum_{\dist(x,\lambda)<R_L/L}(P_L)_{\lambda\lambda},
\]
rotational invariance gives the exact average identity
\[
 \int_{\mathbb S^d}e_L(x)\dd S(x)
       =|\operatorname{cap}(R_L/L)|\,r_L
       \le C_d R_L^d L^{-d}r_L\longrightarrow0.
\]
Select $x_L$ with $e_L(x_L)$ at most this integral divided by $s_d$,
and rotate the entire rule so that $x_L$ is a fixed pole $o$.
For each fixed $R$, let $J_{L,R}$ be the nodes in the cap of radius
$R/L$ about $o$. Its compressed projection has norm at most
$e_L(x_L)=o(1)$, so \eqref{sph:eq:projection} gives
\begin{equation}\label{sph:eq:block}
 a(1-e_L(x_L))I\le(H_L)_{J_{L,R}}\le bI.
\end{equation}
In particular, $c_{L,\lambda}:=\omega_{L,\lambda}k_L$ lies between
$a(1-e_L(x_L))$ and $b$ for every node in a fixed rescaled window.
Together with \eqref{sph:eq:mass}, this bounds the number of nodes
in each such window independently of $L$.

Dividing each block by its diagonal square roots gives the normalized
Gram lower bound $a(1-e_L(x_L))/b$. For two distinct nodes its least
eigenvalue is $1-|K_L(\lambda,\lambda')|/k_L$. By
\cref{sph:prop:kernel}, any two limiting rescaled nodes $u,v$ satisfy
\[
 |k_B(u-v)|/q_d\le1-a/b.
\]
Continuity at zero therefore excludes collisions and gives a positive
separation independent of the window. Extracting on successively larger
windows produces a separated set $\Gamma\subset\R^d$ and coefficients
$c_\gamma\in[a,b]$ such that
\[
 \nu_L\longrightarrow
       \nu=q_d^{-1}\sum_{\gamma\in\Gamma}c_\gamma\delta_\gamma
       \quad\text{vaguely}.
\]
This extraction uses the local counting bounds and then convergence of
the coefficients at each isolated limit point.

By \cref{sph:lem:sampling} and Plancherel, the system
\[
 \left\{\sqrt{c_\gamma}\,|B^d|^{-1/2}e^{-i\gamma\cdot\xi}:
                   \gamma\in\Gamma\right\}
\]
is a frame for $L^2(B^d)$ with bounds $a,b$, in particular complete.
For any finite subset of $\Gamma$, the corresponding blocks of $H_L$
converge, by \cref{sph:prop:kernel}, to its weighted exponential Gram
matrix. Equation~\eqref{sph:eq:block} puts its eigenvalues in $[a,b]$.
Thus the same system is a Riesz sequence and hence a Riesz basis.
Removing the factors $\sqrt{c_\gamma}$ gives bounds $a/b,b/a$.
Every extraction took place within the originally prescribed bandwidth
sequence, as required.
\end{proof}

\begin{corollary}[Positive relative redundancy at fixed conditioning]
\label{sph:cor:gap}
Let $d\ge2$ and $1\le\kappa<\infty$. There are
$\delta(d,\kappa)>0$ and $L_0(d,\kappa)$ such that every positive-weight
sampling rule on $V_L$ with upper-to-lower squared-norm MZ bound ratio
at most $\kappa$ satisfies
\[
 |\Lambda_L|\ge(1+\delta(d,\kappa))N_L,
                 \qquad L\ge L_0(d,\kappa).
\]
This conclusion uses the nonexistence of exponential Riesz bases on
$B^d$ proved independently in \cite[Theorem~1.1]{Wan2026} and
\cite[Theorem~1.1]{Ortega2026}.
\end{corollary}

\begin{proof}
If the conclusion failed, there would be $L_j\to\infty$ and rules with
$1\le M_{L_j}/N_{L_j}<1+1/j$ and bound ratio at most $\kappa$.
Rescaling their weights gives fixed bounds $1,\kappa$.
\cref{sph:thm:transfer} would then produce an exponential Riesz basis
of $L^2(B^d)$, contradicting the cited theorem.
\end{proof}

\subsection{Relation to spherical sampling constructions}

Marzo's density inequalities \cite{Marzo2007} give the critical threshold.
Ehler and Gr\"ochenig \cite{EhlerGrochenig2026} give explicit upper and
lower bounds on the mass of a sampling measure in geodesic balls at
scale $L^{-1}$. For discrete measures with arbitrary positive weights,
\cref{sph:cor:gap} bounds the number of distinct nodes at a fixed
stability ratio without an a priori separation assumption. The proof supplies no
effective value or optimal dependence for $\delta(d,\kappa)$.

Marzo--Ortega-Cerd\`a \cite[Theorem~2.1]{MarzoOrtega2008} show that
Fekete points of degree $\lfloor(1+\eta)L\rfloor$ form an $L^p$-MZ
array for $V_L$, for every fixed $\eta>0$ and $1\le p\le\infty$.
Their cardinality ratio tends to $(1+\eta)^d$, while the stability
bounds may depend on $\eta$. This is compatible with the corollary:
arbitrarily small fixed positive redundancy is possible, but a sequence
of redundancies tending to zero cannot retain a common stability bound.

The restriction $d\ge2$ in \cref{sph:cor:gap} is essential. On $\mathbb S^1$, the space of
trigonometric polynomials of degree at most $L$ has dimension $2L+1$,
and the equally spaced grid satisfies the exact identity
\[
 \int_0^{2\pi}|p(e^{it})|^2\dd t
 =\frac{2\pi}{2L+1}\sum_{j=0}^{2L}
       |p(e^{2\pi i j/(2L+1)})|^2.
\]
Indeed, the grid integrates all Fourier modes of $|p|^2$, whose
nonzero frequencies have magnitude at most $2L$. Its local spectrum
is the interval $[-1,1]$, which admits an exponential orthonormal basis.
For $d\ge2$, constant Christoffel density still prescribes a uniform
distribution of local degrees of freedom, but the ball-shaped local
spectrum obstructs uniformly stable sampling with vanishing relative
redundancy.

\subsection{Spherical designs and packing in the tangent dimension}
\label{sph:sec:design-count}

Let $N(d,s)$ be the minimum cardinality of an equal-weight spherical
$s$-design on $\mathbb S^d$, and write $\Delta_d^{\mathrm{LP}}=
\Delta_{\mathrm{LP}}(B^d)$. For the ball, rotation averaging preserves
admissibility and both values at the origin in \eqref{piv:lpdef}, so
this is the radial Euclidean packing bound used by Gorbachev.
His planar estimate for $\mathbb S^2$ uses explicit packing certificates
\cite[Sections~2--3]{Gorbachev2014}. In his later comparison of the
design and packing linear programs \cite[slides~14--18]{Gorbachev2026Seminar},
the general estimate uses the Euclidean program in dimension $d+1$.
For $d=1,3$ he obtains the reciprocal of the packing LP bound in
dimension $d$ itself. Passing to the tangent sampling measure gives
this same-dimensional bound on the actual node count for every fixed $d$.

\begin{proposition}\label{sph:prop:design-count}
Let $d\ge1$. Suppose a sequence of positive weighted arrays on
$\mathbb S^d$, with $M_L$ distinct nodes, is asymptotically isometric
on $V_L$ along degrees $L\to\infty$. Then
\begin{equation}\label{sph:eq:lp-count}
 \liminf_{L\to\infty}\frac{M_L}{N_L}
 \ \ge\ \frac{1}{\Delta_d^{\mathrm{LP}}}.
\end{equation}
In particular, in Gorbachev's normalization,
\begin{equation}\label{sph:eq:design-ninf}
 N_{\inf}(d):=\liminf_{s\to\infty}\frac{2^d N(d,s)}{\dim V_s}
 \ \ge\ \frac{1}{\Delta_d^{\mathrm{LP}}}.
\end{equation}
The same inequality holds when $N(d,s)$ is replaced by the minimum
number of nodes for positive spherical cubature of precision $s$.
\end{proposition}

\begin{proof}
For $d\ge2$, express the tangent spectrum in cycles per unit length as
$K=(2\pi)^{-1}B^d$. Fix $R\ge1$, write $Q_R=(-R/2,R/2)^d$, and use
an orthonormal frame at $x_0\in\mathbb S^d$ to form the window
$W_{L,R}=\exp_{x_0}(L^{-1}Q_R)$. With normalized Haar measure on
$SO(d+1)$, each node lies in a random rotation of $W_{L,R}$ with
probability $|W_{L,R}|/s_d$. Consequently,
\[
 \min_{g\in SO(d+1)}\#(\Lambda_L\cap gW_{L,R})
 \ \le\ M_L |W_{L,R}|/s_d,\qquad
 |W_{L,R}|=(R/L)^d\bigl(1+O_d((R/L)^2)\bigr),
\]
using the volume element in geodesic normal coordinates for fixed $R$ and large $L$.
Take a subsequence realizing the lower limit of the minimum counts
and use minimizing rotations to choose the centers and tangent frames.
By \cref{sph:lem:sampling}, a further subsequence of the rescaled measures
converges vaguely to a Parseval measure $\sigma$ on $\PW(B^d)$, since the
two sampling bounds tend to one. Each finite collection of distinct
support points in the open cube $Q_R$ has disjoint neighborhoods inside
$Q_R$ of positive $\sigma$-mass. Vague convergence therefore bounds
their number by the lower limit of the window counts.

For $F\in\mathcal F(K)$ and $0<t<1$, \cref{piv:window} now gives
\[
 \liminf_{L\to\infty}\frac{M_L}{s_d L^d}
 \ \ge\ t^d\frac{\fc{F}(0)}{F(0)}
        -C_{K,F,t}\frac{(R+1)^{d-1}}{R^d}.
\]
Let $R\to\infty$, then $t\uparrow1$, and take the supremum over $F$.
The right side becomes $|K|/\Delta_d^{\mathrm{LP}}$.
Since $N_L/(s_dL^d)\to |K|$ by \eqref{sph:eq:diagonal}, this proves
\eqref{sph:eq:lp-count}. For $d=1$, the interval tiles by a lattice,
so $\Delta_1^{\mathrm{LP}}=1$ by \cref{piv:continuity}; the lower
sampling bound implies $M_L\ge N_L$ for all sufficiently large $L$.

Finally, take $L=\lfloor s/2\rfloor$. Every positive cubature rule of
precision $s$ integrates $|p|^2$ exactly for $p\in V_L$, and hence is
isometric on $V_L$ after normalizing its weights to surface area.
Since $2^d N_{\lfloor s/2\rfloor}/\dim V_s\to1$, \eqref{sph:eq:lp-count}
gives \eqref{sph:eq:design-ninf} for both node counts.
\end{proof}

Gorbachev's quantity $N^{\mathrm{LP}}(d,s)$ is the optimum of the
spherical design linear program, with $N^{\mathrm{LP}}(d,s)\le N(d,s)$.
Proposition~\ref{sph:prop:design-count} uses the nodes of an actual
array. Whether the spherical LP optimum itself admits the same-dimensional
estimate in every dimension requires a separate argument.

\clearpage
\phantomsection

\end{document}